\documentclass[11pt]{article}
\usepackage{geometry}
\usepackage{graphicx}
\usepackage{tikz}
\usepackage{caption}
\usepackage{subcaption}
\usepackage{amssymb}
\usepackage{csquotes}
\usepackage{amsmath}
\usepackage{amsthm}
\usepackage{mathtools}
\usepackage{float}
\usepackage{longtable}
\usepackage{physics}
\usepackage{enumitem}
\usepackage{silence}
\usepackage[hidelinks]{hyperref}
\MakeOuterQuote{"}
\newtheorem{theorem}{Theorem}[section]
\newtheorem{definition}[theorem]{Definition}
\newtheorem{proposition}[theorem]{Proposition}
\newtheorem{lemma}[theorem]{Lemma}
\newtheorem{corollary}[theorem]{Corollary}

\newcommand{\R}{\mathbb{R}}
\newcommand{\C}{\mathbb{C}}
\newcommand{\Z}{\mathbb{Z}}
\newcommand{\N}{\mathbb{N}}

\newcommand{\KBS}{K_{BS}}

\newcommand{\er}{\textbf{e}_r}

\newcommand{\DT}{\mathcal{D}_T}

\newcommand{\Ren}{\Upsilon}
\newcommand{\Circ}{\Gamma}
\newcommand{\jap}[1]{\langle #1\rangle}
\newcommand{\Ip}[2]{I_{#1}(#2)}
\newcommand{\pout}{p^{\mathrm{out}}}

\newcommand{\Qm}{Q_m}
\newcommand{\flow}{\mathcal{X}}
\newcommand{\Orot}{\mathsf{O}}

\DeclareMathOperator{\divop}{div}
\DeclareMathOperator{\supp}{supp}

\DeclareMathOperator{\wind}{wind}
\DeclareMathOperator{\RN}{RN}

\newcommand{\xiL}{\xi_L}
\newcommand{\xia}{\xi_a}
\newcommand{\xib}{\xi_b}
\newcommand{\xiR}{\xi_R}
\newcommand{\prof}{\varpi}

\newcommand{\Gone}{\overline{G}}
\newcommand{\Gmax}{G_{\max}}

\newcommand{\NN}{\mathsf{N}}
\newcommand{\NNm}{\mathsf{N}^{-}}
\newcommand{\Nflat}{N^{\flat}}

\newcommand{\Lin}{L}

\newcommand{\kcoef}{k}

\newcommand{\Evans}{\mathcal{E}}

\newcommand{\Mres}{M}
\newcommand{\lamst}{\lambda_*}
\newcommand{\lamref}{\lambda_{\mathrm{ref}}}
\newcommand{\phist}{\phi_*}
\newcommand{\vphist}{\varphi_*}
\newcommand{\ellst}{\ell_*}
\newcommand{\Riesz}{\Pi}

\newcommand{\Thl}[1]{\Theta_{#1}}
\newcommand{\epsst}{\epsilon_*}

\newcommand{\Rect}{\mathcal{R}}
\newcommand{\Cell}{\mathcal{C}}
\newcommand{\sigb}{\sigma_\flat}
\newcommand{\Lam}{\Lambda_{\mathrm{en}}}

\newcommand{\sbar}{\bar s}
\newcommand{\Cmom}{C_{\mathrm{mom}}}
\newcommand{\Cstar}{C_{*}}
\newcommand{\Cvel}{C_{\mathrm{vel}}}
\newcommand{\Crelay}{C_{\mathrm R}}
\newcommand{\Cerr}{C_{\mathrm{err}}}
\newcommand{\maj}{\mathfrak{a}}
\newcommand{\majst}{\mathfrak{a}_*}
\newcommand{\wlin}{w^{\mathrm{lin}}}
\newcommand{\hin}{h^{\mathrm{in}}}
\newcommand{\hout}{h^{\mathrm{out}}}
\newcommand{\etin}{\eta^{\mathrm{in}}}
\newcommand{\etout}{\eta^{\mathrm{out}}}
\newcommand{\KO}{K_{O}}
\newcommand{\Pprod}{\mathsf{P}}
\newcommand{\Ptil}{\tilde{\mathsf{P}}}
\newcommand{\CDuh}{C^{\mathrm{Duh}}_\mu}
\newcommand{\Cls}{C_{\ell_*}}
\newcommand{\Cp}{C_{\mathrm{p}}}
\newcommand{\Cnab}{C_{\nabla}}
\newcommand{\cout}{c_{\mathrm{out}}}
\newcommand{\Abar}{\bar A}
\newcommand{\Klo}{K_{\mathrm{lo}}}
\newcommand{\Kup}{K_{\mathrm{up}}}
\newcommand{\Cseed}{C_{\mathrm{seed}}}
\newcommand{\Tpin}{T_{\mathrm{pin}}}
\newcommand{\rT}{\mathsf{r}}
\newcommand{\Tz}{\mathcal{T}_0}
\newcommand{\Tstar}{T_\star}
\newcommand{\Tend}{T_{\mathrm{end}}}
\newcommand{\tsh}{t^\sharp}

\newcommand{\ksh}{\kappa^\sharp}
\newcommand{\kseed}{\kappa_{\mathrm{seed}}}
\newcommand{\kL}{\kappa_{\mathrm{L}}}
\newcommand{\tauK}{\tau_{\mathrm{K}}}
\newcommand{\Ain}{A^{\mathrm{in}}}
\newcommand{\cin}{c^{\mathrm{in}}}
\newcommand{\qs}{q_{\mathrm{s}}}
\newcommand{\Errg}{\mathrm{Err}_g}
\newcommand{\Errpre}{\mathrm{Err}_{\mathrm{pre}}}
\newcommand{\err}{\mathrm{err}}

\newcommand{\dH}{\delta_0^{\mathrm{H}}}

\newcommand{\rhohat}{\hat\rho}
\newcommand{\Lconf}{\Lambda^{\mathrm{conf}}}
\newcommand{\Phase}{\mathcal{A}}
\newcommand{\Pin}{\mathcal{P}}
\newcommand{\Iint}{\mathcal{I}}

\newcommand{\Psione}{\Psi_1}
\newcommand{\chiseed}{\chi_{\mathrm{s}}}
\newcommand{\chiout}{\chi_{\mathrm{o}}}

\newcommand{\kst}{k_*}
\newcommand{\NT}{N_{\mathrm{T}}}
\newcommand{\epsend}{\varepsilon_{\mathrm{end}}}

\title{Non-Uniqueness of Weak Solutions to the Unforced\\
Two-Dimensional Euler Equations\\ with Vorticity in $L^1\cap L^p$,\\
II: The Nonlinear Construction}
\hypersetup{pdftitle={Non-Uniqueness of Weak Solutions to the Unforced
Two-Dimensional Euler Equations with Vorticity in L1 and Lp,
II: The Nonlinear Construction},
pdfauthor={Matei P. Coiculescu}}
\author{Matei P. Coiculescu\thanks{
Department of Mathematics, New York University,
Email: \texttt{m.coiculescu@nyu.edu}
}}
\date{November 6, 2026}

\begin{document}
\maketitle

\begin{abstract}
This paper and its companion Part I prove that the incompressible Euler
equations in the plane, without external force, admit two distinct renormalized
weak solutions with the same initial vorticity $f$, where $f$ is radial and
nonnegative and lies in $(L^1\cap L^q)(\mathbb{R}^2)$ for every $q<\infty$ but is
not bounded. One solution is radially symmetric and stationary. The other is
$3$-fold rotationally symmetric, is equimeasurable with $f$ at every time, and
departs from $f$ in finite time; the two solutions are distinct in
$L^\infty_t(L^1\cap L^q)_x$ on $[0,T]\times\mathbb{R}^2$ for every
$q\in(4/3,\infty)$ and every sufficiently large $T$. The nonstationary solution
belongs to $C([0,\infty);L^q(\mathbb{R}^2))$ for every $q\in[1,\infty)$.
Part I establishes a compactness criterion for non-uniqueness and conditional
linear instability results for a vortex ring. Part II establishes the hypotheses
of those results and proves the nonlinear non-uniqueness theorem.
\end{abstract}

\section{Introduction}

\subsection{The Main Results}

We consider the incompressible Euler equations in the plane, in vorticity form:
\begin{equation}
\label{eq:euler}
\partial_t \omega + v\cdot\nabla \omega = 0,\qquad
v = \KBS\ast \omega,\qquad
\omega(t=0) = \omega_0 .
\end{equation}
Here $\KBS(x)=x^\perp/(2\pi|x|^2)$, $x^\perp=(-x_2,x_1)$, denotes the
Biot-Savart kernel, so that $v$ is divergence free and
$\partial_1v_2-\partial_2v_1=\omega$. There is no external force, and we call
$\omega_0$ the datum.

This paper is the second of two. Together with Part~I~\cite{PartI}, it proves
that, without force in the plane, uniqueness fails among renormalized weak
solutions with vorticity in $L^\infty_t(L^1\cap L^q)_x$, for every
$q\in(4/3,\infty)$ at once and for an explicit radial datum. Numbers with the
prefix I refer to statements, sections and equations of Part~I. Section~1 of
Part~I states the question and discusses the literature, including Vishik's
non-uniqueness result with a radial force \cite{Vishik1,Vishik2} (see the
account of \cite{ABCDGJK}) and the results obtained without force by convex
integration. It also describes the construction heuristically.

For $T>0$ and $1\le q\le\infty$ we write $X_q(T)$ for the space
$L^\infty_t(L^1\cap L^q)_x$ of functions on $(0,T)\times\R^2$. Renormalized weak solutions in $X_q(T)$, $4/3<q<\infty$ (after
DiPerna and Lions \cite{DiPernaLions}), and equimeasurability are defined in
Definitions~I.\ref{I-def:weak} and~I.\ref{I-def:equimeasurable}. The ring $\Phi$ of
Definition~I.\ref{I-def:ring} is a closed formula in $\log|y|$ times a fixed smooth
cutoff; it is radial, $\Phi\in C_c^\infty(\R^2)$, $0\le\Phi<1$,
$\Phi\not\equiv0$, and $\supp\Phi\subset\{a\le|y|\le1\}$ for an explicit
$a\in(0,1)$ (Lemma~I.\ref{I-lem:ring-facts}). For $A_1>0$, $\mu\in(1,2]$ and
$\sbar\in(0,a/16]$ we put
\begin{equation}
\label{eq:datum}
\begin{gathered}
f=\sum_{k\ge1}A_k\,\Phi(\cdot/r_k),\qquad A_k=A_1\mu^{k-1},\\
r_1=1,\qquad r_{k+1}=s_kr_k,\qquad s_k=\sbar(2\mu)^{-k}\qquad(k\ge1).
\end{gathered}
\end{equation}
This is Definition~\ref{def:datum}. We call $A_k\Phi(\cdot/r_k)$ ring $k$ of
$f$. It is supported in the annulus $\{ar_k\le|x|\le r_k\}$, and these annuli
are pairwise disjoint (Proposition~\ref{prop:integrability}). The main theorem
of the two papers is the following; it is Theorem~I.\ref{I-thm:main}, restated.

\begin{theorem}
\label{thm:main}
Let $\Phi$ and $a$ be as in Definition~I.\ref{I-def:ring}. There are
$\mu\in(1,\frac{17}{16}]$ and $\sbar\in(0,a/16]$ such that, for every $A_1>0$,
the function $f$ of Equation~\eqref{eq:datum} has the following properties. It
is radial and nonnegative, $f\in (L^1\cap L^q)(\R^2)$ for every $q\in[1,\infty)$,
and $f\notin L^\infty(\R^2)$. There is $\Tstar\in(0,\infty)$ such that, for
every $T>\Tstar$ and every $q\in(4/3,\infty)$, Equation~\eqref{eq:euler} with
$\omega_0=f$ has two distinct renormalized weak solutions in $X_q(T)$, in the
sense of Definition~I.\ref{I-def:weak}:
\begin{enumerate}[label=\textup{(\roman*)}]
\item the stationary solution $\omega(t)=f$, $t\in[0,T]$;
\item a solution $\omega_\infty$ such that
  $\omega_\infty\in C([0,\infty);L^p(\R^2))$ for every $p\in[1,\infty)$ and,
  for every $t\ge0$, $\omega_\infty(t)$ is equimeasurable with $f$ and
  $3$-fold rotationally symmetric.
\end{enumerate}
More precisely, there is $t_0\in[0,\Tstar)$ such that $\omega_\infty$ and $f$
are distinct elements of $X_q(T')$ for every $T'>t_0$ and every
$q\in(4/3,\infty)$.
\end{theorem}

The time translation symmetry of the Euler equations allows us to take the
last time up to which $\omega_\infty$ agrees with $f$ as the new time origin.
This gives non-uniqueness on every interval $[0,T]$ with $T>0$.

We do not treat $q\le4/3$. No symmetry is
imposed on the class; the $3$-fold rotational symmetry of $\omega_\infty$ comes from the
construction. This symmetry removes the angular modes $\ell\notin3\Z$ from
the study of the linearization of the ring, among them the modes $\ell=\pm1$
that contain its translations. The remaining modes $\ell\in3\Z$ have spectra
controlled in Part~I under Hypothesis~I.\ref{I-hyp:N}.

The proof of Theorem~\ref{thm:main} is divided between Part~I and this paper.
Part~I proves a non-uniqueness criterion (Section~I.\ref{I-sec:compact}) and the
linear analysis of the ring (Section~I.\ref{I-sec:ring}), and it uses no machine
computation. The main conclusions of its linear analysis assume three explicit
statements, collected in Section~I.\ref{I-sec:ring} as
Hypothesis~I.\ref{I-hyp:N}. The statements \textup{(N1)} and \textup{(N2)} locate
the zeros in the mode $\ell=3$ of the Evans function of the ring, an explicit
holomorphic function, for every $G\in[0,\frac1{50}]$, where $2\pi G$ is the
circulation of a point vortex at the origin (in the reference units of
Definition~I.\ref{I-def:ring}; this is $G\in[0,\Gone]$ in ring units). The
statement \textup{(N3)} is an inequality between explicit numbers controlling
the modes $\ell=\pm6$ and $|\ell|\ge9$. This paper proves that
\textup{(N1)}--\textup{(N3)} imply Theorem~\ref{thm:main}, by a nonlinear
construction followed by the criterion of Part~I, and it proves
\textup{(N1)}--\textup{(N3)}. The implication is the following reduction,
which we prove in Section~\ref{sec:bootstrap}.

\begin{proposition}
\label{prop:reduction}
Suppose that \textup{(N1)}--\textup{(N3)} of Hypothesis~I.\ref{I-hyp:N} hold. Then
Theorem~\ref{thm:main} holds.
\end{proposition}

The proof of Proposition~\ref{prop:reduction} uses no machine computation.
Appendix~\ref{app:verify} proves \textup{(N3)} by hand and verifies
\textup{(N1)} and \textup{(N2)} with rigorous interval arithmetic, the only
machine-computed part of the proof (code: the package
\texttt{euler-ring-certificate}, version~1.0.0, in the public repository
\url{https://github.com/mcoiculescu/NonUniqueEuler};
Section~\ref{app:reproduce}). It thus proves the following theorem.

\begin{theorem}
\label{thm:N-holds}
Hypothesis~I.\ref{I-hyp:N} holds: the statements \textup{(N1)}, \textup{(N2)} and
\textup{(N3)} are true.
\end{theorem}

\subsection{Results of Part I}
\label{sec:II-partI}

Part~I's introduction summarizes what this paper uses in two theorems: the
conditional non-uniqueness theorem (Theorem~I.\ref{I-thm:conditional}, the form
of Corollary~I.\ref{I-cor:defect} given by Remark~I.\ref{I-rem:conditional}) and
the linear instability of the ring under \textup{(N1)}--\textup{(N3)}
(Theorem~I.\ref{I-thm:linear}).

In detail, we use the following definitions and results of Part~I, which we
cite by their numbers there and do not restate. From
Section~I.\ref{I-sec:prelim}, we use the properties of the Biot-Savart kernel
(Lemma~I.\ref{I-lem:BS-basic}) and the multipole moments, equimeasurability and
renormalized weak solutions (Definitions~I.\ref{I-def:multipole},
I.\ref{I-def:equimeasurable} and~I.\ref{I-def:weak}, with
Equations~(I.\ref{I-eq:weak}) and~(I.\ref{I-eq:renorm})). From
Section~I.\ref{I-sec:compact}, we use the fact that a radial function in
$(L^1\cap L^q)(\R^2)$, $4/3<q<\infty$, is a stationary renormalized weak
solution (Lemma~I.\ref{I-lem:radial-stationary}), the $2$-admissible sequences
(Definition~I.\ref{I-def:adm}), and the non-uniqueness criterion of
Corollary~I.\ref{I-cor:defect}, which rests on Theorem~I.\ref{I-thm:criterion}
and which we apply in the proof of Proposition~\ref{prop:reduction}. From
Section~I.\ref{I-sec:ring}, we use the ring and the annuli $N$, $N'$
(Definition~I.\ref{I-def:ring}, Lemma~I.\ref{I-lem:ring-facts}); the mode
settings, the decomposition into angular modes, the Evans function and its
simple zeros, and the semigroup (Definition~I.\ref{I-def:sector-setting},
Lemmas~I.\ref{I-lem:sector}, I.\ref{I-lem:evans-basic}
and~I.\ref{I-lem:semigroup-construct}, Theorem~I.\ref{I-thm:evans-simple}); and
Hypothesis~I.\ref{I-hyp:N}. Under \textup{(N1)}--\textup{(N3)} we use, from the
same section, the spectrum of the linearization (Proposition~I.\ref{I-prop:SU}),
the eigendata (Definition~I.\ref{I-def:eigendata}), the continuity in $G$ of
$\lamst(G)$ and of the sizes $|\pout_*(G)|$, $|\beta(G)|$,
$\|\ellst(G)\|_{(X^0)^*}$ and $\|\phist(G)\|_{H^k}$ of the eigendata, the
semigroup bounds and the uniformity of the constants in $G$
(Propositions~I.\ref{I-prop:U} and~I.\ref{I-prop:semigroup},
Lemma~I.\ref{I-lem:uniform}).

\subsection{Strategy}

We describe the construction with references to the statements of this paper
that carry it out. The background is the datum $f$, which is radial, hence
stationary (Lemma~I.\ref{I-lem:radial-stationary}). Its inner rings have larger
amplitudes and rotate faster, and the task is to
carry an instability from an inner ring outward through every ring in a time
independent of the depth, uniformly in the number of interacting rings.

The linear input is the instability of one ring. The linearization of
Equation~\eqref{eq:euler} at $\Phi$ plus a point vortex of small circulation
$2\pi G$ at the origin, which stands for everything inside the ring, is the
operator $\Lin(G)$ of Theorem~I.\ref{I-thm:linear}. On $3$-fold rotationally
symmetric perturbations it has an eigenvalue $\lamst(G)$ of the mode $\ell=3$,
algebraically simple there, with real part $\sigma(G)>0$, and its conjugate
$\overline{\lamst(G)}$ in the mode $\ell=-3$; the rest of the spectrum has
smaller real part (Proposition~I.\ref{I-prop:SU}). Ring $k$ of $f$ evolves on
the time scale $1/A_k$ and feels the inner rings mostly through their
circulation $2\pi G_kA_kr_k^2$, with $G_k\ge0$ small when $\sbar$ is small
(Lemma~\ref{lem:inner-circ}). In ring-$k$ units (Section~\ref{sec:datum}), the linearization at
ring $k$ plus this point vortex is therefore $\Lin(G_k)$.

The approximate datum $f^{(n)}$ keeps the first $n$ rings of $f$ and replaces
ring $n$, the innermost of them, by a rearrangement of it whose coefficient
along the unstable eigenfunction has a small size of order $\kseed$, fixed
independently of $n$ (Definition~\ref{def:approx}, Lemma~\ref{lem:seed}). As
$r_n$ decreases faster than $A_n$ increases, $f^{(n)}\to f$ in
$(L^1\cap L^q)(\R^2)$ for every $q<\infty$ (Lemma~\ref{lem:adm}). Each
$f^{(n)}$ is smooth with compact support, so it has a unique global smooth
solution $\omega^{(n)}$ (Theorem~\ref{thm:smooth}), which is $3$-fold
rotationally symmetric (Lemma~\ref{lem:adm}).

The scalar linear model
\[
b'(t)=\gamma b(t)+d(t),\qquad
b(t)=e^{\gamma t}\Big(b(0)+\int_0^t e^{-\gamma s}d(s)\dd s\Big)
\]
explains the terminology used below. The evolution $e^{\gamma t}b(0)$ is the
\emph{free growth} when $\Re\gamma>0$, and the additive term $d(t)$ is the
\emph{drive}. In the ring interaction, this term comes from the velocity
induced by the other rings. Factoring out $e^{\gamma t}$ leaves a normalized
coefficient whose change is the weighted integral of the drive.

In $\omega^{(n)}$ the instability passes outward like a relay. We say that
ring $k$ saturates when its free-growth envelope, growing at the rate
$A_k\sigma(G_k)$ per unit physical time, reaches the fixed size $\delta_0$.
The 3rd angular mode of the velocity induced by ring $k$ and the rings
inside it advects the vorticity of ring $k-1$ and drives its unstable
coefficient. In ring-$(k-1)$ units, this coupling carries the factor
$\mu s_{k-1}^5$, where $s_{k-1}=r_k/r_{k-1}$ and
$\mu=A_k/A_{k-1}$ (Lemma~\ref{lem:coupling}). In ring-$k$ time $\tau=A_kt$,
this velocity component grows at rate $\sigma(G_k)$. Ring $k-1$ grows freely
at physical rate $A_{k-1}\sigma(G_{k-1})$, so its rate in $\tau$ is
$(A_{k-1}/A_k)\sigma(G_{k-1})=\sigma(G_{k-1})/\mu$.
Once $\mu$ is fixed, we separate the rings to make the central circulations
small. By continuity of $\sigma(G)$ at $G=0$, the change in each growth rate is
then only a small fraction of the gap at $G_k=G_{k-1}=0$
(Lemma~\ref{lem:detuning}). Thus
\[
\sigma(G_k)-\frac{\sigma(G_{k-1})}{\mu}
\asymp\sigma(0)\frac{\mu-1}{\mu}>0.
\]
Here $\sigma(0)=\Re\lamst(0)$ is the instability growth rate of the reference
ring without a central point vortex. After the free linear evolution of ring
$k-1$ is factored out, the driven part of its unstable coefficient is a
truncated Laplace transform of the moment of the 3rd angular mode of the vorticity in ring
$k$ and the rings inside it (Lemma~\ref{lem:coupling},
Equation~\eqref{eq:kick-A}). The positive difference of growth rates makes
this integral dominated by the contribution near the saturation of ring $k$.
The resulting amplitude remains essentially stationary during the subsequent
handoff to ring $k-1$ (Proposition~\ref{prop:kick} and its proof).

As $\mu$ is fixed, ring $k-1$ therefore starts with a coefficient of order
$s_{k-1}^5$ (the hand-off in Proposition~\ref{prop:relay-step}), and it
saturates after a physical time of about
$5\log(1/s_{k-1})/(\sigma(G_{k-1})A_{k-1})$. Since $\log(1/s_k)$ grows
linearly in $k$ and $A_k$ geometrically, the bounds on the relay times are
summable over infinitely many rings. Thus the total time for the relay remains
bounded even as the number of rings tends to infinity: every finite relay
finishes by a time $\Tstar$ independent of $n$ (Lemma~\ref{lem:time}).

The timing of the relay depends on the solution itself. At each handoff we
allow time for the inner ring's contribution to accumulate, then read off the
next ring's starting amplitude from $\omega^{(n)}$. The clock of
Definition~\ref{def:clock} uses this amplitude to predict when free growth
will bring the ring to $\delta_0$ and to schedule the next handoff.

Saturated rings continue to evolve and interact. During each active interval,
we follow the outermost activated ring and control the combined field of the
material inside it. The separation of spatial scales keeps this material
confined well inside the active ring (Proposition~\ref{prop:confine}). Its
circulation field remains in the linearization. After the handoff, the
polynomial bound on its multipole moments, together with the exponential
weight in the Duhamel integral, makes its further contribution to the
normalized unstable coefficient small (Lemma~\ref{lem:drift}). Meanwhile,
the outer quiet rings stay close to their radial profiles until activation,
with controlled feedback on the active ring (Theorem~\ref{thm:quiet}).

The active ring follows its unstable eigenfunction (Theorem~\ref{thm:active}).
These bounds are verified jointly, for all rings and all times until the outermost ring
saturates, by one bootstrap with constants independent of $n$ (proof of
Proposition~\ref{prop:bootstrap}).

When the outermost ring saturates, at $T^{(n)}\le\Tstar$, its moment
$\int\omega^{(n)}(T^{(n)})\Psione$, where $\Psione(x)=\chiout(|x|)\,\bar x^{\,3}$
with a fixed real cutoff $\chiout\in C_c^\infty((a/4,3))$, has modulus at least
$c_0>0$ independent of $n$ (proof of Proposition~\ref{prop:reduction}), while
it vanishes for $f$. This lower bound is the moment hypothesis of
Theorem~I.\ref{I-thm:conditional}. The $\omega^{(n)}$ are bounded in
$(L^1\cap L^2)(\R^2)$ uniformly in time, and compactness of renormalized
solutions in $L^2$ (Section~I.\ref{I-sec:compact}) gives a subsequence
converging to a renormalized weak solution $\omega_\infty$ with datum $f$. The
moment bound passes to the limit and separates $\omega_\infty$ from $f$
(Corollary~I.\ref{I-cor:defect}).

Section~\ref{sec:II-prelim} collects notation and preliminaries.
Section~\ref{sec:single} gives the quiet, active and seeded ring estimates.
Section~\ref{sec:datum} defines the datum, its approximants and the clock,
and Section~\ref{sec:several} controls the fields between rings and the
confinement of the inner material. Section~\ref{sec:kick} computes the driving,
and Section~\ref{sec:relay} proves one step of the relay.
Section~\ref{sec:param} chooses the parameters, and Section~\ref{sec:bootstrap}
iterates the step and completes the proof of Proposition~\ref{prop:reduction}.
Appendices~\ref{app:verify} and~\ref{app:smooth}
prove Theorems~\ref{thm:N-holds} and~\ref{thm:smooth};
Appendix~\ref{sec:construction-details} gives the exact parameter conditions
and quantitative verifications used in Sections~\ref{sec:relay}--\ref{sec:bootstrap}.

The main idea of a linear instability in physical coordinates driving the
non-uniqueness mechanism, the use of non-monotone unstable radial vortices and of
a multiscale datum, and the overall strategy of the proof are due to the author.
Within this strategy, the detailed proofs and estimates, the code for the
computer-assisted verification and the text were largely produced with
Claude Opus 5.5 (Anthropic), under the direction of the author. The author has
checked all proofs by hand and takes full responsibility for the paper.

\paragraph{Acknowledgments.} I am supported at New York University by a Junior
Fellowship of the Simons Society of Fellows, funded by Simons Foundation
International and administered by the Simons Foundation. I would also like to thank
Tristan Buckmaster for generously offering his time and advice on this and many
other projects.

\subsection{Preliminaries}
\label{sec:II-prelim}

This subsection collects the notation that we take from Part~I and the
preliminaries used only in this paper. We use all notation and conventions of
Section~I.\ref{I-sec:prelim}, including the convention on constants and
$O(\cdot)$. From
Section~I.\ref{I-sec:ring} we use the convention $m=3$, the ring units and the
number $\Gone$ of Definition~I.\ref{I-def:ring}, the annuli $N$ and $N'$ of
Section~I.\ref{I-sec:ring-def}, the spaces $X^s$ and the operators $\Lin(G)$ on
them, fixed before Proposition~I.\ref{I-prop:SU}, and the constants $c_N$, $C_N$,
$\cout$, $\|\ellst\|$, $\Cls$, $\Cp$ and $\Cnab(k)$ of
Lemma~I.\ref{I-lem:uniform}. The results and definitions of
Part~I that we use are listed in Section~\ref{sec:II-partI}.

In ring units we also use the closed annuli
\[
\NNm=\{\tfrac78a\le|y|\le\tfrac54\},\qquad \NN=\{\tfrac34a\le|y|\le\tfrac32\},\qquad
\Nflat=\{\tfrac58a\le|y|\le\tfrac74\}.
\]
They mark where the support of a ring lies at the stages of the construction:
in $\NNm$ for a quiet ring and at the start of the active phase
(Theorems~\ref{thm:quiet} and~\ref{thm:active}), in $\NN$ during the active
phase, and in $\Nflat$ under the bootstrap assumption of
Section~\ref{sec:bootstrap}, which the bootstrap improves to $\NN$. All of them
lie in $N$, the common support of the spaces $X^s$ of Part~I, and the strict
inclusion $\NN\subset\operatorname{int}\Nflat$ leaves the room that the bootstrap
needs.
We will apply the statement about $S'$ below to a transported ring support,
to keep it in the next larger annulus when its radial displacement is at most
$a/16$.

\begin{lemma}
\label{lem:annuli}
Let $a\in(0,1)$ and define the annuli above, and $N$, $N'$, with this $a$. Then
$\{a\le|y|\le1\}\subset\NNm\subset\NN\subset\operatorname{int}\Nflat\subset\Nflat\subset N\subset N'$.
Let $S\subset\{a\le|y|\le1\}$ \textup{(}respectively $S\subset\NNm$\textup{)},
and let $S'$ be a set such that every $y'\in S'$ satisfies
$\big||y'|-|y|\big|\le a/16$ for some $y\in S$. Then $S'\subset\NNm$
\textup{(}respectively $S'\subset\NN$\textup{)}.
\end{lemma}
\begin{proof}
The inner radii $a$, $\frac78a$, $\frac34a$, $\frac58a$, $\frac12a$, $\frac14a$
strictly decrease and the outer radii $1$, $\frac54$, $\frac32$, $\frac74$, $2$, $3$
strictly increase, and the last two radii of each list are those of $N$ and $N'$;
this gives the inclusion chain. If $a\le|y|\le1$ and
$||y'|-|y||\le a/16$, then $\frac{15}{16}a\le|y'|\le1+\frac a{16}$, and
$\frac{15}{16}a\ge\frac78a$ and $1+\frac a{16}\le\frac54$ because $a\le4$. If
instead $y\in\NNm$, then $\frac{13}{16}a\le|y'|\le\frac54+\frac a{16}$, and
$\frac{13}{16}a\ge\frac34a$ and $\frac54+\frac a{16}\le\frac32$, for the same
reason.
\end{proof}

Away from the support of a vorticity, complex variables turn the Biot-Savart
law into a Cauchy integral, which can be expanded in powers of the complex
variable.
In these complex-variable calculations we use a single spatial variable
$z\in\C$. For a real vector field $V=(V_1,V_2)$, we write
$V(z)=V(\Re z,\Im z)$, and use the same convention for scalar functions.
The notation $\nabla^qV(z)$ always means the tensor of ordinary Cartesian
derivatives evaluated at $(\Re z,\Im z)$; it does not denote complex
differentiation. We call $F(z)=V_1(z)-iV_2(z)$ the \textbf{complex form} of
$V$.
Since $(\KBS)_1(z)-i(\KBS)_2(z)=1/(2\pi iz)$,
the velocity of every $g\in L^1(\R^2)$ has the complex form
\begin{equation}
\label{eq:BS-complex}
(\KBS\ast g)_1(z)-i(\KBS\ast g)_2(z)=\frac1{2\pi i}\int\frac{g(\zeta)}{z-\zeta}\dd A(\zeta)
\qquad(z\in\C\setminus\supp g),
\end{equation}
where $\dd A$ is the area element, the integral converges absolutely and, by differentiation under the
integral sign, defines a holomorphic function of $z\in\C\setminus\supp g$. Let
$V$ be a real vector field on an open set $U$ whose complex form $F$ is
holomorphic on $U$. Writing $z=x_1+ix_2$ only to name its two real
coordinates, recall the Wirtinger derivatives
\[
\partial_z=\tfrac12(\partial_{x_1}-i\partial_{x_2}),\qquad
\partial_{\bar z}=\tfrac12(\partial_{x_1}+i\partial_{x_2}).
\]
Since $F$ is holomorphic, $\partial_{\bar z}F=0$, and we write
$F'=dF/dz=\partial_zF$ for its derivative in the complex variable $z$, and
$F^{(q)}=\partial_z^qF$ for its $q$-th complex derivative.
Then $\partial_{x_1}F=F'$ and $\partial_{x_2}F=iF'$, so every
partial derivative of order $q$ of $F$ is $i^{\gamma_2}F^{(q)}$ for the
corresponding multi-index $\gamma$. Fix one of the $2^q$ ways of applying $q$
successive derivatives, choosing $\partial_{x_1}$ or $\partial_{x_2}$ at each
step. Since $V_1=\Re F$ and $V_2=-\Im F$, the corresponding derivatives of
$V_1$ and $V_2$ at $z$ are the real part and minus the imaginary part of the
corresponding derivative of $F$ at $z$, a number of modulus $|F^{(q)}(z)|$; so this
pair of entries of $\nabla^qV(z)$ has Euclidean norm $|F^{(q)}(z)|$. Summing the
squares over these $2^q$ choices gives
\begin{equation}
\label{eq:BS-complex-deriv}
|\nabla^qV(z)|=2^{q/2}\,|F^{(q)}(z)|\qquad(z\in U,\ q\ge0).
\end{equation}
The point of this discussion is that, for these holomorphic fields, complex
derivative bounds and ordinary Cartesian derivative bounds are equivalent,
with constants depending only on the derivative order $q$.

Lemmas~\ref{lem:hollow} and~\ref{lem:ext-multipole} bound the velocity induced by a rotationally symmetric
vorticity away from its support. Both assume that $g$ is real-valued and are
applied only to real-valued vorticities or vorticity perturbations.
In Lemma~\ref{lem:hollow}, $\bar g$ denotes the angular average of $g$, not
its complex conjugate, which would equal $g$; bars on complex coordinates and
moments continue to denote complex conjugation.

Lemma~\ref{lem:hollow} will control the
velocity that the outer rings induce on an inner ring: the angular average
$\bar g$ produces no velocity in the central hole, so only the nonradial
part $g-\bar g$ contributes there, and rotational symmetry makes its
velocity small near the origin. Lemma~\ref{lem:ext-multipole} gives the
corresponding exterior estimate; Section~\ref{sec:several} applies both.

\begin{lemma}
\label{lem:hollow}
Let $m_0\ge2$ be an integer and $R>0$. Let $g\in (L^1\cap L^\infty)(\R^2)$ be
real-valued and $m_0$-fold rotationally symmetric with $\supp g\subset\{|z|\ge R\}$, and let
$\bar g(z)=\frac1{2\pi}\int_0^{2\pi}g(\Orot_\phi z)\dd\phi$ be its angular
average. Then $\KBS\ast\bar g=0$ on $B_R$, and for $|z|\le R/2$,
\[
|\KBS\ast g(z)|\le\tfrac4{3\pi}\|g\|_{L^1}|z|^{m_0-1}R^{-m_0}.
\]
\end{lemma}

\begin{proof}
Let $|z|<R$, so that $|z-\zeta|\ge R-|z|>0$ for $\zeta\in\supp g$, and let $u=\KBS\ast g$. For $|\zeta|\ge R$ the series
$(z-\zeta)^{-1}=-\zeta^{-1}\sum_{l\ge0}(z/\zeta)^l$ converges uniformly in $\zeta$,
so Equation~\eqref{eq:BS-complex} gives, term by term,
$u_1-iu_2=\frac{i}{2\pi}\sum_{l\ge0}c_l(g)z^l$ with
$c_l(g)=\int g(\zeta)\zeta^{-l-1}\dd A(\zeta)$. The substitution $\zeta\mapsto e^{2\pi i/m_0}\zeta$ and the $m_0$-fold rotational
symmetry of $g$ show $c_l(g)=e^{-2\pi
i(l+1)/m_0}c_l(g)$, so $c_l(g)=0$ unless $m_0\mid(l+1)$. The angular average
$\bar g$ is radial and, like $g$, supported in $\{|z|\ge R\}$, so the same
expansion applies to it, and integration in the angle
shows $c_l(\bar g)=0$ for $l\ge0$; this proves $\KBS\ast\bar g=0$ on $B_R$. It
also gives $c_l(g)=c_l(g-\bar g)$, so $|c_l(g)|\le\|g-\bar g\|_{L^1}R^{-l-1}$,
as $|\zeta|\ge R$ on the support of $g-\bar g$.
Since $g$ is real, $|u|=|u_1-iu_2|$. Only the terms with $l=(j+1)m_0-1$, $j\ge0$,
remain, so for $|z|\le R/2$ these bounds and the sum of a geometric series
give
\[
\begin{aligned}
|u(z)|&\le\frac{\|g-\bar g\|_{L^1}}{2\pi}\sum_{j\ge0}|z|^{(j+1)m_0-1}R^{-(j+1)m_0}\\
&=\frac{\|g-\bar g\|_{L^1}}{2\pi}|z|^{m_0-1}R^{-m_0}
\sum_{j\ge0}\Big(\frac{|z|}{R}\Big)^{jm_0}\\
&\le\frac{\|g-\bar g\|_{L^1}}{2\pi(1-2^{-m_0})}|z|^{m_0-1}R^{-m_0},
\end{aligned}
\]
and $1/(1-2^{-m_0})\le\frac43$ because $m_0\ge2$. Finally
$\|g-\bar g\|_{L^1}\le2\|g\|_{L^1}$, because averaging over rotations does
not increase the $L^1$ norm.
\end{proof}

For an integer $m\ge1$ and $g$ as in Definition~I.\ref{I-def:multipole}, let
$\pout_m[g]$ and $\Qm[g]$ be the multipole moments defined there.
The following lemma gives the angular modes of the stream function of the
inner rings outside the disk containing their support, where the outer rings
are located; Section~\ref{sec:kick} uses this formula.

\begin{lemma}
\label{lem:exterior-mode}
Let $m\ge1$ be an integer and $\varrho>0$. Let $g\in L^1(\R^2)$ be
real-valued with compact support in $\overline B_{\varrho}\setminus\{0\}$,
and put $\psi=(2\pi)^{-1}\log|\cdot|\ast g$. Then the angular mode $m$ of
$\psi$ is
\[
\psi_m(r)=-\frac1{2m}\pout_m[g]\,r^{-m}\qquad(r>\varrho).
\]
\end{lemma}

\begin{proof}
For complex numbers $z,w$ with $|w|<|z|$, we have
\[
\log|z-w|=\log|z|-\Re\sum_{k\ge1}\frac1k\Big(\frac wz\Big)^k
.
\]
For fixed $|z|>\varrho$, the series converges uniformly in
$w\in\overline B_\varrho$, so it may be integrated term by term against $g$.
Writing $z=re^{i\theta}$, the angular mode $m$ of
$\Re((w/z)^k)$ vanishes unless $k=m$, when it equals
$\tfrac12\bar w^{\,m}r^{-m}$. Hence
\[
\psi_m(r)=-\frac{r^{-m}}{4\pi m}\int g(w)\bar w^{\,m}\dd A(w)
=-\frac1{2m}\pout_m[g]\,r^{-m}.
\qedhere
\]
\end{proof}

In Lemma~\ref{lem:ext-multipole} below, we estimate the exterior velocity after
subtracting its rotationally invariant circulation field, or \emph{monopole},
$\frac{\Gamma}{2\pi}\frac{(-\Im z,\Re z)}{|z|^2}$, where
$\Gamma=\int g$ is the total circulation, or zeroth vorticity moment.
This is the velocity of a point vortex at the origin with circulation
$\Gamma$. Subtracting it isolates the nonradial velocity, which decays faster
with distance and transmits the instability to the outer rings.

\begin{lemma}
\label{lem:ext-multipole}
Let $m\ge1$ and $q\ge0$ be integers, and let $\varrho>0$. Let $g\in L^1(\R^2)$ be
real-valued and $m$-fold rotationally symmetric with $\supp g\subset\overline{B}_\varrho$, and let $E$ be the
velocity field $\KBS\ast g$ minus its monopole
$\frac{1}{2\pi}\big(\int g\big)\frac{(-\Im z,\Re z)}{|z|^2}$. Then for $|z|\ge2\varrho$,
\[
|\nabla^qE(z)|\le C_{m,q}\,\|g\|_{L^1}\,\varrho^m|z|^{-m-1-q},\qquad
C_{m,q}=\frac{2^q}{2\pi}\sup_{0<\rho\le1/2}\rho^{-m}\sum_{l\in m\N}(l+q)^q\rho^{l}<\infty .
\]
\end{lemma}

\begin{proof}
Let $|z|>\varrho$. For $|\zeta|\le\varrho$ the series
$(z-\zeta)^{-1}=\sum_{l\ge0}\zeta^lz^{-l-1}$ converges uniformly in $\zeta$, so
Equation~\eqref{eq:BS-complex} gives
$(\KBS\ast g)_1(z)-i(\KBS\ast g)_2(z)=\frac1{2\pi i}\sum_{l\ge0}Q_lz^{-l-1}$ with
$Q_l=\int g(\zeta)\zeta^l\dd A(\zeta)$. The term $l=0$, which is $(\int g)/(2\pi iz)$, is the
monopole, and $Q_l=0$ unless $m\mid l$, because the substitution $\zeta\mapsto
e^{2\pi i/m}\zeta$ and the $m$-fold rotational symmetry of $g$ give $Q_l=e^{2\pi il/m}Q_l$. Hence the complex form of $E$ is
$F(z)=\frac1{2\pi i}\sum_{l\in m\N}Q_lz^{-l-1}$, where
$|Q_l|\le\|g\|_{L^1}\varrho^l$ because $|\zeta|\le\varrho$ on $\supp g$. The series converges locally uniformly on
$\{|z|>\varrho\}$, so $F$ is holomorphic there, may be differentiated term by
term, and, as $E$ is real,
Equation~\eqref{eq:BS-complex-deriv} gives
$|\nabla^qE(z)|=2^{q/2}|F^{(q)}(z)|\le2^q|F^{(q)}(z)|$. With
$|\partial_z^qz^{-l-1}|=(l+1)\cdots(l+q)|z|^{-l-1-q}\le(l+q)^q|z|^{-l-1-q}$ and
$\rho=\varrho/|z|\le\frac12$,
\[
|F^{(q)}(z)|\le\frac{\|g\|_{L^1}}{2\pi}|z|^{-1-q}\sum_{l\in m\N}(l+q)^q\rho^l
\le\frac{C_{m,q}}{2^q}\|g\|_{L^1}\rho^m|z|^{-1-q},
\]
and $\rho^m|z|^{-1-q}=\varrho^m|z|^{-m-1-q}$. The constant is finite, because
$\rho^{-m}\sum_{l\in m\N}(l+q)^q\rho^l=\sum_{j\ge1}(mj+q)^q\rho^{m(j-1)}$ is
nondecreasing in $\rho$ and finite at $\rho=\frac12$.
\end{proof}

The approximate solutions $\omega^{(n)}$ of Definition~\ref{def:approx} are the
global smooth solutions with the smooth, compactly supported data $f^{(n)}$,
given by the following classical theorem, which is proved in
Appendix~\ref{app:smooth}.

\begin{theorem}[cf. \cite{Wolibner}, \cite{Holder}, \cite{Kato1967}]
\label{thm:smooth}
Let $\omega_0\in C_c^\infty(\R^2)$ be real with $\supp\omega_0\subset B_{R_0}$.
There is a unique $\omega\in C^\infty(\R\times\R^2)$ such that:
\begin{enumerate}[label=\textup{(\arabic*)}]
\item $\supp\omega(t)\subset B_{R(t)}$ for every $t\in\R$, with
  $R(t)=R_0+(2\|\omega_0\|_{L^\infty}\|\omega_0\|_{L^1}/\pi)^{1/2}|t|$;
\item $v:=\KBS\ast\omega\in C^\infty(\R\times\R^2)$, $v$ and every derivative of
  $v$ are bounded on $[-T,T]\times\R^2$ for every $T>0$, and $\divop v=0$;
\item $\partial_t\omega+v\cdot\nabla\omega=0$ on $\R\times\R^2$ and $\omega(0)=\omega_0$.
\end{enumerate}
Moreover:
\begin{enumerate}[label=\textup{(\arabic*)},resume]
\item the flow $\flow$ of $v$, $\partial_t\flow(t,y)=v(t,\flow(t,y))$,
  $\flow(0,y)=y$, is defined for all $t\in\R$; each $\flow_t=\flow(t,\cdot)$ is
  a $C^\infty$ diffeomorphism of $\R^2$ with $\det D\flow_t\equiv1$; the maps
  $\flow$ and $(t,x)\mapsto\flow_t^{-1}(x)$ lie in $C^\infty(\R\times\R^2)$; and
  $\omega(t)=\omega_0\circ\flow_t^{-1}$; consequently
  $\int F(\omega(t))\dd x=\int F(\omega_0)\dd x$ for every continuous $F$ with
  $F(0)=0$;
\item if $I\ni0$ is an interval and $\tilde\omega\in C^1(I\times\R^2)$ satisfies
  $\tilde\omega(0)=\omega_0$, $\partial_t\tilde\omega+(\KBS\ast\tilde\omega)\cdot\nabla\tilde\omega=0$
  on $I\times\R^2$, and, for each compact subinterval $I'\subset I$,
  $\supp\tilde\omega(t)\subset B_{\tilde R}$ for $t\in I'$ with $\tilde R$
  depending on $I'$, then $\tilde\omega=\omega$ on $I\times\R^2$;
\item if $\omega_0\circ\Orot_\theta=\omega_0$ for a rotation $\Orot_\theta\in SO(2)$,
  then $\omega(t)\circ\Orot_\theta=\omega(t)$ for every $t\in\R$;
\item for every $T>0$ and every $q\in(4/3,\infty)$, the restriction of $\omega$
  to $[0,T)$ lies in $L^\infty_t(L^1\cap L^\infty)_x$ and is a renormalized
  weak solution in $X_q(T)$ with datum $\omega_0$, in the sense of
  Definition~I.\ref{I-def:weak}.
\end{enumerate}
\end{theorem}

The energy estimates of Section~\ref{sec:single}, the construction of
Appendix~\ref{app:smooth} and the proof of Lemma~\ref{lem:Hk-calculus} below use
the following Sobolev inequalities.

\begin{lemma}
\label{lem:embed}
The following embeddings hold on $\R^2$, with explicit constants.
\begin{enumerate}[label=\textup{(E\arabic*)}]
\item For $g\in H^1$, $\|g\|_{L^4}\le\|g\|_{H^1}$.
\item For $g\in H^2$, $g$ has a bounded continuous
  representative, and $\|g\|_{L^\infty}\le\|g\|_{H^2}$.
\item For an integer $k\ge0$ and $g\in H^{2k}$,
  $\|g\|_{H^k}^2\le3^k\|g\|_{L^2}\|g\|_{H^{2k}}$.
\end{enumerate}
\end{lemma}

\begin{proof}
(E1) Let first $g\in C_c^\infty(\R^2)$. Since
$|g(x)|^2=\int_{-\infty}^{x_1}\partial_1|g|^2(t,x_2)\dd t$, we have
\[
|g(x)|^2\le\int_\R|\partial_1|g|^2|(t,x_2)\dd t=:a(x_2),
\]
and likewise
$|g(x)|^2\le\int_\R|\partial_2|g|^2|(x_1,t)\dd t=:b(x_1)$. The product of
these two bounds, integrated over $x$, gives
$\int|g|^4\le\int a\int b$. Since $|\partial_j|g|^2|\le2|g||\partial_jg|$,
the Cauchy-Schwarz inequality, and then the inequality between the arithmetic
and geometric means, give
$\int|g|^4\le4\|g\|_{L^2}^2\allowbreak\|\partial_1g\|_{L^2}\|\partial_2g\|_{L^2}\le2\|g\|_{L^2}^2\|\nabla g\|_{L^2}^2$.
By density of $C_c^\infty$ in $H^1$ the inequality holds on $H^1$. Then, by the
same inequality between the means,
$\|g\|_{L^4}\le2^{1/4}\|g\|_{L^2}^{1/2}\|\nabla g\|_{L^2}^{1/2}\le2^{-1/4}\|g\|_{H^1}$.

(E2) Since $\Delta g=\partial_1^2g+\partial_2^2g$, the triangle inequality
and then the arithmetic mean--quadratic mean (AM--QM) inequality give
\[
\begin{aligned}
\|g-\Delta g\|_{L^2}
&\le\|g\|_{L^2}+\|\partial_1^2g\|_{L^2}+\|\partial_2^2g\|_{L^2}\\
&\le\sqrt3\big(\|g\|_{L^2}^2+\|\partial_1^2g\|_{L^2}^2
+\|\partial_2^2g\|_{L^2}^2\big)^{1/2}
\le\sqrt3\,\|g\|_{H^2},
\end{aligned}
\]
since the three squared norms are terms in $\|g\|_{H^2}^2$.
Our Fourier normalization gives $\|\hat h\|_{L^2}=2\pi\|h\|_{L^2}$ by
Plancherel's theorem. Thus, using $\widehat{\Delta g}=-|\zeta|^2\hat g$,
\[
\|(1+|\zeta|^2)\hat g\|_{L^2}=2\pi\|g-\Delta g\|_{L^2}
\le2\sqrt3\,\pi\|g\|_{H^2}.
\]
As $\int_{\R^2}(1+|\zeta|^2)^{-2}\dd\zeta=\pi$, Cauchy-Schwarz gives
$\|\hat g\|_{L^1}\le\sqrt\pi\,\|(1+|\zeta|^2)\hat g\|_{L^2}$. Thus $\hat g\in L^1$, the
Fourier inversion formula gives a bounded continuous representative, and
$\|g\|_{L^\infty}\le(2\pi)^{-2}\|\hat g\|_{L^1}\le\frac{\sqrt3}{2\sqrt\pi}\|g\|_{H^2}\le\|g\|_{H^2}$.

(E3) Put $P_k(\zeta)=\sum_{|\gamma|\le k}\zeta^{2\gamma}$. Every monomial
$\zeta^{2\gamma}$ with $|\gamma|\le k$ occurs in the expansion of
$(1+\zeta_1^2+\zeta_2^2)^k$ exactly once, with a multinomial coefficient between
$1$ and $3^k$; hence $P_k\le(1+|\zeta|^2)^k\le3^kP_k$. By Plancherel's theorem and
Cauchy-Schwarz,
\[
\|g\|_{H^k}^2\le\frac{\int(1+|\zeta|^2)^k|\hat g|^2\dd\zeta}{(2\pi)^2}
\le\frac{\|\hat g\|_{L^2}\|(1+|\zeta|^2)^{k}\hat g\|_{L^2}}{(2\pi)^2}
\le3^k\|g\|_{L^2}\|g\|_{H^{2k}} .\qedhere
\]
\end{proof}

In Section~\ref{sec:single} the perturbation of a ring is supported in a fixed
annulus, and the field that transports it is controlled only on an open
neighborhood of that annulus. The following calculus in $H^k$ is stated for
such data: $w$ is supported in a compact set $K$, and the field $U$ need only be
given on an open set $\mathcal O\supset K$.

\begin{lemma}
\label{lem:Hk-calculus}
Let $k\ge3$ be an integer, $K\subset\R^2$ compact, $\mathcal O\supset K$ open and
bounded, and $\chi_K\in C_c^\infty(\mathcal O)$ with $\chi_K=1$ on a neighborhood of
$K$. There is $C=C(k,\chi_K)$ such that:
\begin{enumerate}[label=\textup{(\arabic*)}]
\item for $w\in H^{k+1}(\R^2)$ with $\supp w\subset K$, a vector field
  $U\in W^{1,\infty}(\mathcal O)\cap H^k(\mathcal O)$, and every multi-index $\alpha$ with
  $|\alpha|\le k$,
  \[
  \|[\partial^\alpha,U\cdot\nabla]w\|_{L^2}\le C\big(\|U\|_{W^{1,\infty}(\mathcal O)}+\|U\|_{H^k(\mathcal O)}\big)\|w\|_{H^k};
  \]
  if moreover $U$ is real and $\divop U=0$ on $\mathcal O$, then
  $|\Re\langle\partial^\alpha(U\cdot\nabla w),\partial^\alpha w\rangle_{L^2}|\le
  C\big(\|U\|_{W^{1,\infty}(\mathcal O)}+\|U\|_{H^k(\mathcal O)}\big)\|w\|_{H^k}^2$, and for
  real $w$ the inner product is real;
\item $\|gh\|_{H^j}\le C\|g\|_{H^j}\|h\|_{H^j}$ for every integer $j$ with
  $2\le j\le k$ and all $g,h\in H^j(\R^2)$;
\item for every integer $j\ge0$, every $\rho>0$ and every $\chi\in C_c^\infty(B_\rho)$
  there is $C'=C'(j,\chi,\rho)$ such that
  $\|\chi\,\KBS\ast w\|_{H^{j+1}}\le C'\|w\|_{H^j}$ for every $w\in H^j(\R^2)$ with
  $\supp w\subset B_\rho$.
\end{enumerate}
Here $[\partial^\alpha,U\cdot\nabla]w=\partial^\alpha(U\cdot\nabla w)-U\cdot\nabla\partial^\alpha w$,
with both terms extended by $0$ outside $\mathcal O$.
\end{lemma}

\begin{proof}
We prove (2) first. Let $g,h\in C_c^\infty(\R^2)$, $|\gamma|\le j$ and
$\beta\le\gamma$. If $\beta=0$ or $\beta=\gamma$, then
$\|\partial^\beta g\,\partial^{\gamma-\beta}h\|_{L^2}\le\|g\|_{H^j}\|h\|_{H^j}$ by
Lemma~\ref{lem:embed}\textup{(E2)}, since $j\ge2$. If
$1\le|\beta|\le|\gamma|-1$, then by H\"older's inequality and
Lemma~\ref{lem:embed}\textup{(E1)},
$\|\partial^\beta g\,\partial^{\gamma-\beta}h\|_{L^2}\le\|\partial^\beta g\|_{H^1}\|\partial^{\gamma-\beta}h\|_{H^1}
\le\|g\|_{H^{|\beta|+1}}\|h\|_{H^{|\gamma|-|\beta|+1}}\le\|g\|_{H^j}\|h\|_{H^j}$.
The Leibniz formula
$\partial^\gamma(gh)=\sum_{\beta\le\gamma}\binom\gamma\beta\partial^\beta g\,\partial^{\gamma-\beta}h$
and $\sum_{\beta\le\gamma}\binom\gamma\beta=2^{|\gamma|}$ give (2) for smooth $g,h$,
with $C=C(k)$. For $g,h\in H^j$ take $g_n,h_n\in C_c^\infty$ converging to
$g,h$ in $H^j$; by the bound just proved, $g_nh_n$ is Cauchy in $H^j$, and by
(E2) it converges to $gh$ in $L^2$, so $gh\in H^j$ and (2) holds.

We now prove (1). Put $\tilde U=\chi_KU$, extended by $0$ outside $\mathcal O$; by the
Leibniz rule, $\|\tilde U\|_{W^{1,\infty}(\R^2)}\le C\|U\|_{W^{1,\infty}(\mathcal O)}$ and
$\|\tilde U\|_{H^k(\R^2)}\le C\|U\|_{H^k(\mathcal O)}$, with $C=C(k,\chi_K)$. Since
$\chi_K=1$ on a neighborhood $V$ of $K$ and $w=0$ outside $K$,
$[\partial^\alpha,U\cdot\nabla]w=[\partial^\alpha,\tilde U\cdot\nabla]w$, and
\[
[\partial^\alpha,\tilde U\cdot\nabla]w=\sum_{0<\beta\le\alpha}\binom\alpha\beta\,
\partial^\beta\tilde U\cdot\nabla\partial^{\alpha-\beta}w :
\]
with $\tilde U_n\in C_c^\infty$ in place of
$\tilde U$, this identity holds in $\mathcal D'$ by the Leibniz formula. Since
$\tilde U\in H^k(\R^2)$ by the bound above, there are such $\tilde U_n$ with
$\tilde U_n\to\tilde U$ in $H^k$. As $w\in H^{k+1}$, we have
$\nabla\partial^{\alpha-\beta}w\in L^2$ for $\beta\le\alpha$; so each product in
the identity, as
well as $\tilde U_n\cdot\nabla w$ and $\tilde U_n\cdot\nabla\partial^\alpha w$,
converges in $L^1$, and the identity passes to the limit. We bound the terms on
the right according to $|\beta|$. Here $l\in\{1,2\}$ indexes the Cartesian
components $U_l$ and $\tilde U_l$, whereas $n$ above indexes the approximating
sequence $\tilde U_n$. Write each term as
$\sum_l\partial^\beta\tilde U_l\,\partial_l\partial^{\alpha-\beta}w$. If $|\beta|=1$,
its $L^2$ norm is at most $2\|\nabla\tilde U\|_{L^\infty}\|w\|_{H^k}$. If
$|\beta|=|\alpha|\ge2$, it is at most
$\sum_l\|\partial^\alpha\tilde U_l\|_{L^2}\|\partial_lw\|_{L^\infty}\le2\|\tilde U\|_{H^k}\|w\|_{H^3}
\le2\|\tilde U\|_{H^k}\|w\|_{H^k}$, by (E2) and $k\ge3$. If
$2\le|\beta|\le|\alpha|-1$, then by H\"older's inequality and (E1) it is at most
$\sum_l\|\partial^\beta\tilde U_l\|_{H^1}\|\partial_l\partial^{\alpha-\beta}w\|_{H^1}
\le2\|\tilde U\|_{H^{|\beta|+1}}\|w\|_{H^{|\alpha|-|\beta|+2}}\le2\|\tilde U\|_{H^k}\|w\|_{H^k}$,
because $|\beta|+1\le k$ and $|\alpha|-|\beta|+2\le k$. Summing over $\beta$ gives
the commutator bound. For the second claim, put $g=\partial^\alpha w\in H^1$,
which vanishes outside $K$. Then
$\partial^\alpha(U\cdot\nabla w)=[\partial^\alpha,\tilde U\cdot\nabla]w+\tilde U\cdot\nabla g$
and
$\Re\langle\tilde U\cdot\nabla g,g\rangle_{L^2}=\frac12\int\tilde U\cdot\nabla|g|^2\dd x
=-\frac12\int(\divop\tilde U)|g|^2\dd x=0$:
the integration by parts holds because $|g|^2\in W^{1,1}$ and $\tilde U$ is
Lipschitz with compact support, and the last integral vanishes because
$\divop\tilde U=\chi_K\divop U+U\cdot\nabla\chi_K=0$ on $V$ and $g=0$ outside $K$.
Hence the real part is at most
$\|[\partial^\alpha,\tilde U\cdot\nabla]w\|_{L^2}\|g\|_{L^2}$, by the
Cauchy-Schwarz inequality, and the commutator bound gives the claim, because
the $L^2$ norm of $g$ is at most the $H^k$ norm of $w$.

We now prove (3). Let first $w\in C_c^\infty(B_\rho)$ and $u=\KBS\ast w$. By
Lemma~I.\ref{I-lem:BS-basic}(d), and (c) with $R=\rho$,
$\|\nabla u\|_{H^j}=\|w\|_{H^j}$ and $\|u\|_{L^2(B_\rho)}\le2\rho\|w\|_{L^2}$.
For $|\gamma|\le j+1$, in the Leibniz formula
$\partial^\gamma(\chi u)=\sum_{\beta\le\gamma}\binom\gamma\beta\partial^{\gamma-\beta}\chi\,\partial^\beta u$
the term $\beta=0$ has $L^2$ norm at most
$\|\partial^\gamma\chi\|_{L^\infty}\|u\|_{L^2(B_\rho)}$, because $\chi$ vanishes
outside $B_\rho$, and every other term at most
$\|\partial^{\gamma-\beta}\chi\|_{L^\infty}\|\nabla u\|_{H^j}$. This gives (3) for
smooth $w$. For general $w$, $\supp w$ is a compact subset of $B_\rho$, so
mollification gives $w_n\in C_c^\infty(B_\rho)$ with $w_n\to w$ in $H^j$; then
$\chi\KBS\ast w_n$ is Cauchy in $H^{j+1}$ and, by Lemma~I.\ref{I-lem:BS-basic}(c),
converges to $\chi\KBS\ast w$ in $L^2$, so (3) holds for $w$.
\end{proof}

\section{Nonlinear Estimates for One Ring}
\label{sec:single}

In this section we study, in ring units, one ring of the approximate solutions: a
ring close to $\Phi$, whose deviation from $\Phi$ is driven by its own
nonlinearity and by the field $h$ of the vorticity inside and outside it, apart
from the point-vortex part of the inner vorticity, which enters the linearization
$\Lin(G)$ (Lemma~\ref{lem:reduction}). We prove estimates for a quiet ring
(Theorem~\ref{thm:quiet}) and an active ring (Theorem~\ref{thm:active}), and
construct the seeded innermost ring by rearranging $\Phi$
(Lemma~\ref{lem:seed}).
Sections~\ref{sec:relay} and~\ref{sec:bootstrap} apply Theorems~\ref{thm:quiet}
and~\ref{thm:active} to every ring in its own units.

Throughout this section we assume \textup{(N1)}--\textup{(N3)}, we work in ring
units, and we use the objects of Definition~I.\ref{I-def:eigendata}: the unstable
eigenvalue $\lamst(G)$ of $\Lin(G)$ in the mode $3$ and its eigenfunction
$\phist(G)$, normalized in $L^2(\R^2)$; the functional $\ellst(G)$ that gives the
coefficient of a function along $\phist(G)$; the Riesz projection $\Riesz(G)$
onto the span of $\phist(G)$ and $\overline{\phist(G)}$; and the growth rates
$\sigma(G)=\Re\lamst(G)$ and $\sigb=\min_{G\in[0,\Gone]}\sigma(G)$. When $G$ is fixed
we write $\lamst=\lamst(G)$, $\phist=\phist(G)$, $\ellst=\ellst(G)$,
$\Riesz=\Riesz(G)$. The proofs of
Theorems~\ref{thm:quiet} and~\ref{thm:active} follow Grenier \cite{Grenier}: an
approximate solution of high order is built from the linear semigroup, and the order is so large that a crude energy estimate for the
remainder suffices.

We use the norms of $H^k$, $W^{1,\infty}$ and $C^k(N')$ fixed in
Section~I.\ref{I-sec:prelim}, and the embedding $\|f\|_{L^\infty}\le\|f\|_{H^2}$ of
Lemma~\ref{lem:embed}\textup{(E2)}. A vector field $h$ on a neighborhood of $N'$
is \textit{$3$-fold rotationally symmetric} if $h(\Orot x)=\Orot h(x)$.

The equation for the deviation of a ring from $\Phi$
(Lemma~\ref{lem:reduction}\textup{(3)} below) assumes that the ring stays in
the annulus $N$, and the continuity arguments of this section keep it there by
bounding how far its material moves radially. We write $\omega_{\mathrm r}(0)$
for the initial vorticity of the selected ring, with $\mathrm r$ standing for
``ring''. Let $\flow$ be the flow
of a solution of Theorem~\ref{thm:smooth}, let $\omega_{\mathrm r}(0)$ have compact
support, let $\omega_{\mathrm r}(t)=\omega_{\mathrm r}(0)\circ\flow_t^{-1}$, and let
$\delta_{\mathrm r}(t)$ be the supremum of $\big||\flow_\tau(y)|-|y|\big|$ over
$\tau\in[0,t]$ and $y\in\supp\omega_{\mathrm r}(0)$ ($0$ if the support is empty);
it is a maximum if the support is nonempty, it is continuous and
nondecreasing, and $\delta_{\mathrm r}(0)=0$, because $(\tau,y)\mapsto\flow_\tau(y)$
is continuous. We say that \textit{every trajectory of $\omega_{\mathrm r}$ moves
radially by at most $\delta$ on $[0,t]$} if $\delta_{\mathrm r}(t)\le\delta$. Since
$\supp\omega_{\mathrm r}(\tau)=\flow_\tau(\supp\omega_{\mathrm r}(0))$, the support
of $\omega_{\mathrm r}(\tau)$, $\tau\le t$, then lies within radial distance $\delta$
of $\supp\omega_{\mathrm r}(0)$.

For $G\in[0,\Gone]$ and an integer $k\ge0$, $e^{t\Lin(G)}$ denotes the semigroup
of Lemma~I.\ref{I-lem:semigroup-construct} for the ring in ring units. It preserves
each mode, hence the subspace $X^k$ of the space $Y^k$ of that lemma
consisting of the functions whose mode components vanish outside $3\Z$. It maps
real-valued functions to real-valued functions, since it is the sum of a
Dyson-Phillips series built from the two operators of that lemma, which both do:
the transport group $e^{t\Lin^{\rm tr}}$, which rotates each circle $\{|y|=r\}$ by
the angle $\Omega(r;G)t$, and the coupling operator $\mathcal K$, here
$\mathcal Kw=-(\KBS\ast w)\cdot\nabla\Phi$. The approximate solutions of this
section are Duhamel integrals of this semigroup. The next lemma shows that they are
solutions of the inhomogeneous linear equation, differentiable in time as functions with
values in $X^{k-1}$ when the data lie in $X^k$.

\begin{lemma}
\label{lem:duhamel}
Let $G\in[0,\Gone]$, $T>0$, and let $k\ge1$ be an integer. Let $f_0\in X^k$ and
$g\in C([0,T];X^k)$, and put
$u(t)=e^{t\Lin(G)}f_0+\int_0^te^{(t-\tau)\Lin(G)}g(\tau)\dd\tau$. Then
$u\in C([0,T];X^k)\cap C^1([0,T];X^{k-1})$, and $\partial_tu=\Lin(G)u+g$ on
$[0,T]$. If $f_0$ and $g$ are real-valued, so is $u$.
\end{lemma}
\begin{proof}
We first prove the continuity of $u$ in $X^k$, and then its differentiability in
$X^{k-1}$, where the generator of the semigroup acts on $X^k$ as $\Lin(G)$.
The semigroup $S(t)=e^{t\Lin(G)}$ is strongly continuous on $X^k$ and, by the
uniform boundedness principle, bounded on $[0,T]$. So $\tau\mapsto S(t-\tau)g(\tau)$
is continuous into $X^k$: the argument given below for the continuity of
$(s,\tau)\mapsto S(s-\tau)H(\tau)$ in $X^{k-1}$ uses only the strong continuity of
$S$ and a bound on $\|S(t)\|$ for $t$ in a compact interval, and both hold in
$X^k$. Hence
$\int_0^tS(t-\tau)g(\tau)\dd\tau=\int_0^tS(\tau)g(t-\tau)\dd\tau$
is continuous in $t$ by dominated convergence, and $u\in C([0,T];X^k)$. Moreover
$X^k\hookrightarrow X^{k-1}$, and $S$ is a $C_0$-semigroup on $X^{k-1}$, since it
is one on $Y^{k-1}$ and $X^{k-1}=Y^{k-1}_{3\Z}$ is a closed invariant subspace,
both by Lemma~I.\ref{I-lem:semigroup-construct}\textup{(2)}. By
Lemma~I.\ref{I-lem:semigroup-construct}\textup{(2)}, the generator $A$ of $S$ on
$X^{k-1}$ (the part in $X^{k-1}$ of the generator on $Y^{k-1}$ of that lemma)
satisfies $X^k\subset D(A)$, $Ay=\Lin(G)y$ and $\|Ay\|_{X^{k-1}}\le C\|y\|_{X^k}$ for
$y\in X^k$; here $\Lin(G)$ maps $X^k$ into $X^{k-1}$ because it preserves modes
(Lemma~I.\ref{I-lem:sector}\textup{(2)}). Below, $\|\cdot\|$ is the norm of $X^{k-1}$ and
the operator norm on $X^{k-1}$.

All integrals below are Riemann integrals of continuous functions with values
in a Banach space; Fubini's theorem holds for them, because it holds after
applying any continuous linear functional. We extend $g$ by $g(t)=g(T)$ for
$T<t\le2T$. Since each $t\mapsto S(t)x$, $x\in X^{k-1}$, is continuous, hence bounded on $[0,2T]$, the
uniform boundedness principle gives $M_T:=\sup_{0\le t\le2T}\|S(t)\|<\infty$.

Let $x\in D(A)$. For $t\ge0$,
$\epsilon^{-1}(S(\epsilon)-I)S(t)x=S(t)\epsilon^{-1}(S(\epsilon)x-x)\to S(t)Ax$, so $S(t)x\in D(A)$,
$AS(t)x=S(t)Ax$, and $S(t)Ax$ is the right derivative of $t\mapsto S(t)x$. For
$0<\epsilon\le t\le2T$,
$\epsilon^{-1}(S(t)x-S(t-\epsilon)x)-S(t)Ax=S(t-\epsilon)\big(\epsilon^{-1}(S(\epsilon)x-x)-Ax\big)+\big(S(t-\epsilon)-S(t)\big)Ax\to0$,
because $\|S(t-\epsilon)\|\le M_T$. Thus $t\mapsto S(t)x$ lies in $C^1([0,2T];X^{k-1})$ with
derivative $S(t)Ax$, and
\begin{equation}
\label{eq:semigroup-FTC}
S(r)x-x=\int_0^rS(s)Ax\dd s\qquad(x\in D(A),\ 0\le r\le2T).
\end{equation}

Let $H\in C([0,2T];X^{k-1})$. The identity
$S(s'-\tau')H(\tau')-S(s-\tau)H(\tau)=S(s'-\tau')(H(\tau')-H(\tau))
+(S(s'-\tau')-S(s-\tau))H(\tau)$, whose first term has norm at most
$M_T\|H(\tau')-H(\tau)\|$ and whose second tends to $0$ by strong continuity,
shows that $(s,\tau)\mapsto S(s-\tau)H(\tau)$ is continuous on
$\{0\le\tau\le s\le2T\}$. Moreover
$s\mapsto\int_0^sS(s-\tau)H(\tau)\dd\tau$ is continuous on $[0,2T]$: for
$0\le s<s'\le2T$ the difference of its values is
$\int_s^{s'}S(s'-\tau)H(\tau)\dd\tau+\int_0^sS(s-\tau)(S(s'-s)-I)H(\tau)\dd\tau$,
of norm at most $M_T\sup\|H\|(s'-s)+2TM_T\sup_\tau\|(S(s'-s)-I)H(\tau)\|$; and the
operators $S(\epsilon)-I$, $0\le\epsilon\le2T$, are bounded by $M_T+1$ and tend strongly to
$0$ as $\epsilon\to0$, hence uniformly on the compact set $H([0,2T])$.

Now put $V(t)=\int_0^tS(t-\tau)g(\tau)\dd\tau$ and
$W(t)=\int_0^tS(t-\tau)Ag(\tau)\dd\tau$. Since
$\|A(g(\tau)-g(\tau'))\|\le C\|g(\tau)-g(\tau')\|_{X^k}$, $Ag$ is continuous
into $X^{k-1}$, so $W$ is continuous by the previous paragraph with $H=Ag$. By
Equation~\eqref{eq:semigroup-FTC} with $x=g(\tau)\in X^k\subset D(A)$ and Fubini's theorem,
\[
V(t)=\int_0^t\Big(g(\tau)+\int_\tau^tS(s-\tau)Ag(\tau)\dd s\Big)\dd\tau
=\int_0^t\big(g(s)+W(s)\big)\dd s ,
\]
so $V\in C^1([0,2T];X^{k-1})$ and $V'=g+W$. For $t\in[0,T]$ and $0<\epsilon\le T$ we have
$V(t+\epsilon)=S(\epsilon)V(t)+\int_t^{t+\epsilon}S(t+\epsilon-\tau)g(\tau)\dd\tau$, and $\epsilon^{-1}$ times the
last integral tends to $g(t)$ as $\epsilon\to0$, by the continuity of
$(s,\tau)\mapsto S(s-\tau)g(\tau)$ at $(t,t)$. Hence
$\epsilon^{-1}(S(\epsilon)V(t)-V(t))\to V'(t)-g(t)=W(t)$ in $X^{k-1}$.
By definition of the generator $A$, existence of this right difference-quotient
limit means that $V(t)\in D(A)$, and the limit is $AV(t)$; hence $AV(t)=W(t)$.
Now apply the semigroup differentiation argument immediately preceding
Equation~\eqref{eq:semigroup-FTC} to $x=f_0$, which is admissible since
$f_0\in X^k\subset D(A)$. It gives $S(t)f_0\in D(A)$ and
$\frac{\dd}{\dd t}S(t)f_0=AS(t)f_0=S(t)Af_0$.
Combining these facts, $u=S(\cdot)f_0+V$ lies in
$C^1([0,T];X^{k-1})$, $u(t)\in D(A)$ with $Au(t)=S(t)Af_0+W(t)$, and
$\partial_tu=S(\cdot)Af_0+g+W=Au+g=\Lin(G)u+g$, since $u(t)\in X^k$. Finally, $S(t)$
preserves real-valued functions.
\end{proof}

\emph{Remark.} The function $u(t)$ in Lemma~\ref{lem:duhamel} is the
\emph{first Duhamel iterate}, combining the free evolution of $f_0$ with the
response to the prescribed driving term $g$.

\subsection{The Perturbation Equation}

We start with the global smooth solution $\omega$ given by
Theorem~\ref{thm:smooth} for a real-valued, $3$-fold rotationally symmetric
initial vorticity $\omega_0\in C_c^\infty(\R^2)$. In the applications,
$\omega$ is a rescaled, and possibly time-translated, version of the smooth
approximate solution $\omega^{(n)}$ of Definition~\ref{def:approx}.
We split its initial datum as
$\omega_0=\omega_{\mathrm r}(0)+\omega_{\mathrm{inn}}(0)+\omega_{\mathrm{ou}}(0)$,
corresponding to the selected ring, the vorticity inside it and the vorticity
outside it, and transport all three components by the full flow.

Our principal analysis concerns $\omega_{\mathrm r}$, through its deviation
$w=\omega_{\mathrm r}-\Phi$ from the reference ring. The inner vorticity's
circulation $2\pi G$ enters the linearization $\Lin(G)$ as a point vortex at
the origin. The rest of its velocity field, together with the velocity field
of $\omega_{\mathrm{ou}}$, forms the field $h$ in the perturbation equation.
The later estimates use the separation between the rings to control $h$ as
a small perturbation on the selected ring; the inner and outer vorticities
themselves need not be small, and $h$ retains its role in driving the
instability.

In the approximate data, the initial ring supports are widely separated.
The components $\omega_\bullet(t)$ below follow the material transported
from these initial supports. Since the full flow is a diffeomorphism at every
time (Theorem~\ref{thm:smooth}\textup{(4)}), their transported supports remain
pairwise disjoint for all time, and the decomposition remains well-defined.
The quantitative radial separation required below is an additional
hypothesis on $[0,T]$, as is the condition that the selected ring stays in
$N$. The later confinement and bootstrap estimates verify these conditions
on the time windows needed for the relay and the non-uniqueness argument.

\begin{lemma}
\label{lem:reduction}
Let $\omega$ be a solution given by Theorem~\ref{thm:smooth} whose datum
$\omega(0)$ is $3$-fold rotationally symmetric, and let $\flow_t$ be its flow. Let
$\omega(0)=\omega_{\mathrm r}(0)+\omega_{\mathrm{inn}}(0)+\omega_{\mathrm{ou}}(0)$ with
pairwise disjoint supports, $\omega_{\mathrm r}(0)$, $\omega_{\mathrm{inn}}(0)$,
$\omega_{\mathrm{ou}}(0)$ each $3$-fold rotationally symmetric, and put
$\omega_\bullet(t)=\omega_\bullet(0)\circ\flow_t^{-1}$ for
$\bullet\in\{\mathrm r,\mathrm{inn},\mathrm{ou}\}$. Put
$G=\frac1{2\pi}\int\omega_{\mathrm{inn}}(0)$. Let $T>0$ and assume that for every
$t\in[0,T]$, $\supp\omega_{\mathrm{inn}}(t)\subset\overline B_{a/8}$ and
$\supp\omega_{\mathrm{ou}}(t)\cap B_6=\emptyset$. Put
\[
h(t)=\KBS\ast\big(\omega_{\mathrm{inn}}(t)+\omega_{\mathrm{ou}}(t)\big)-G\,\frac{x^\perp}{|x|^2}
\qquad\text{on a neighborhood of }N'.
\]
Then:
\begin{enumerate}[label=\textup{(\arabic*)}]
\item for every $t\in[0,T]$, $h(t)$ is smooth, divergence-free and
  $3$-fold rotationally symmetric on a neighborhood of $N'$; and $t\mapsto h(t)\in C^k(N')$ is
  continuous on $[0,T]$ for every $k$;
\item $\omega_{\mathrm r}\in C^\infty(\R\times\R^2)$, and on $[0,T]\times N'$ the
  velocity $v=\KBS\ast\omega$ is
  \[
  v=\Omega(|x|;G)\,x^\perp+h+\KBS\ast(\omega_{\mathrm r}-\Phi),
  \]
  with $\Omega$ of Definition~I.\ref{I-def:sector-setting} for $(\Phi,G)$; consequently,
  if $t\in[0,T]$ and $\flow_\tau(y)\in N'$ for every $\tau\in[0,t]$, then
  $|\flow_t(y)|-|y|=\int_0^t\big((h+\KBS\ast(\omega_{\mathrm r}-\Phi))\cdot\er\big)(\tau,\flow_\tau(y))\dd\tau$;
\item if $\supp\omega_{\mathrm r}(t)\subset N$ for $t\in[0,T]$, then the function
  $w=\omega_{\mathrm r}-\Phi$ is real-valued, it lies in $C^1([0,T];X^k)$ for
  every $k$, and
  \begin{equation}
  \label{eq:ring}
  \partial_tw=\Lin(G)w+B(w,w)-h\cdot\nabla w-h\cdot\nabla\Phi,\qquad B(v,w)=-(\KBS\ast v)\cdot\nabla w,
  \end{equation}
  on $[0,T]\times\R^2$, where $h\cdot\nabla w$ and $h\cdot\nabla\Phi$ denote the
  functions equal to these products on $N'$ and to $0$ elsewhere.
\end{enumerate}
\end{lemma}
\begin{proof}
Put $K_\bullet=\supp\omega_\bullet(0)$. These compact sets are pairwise
disjoint, so each $\omega_\bullet(0)$ agrees with $\omega(0)$ near $K_\bullet$ and
vanishes near the other two; hence $\omega_\bullet(0)\in C_c^\infty(\R^2)$. The
map $(t,x)\mapsto\flow_t^{-1}(x)$ is the value at time $0$ of the solution of
$\dot z=v(\tau,z)$ with $z(t)=x$. Since $v$ is smooth and bounded on
$[-T',T']\times\R^2$ for every $T'>0$ (Theorem~\ref{thm:smooth}\textup{(2)}),
Lemma~\ref{lem:flow}\textup{(a)}, \textup{(b)} (Appendix~\ref{app:smooth}) shows that this map
is defined and
$C^\infty$ on $\R\times\R^2$. Hence $\omega_\bullet\in C^\infty(\R\times\R^2)$,
$\supp\omega_\bullet(t)=\flow_t(K_\bullet)$, and $\int\omega_{\mathrm{inn}}(t)=2\pi G$
because $\flow_t$ preserves area. By Theorem~\ref{thm:smooth}\textup{(6)},
$\omega(t)$ is $3$-fold rotationally symmetric. Since $\KBS(\Orot x)=\Orot\KBS(x)$, we get
$(\KBS\ast g)(\Orot x)=\Orot(\KBS\ast g)(x)$ for every $3$-fold rotationally symmetric
$g\in (L^1\cap L^\infty)(\R^2)$; so $v(t,\Orot x)=\Orot v(t,x)$, uniqueness for the flow
gives $\flow_t(\Orot y)=\Orot\flow_t(y)$, and each $\omega_\bullet(t)$ is $3$-fold
rotationally symmetric.

\emph{Proof of (1).} By continuity of the flow on compact sets, and since
$\flow_t(K_{\mathrm{inn}})=\supp\omega_{\mathrm{inn}}(t)\subset\overline B_{a/8}$,
$\rho_{\mathrm{inn}}:=\max\{|\flow_t(y)|:t\in[0,T],\,y\in K_{\mathrm{inn}}\}\le a/8<a/4$
($\rho_{\mathrm{inn}}:=0$ if $K_{\mathrm{inn}}=\emptyset$), and $|\flow_t(y)|\ge6$ for
$t\in[0,T]$ and $y\in K_{\mathrm{ou}}$. So the open annulus
$U=\{\frac12(\rho_{\mathrm{inn}}+\frac a4)<|x|<4\}\supset N'$ has distance at least
$\min(\frac12(\frac a4-\rho_{\mathrm{inn}}),2)>0$ from
$\supp(\omega_{\mathrm{inn}}(t)+\omega_{\mathrm{ou}}(t))$ for $t\in[0,T]$. Since
$|\partial^\beta\KBS(z)|\le C_\beta|z|^{-1-|\beta|}$, and $\omega_{\mathrm{inn}}$,
$\omega_{\mathrm{ou}}$ are continuous with supports in a fixed compact set for
$t\in[0,T]$ (Theorem~\ref{thm:smooth}\textup{(1)}), we may differentiate under the
integral sign: $h(t)\in C^\infty(U)$, every $\partial^\beta_xh$ is continuous on
$[0,T]\times U$, and in particular $t\mapsto h(t)\in C^k(N')$ is continuous. The
fields $x\mapsto\KBS(x-y)=\nabla^\perp_x\big(\frac1{2\pi}\log|x-y|\big)$ and
$x^\perp/|x|^2$ are divergence-free away from $y$ and from $0$; hence
$h(t)$ is divergence-free on $U$, and it is $3$-fold rotationally symmetric by the rotation identities above and
$(\Orot x)^\perp=\Orot x^\perp$.

\emph{Proof of (2).} The radial Biot-Savart formula of
Lemma~I.\ref{I-lem:radial-stationary} applies to $\Phi$ and holds pointwise
away from $0$ by smoothness. Thus, by the definition of $h$, on $[0,T]\times U$,
\[
v=\KBS\ast\Phi+\KBS\ast(\omega_{\mathrm r}-\Phi)+\KBS\ast(\omega_{\mathrm{inn}}+\omega_{\mathrm{ou}})
=\Big(\frac{\Circ_\Phi(|x|)}{2\pi|x|^2}+\frac G{|x|^2}\Big)x^\perp+\KBS\ast(\omega_{\mathrm r}-\Phi)+h,
\]
and the bracket is $\Omega(|x|;G)$. Since $\partial_\tau\flow_\tau(y)=v(\tau,\flow_\tau(y))$
and $x^\perp\cdot\er=0$, we get $\frac{\dd}{\dd\tau}|\flow_\tau(y)|=
\big((h+\KBS\ast(\omega_{\mathrm r}-\Phi))\cdot\er\big)(\tau,\flow_\tau(y))$ while
$\flow_\tau(y)\in N'$, which gives the integral formula.

\emph{Proof of (3).} Since $\supp\Phi\subset\{a\le|y|\le1\}\subset N$
(Lemma~I.\ref{I-lem:ring-facts}), $w(t)$ is a smooth, real, $3$-fold rotationally symmetric function
supported in $N$. The derivatives $\partial^\alpha_xw$ and
$\partial_t\partial^\alpha_xw$ are uniformly continuous on the compact set
$[0,T]\times N$ and, for each $t\in[0,T]$, vanish for $x\notin N$.
Thus the difference quotients of
$t\mapsto w(t)$ converge in every $H^k$, and $w\in C^1([0,T];X^k)$.

The ring component is transported by the full velocity $v$. Differentiating
its defining identity along a trajectory gives
\[
0=\frac{\dd}{\dd t}\omega_{\mathrm r}(t,\flow_t(y))
=\big(\partial_t\omega_{\mathrm r}+v\cdot\nabla\omega_{\mathrm r}\big)(t,\flow_t(y)).
\]
Since $\flow_t$ is onto, $\partial_t\omega_{\mathrm r}
+v\cdot\nabla\omega_{\mathrm r}=0$ at every spatial point.
Now $\omega_{\mathrm r}=\Phi+w$ and $\Phi$ is independent of time, so
$\partial_t\omega_{\mathrm r}=\partial_tw$ and
$\nabla\omega_{\mathrm r}=\nabla\Phi+\nabla w$.
Substituting the velocity decomposition from (2), on $[0,T]\times N'$ we obtain
\[
\partial_tw=-\big(\Omega(|x|;G)x^\perp+\KBS\ast w+h\big)
\cdot(\nabla\Phi+\nabla w).
\]
Because $\Phi$ is radial,
$x^\perp\cdot\nabla\Phi=\partial_\theta\Phi=0$, whereas
$x^\perp\cdot\nabla w=\partial_\theta w$. Expanding the product therefore gives
\[
\begin{aligned}
\partial_tw
&=-\Omega(\cdot;G)\,\partial_\theta w-(\KBS\ast w)\cdot\nabla\Phi\\
&\quad-(\KBS\ast w)\cdot\nabla w-h\cdot\nabla w-h\cdot\nabla\Phi.
\end{aligned}
\]
The first line on the right is $\Lin(G)w$ by
Definition~I.\ref{I-def:sector-setting}, with background $\Phi$.
The self-interaction term $-(\KBS\ast w)\cdot\nabla w$ is $B(w,w)$, and
the remaining two terms describe advection of $w$ and $\Phi$ by $h$.
This is Equation~\eqref{eq:ring} on $N'$. The spatial support properties
above extend it to $\R^2$, with the products involving $h$ extended by zero
as in the statement.
\end{proof}

By Theorem~\ref{thm:smooth}\textup{(1)}, $\omega_{\mathrm r}$ is smooth with
support in a fixed compact set on $[0,T]$, so
$t\mapsto\|\omega_{\mathrm r}(t)-\Phi\|_{H^k}$ is continuous for every $k$.
We fix the cutoffs used in the estimates below. The functions we estimate are
supported in $N$, and the velocities that act on them are needed only on a
neighborhood of $N$. We therefore multiply the Biot-Savart velocity by a compactly
supported cutoff $\chi_B$ that does not change it on $N'$, so that
Lemma~\ref{lem:Hk-calculus}\textup{(3)} applies, and we use a cutoff $\chi_N$ around $N$ in Lemma~\ref{lem:Hk-calculus}\textup{(1)}
and \textup{(2)}. Let $N_1=\{\frac a4<|y|<3\}$ be
the interior of $N'$. Fix $\chi_B\in C_c^\infty(B_4)$ with $\chi_B=1$ on $B_{7/2}$,
and $\chi_N\in C_c^\infty(N_1)$ with $\chi_N=1$ on a neighborhood of $N$.
The constants in the applications of Lemma~\ref{lem:Hk-calculus} below
depend only on the derivative order and $a$.

In the analytic estimates that follow, $A\lesssim B$ means $A\le CB$ for
an auxiliary constant whose dependencies are specified in the argument.
Constants are uniform in $G\in[0,\Gone]$, the number of rings, their indices
and $A_1$, unless stated otherwise. Named constants are fixed before the
parameters whose choice depends on them.

\begin{lemma}
\label{lem:bilinear}
Let $k\ge2$ be an integer. There exist a constant $C_{\mathrm{bil}}(k)$,
depending only on $k$ and $a$, and an absolute constant
$C_{\mathrm{BS}}$, such that for $v,w\in X^{k+1}$ and every vector field
$h\in C^k(N';\R^2)$:
\begin{enumerate}[label=\textup{(\arabic*)}]
\item $\|B(v,w)\|_{H^k}\le C_{\mathrm{bil}}(k)\|v\|_{H^k}\|w\|_{H^{k+1}}$;
\item $\|h\cdot\nabla w\|_{H^k}\le C_{\mathrm{bil}}(k)\|h\|_{C^k(N')}\|w\|_{H^{k+1}}$;
\item $\|h\cdot\nabla\Phi\|_{H^k}\le\Cnab\|h\|_{C^k(N')}$, with $\Cnab=\Cnab(k)$ of
  Lemma~I.\ref{I-lem:uniform};
\item $\sup_{|x|\le3}|(\KBS\ast v)(x)|\le C_{\mathrm{BS}}\|v\|_{H^1}$.
\end{enumerate}
The functions in \textup{(1)}--\textup{(3)} are supported in $N$, and they lie in
$X^k$ if $h$ is $3$-fold rotationally symmetric.
\end{lemma}
\begin{proof}
Since $\nabla w$ is supported in $N$,
$B(v,w)=-\sum_{i=1,2}(\chi_B\KBS\ast v)_i\,\partial_iw$. By
Lemma~\ref{lem:Hk-calculus}\textup{(2)} (with $\max(k,3)$ in place of its $k$, and
$j=k$) and \textup{(3)},
\[
\|B(v,w)\|_{H^k}\lesssim
\|\chi_B\KBS\ast v\|_{H^k}\|w\|_{H^{k+1}}
\lesssim\|v\|_{H^{k-1}}\|w\|_{H^{k+1}},
\]
with constants depending only on $k$ and $a$, which gives (1).
For (2), write $h\cdot\nabla w=(\chi_Nh)\cdot\nabla w$ and apply the
Leibniz rule. Each term is bounded by a constant depending only on $k$ and
$\chi_N$ times $\|h\|_{C^k(N')}\|w\|_{H^{k+1}}$.
Item (3) is the
last inequality of Lemma~I.\ref{I-lem:uniform}. For (4), $\KBS\ast v=\chi_B\KBS\ast v$ on
$B_{7/2}$, and Lemma~\ref{lem:embed}\textup{(E2)}, applied to each component, and
Lemma~\ref{lem:Hk-calculus}\textup{(3)} give
$\|\chi_B\KBS\ast v\|_{L^\infty}\lesssim\|v\|_{H^1}$ with an absolute constant,
which gives (4), since $\{|x|\le3\}\subset B_{7/2}$. The functions in
(1)--(3) vanish outside $\supp\nabla w\cup\supp\nabla\Phi\subset N$. If $v$, $w$ and
$h$ are $3$-fold rotationally symmetric, the vector fields $\KBS\ast v$,
$\nabla w$ and $\nabla\Phi$ satisfy $V(\Orot x)=\Orot V(x)$; for
$\KBS\ast v$ this was proved in Lemma~\ref{lem:reduction}.
Their scalar products are therefore $3$-fold rotationally symmetric.
\end{proof}

Lemma~\ref{lem:energy} compares the deviation $w$ of a ring with an
approximate solution $w_{\mathrm{app}}$ by an energy estimate in $H^s$.
In the proofs of Theorem~\ref{thm:quiet} and Lemma~\ref{lem:window},
$w_{\mathrm{app}}$ will be a finite sum of iterated Duhamel terms, starting
with the linear response and adding nonlinear corrections as in
Equation~\eqref{eq:hier} below. The estimate controls the error of this
approximation by its initial error and its residual.

\begin{lemma}
\label{lem:energy}
Let $s\ge3$ be an integer. Fix a sufficiently large rate $\Lam(s)>0$,
depending only on $\Phi$, $\Gone$ and $s$, and assume in addition that
\begin{equation}
\label{eq:Lambda}
\Lam(s)\ge2\max_{G\in[0,\Gone]}\sigma(G).
\end{equation}
Let $G\in[0,\Gone]$
and $t>0$, and let the assumptions of Lemma~\ref{lem:reduction}\textup{(3)} hold on
$[0,t]$, with $w=\omega_{\mathrm r}-\Phi$. Let
$w_{\mathrm{app}}\in C([0,t];X^{s+1})\cap C^1([0,t];X^s)$ be real-valued, and put
$v=w-w_{\mathrm{app}}$ and
\[
R_{\mathrm{app}}=\partial_tw_{\mathrm{app}}-\big(\Lin(G)w_{\mathrm{app}}+B(w_{\mathrm{app}},w_{\mathrm{app}})-h\cdot\nabla w_{\mathrm{app}}-h\cdot\nabla\Phi\big).
\]
If throughout $[0,t]$,
\[
\|w\|_{H^s}\le1,\qquad
\|w_{\mathrm{app}}\|_{H^{s+1}}\le1,\qquad
\|h\|_{C^{s+1}(N')}\le1,
\]
then for
$\tau\in[0,t]$,
\[
\|v(\tau)\|_{H^s}\le e^{\Lam(s)\tau}\|v(0)\|_{H^s}+\int_0^\tau e^{\Lam(s)(\tau-\tau')}\|R_{\mathrm{app}}(\tau')\|_{H^s}\dd\tau' .
\]
\end{lemma}
\begin{proof}
On $N_1$ put $U=\Omega(|x|;G)x^\perp+h+\KBS\ast w$, and for $f$ supported in $N$ let
$U\cdot\nabla f$ mean $(\chi_NU)\cdot\nabla f$.
Recall that $\Lin(G)w=-\Omega(|x|;G)\partial_\theta w-(\KBS\ast w)\cdot\nabla\Phi$
and $\partial_\theta w=x^\perp\cdot\nabla w$. Substituting this and
$B(w,w)=-(\KBS\ast w)\cdot\nabla w$ into the equation of
Lemma~\ref{lem:reduction}\textup{(3)} gives
\[
\begin{aligned}
\partial_tw
&=-\Omega(|x|;G)x^\perp\cdot\nabla w-(\KBS\ast w)\cdot\nabla\Phi
  -(\KBS\ast w)\cdot\nabla w-h\cdot\nabla w-h\cdot\nabla\Phi\\
&=-U\cdot\nabla w-(\KBS\ast w+h)\cdot\nabla\Phi.
\end{aligned}
\]
Moving the transport term to the left yields
$\partial_tw+U\cdot\nabla w=-(\KBS\ast w)\cdot\nabla\Phi-h\cdot\nabla\Phi$.
For the approximate solution, put
$U_{\mathrm{app}}=\Omega(|x|;G)x^\perp+h+\KBS\ast w_{\mathrm{app}}$.
The definition of $R_{\mathrm{app}}$ gives
\[
\partial_tw_{\mathrm{app}}+U_{\mathrm{app}}\cdot\nabla w_{\mathrm{app}}
=-(\KBS\ast w_{\mathrm{app}})\cdot\nabla\Phi-h\cdot\nabla\Phi+R_{\mathrm{app}}.
\]
Since $U-U_{\mathrm{app}}=\KBS\ast v$ and
\[
U\cdot\nabla w-U_{\mathrm{app}}\cdot\nabla w_{\mathrm{app}}
=U\cdot\nabla v+(\KBS\ast v)\cdot\nabla w_{\mathrm{app}},
\]
subtracting the equation for $w_{\mathrm{app}}$ from that for $w$ gives
\[
\partial_tv+U\cdot\nabla v=-(\KBS\ast v)\cdot\nabla(\Phi+w_{\mathrm{app}})-R_{\mathrm{app}},
\]
where every term lies in $C([0,t];H^s)$, and
$v\in C^1([0,t];H^s)\cap C([0,t];H^{s+1})$ is real and supported in $N$.

The field $U$ is divergence-free on $N_1$: $\divop h=0$ by
Lemma~\ref{lem:reduction}\textup{(1)}, $\divop(\KBS\ast w)=0$, and
$\divop(f(|x|)x^\perp)=f'(|x|)\frac{x}{|x|}\cdot x^\perp=0$. By
Lemma~\ref{lem:Hk-calculus}\textup{(3)}, Sobolev embedding and the hypotheses give
\[
\|U\|_{W^{1,\infty}(N_1)}+\|U\|_{H^s(N_1)}\lesssim1.
\]
Here the implicit constant depends only on $\Phi$, $\Gone$ and $s$:
$\Omega(r;G)=\Omega(r;0)+Gr^{-2}$ is uniformly smooth on
$[a/4,3]\times[0,\Gone]$. Since $U$ and $v$ are real and $\divop U=0$ on
$N_1$, Lemma~\ref{lem:Hk-calculus}\textup{(1)} bounds the transport contribution
by $C\|v\|_{H^s}^2$. Parts \textup{(2)} and \textup{(3)} of that lemma also give
\[
\|(\KBS\ast v)\cdot\nabla(\Phi+w_{\mathrm{app}})\|_{H^s}
\lesssim\|v\|_{H^s}
\big(\|\Phi\|_{H^{s+1}}+\|w_{\mathrm{app}}\|_{H^{s+1}}\big)
\lesssim\|v\|_{H^s}.
\]
Apply $\partial^\alpha$ to the equation of $v$, take the $L^2$ inner product
with $\partial^\alpha v$, and sum over $|\alpha|\le s$. Choosing $\Lam(s)$
large enough to dominate these bounds gives
$\frac12\frac{\dd}{\dd\tau}\|v\|_{H^s}^2\le\Lam(s)\|v\|_{H^s}^2+\|R_{\mathrm{app}}\|_{H^s}\|v\|_{H^s}$.
To divide by $\|v\|_{H^s}$, which may vanish, we regularize it. For $\epsilon>0$
the function $y_\epsilon=(\|v\|_{H^s}^2+\epsilon^2)^{1/2}$ is $C^1$, as
$v\in C^1([0,t];H^s)$ and $y_\epsilon\ge\epsilon$. Since
$y_\epsilon y_\epsilon'=\frac12\frac{\dd}{\dd\tau}\|v\|_{H^s}^2$ and
$\|v\|_{H^s}\le y_\epsilon$, it satisfies
$y_\epsilon'\le\Lam(s)y_\epsilon+\|R_{\mathrm{app}}\|_{H^s}$, so
$(e^{-\Lam(s)\tau}y_\epsilon)'\le e^{-\Lam(s)\tau}\|R_{\mathrm{app}}\|_{H^s}$.
Integrating from $0$ to $\tau$, and letting $\epsilon\to0$ with
$\|v\|_{H^s}\le y_\epsilon$ and $y_\epsilon(0)\to\|v(0)\|_{H^s}$, proves the lemma.
\end{proof}

\emph{Approximate solutions.} The semigroup estimates of Part~I give precise
control of the linear evolution and its response to a source. Following
Grenier, we use this control at each order, approximating the deviation $w$
by the first $J$ terms $w^{[1]},\dots,w^{[J]}$ of its formal expansion in powers of
the size of the data, in which $h$ and the initial value count as terms of first
order: $w^{[1]}$ solves the linear equation driven by $-h\cdot\nabla\Phi$, with an
initial value $w_0$ that is $0$ or $w(0)$ in our applications, and $w^{[j]}$,
$j\ge2$, solves the linear equation with initial value $0$, driven by the terms of
order $j$ of $B(w,w)-h\cdot\nabla w$. For example, the order-two term of this
nonlinear expression is $B(w^{[1]},w^{[1]})-h\cdot\nabla w^{[1]}$.
Spatial differentiation preserves the formal order: both $w^{[j]}$ and
$\nabla w^{[j]}$ have order $j$. Multiplying by the first-order field $h$
therefore makes $h\cdot\nabla w^{[j]}$ of order $j+1$. The background $\Phi$
is fixed and has order zero, so $h\cdot\nabla\Phi$ has order one and enters
the equation for $w^{[1]}$. The residual of the sum consists of terms of
order larger than $J$ (Lemma~\ref{lem:hier-res} below), and
Lemma~\ref{lem:energy} controls the difference
between $w$ and the sum. Let the assumptions of Lemma~\ref{lem:reduction} hold
on $[0,T]$, assume $G\in[0,\Gone]$, and put $\Lin=\Lin(G)$ and $F=-h\cdot\nabla\Phi$.
By Lemma~\ref{lem:reduction}\textup{(1)} and Lemma~\ref{lem:bilinear}, $F$ is real
and lies in $C([0,T];X^k)$ for every $k$; so, by Lemma~\ref{lem:duhamel}, for every
real $w_0\in\bigcap_kX^k$ the function
$e^{t\Lin}w_0+\int_0^te^{(t-\tau)\Lin}F(\tau)\dd\tau$ is real, lies in
$C([0,T];X^k)\cap C^1([0,T];X^{k-1})$ for every $k\ge1$, and solves
$\partial_tu=\Lin u+F$. If $\supp\omega_{\mathrm r}(0)\subset N$, this applies to
$w_0=w(0)=\omega_{\mathrm r}(0)-\Phi$, which is then smooth
(Lemma~\ref{lem:reduction}\textup{(2)}), real, $3$-fold rotationally symmetric and supported in $N$,
hence lies in every $X^k$. Let $w^{[1]}$ be such a function, let $J\ge1$ be an integer,
and for $2\le j\le J$ define recursively
\begin{equation}
\label{eq:hier}
w^{[j]}(t)=\int_0^te^{(t-\tau)\Lin}g^{[j]}(\tau)\dd\tau,\qquad
g^{[j]}=\sum_{i+l=j}B(w^{[i]},w^{[l]})-h\cdot\nabla w^{[j-1]},
\end{equation}
where here and below $i,l\ge1$ in the sums. By induction on $j$,
Lemma~\ref{lem:bilinear} (whose bounds, with the continuity of $t\mapsto h(t)$ in
every $C^k(N')$, give continuity in time) and Lemma~\ref{lem:duhamel} show that
$g^{[j]}$ and $w^{[j]}$ are real, $g^{[j]}\in C([0,T];X^k)$ and
$w^{[j]}\in C([0,T];X^k)\cap C^1([0,T];X^{k-1})$ for every $k\ge1$, and
$\partial_tw^{[j]}=\Lin w^{[j]}+g^{[j]}$.

\begin{lemma}[Residual of the approximate solution]
\label{lem:hier-res}
Let $w^{[1]},\dots,w^{[J]}$ be constructed above, and put
$w_{\mathrm{app}}=\sum_{j=1}^Jw^{[j]}$. Then $w_{\mathrm{app}}(0)=w_0$, and its
residual, as defined in Lemma~\ref{lem:energy}, is
\begin{equation}
\label{eq:hier-res}
R_{\mathrm{app}}=-\sum_{i,l\le J,\,i+l>J}B(w^{[i]},w^{[l]})+h\cdot\nabla w^{[J]} .
\end{equation}
Every term on the right has formal order larger than $J$.
\end{lemma}
\begin{proof}
The first-order equation is
$\partial_tw^{[1]}-\Lin w^{[1]}=F=-h\cdot\nabla\Phi$.
Summing it with the equations for $2\le j\le J$ gives
\[
\partial_tw_{\mathrm{app}}-\Lin w_{\mathrm{app}}
=F+\sum_{\substack{1\le i,l\le J\\i+l\le J}}B(w^{[i]},w^{[l]})
  -h\cdot\nabla\sum_{j=1}^{J-1}w^{[j]}.
\]
By bilinearity,
$B(w_{\mathrm{app}},w_{\mathrm{app}})=\sum_{1\le i,l\le J}B(w^{[i]},w^{[l]})$,
and $\sum_{j=1}^{J-1}w^{[j]}=w_{\mathrm{app}}-w^{[J]}$. Substituting these
identities into the definition of the residual yields
\[
R_{\mathrm{app}}=
\underbrace{F+h\cdot\nabla\Phi}_{=\,0}
-\sum_{i,l\le J,\,i+l>J}B(w^{[i]},w^{[l]})+h\cdot\nabla w^{[J]},
\]
which proves Equation~\eqref{eq:hier-res}. Thus the background-driving term
cancels because it is already included in the first-order equation.
The bilinear terms have order $i+l>J$, and the last term has order $J+1$.
Finally, $w^{[1]}(0)=w_0$ and $w^{[j]}(0)=0$ for $j\ge2$ give the initial value.
\end{proof}

\subsection{A Quiet Ring}

A ring is quiet while the instability has not reached it: it starts exactly at
$\Phi$, so its deviation is created by the field $h$ alone, and
Theorem~\ref{thm:quiet} shows that it stays small. We bound a norm of the field by a
function $\eta$ of time, and its effect by the Duhamel bound of $\eta$, defined
next. Since the semigroup bound of Proposition~I.\ref{I-prop:semigroup} grows like
$e^{\sigma(G)t}$, this integral controls, up to a constant factor, the linear
response to a driving term of size $\eta$.

\begin{definition}
\label{def:majorant}
Let $G\in[0,\Gone]$, $T>0$, and let $\eta\colon[0,T]\to[0,\infty)$ be integrable.
The \textbf{Duhamel bound} of $\eta$ is
$\maj(t)=\int_0^te^{\sigma(G)(t-\tau)}\eta(\tau)\dd\tau$, $t\in[0,T]$.
\end{definition}

\begin{lemma}
\label{lem:majorant}
Let $G$, $T$, $\eta$ and $\maj$ be as in Definition~\ref{def:majorant}. Then for
$0\le\tau\le t\le T$:
\textup{(M1)} $e^{\sigma(G)(t-\tau)}\maj(\tau)\le\maj(t)$;
\textup{(M2)} $\int_\tau^t\eta\le\maj(t)$;
\textup{(M3)} $\int_0^t\maj\le\maj(t)/\sigma(G)$.
\end{lemma}
\begin{proof}
Write $\sigma=\sigma(G)$, which is positive by Proposition~I.\ref{I-prop:SU}. We have
$\maj(t)=e^{\sigma(t-\tau)}\maj(\tau)+\int_\tau^te^{\sigma(t-v)}\eta(v)\dd v$, and the
integral is nonnegative, which gives (M1), and at least $\int_\tau^t\eta$, which
gives (M2). By (M1), $\maj(\tau)\le e^{-\sigma(t-\tau)}\maj(t)$; integrating over
$\tau\in[0,t]$ gives (M3).
\end{proof}

Recall that $\sigb=\min_{G\in[0,\Gone]}\sigma(G)>0$ is the uniform lower
bound on the unstable growth rate. In the next theorem, $Q$ stands for
``quiet'': $N_Q$ is the order of the approximate solution and $C_Q$ is the
constant in the bounds for a quiet ring; $\majst$ is the smallness threshold for the
Duhamel bound. The subscript $h$ labels the regularity indices used in the
field estimates: the perturbation is controlled in $H^{s_h}$ and the field
$h$ is assumed bounded in $C^{M_h}$. Taking $M_h=s_h+N_Q+1$ supplies the derivatives needed
by the successive Duhamel terms and the energy estimate. Section~\ref{sec:relay}
uses these bounds to control the velocity fields induced by quiet rings,
and Section~\ref{sec:bootstrap} keeps those rings close to $\Phi$ until activation.

The condition $N_Q\sigb\ge\Lam(s_h)$ ensures that truncation at order $N_Q$
overcomes the growth allowed by the energy estimate. Indeed, the residual
contains a factor $\maj^{N_Q}$, and (M1) gives, for $0\le\tau\le t$,
\[
e^{\Lam(s_h)(t-\tau)}\maj(\tau)^{N_Q}
\le e^{(\Lam(s_h)-N_Q\sigma(G))(t-\tau)}\maj(t)^{N_Q}
\le\maj(t)^{N_Q}.
\]
Thus the remainder retains the higher-order smallness needed to compare
$w$ with its linear response.

\begin{theorem}
\label{thm:quiet}
Assume \textup{(N1)}--\textup{(N3)}. Let $s_h\ge3$ and $N_Q\ge2$ be integers with
$N_Q\sigb\ge\Lam(s_h)$, and put $M_h=s_h+N_Q+1$. There are $\majst\in(0,1]$ and
$C_Q<\infty$, depending only on $(\Phi,\Gone,s_h,N_Q)$, such that the following
holds. Let $\omega$,
$\omega_{\mathrm r},\omega_{\mathrm{inn}},\omega_{\mathrm{ou}}$, $G$ and $h$ be as in
Lemma~\ref{lem:reduction} on $[0,T]$, with $\omega_{\mathrm r}(0)=\Phi$ and
$G\in[0,\Gone]$, and put $w=\omega_{\mathrm r}-\Phi$. Let $\eta$ be measurable with
$\|h(t)\|_{C^{M_h}(N')}\le\eta(t)\le1$ on $[0,T]$, let $\maj$ be its Duhamel
bound, and assume $\maj\le\majst$ on $[0,T]$. Then for every $t\in[0,T]$:
\begin{enumerate}[label=\textup{(\alph*)}]
\item $\supp\omega_{\mathrm r}(t)\subset\NNm$;
\item $\|w(t)-\wlin(t)\|_{H^{s_h}}\le C_Q\maj(t)^2$, where
  $\wlin(t)=\int_0^te^{(t-\tau)\Lin(G)}\big(-h(\tau)\cdot\nabla\Phi\big)\dd\tau$
  is the \textbf{linear response};
\item $\|w(t)\|_{H^{s_h}}\le C_Q\maj(t)$.
\end{enumerate}
\end{theorem}
\begin{proof}
Throughout, $\sigma=\sigma(G)\ge\sigb$ and $\Lin=\Lin(G)$. We compare $w$ with the
approximate solution of Equation~\eqref{eq:hier} with $J=N_Q$ and $w^{[1]}=\wlin$,
bound its terms by powers of $\maj$, and close with Lemma~\ref{lem:energy} in a
continuity argument on the radial displacement and on $\|w\|_{H^{s_h}}$.

\emph{Constants.} Fix the semigroup constants of
Proposition~I.\ref{I-prop:semigroup} at the rate
$\frac12(\frac9{200}+\sigb)\in(\frac9{200},\sigb)$, and the constants of
Lemma~\ref{lem:bilinear} through order $M_h$. They are finite and uniform
in $G$. All implicit constants in this proof depend only on
$(\Phi,\Gone,s_h,N_Q)$. Put $k_j=M_h-j$ for $1\le j\le N_Q$, so that
$k_{N_Q}=s_h+1$ and $k_j+1=k_{j-1}$. The finite induction below gives the
hierarchy, residual and energy bounds with these dependencies. Choose $C_Q$
large enough to dominate the coefficients in the resulting norm bounds,
and then choose $\majst>0$ sufficiently small that
\[
\majst\le1,\qquad C_Q\majst\le\tfrac12,\qquad
\big(1+C_{\mathrm{BS}}C_Q/\sigb\big)\majst\le a/16.
\]

\emph{Approximate solution.} Let $w^{[j]}$, $1\le j\le N_Q$, be the functions of
Equation~\eqref{eq:hier} with $J=N_Q$ and $w^{[1]}=\wlin$ (the case $w_0=0$), so that
they are regular by Lemmas~\ref{lem:duhamel} and \ref{lem:bilinear}. By
Lemma~\ref{lem:bilinear}\textup{(3)},
$\|F(\tau)\|_{H^{M_h}}\lesssim\eta(\tau)$. We claim that
\begin{equation}
\label{eq:Q-hier}
\|w^{[j]}(t)\|_{H^{k_j}}\lesssim\maj(t)^j\qquad(t\in[0,T],\ 1\le j\le N_Q).
\end{equation}
For $j=1$,
$\|\wlin(t)\|_{H^{k_1}}\lesssim\int_0^te^{\sigma(t-\tau)}\eta(\tau)\dd\tau\le\maj(t)$.
Let $j\ge2$ and assume Equation~\eqref{eq:Q-hier} for all indices less than $j$.
Since $k_j+1=k_{j-1}\le k_l$ for $l\le j-1$, Lemma~\ref{lem:bilinear} gives
$\|g^{[j]}(\tau)\|_{H^{k_j}}\lesssim\maj(\tau)^j+\eta(\tau)\maj(\tau)^{j-1}$.
By (M1), $e^{\sigma(t-\tau)}\maj(\tau)^j\le e^{-(j-1)\sigma(t-\tau)}\maj(t)^j$, whose
integral over $[0,t]$ is at most $\maj(t)^j/((j-1)\sigb)$; and
$e^{\sigma(t-\tau)}\eta(\tau)\maj(\tau)^{j-1}\le\eta(\tau)\maj(t)^{j-1}$, whose integral
is at most $\maj(t)^j$ by (M2). With $\|e^{(t-\tau)\Lin}\|\lesssim e^{\sigma(t-\tau)}$ this
gives Equation~\eqref{eq:Q-hier} for $j$. Since $k_j\ge s_h+1$ and $\maj\le1$,
Equations~\eqref{eq:hier-res} and \eqref{eq:Q-hier} and Lemma~\ref{lem:bilinear} give
\begin{equation}
\label{eq:Q-res}
\|R_{\mathrm{app}}(\tau)\|_{H^{s_h}}\lesssim\big(\maj(\tau)^{N_Q+1}+\eta(\tau)\maj(\tau)^{N_Q}\big),
\end{equation}
and $\|w_{\mathrm{app}}\|_{H^{s_h+1}}\lesssim\maj$ on $[0,T]$. Our choice of
$C_Q$ gives $\|w_{\mathrm{app}}\|_{H^{s_h+1}}\le C_Q\majst\le\frac12$ there.

\emph{Continuity argument.} Let $t_*$ be the largest $t\in[0,T]$ such that
$\delta_{\mathrm r}(t)\le a/8$ and $\|w(\tau)\|_{H^{s_h}}\le1$ for $\tau\in[0,t]$; it
exists and is positive since $\delta_{\mathrm r}(0)=0$, $w(0)=0$, $\delta_{\mathrm r}$
is continuous (see its definition at the start of this section), and
$\|w\|_{H^{s_h}}$ is continuous by Theorem~\ref{thm:smooth}\textup{(1)}. Because $\supp\Phi\subset\{a\le|y|\le1\}$ and
$1+\frac a8\le\frac54$, the support of $\omega_{\mathrm r}(\tau)$ lies in
$\NNm\subset N$ for $\tau\le t_*$. So Lemma~\ref{lem:reduction}\textup{(3)} and
Lemma~\ref{lem:energy} with $s=s_h$ apply on $[0,t_*]$; the conditions
$\|w_{\mathrm{app}}\|_{H^{s_h+1}}\le1$ and $\|h\|_{C^{s_h+1}(N')}\le\eta\le1$ hold.
Since $v=w-w_{\mathrm{app}}$ has $v(0)=0$, Equation~\eqref{eq:Q-res} gives, for
$t\le t_*$,
\[
\|v(t)\|_{H^{s_h}}\lesssim\int_0^te^{\Lam(s_h)(t-\tau)}\maj(\tau)^{N_Q}\big(\maj(\tau)+\eta(\tau)\big)\dd\tau .
\]
For $0\le\tau\le t$, (M1) and $N_Q\sigma\ge N_Q\sigb\ge\Lam(s_h)$ give
$e^{\Lam(s_h)(t-\tau)}\maj(\tau)^{N_Q}\le e^{(\Lam(s_h)-N_Q\sigma)(t-\tau)}\maj(t)^{N_Q}$,
which is at most $\maj(t)^{N_Q}$. By (M2) and (M3) we conclude, for $t\le t_*$,
\begin{equation}
\label{eq:Q-v}
\|v(t)\|_{H^{s_h}}\lesssim\maj(t)^{N_Q}\int_0^t(\maj+\eta)\lesssim\big(1+\sigb^{-1}\big)\maj(t)^{N_Q+1}.
\end{equation}
Since $w=\wlin+\sum_{j\ge2}w^{[j]}+v$ and $\maj\le1$, Equations~\eqref{eq:Q-hier}
and \eqref{eq:Q-v} give, for $t\le t_*$,
\begin{equation}
\label{eq:Q-b}
\|w(t)-\wlin(t)\|_{H^{s_h}}\le C_Q\maj(t)^2,\qquad\|w(t)\|_{H^{s_h}}\le C_Q\maj(t)\le\tfrac12 .
\end{equation}
For $y\in\supp\Phi$ and $0\le u\le t\le t_*$, the point $\flow_u(y)$ lies in $N$.
The background velocity $\Omega(|x|;G)x^\perp$ is tangential, so
Lemma~\ref{lem:reduction}\textup{(2)} gives
\[
\frac{\dd}{\dd u}|\flow_u(y)|
=\big(h+\KBS\ast w\big)(u,\flow_u(y))\cdot\er(\flow_u(y)).
\]
Since $\er$ is a unit vector and $\flow_u(y)\in N\subset B_3$, the absolute
value of the right side is at most
$\|h(u)\|_{C^0(N')}+C_{\mathrm{BS}}\|w(u)\|_{H^1}
\le\eta(u)+C_{\mathrm{BS}}\|w(u)\|_{H^{s_h}}$,
by Lemma~\ref{lem:bilinear}\textup{(4)} and $s_h\ge3$.
Integrating from $0$ to $\tau\le t$ and enlarging the integral to $[0,t]$
gives the first inequality below; Equation~\eqref{eq:Q-b}, (M2) and (M3)
then give the remaining bounds:
\[
\big||\flow_\tau(y)|-|y|\big|\le\int_0^t\big(\eta+C_{\mathrm{BS}}\|w\|_{H^{s_h}}\big)\le\maj(t)+\frac{C_{\mathrm{BS}}C_Q}{\sigb}\maj(t).
\]
So $\delta_{\mathrm r}\le(1+C_{\mathrm{BS}}C_Q/\sigb)\maj\le a/16$ and $\|w\|_{H^{s_h}}\le\frac12$ on
$[0,t_*]$. If $t_*<T$, continuity would give a larger $t$ with
$\delta_{\mathrm r}(t)\le a/8$ and $\|w\|_{H^{s_h}}\le1$ on $[0,t]$; hence $t_*=T$.
Then (a) follows from $\delta_{\mathrm r}\le(1+C_{\mathrm{BS}}C_Q/\sigb)\maj$ and Lemma~\ref{lem:annuli},
and (b), (c) are
Equation~\eqref{eq:Q-b}.
\end{proof}

\subsection{An Active Ring}

A ring is active once its deviation from $\Phi$ has a component along the
unstable eigenfunction that grows at the rate $\sigma(G)$. Theorem~\ref{thm:active}
shows that the deviation then follows this component, up to a transient and
controlled errors, until its envelope $\kappa e^{\sigma(G)t}$ reaches a small
size $\delta_0$. The transient reflects the slower linear evolution outside
the unstable eigenspace: its contribution becomes small relative to the
growing unstable component.
Corollary~\ref{cor:active-moment} gives the outer moment of the ring, through
which it drives the rings outside it. The unstable component of the linear
response is explicit, by Proposition~I.\ref{I-prop:semigroup}; the next lemma
writes it with a running coefficient $c(t)$, which accounts for the driving by
the velocity field $h$.

\begin{lemma}
\label{lem:running}
Let $G\in[0,\Gone]$, $\bar t>0$, $w_0\in X^0$ real, and let $h$ be continuous on
$[0,\bar t]$ with values in $C^0(N';\R^2)$ and $3$-fold rotationally symmetric. Put
$F(t)=-h(t)\cdot\nabla\Phi$ and
\begin{equation}
\label{eq:running}
c(t)=\ellst(w_0)+\int_0^te^{-\lamst v}\,\ellst(F(v))\dd v\qquad(0\le t\le\bar t),
\end{equation}
the \textbf{running coefficient}. Then
$\wlin_0(t):=e^{t\Lin(G)}w_0+\int_0^te^{(t-v)\Lin(G)}F(v)\dd v$ satisfies
$\ellst(\wlin_0(t))=c(t)e^{\lamst t}$ and
$\Riesz\,\wlin_0(t)=2\Re\big(c(t)e^{\lamst t}\phist\big)$ for $0\le t\le\bar t$.
\end{lemma}
\begin{proof}
The function $F$ is continuous on $[0,\bar t]$ with values in $X^0$, since
$\|h(t)\cdot\nabla\Phi-h(t')\cdot\nabla\Phi\|_{L^2}\le\|h(t)-h(t')\|_{C^0(N')}\|\nabla\Phi\|_{L^2}$,
and $F(t)$ is real, supported in $N$ and $3$-fold rotationally symmetric. The functional
$\ellst$ is bounded on $X^0$, so it commutes with the integral, and
Proposition~I.\ref{I-prop:semigroup} gives $\ellst(e^{\tau\Lin}f)=e^{\lamst\tau}\ellst(f)$
for $f\in X^0$. Hence
$\ellst(\wlin_0(t))=e^{\lamst t}\ellst(w_0)+\int_0^te^{\lamst(t-v)}\ellst(F(v))\dd v=c(t)e^{\lamst t}$.
Since $\wlin_0(t)$ is real, Definition~I.\ref{I-def:eigendata} gives
$\Riesz\,\wlin_0(t)=2\Re\big(\ellst(\wlin_0(t))\phist\big)=2\Re\big(c(t)e^{\lamst t}\phist\big)$.
\end{proof}

Recall that the Japanese bracket is $\jap{t}=(1+t^2)^{1/2}$, so
$\jap{t}^p$ is a polynomial weight that is nonzero at $t=0$.
On the range of $1-\Riesz$ the semigroup grows at most at a rate
$\sigma'<\sigb\le\sigma(G)$ (Proposition~I.\ref{I-prop:semigroup}).
The next lemma bounds this part of the linear response while keeping its
two sources of growth separate: the initial deviation and the polynomial
part of $h$ give a bound proportional to $\jap{t}^pe^{\sigma't}$, whereas
the part of $h$ growing like $e^{\sigma t}$ produces a response at that same
rate. The former becomes small relative to the growing unstable component;
the latter is controlled later by choosing its prefactor $\eta_b$ small.

\begin{lemma}
\label{lem:linear-part}
Assume \textup{(N1)}--\textup{(N3)}. Let $\sigma'\in(\frac9{200},\sigb)$, and let
$p\ge0$ and $k\ge2$ be integers. There is a constant $C_{\mathrm{lin}}(k)\ge1$,
depending only on $\Phi$, $\Gone$, $k$ and $\sigma'$, with the following property. Let $G\in[0,\Gone]$,
$\sigma=\sigma(G)$, $\lambda=\lamst(G)$, $\bar t>0$, $A\ge0$, $\kappa\ge0$ and
$\eta_b\ge0$. Let $w_0\in X^k$ be real with $\|w_0\|_{H^k}\le A$, and let $h$ be
continuous on $[0,\bar t]$ with values in $C^k(N';\R^2)$, $3$-fold rotationally symmetric, with
$\|h(t)\|_{C^k(N')}\le A\jap{t}^p+\eta_b\kappa e^{\sigma t}$. Let $c$ and $\wlin_0$ be as
in Lemma~\ref{lem:running}. Then for $t\in[0,\bar t]$,
\[
\begin{gathered}
\wlin_0(t)=2\Re\big(c(t)e^{\lambda t}\phist\big)+(1-\Riesz)\wlin_0(t),\\
\|(1-\Riesz)\wlin_0(t)\|_{H^k}\le C_{\mathrm{lin}}(k)\big(A\jap{t}^pe^{\sigma't}+\eta_b\kappa e^{\sigma t}\big).
\end{gathered}
\]
\end{lemma}
\begin{proof}
The decomposition is Lemma~\ref{lem:running}. For real $f\in X^k$,
Proposition~I.\ref{I-prop:semigroup} gives
$\Riesz e^{\tau\Lin}f=2\Re\big(\ellst(e^{\tau\Lin}f)\phist\big)=2\Re\big(e^{\lambda\tau}\ellst(f)\phist\big)=e^{\tau\Lin}\Riesz f$.
Since $\Riesz$ is bounded on $X^k$, it commutes with the integral, and
$(1-\Riesz)\wlin_0(t)=e^{t\Lin}(1-\Riesz)w_0+\int_0^te^{(t-v)\Lin}(1-\Riesz)F(v)\dd v$.
By Proposition~I.\ref{I-prop:semigroup} and the bound
$\|h\cdot\nabla\Phi\|_{H^k}\le\Cnab(k)\|h\|_{C^k(N')}$ of Lemma~I.\ref{I-lem:uniform}, its
norm in $H^k$ is at most
$C_S(k,\sigma')\big(e^{\sigma't}A+\Cnab(k)\int_0^te^{\sigma'(t-v)}(A\jap{v}^p+\eta_b\kappa e^{\sigma v})\dd v\big)$.
Finally
$\int_0^te^{\sigma'(t-v)}\jap{v}^p\dd v\le\jap{t}^pe^{\sigma't}/\sigma'$ and
$\int_0^te^{\sigma'(t-v)}e^{\sigma v}\dd v\le e^{\sigma t}/(\sigma-\sigma')\le e^{\sigma t}/(\sigb-\sigma')$.
\end{proof}

For $x\in(0,1]$ we write $\varsigma(x)=1+\sigb^{-1}\log(1/x)$. Since
$\sigma(G)\ge\sigb$, this number bounds $1+t$ for $t\le\sigma(G)^{-1}\log(1/x)$,
the time in which a coefficient of size $x$ that grows at the rate $\sigma(G)$
reaches size $1$. The cubic time weight in the active-ring estimates below matches the
mode-$3$ moment bound for the inner vorticity: this moment is bounded by a
constant times the cube of its confinement radius, which grows at most linearly in time
(Proposition~\ref{prop:confine}). For an integer $s\ge3$ and a number
$\sigma'\in(\frac9{200},\sigb)$ we put
\begin{equation}
\label{eq:W-const}
N_G=\min\big\{n\in\Z:n\ge2,\ n\ge\Lam(s)/\sigb+1\big\},\qquad
M=s+N_G+2,\qquad N_W=5N_G+1.
\end{equation}
The semigroup bounds, Lemmas~\ref{lem:bilinear} and~\ref{lem:linear-part}, and
the eigenfunction bounds of Lemma~I.\ref{I-lem:uniform} have constants
uniform in $G$ through these derivative orders. Also
$\sup_{t\ge0}\jap{t}^3e^{-(\sigb-\sigma')t}<\infty$ because
$\sigb>\sigma'$. These bounds depend only on $(\Phi,\Gone,s,\sigma')$. As in
Theorem~\ref{thm:quiet}, the order $N_G$ of the approximate solution is chosen
large, so that its residual, which is one order higher, is small enough on the
window of Lemma~\ref{lem:window} to absorb the growth $e^{\Lam(s)t}$ allowed by
Lemma~\ref{lem:energy} (Equation~\eqref{eq:W-exp} below).

Lemma~\ref{lem:window} is the main step in the proof of Theorem~\ref{thm:active}.
It treats an active ring whose initial deviation and running coefficient are
bounded by a fixed multiple $C_\natural$ of a number $\kappa_1$, and whose field is
bounded accordingly. It does so on the window of time in which a component of
initial size $\kappa_1$, growing at the rate $\sigma(G)$, stays below $\delta_0$.

\begin{lemma}
\label{lem:window}
Assume \textup{(N1)}--\textup{(N3)}. Let $s\ge3$ be an integer and
$\sigma'\in(\frac9{200},\sigb)$, let $N_G$, $M$ and $N_W$ be as in
Equation~\eqref{eq:W-const}, and let $C_\natural\ge1$. There are $\delta_W\in(0,\frac12]$, $\kappa_W\in(0,1]$ and
$C_W\ge1$, depending only on $(\Phi,\Gone,s,\sigma',C_\natural)$, such that the
following holds. Let $G\in[0,\Gone]$, $\sigma=\sigma(G)$, $\lambda=\lamst(G)$,
$\kappa_1\in(0,\kappa_W]$, $\delta_0\in(0,\delta_W]$, $\eta_b\in[0,1]$,
$\hat T_1=\sigma^{-1}\log(\delta_0/\kappa_1)$ and $\bar t\in(0,\hat T_1]$. Let
$\omega$, $\omega_{\mathrm r},\omega_{\mathrm{inn}},\omega_{\mathrm{ou}}$ and $h$ be as in
Lemma~\ref{lem:reduction} on $[0,\bar t]$ with this $G$, put
$w=\omega_{\mathrm r}-\Phi$, and let $c$ be the running coefficient of Equation~\eqref{eq:running}
with $w_0=w(0)$. Assume that $\|w(0)\|_{H^{M-2}}\le C_\natural\kappa_1$, that
$\supp(\Phi+w(0))\subset\{(\frac78-\frac1{64})a\le|y|\le\frac54+\frac a{64}\}$, and
that for $t\in[0,\bar t]$,
\[
\|h(t)\|_{C^M(N')}\le C_\natural\kappa_1\jap{t}^3+\eta_b\kappa_1e^{\sigma t},\qquad
|c(t)|\le C_\natural\kappa_1 .
\]
Put $\hat a_1(t)=\kappa_1e^{\sigma t}$ and $\rho_1(t)=\jap{t}^3e^{-(\sigma-\sigma')t}$.
Then for every $t\in[0,\bar t]$, every trajectory of $\omega_{\mathrm r}$ moves
radially by at most $C_W\delta_0+C_W\kappa_1\varsigma(\kappa_1)^{N_W}\le a/32$ on
$[0,t]$, and
\[
\big\|w(t)-2\Re\big(c(t)e^{\lambda t}\phist\big)\big\|_{H^s}\le C_W\hat a_1(t)\big[\rho_1(t)+\hat a_1(t)+\kappa_1\varsigma(\kappa_1)^{N_W}+\eta_b\big].
\]
\end{lemma}
\begin{proof}
Write $\Lam=\Lam(s)$ and $\Lin=\Lin(G)$. Since $N_G\ge\Lam/\sigb+1$ and
$\Lam\ge2\sigma\ge2\sigb$ (Lemma~\ref{lem:energy}),
\begin{equation}
\label{eq:W-exp}
(N_G+1)\sigma-\Lam\ge2\sigma,\qquad\Lam-\sigma\ge\sigma\ge\sigb,\qquad N_G+1-\Lam/\sigma\ge2 .
\end{equation}

\emph{Constants.} All implicit constants in this proof depend only on
$(\Phi,\Gone,s,\sigma',C_\natural)$. The finite hierarchy, residual and energy
estimates below produce constants with these dependencies. Choose
$C_W\ge C_\natural$ large enough to dominate their coefficients and those in
the final norm and displacement estimates. Then choose
$\delta_W\in(0,1/2]$ so that $C_W\delta_W\le a/64$, and
$\kappa_W\in(0,1]$ so that
$C_Wx\varsigma(x)^{N_W}\le a/64$ for every $0<x\le\kappa_W$.
To justify the latter choice, recall that $N_W=5N_G+1$ and $\sigb>0$ are
already fixed. With $L=\log(1/x)$,
\[
x\varsigma(x)^{N_W}=e^{-L}(1+L/\sigb)^{N_W}.
\]
Its logarithmic derivative with respect to $L$ is
\[
-1+\frac{N_W}{\sigb+L}\le-\frac12
\qquad\text{whenever }\sigb+L\ge2N_W.
\]
Thus the quantity decreases to zero as $L\to\infty$. Choose $\kappa_W$
small enough that $\sigb+\log(1/\kappa_W)\ge2N_W$ and
$C_W\kappa_W\varsigma(\kappa_W)^{N_W}\le a/64$; monotonicity then gives
the required inequality for every $0<x\le\kappa_W$.

Now fix the data, and put $\varepsilon_W=C_\natural\kappa_1$,
$\mathfrak b(t)=\hat a_1(t)+\varepsilon_W\jap{t}^{5}$ and $k_j=M-1-j$ for
$1\le j\le N_G$, so $k_1=M-2$, $k_{N_G}=s+1$, $k_j+1=k_{j-1}$. The function
$\mathfrak b$ plays the role of $\maj$ in the proof of Theorem~\ref{thm:quiet}: the
terms $w^{[j]}$ below are bounded by multiples of $\hat a_1\mathfrak b^{j-1}$. For
$t\le\bar t\le\hat T_1\le\sigma^{-1}\log(1/\kappa_1)$ we have $\jap{t}\le1+t\le\varsigma(\kappa_1)$,
$\hat a_1(t)\le\delta_0$ and $\kappa_1\le\hat a_1(t)$. By the choice of $\kappa_W$ and
$\delta_W$, for $t\in[0,\bar t]$,
\begin{equation}
\label{eq:W-small}
\mathfrak b(t)\le\delta_0+C_\natural\kappa_1\varsigma(\kappa_1)^{5}\le1,\qquad
\|h(t)\|_{C^{s+1}(N')}\le C_\natural\kappa_1\varsigma(\kappa_1)^3+\delta_0\le1 .
\end{equation}

\emph{Approximate solution.} Let $w^{[j]}$ be the functions of
Equation~\eqref{eq:hier} with $J=N_G$ and $w^{[1]}=\wlin_0$, the function of
Lemma~\ref{lem:running} with $w_0=w(0)$. This choice of $w^{[1]}$ is admissible in
Equation~\eqref{eq:hier} because $\supp(\Phi+w(0))\subset N$, and the functions
$w^{[j]}$ are regular by Lemmas~\ref{lem:duhamel} and \ref{lem:bilinear}. We
claim that
\begin{equation}
\label{eq:W-hier}
\|w^{[j]}(t)\|_{H^{k_j}}\lesssim\hat a_1(t)\,\mathfrak b(t)^{j-1}\qquad(t\in[0,\bar t],\ 1\le j\le N_G).
\end{equation}
For $j=1$, Lemma~\ref{lem:linear-part} with $k=M-2$, $A=C_\natural\kappa_1$ and
$\kappa=\kappa_1$ gives
$\|(1-\Riesz)w^{[1]}(t)\|_{H^{M-2}}\lesssim\hat a_1(t)\big(C_\natural\rho_1(t)+\eta_b\big)$,
since $C_\natural\kappa_1\jap{t}^3e^{\sigma't}=C_\natural\hat a_1(t)\rho_1(t)$; and
$\|2\Re(c(t)e^{\lambda t}\phist)\|_{H^{M-2}}\lesssim C_\natural\hat a_1(t)$. As
$\rho_1$ is uniformly bounded and $\eta_b\le1$, this is Equation~\eqref{eq:W-hier} for $j=1$.
Let $j\ge2$. By Lemma~\ref{lem:bilinear}, Equation~\eqref{eq:W-hier} for the
indices less than $j$, and
$\|h(v)\|_{C^{k_j}(N')}\le\varepsilon_W\jap{v}^3+\hat a_1(v)$,
\[
\|g^{[j]}(v)\|_{H^{k_j}}\lesssim\hat a_1(v)\big(\hat a_1(v)+\varepsilon_W\jap{v}^3\big)\mathfrak b(v)^{j-2}.
\]
Since $e^{\sigma(t-v)}\hat a_1(v)=\hat a_1(t)$ and
$\|e^{(t-v)\Lin}\|\lesssim e^{\sigma(t-v)}$, it suffices to show that
\[
\int_0^t\big(\hat a_1+\varepsilon_W\jap{v}^3\big)\mathfrak b^{j-2}\dd v\le2^{j-1}\max(1,\sigb^{-1})\,\mathfrak b(t)^{j-1}.
\]
Expand $\mathfrak b^{j-2}=\sum_{i+q=j-2}\binom{j-2}i\hat a_1^i(\varepsilon_W\jap{v}^{5})^q$,
and use that $\hat a_1$ and $\jap{\cdot}$ are nondecreasing. Then
\[
\begin{gathered}
\int_0^t\hat a_1^{i+1}(\varepsilon_W\jap{v}^{5})^q\dd v\le\frac{\hat a_1(t)^{i+1}(\varepsilon_W\jap{t}^{5})^q}{(i+1)\sigma},\\
\int_0^t\varepsilon_W\jap{v}^3\hat a_1^i(\varepsilon_W\jap{v}^{5})^q\dd v\le\frac{\hat a_1(t)^i(\varepsilon_W\jap{t}^{5})^{q+1}}{i\sigma}\quad(i\ge1),\\
\int_0^t\varepsilon_W^{j-1}\jap{v}^{3+5(j-2)}\dd v\le(\varepsilon_W\jap{t}^{5})^{j-1},
\end{gathered}
\]
where the second line uses $\jap{t}^3\le\jap{t}^{5}$ and the third uses
$t\le\jap{t}$. Each right side is a monomial of degree $j-1$ in
$\hat a_1(t)$ and $\varepsilon_W\jap{t}^{5}$, hence at most
$\max(1,\sigb^{-1})\mathfrak b(t)^{j-1}$, and there are $2\cdot2^{j-2}$ terms counted with
the binomial weights. This proves Equation~\eqref{eq:W-hier}.

The residual is given by Equation~\eqref{eq:hier-res} with $J=N_G$. Since
$\mathfrak b\le1$ and $\hat a_1\le\mathfrak b$, we have
$\hat a_1^2\mathfrak b^{i+l-2}\le\hat a_1\mathfrak b^{N_G}$ for $i+l>N_G$, and
$\|h\|_{C^s(N')}\le\varepsilon_W\jap{t}^3+\hat a_1\le2\mathfrak b$. So
\begin{equation}
\label{eq:W-res}
\|R_{\mathrm{app}}(t)\|_{H^s}\lesssim\hat a_1(t)\mathfrak b(t)^{N_G},\qquad
\|w_{\mathrm{app}}(t)\|_{H^{s+1}}\lesssim\hat a_1(t)\le\delta_0 .
\end{equation}
Our choice of $C_W$ gives $\|w_{\mathrm{app}}\|_{H^{s+1}}\le C_W\delta_0\le1$.

\emph{Continuity argument.} Let $t_*$ be the largest $t\in[0,\bar t]$ with
$\delta_{\mathrm r}(t)\le a/16$ and $\|w(v)\|_{H^s}\le1$ for $v\in[0,t]$; at $t=0$
both hold, since $\delta_{\mathrm r}(0)=0$ and, as $s\le M-2$,
$\|w(0)\|_{H^s}\le C_\natural\kappa_1\le a/64$ (recall $C_W\ge C_\natural$, $\varsigma\ge1$
and the choice of $\kappa_W$). For $t\le t_*$ the
support of $\omega_{\mathrm r}(t)$ lies in
$\{\frac{51}{64}a\le|y|\le\frac54+\frac5{64}a\}\subset\Nflat\subset N$
(Lemma~\ref{lem:annuli}). So
Lemma~\ref{lem:reduction}\textup{(3)} and Lemma~\ref{lem:energy} apply on $[0,t_*]$,
by Equations~\eqref{eq:W-small} and \eqref{eq:W-res}. As $v=w-w_{\mathrm{app}}$
vanishes at $t=0$, and
$\mathfrak b^{N_G}\le2^{N_G-1}\big(\hat a_1^{N_G}+(\varepsilon_W\jap{\cdot}^{5})^{N_G}\big)$,
\[
\begin{aligned}
\|v(t)\|_{H^s}\lesssim\Big(&\int_0^te^{\Lam(t-v)}\hat a_1(v)^{N_G+1}\dd v\\
&+(\varepsilon_W\jap{t}^{5})^{N_G}\int_0^te^{\Lam(t-v)}\hat a_1(v)\dd v\Big).
\end{aligned}
\]
By Equation~\eqref{eq:W-exp}, the first integral is at most
$\hat a_1(t)^{N_G+1}/((N_G+1)\sigma-\Lam)\le\hat a_1(t)^{N_G+1}/(2\sigb)$, and the second is
at most $\kappa_1e^{\Lam t}/(\Lam-\sigma)\le\kappa_1e^{\Lam t}/\sigb$. For $t\le\hat T_1$,
$\kappa_1e^{\Lam t}\le\kappa_1^{1-\Lam/\sigma}\delta_0^{\Lam/\sigma}$. Hence, for $t\le t_*$,
\begin{equation}
\label{eq:W-v}
\|v(t)\|_{H^s}\lesssim\big(\hat a_1(t)^{N_G+1}+\Theta_W(t)\big),\qquad
\Theta_W(t)=(\varepsilon_W\jap{t}^{5})^{N_G}\kappa_1^{1-\Lam/\sigma}\delta_0^{\Lam/\sigma},
\end{equation}
and, by Equation~\eqref{eq:W-exp}, $\kappa_1\le1$ and $\delta_0\le1$,
\begin{equation}
\label{eq:W-theta}
\Theta_W(t)\le C_\natural^{N_G}\varsigma(\kappa_1)^{5N_G}\kappa_1^2\le C_\natural^{N_G}\kappa_1\varsigma(\kappa_1)^{N_W}\cdot\kappa_1 .
\end{equation}

We collect the bounds for $t\le t_*$. The leading approximation is the
unstable part of the linear response $w^{[1]}=\wlin_0$:
Lemma~\ref{lem:running} gives $\Riesz w^{[1]}=2\Re(ce^{\lambda t}\phist)$.
Splitting $w^{[1]}=\Riesz w^{[1]}+(1-\Riesz)w^{[1]}$ in
$w=\sum_{j=1}^{N_G}w^{[j]}+v$ therefore yields
\[
w-2\Re(ce^{\lambda t}\phist)
=(1-\Riesz)w^{[1]}+\sum_{j=2}^{N_G}w^{[j]}+v.
\]
The higher-order terms $w^{[j]}$, $j\ge2$, retain both their unstable and
complementary components; their full $H^s$ norms are already small by
Equation~\eqref{eq:W-hier}, so both components are included in the error.
More precisely, $\mathfrak b\le1$ gives $\mathfrak b^{j-1}\le\mathfrak b$ for
$j\ge2$, so their sum is bounded by a constant times $\hat a_1\mathfrak b$.
Substituting $\mathfrak b=\hat a_1+C_\natural\kappa_1\jap{t}^5$ and
$\jap{t}\le\varsigma(\kappa_1)$ produces the term
$C_\natural\kappa_1\varsigma(\kappa_1)^5$ after factoring out $\hat a_1$.
By the bound on
$(1-\Riesz)w^{[1]}$ above, Equation~\eqref{eq:W-hier} with $\mathfrak b\le1$,
Equations~\eqref{eq:W-v} and \eqref{eq:W-theta}, $\hat a_1^{N_G+1}\le\hat a_1^2$ and $\Theta_W\le\hat a_1\Theta_W/\kappa_1$,
\[
\begin{aligned}
\big\|w-2\Re(ce^{\lambda t}\phist)\big\|_{H^s}\lesssim{}&\hat a_1(C_\natural\rho_1+\eta_b)\\
&+\hat a_1\big(2\hat a_1+C_\natural\kappa_1\varsigma(\kappa_1)^{5}+C_\natural^{N_G}\kappa_1\varsigma(\kappa_1)^{N_W}\big),
\end{aligned}
\]
which, after including the implicit constant in the choice of $C_W$, is at most
$C_W\hat a_1\big[\rho_1+\hat a_1+\kappa_1\varsigma(\kappa_1)^{N_W}+\eta_b\big]$.
Similarly,
$\|w(t)\|_{H^s}\le C\big(\hat a_1(t)+\Theta_W(t)\big)\le C_W\big(\delta_0+\kappa_1\varsigma(\kappa_1)^{N_W}\big)\le\tfrac a{32}$.
For the displacement, Lemma~\ref{lem:reduction}\textup{(2)} and
Lemma~\ref{lem:bilinear}\textup{(4)} give, for $t\le t_*$,
\[
\begin{aligned}
\delta_{\mathrm r}(t)&\le\int_0^t\big(\|h\|_{C^0(N')}+C_{\mathrm{BS}}\|w\|_{H^s}\big)\\
&\lesssim C_\natural\kappa_1\varsigma(\kappa_1)^{4}+\frac{\delta_0}{\sigb}+C_{\mathrm{BS}}\Big(\frac{\delta_0}{\sigb}+C_\natural^{N_G}\varsigma(\kappa_1)^{N_W}\kappa_1^2\Big),
\end{aligned}
\]
where we used $\|w\|_{H^s}\le\sum_j\|w^{[j]}\|_{H^s}+\|v\|_{H^s}\lesssim\hat a_1+\Theta_W$,
$\int_0^t\hat a_1\le\hat a_1(t)/\sigma\le\delta_0/\sigb$, and
$\int_0^t\Theta_W\le\bar t\,C_\natural^{N_G}\varsigma(\kappa_1)^{5N_G}\kappa_1^2$ with
$\bar t\le\varsigma(\kappa_1)$. After including the implicit constant, our choice
of $C_W$ bounds the right side by
$C_W\delta_0+C_W\kappa_1\varsigma(\kappa_1)^{N_W}\le a/32$. If $t_*<\bar t$, continuity
would give a larger $t$ with $\delta_{\mathrm r}(t)\le a/16$ and $\|w\|_{H^s}\le1$ on $[0,t]$, a
contradiction. So $t_*=\bar t$, and the bounds just proved hold on $[0,\bar t]$.
\end{proof}

\begin{theorem}
\label{thm:active}
Assume \textup{(N1)}--\textup{(N3)}. Fix an integer $s\ge3$ and a number
$\sigma'\in(\frac9{200},\sigb)$. Let $N_G$ be the least integer
with $N_G\ge2$ and $N_G\ge\Lam(s)/\sigb+1$, and put $M=s+N_G+2$ (so $M\ge s+4$).
There are $\dH\in(0,1)$, $C_2\ge1$ and $C_3>0$, depending only on
$(\Phi,\Gone,s,\sigma')$, and for every $C_0\ge1$ a number $\kappa_0(C_0)>0$,
depending only on $C_0$ and $(\Phi,\Gone,s,\sigma')$, such that the following holds.

Let $G\in[0,\Gone]$, $\sigma=\sigma(G)$, $\lambda=\lamst(G)$, $\kappa>0$,
$\eta_b\in[0,1]$, $\delta_0\in(0,\dH]$, $\hat a(t)=\kappa e^{\sigma t}$ and
$\hat T=\sigma^{-1}\log(\delta_0/\kappa)$. Let $t_1=t_1(C_0,G)$ be the least
$t\ge0$ such that $2^{3/2}C_0\jap{t'}^3e^{-(\sigma-\sigma')t'}\le1$ for all $t'\ge t$,
and define the \textbf{transient}
\[
\rhohat(t)=C_0\jap{t}^3e^{-(\sigma-\sigma')t}\ (0\le t<t_1),\qquad
\rhohat(t)=\jap{t-t_1}^3e^{-(\sigma-\sigma')(t-t_1)}\ (t\ge t_1).
\]
Let $\bar t>0$, let $\omega$, $\omega_{\mathrm r},\omega_{\mathrm{inn}},\omega_{\mathrm{ou}}$
and $h$ be as in Lemma~\ref{lem:reduction} on $[0,\bar t]$ with this $G$, put
$w(t)=\omega_{\mathrm r}(t)-\Phi$, and let $c$ be the running coefficient of
Equation~\eqref{eq:running} with $w_0=w(0)$. Assume:
\begin{description}[font=\normalfont]
\item[\textup{(Init)}] $\|w(0)\|_{H^{M-1}}\le C_0\kappa$ and $\supp(\Phi+w(0))\subset\NNm$;
\item[\textup{(Field)}] $h=h_1+h_2$ on $[0,\bar t]$ with
  $\|h_1(t)\|_{C^M(N')}\le C_0\kappa\jap{t}^3$ and
  $\|h_2(t)\|_{C^M(N')}\le\eta_b\,\hat a(t)$;
\item[\textup{(Coef)}] $|c(t)|\le\kappa\big(2+C_0\jap{t}^3e^{-\sigma t}\big)$ on $[0,\bar t]$;
\item[\textup{(Small)}] $\kappa\le\kappa_0(C_0)$.
\end{description}
Then for every $t\in[0,\min(\bar t,\hat T)]$:
\begin{enumerate}[label=\textup{(\alph*)}]
\item $\supp(\Phi+w(t))\subset\NN$;
\item $\big\|w(t)-2\Re\big(c(t)e^{\lambda t}\phist\big)\big\|_{H^s}\le
  C_2\hat a(t)\big[\rhohat(t)+\hat a(t)+\kappa^{1/2}+\eta_b\big]$ and
  $\|w(t)\|_{H^s}\le C_2\hat a(t)(1+\rhohat(t))$.
\end{enumerate}
Moreover, let $\sigma_{\max}=\max_{G\in[0,\Gone]}\sigma(G)$, and let $\bar t_1(C_0)$ be the
least $t\ge0$ such that
\[
2^{3/2}C_0\jap{t'}^3e^{-(\sigb-\sigma')t'}\le1\qquad\text{for all }t'\ge t .
\]
Then $t_1(C_0,G)\le\bar t_1(C_0)$, and for all $t\ge0$,
\begin{equation}
\label{eq:active-transient}
\hat a(t)\rhohat(t)\le C(C_0)\,\kappa\jap{t}^3e^{\sigma't},\qquad
C(C_0)=\max\big(C_0,e^{\sigma_{\max}\bar t_1(C_0)}\big);
\end{equation}
the number $C(C_0)$ depends only on $C_0$, $\sigma'$, $\sigb$ and $\sigma_{\max}$.
\end{theorem}
\begin{proof}
Throughout,
$\Lin=\Lin(G)$, $F=-h\cdot\nabla\Phi$, and $\wlin$ is the function $\wlin_0$ of
Lemma~\ref{lem:running} with $w_0=w(0)$; it is an admissible choice of $w^{[1]}$
in Equation~\eqref{eq:hier}, because $\supp(\Phi+w(0))\subset\NNm\subset N$ by
(Init). The proof has two stages. We call $[0,\min(\bar t,t_1)]$ the
\emph{initial layer}: during this initial interval, the contribution of the
initial data and polynomial driving may still be large relative to
$\kappa e^{\sigma t}$. By the definition of $t_1$, the factor
$C_0\jap{t}^3e^{-(\sigma-\sigma')t}$ controlling this relative size is at most
$2^{-3/2}$ for every $t\ge t_1$.
The initial interval has length at most $\bar t_1(C_0)$, independently of
$\kappa$, and we control it by comparing $w$ with $\wlin$ using
Lemma~\ref{lem:energy}. At $t_1$ the relative contribution is small enough
to restart the problem and apply Lemma~\ref{lem:window} on the remaining interval.

\emph{Constants.} Fix a constant $C\ge1$, depending only on
$(\Phi,\Gone,s,\sigma')$, that bounds the semigroup, eigenfunction and
linear-response estimates used below. Let $\|\ellst\|$ and $\Cp$ be as in
Lemma~I.\ref{I-lem:uniform}. Choose $C_\natural\ge4$ sufficiently large for
the restarted initial-data bound below, in particular $C_\natural\ge8C+1$.
Let $\delta_W,\kappa_W,C_W$ be the constants of Lemma~\ref{lem:window} for
this $C_\natural$, and choose
\[
C_2\ge1,\qquad C_3\ge(1+\Cp)C_2,\qquad
0<\dH\le\delta_W,\quad C_2\dH\le a/32,
\]
with $C_2$ large enough for the norm and displacement estimates below.
These depend only on $(\Phi,\Gone,s,\sigma')$. Now let $C_0\ge1$, write
$\bar t_1=\bar t_1(C_0)$, $\Lambda'=\Lam(M-1)$ and $\Lambda''=\Lam(M-2)$.
On the fixed interval $[0,\bar t_1]$, the estimates below bound
$\|h\|_{C^M}/\kappa$, $\|w\|_{H^{M-1}}/\kappa$,
$\delta_{\mathrm r}/\kappa$, $\|\wlin\|_{H^{M-1}}/\kappa$ and
$\|w-\wlin\|_{H^{M-2}}/\kappa^2$ by constants depending only on $C_0$ and
$(\Phi,\Gone,s,\sigma')$. The implicit constants in these initial-layer
estimates may also depend on $C_0$ and are obtained independently of the
constant now chosen. Choose one fixed $C_{\mathrm{lay}}=C_{\mathrm{lay}}(C_0)\ge C_0+1$ large enough to
dominate all these coefficients.
The smallness of $\kappa$ is imposed by three groups of conditions: (L1) is used on
the initial layer, (L2) to bound the terms that contain the error $v_1=w-\wlin$ of
the initial layer, and (L3) to apply Lemma~\ref{lem:window} after the restart.
Choose $\kappa_0(C_0)\in(0,(a/32)^2]$ sufficiently small that every
$\kappa\in(0,\kappa_0(C_0)]$ satisfies:
\begin{enumerate}[label=\textup{(L\arabic*)}]
\item $C_{\mathrm{lay}}\kappa\le\min(\frac12,\frac a{128},\frac12\kappa^{1/2})$;
\item $(1+\|\ellst\|)C_{\mathrm{lay}}\kappa\le1$ and
  $4(1+C\|\ellst\|)C_{\mathrm{lay}}\kappa^{1/2}\le1$;
\item $e^{\sigma_{\max}\bar t_1}\kappa\le\kappa_W$ and
  $C_We^{\sigma_{\max}\bar t_1}\kappa\,\varsigma(\kappa)^{N_W}\le\frac12\kappa^{1/2}$.
\end{enumerate}
Such a choice is possible because in each inequality the ratio of the left
side to the right side tends to zero as $\kappa\to0^+$. Thus $\kappa_0(C_0)$
depends only on $C_0$ and $(\Phi,\Gone,s,\sigma')$, and (Small) gives
(L1)--(L3) and $\kappa\le(a/32)^2$.

\emph{The time $t_1$.} Since $\sigma\ge\sigb$,
$2^{3/2}C_0\jap{t'}^3e^{-(\sigma-\sigma')t'}\le2^{3/2}C_0\jap{t'}^3e^{-(\sigb-\sigma')t'}$, so
$t_1\le\bar t_1$. For $t\ge t_1$, the definition of $t_1$ with $t'=t$ gives
\begin{equation}
\label{eq:A-t1}
C_0\jap{t}^3e^{-\sigma t}\le2^{-3/2}e^{-\sigma't}\le1,\qquad\text{hence }|c(t)|\le3\kappa\text{ by (Coef)}.
\end{equation}
For $t<t_1$,
$\hat a(t)\rhohat(t)=C_0\kappa\jap{t}^3e^{\sigma't}$. For $t\ge t_1$,
$\hat a(t)\rhohat(t)=\kappa e^{\sigma t_1}\jap{t-t_1}^3e^{\sigma'(t-t_1)}\le e^{\sigma_{\max}\bar t_1}\kappa\jap{t}^3e^{\sigma't}$.
This proves the last assertion of the theorem.

\emph{The initial layer $[0,t_{\mathrm{in}}]$, $t_{\mathrm{in}}=\min(\bar t,t_1)$.}
Here $t\le t_1\le\bar t_1$, $\sigma\le\sigma_{\max}$ and $\eta_b\le1$.
Thus the triangle inequality and (Field) give
\[
\begin{aligned}
\|h(t)\|_{C^M(N')}
&\le C_0\kappa\jap{t}^3+\eta_b\kappa e^{\sigma t}\\
&\le\kappa\big(C_0\jap{t}^3+e^{\sigma_{\max}t}\big)
\le C_{\mathrm{lay}}\kappa\le1,
\end{aligned}
\]
by the choice of $C_{\mathrm{lay}}$ and (L1).

Let $t_*$ be the largest $t\in[0,t_{\mathrm{in}}]$ such that
$\delta_{\mathrm r}(t)\le a/64$ and $\|w(v)\|_{H^{M-1}}\le1$ for $v\in[0,t]$.
At $t=0$, $\delta_{\mathrm r}(0)=0$, and (Init) gives
$\|w(0)\|_{H^{M-1}}\le C_0\kappa<C_{\mathrm{lay}}\kappa\le1$:
indeed, $C_{\mathrm{lay}}\ge C_0+1$, and (L1) gives $C_{\mathrm{lay}}\kappa\le1/2<1$. Both continuity conditions therefore hold
strictly at $t=0$.
For $t\le t_*$ the support of $\omega_{\mathrm r}(t)$ lies in
$\{(\frac78-\frac1{64})a\le|y|\le\frac54+\frac a{64}\}\subset N$
(Lemma~\ref{lem:annuli}), so Lemma~\ref{lem:reduction}\textup{(3)} applies on $[0,t_*]$.

We first obtain a bound of order $\kappa$ for $w$ itself. Taking
$w_{\mathrm{app}}=0$ in Lemma~\ref{lem:energy} makes its error equal to $w$
and its residual equal to $h\cdot\nabla\Phi$. At level $M-1$, the lemma's
hypotheses hold by the continuity assumption on $w$ and the bound on
$\|h\|_{C^M(N')}$ above. Lemma~\ref{lem:bilinear}\textup{(3)} gives
$\|R_{\mathrm{app}}\|_{H^{M-1}}\lesssim\kappa$, with a coefficient depending only on $C_0$ and the fixed data.
Consequently,
\[
\begin{aligned}
\|w(t)\|_{H^{M-1}}
&\lesssim e^{\Lambda't}C_0\kappa
 +\kappa\int_0^te^{\Lambda'(t-v)}\dd v\\
&\lesssim e^{\Lambda'\bar t_1}(C_0+\bar t_1)\kappa
 \le C_{\mathrm{lay}}\kappa\le\tfrac12.
\end{aligned}
\]
The displayed coefficient depends only on $C_0$ and the fixed data and is
included in the choice of $C_{\mathrm{lay}}$; the last inequality is (L1).
Likewise, Lemma~\ref{lem:reduction}\textup{(2)} and
Lemma~\ref{lem:bilinear}\textup{(4)} yield
\[
\delta_{\mathrm r}(t)
\le\int_0^t\big(\|h\|_{C^0(N')}+C_{\mathrm{BS}}\|w\|_{H^1}\big)
\lesssim\bar t_1\kappa\le C_{\mathrm{lay}}\kappa\le a/128.
\]
Here the coefficient follows from the bounds on $h$ and $w$ just proved,
and the last inequality is again (L1).
These strict improvements of both continuity bounds imply $t_*=t_{\mathrm{in}}$.

We next compare $w$ with its linear response. Lemma~\ref{lem:duhamel} gives
$\wlin\in C([0,\bar t];X^{M-1})\cap C^1([0,\bar t];X^{M-2})$, real-valued,
with $\partial_t\wlin=\Lin\wlin+F$ and $\wlin(0)=w(0)$.
In its Duhamel formula, Proposition~I.\ref{I-prop:semigroup} bounds the
semigroup by $Ce^{\sigma t}$ at level $M-1$, while
Lemma~\ref{lem:bilinear}\textup{(3)} gives
$\|F(v)\|_{H^{M-1}}\lesssim\kappa$ with the same allowed dependencies.
Thus, for $t\le t_{\mathrm{in}}$,
\[
\begin{aligned}
\|\wlin(t)\|_{H^{M-1}}
&\lesssim e^{\sigma t}C_0\kappa
 +\kappa\int_0^te^{\sigma(t-v)}\dd v\\
&\lesssim e^{\sigma_{\max}\bar t_1}(C_0+\bar t_1)\kappa
 \le C_{\mathrm{lay}}\kappa\le1.
\end{aligned}
\]
This bound supplies the $H^{M-1}$ hypothesis on the approximate solution
when applying Lemma~\ref{lem:energy} at level $M-2$ with
$w_{\mathrm{app}}=\wlin$. The required bounds on $w$ and $h$ follow from
those already proved, since $H^{M-1}\subset H^{M-2}$ and
$\|h\|_{C^{M-1}(N')}\le\|h\|_{C^M(N')}\le1$.

By Equation~\eqref{eq:hier-res} with $J=1$, the residual is
$R_{\mathrm{app}}=-B(\wlin,\wlin)+h\cdot\nabla\wlin$.
The bilinear estimates of Lemma~\ref{lem:bilinear} give
\[
\begin{aligned}
\|R_{\mathrm{app}}\|_{H^{M-2}}
&\lesssim
 \big(\|\wlin\|_{H^{M-2}}+\|h\|_{C^{M-2}(N')}\big)
 \|\wlin\|_{H^{M-1}}\\
&\lesssim\kappa^2.
\end{aligned}
\]
Here we use the preceding order-$\kappa$ bounds with their coefficients
depending only on $C_0$ and the fixed data, before bounding them by $C_{\mathrm{lay}}$.
Since $v_1=w-\wlin$ has $v_1(0)=0$, Lemma~\ref{lem:energy} now yields
\[
\begin{aligned}
\|v_1(t)\|_{H^{M-2}}
&\lesssim\kappa^2\int_0^te^{\Lambda''(t-v)}\dd v\\
&\lesssim\bar t_1e^{\Lambda''\bar t_1}\kappa^2
 \le C_{\mathrm{lay}}\kappa^2,
\end{aligned}
\]
where the constants in the preceding bounds depend only on $C_0$ and the
fixed data, and $C_{\mathrm{lay}}$ dominates the final coefficient.
Together with the radial-displacement bound proved above, this gives
\begin{equation}
\label{eq:A-layer}
\|v_1(t)\|_{H^{M-2}}\le C_{\mathrm{lay}}\kappa^2,\qquad \delta_{\mathrm r}(t)\le C_{\mathrm{lay}}\kappa\qquad(0\le t\le t_{\mathrm{in}}).
\end{equation}

\emph{Proof of (b) for $t\le\min(t_{\mathrm{in}},\hat T)$.} By
Lemma~\ref{lem:linear-part} with $k=s$ and $A=C_0\kappa$, and since
$C_0\kappa\jap{t}^3e^{\sigma't}\le\hat a(t)\rhohat(t)$ for $t\le t_1$ (with
equality for $t<t_1$, while at $t=t_1$ the left side is at most
$2^{-3/2}\hat a(t_1)$ by the definition of $t_1$ and $\rhohat(t_1)=1$),
\[
\begin{gathered}
\wlin(t)=2\Re\big(c(t)e^{\lambda t}\phist\big)+(1-\Riesz)\wlin(t),\\
\|(1-\Riesz)\wlin(t)\|_{H^s}\le C\hat a(t)\big(\rhohat(t)+\eta_b\big).
\end{gathered}
\]
By Equation~\eqref{eq:A-layer} and (L2),
$\|v_1(t)\|_{H^s}\le C_{\mathrm{lay}}\kappa^2\le\kappa^{1/2}\kappa\le\kappa^{1/2}\hat a(t)$. Hence
$\|w-2\Re(ce^{\lambda t}\phist)\|_{H^s}\le(C+1)\hat a\big[\rhohat+\kappa^{1/2}+\eta_b\big]$.
By (Coef), $|c(t)|e^{\sigma t}\le2\hat a(t)+C_0\kappa\jap{t}^3\le\hat a(t)(2+\rhohat(t))$,
so $\|w(t)\|_{H^s}\le2C(2+\rhohat)\hat a+(C+1)(\rhohat+2)\hat a\le C_2\hat a(1+\rhohat)$.
The displacement is at most $C_{\mathrm{lay}}\kappa\le\frac12\kappa^{1/2}$ by (L1).
This proves (b) when $\min(\bar t,\hat T)\le t_1$.

\emph{Restart.} Assume now that $t_1<t_e:=\min(\bar t,\hat T)$. The function
$\tilde\omega(t')=\omega(t_1+t')$ is the solution of Theorem~\ref{thm:smooth} with the
$3$-fold rotationally symmetric datum $\omega(t_1)$, by Theorem~\ref{thm:smooth}\textup{(5)}, and its
flow is $\tilde\flow_{t'}=\flow_{t_1+t'}\circ\flow_{t_1}^{-1}$. With
$\tilde\omega_\bullet(0)=\omega_\bullet(t_1)$ we get $\tilde\omega_\bullet(t')=\omega_\bullet(t_1+t')$,
the same $G$ (since $\int\omega_{\mathrm{inn}}(t_1)=2\pi G$), and the field
$\tilde h(t')=h(t_1+t')$. So the restarted configuration satisfies the assumptions of
Lemma~\ref{lem:reduction} on $[0,\bar t-t_1]$, and its perturbation is
$\tilde w(t')=w(t_1+t')$. Its running coefficient is
$\tilde c(t')=\ellst(w(t_1))+\int_0^{t'}e^{-\lambda v}\ellst(F(t_1+v))\dd v$. By
Equation~\eqref{eq:running} and Lemma~\ref{lem:running},
$e^{\lambda t_1}c(t_1+t')=\ellst(\wlin(t_1))+\int_0^{t'}e^{-\lambda v}\ellst(F(t_1+v))\dd v$, so
\begin{equation}
\label{eq:A-restart}
\tilde c(t')=e^{\lambda t_1}c(t_1+t')+\ellst(v_1(t_1)).
\end{equation}

We check the assumptions of Lemma~\ref{lem:window} for the restarted problem,
with $C_\natural$ as above, $\kappa_1=\hat a(t_1)=\kappa e^{\sigma t_1}$, the same
$\delta_0$ and $\eta_b$, and $\bar t'=t_e-t_1$. First,
$\kappa\le\kappa_1\le e^{\sigma_{\max}\bar t_1}\kappa\le\kappa_W$ by (L3),
$\delta_0\le\dH\le\delta_W$, and $\sigma^{-1}\log(\delta_0/\kappa_1)=\hat T-t_1\ge\bar t'$.
For the initial norm at the restart, write
\[
w(t_1)=2\Re\big(c(t_1)e^{\lambda t_1}\phist\big)
 +(1-\Riesz)\wlin(t_1)+v_1(t_1).
\]
Equation~\eqref{eq:A-t1} and (Coef) give $|c(t_1)|\le3\kappa$.
The definition of $t_1$ also gives
\[
C_0\kappa\jap{t_1}^3e^{\sigma't_1}
=\kappa_1 C_0\jap{t_1}^3e^{-(\sigma-\sigma')t_1}
\le2^{-3/2}\kappa_1\le\kappa_1,
\qquad \eta_b\kappa e^{\sigma t_1}\le\kappa_1.
\]
We bound the three terms separately. First,
$\|\phist\|_{H^{M-2}}\le C$, so
\[
\big\|2\Re\big(c(t_1)e^{\lambda t_1}\phist\big)\big\|_{H^{M-2}}
\le2C|c(t_1)|e^{\sigma t_1}\le6C\kappa_1.
\]
Next, Lemma~\ref{lem:linear-part} with $k=M-2$ and $A=C_0\kappa$ yields
\[
\|(1-\Riesz)\wlin(t_1)\|_{H^{M-2}}
\le C\big(C_0\kappa\jap{t_1}^3e^{\sigma't_1}
 +\eta_b\kappa e^{\sigma t_1}\big)
\le2C\kappa_1.
\]
Finally, in the restart case the initial interval ends at $t_1$, so
Equation~\eqref{eq:A-layer} applies there. The first inequality of (L2),
$(1+\|\ellst\|)C_{\mathrm{lay}}\kappa\le1$, implies
$C_{\mathrm{lay}}\kappa^2\le\kappa\le\kappa_1$. Thus
$\|v_1(t_1)\|_{H^{M-2}}\le\kappa_1$. Adding the three bounds gives
\[
\|w(t_1)\|_{H^{M-2}}
\le(8C+1)\kappa_1\le C_\natural\kappa_1,
\]
by our choice of $C_\natural$.
Second, by Equation~\eqref{eq:A-layer} the support of $\Phi+w(t_1)$ lies in
$\{(\frac78-\frac1{64})a\le|y|\le\frac54+\frac a{64}\}$. Third, by Peetre's inequality
$\jap{t_1+t'}\le\sqrt2\jap{t_1}\jap{t'}$ and the definition of $t_1$,
\[
C_0\kappa\jap{t_1+t'}^3\le\big[2^{3/2}C_0\jap{t_1}^3e^{-(\sigma-\sigma')t_1}\big]e^{-\sigma't_1}\kappa_1\jap{t'}^3\le\kappa_1\jap{t'}^3,
\]
and $\eta_b\hat a(t_1+t')=\eta_b\kappa_1e^{\sigma t'}$; so
$\|\tilde h(t')\|_{C^M(N')}\le C_\natural\kappa_1\jap{t'}^3+\eta_b\kappa_1e^{\sigma t'}$.
Fourth, for every $t'\in[0,\bar t']$, $t_1+t'\ge t_1$, so (Coef) and
Equation~\eqref{eq:A-t1} give
\[
|c(t_1+t')|\le\kappa\big(2+C_0\jap{t_1+t'}^3e^{-\sigma(t_1+t')}\big)
\le3\kappa.
\]
Since $\|\ellst\|$ bounds the operator norm on $X^0$,
Equation~\eqref{eq:A-layer} gives
\[
|\ellst(v_1(t_1))|
\le\|\ellst\|\,\|v_1(t_1)\|_{L^2}
\le\|\ellst\|C_{\mathrm{lay}}\kappa^2\le\kappa\le\kappa_1,
\]
where the penultimate inequality again follows from the first inequality
of (L2). Substituting these bounds into Equation~\eqref{eq:A-restart}, and
using $|e^{\lambda t_1}|=e^{\sigma t_1}$ and $C_\natural\ge4$, yields
\[
\begin{aligned}
|\tilde c(t')|
&\le e^{\sigma t_1}|c(t_1+t')|+|\ellst(v_1(t_1))|\\
&\le3\kappa e^{\sigma t_1}+\kappa_1
=4\kappa_1\le C_\natural\kappa_1.
\end{aligned}
\]

\emph{Proof of (b) for $t\in[t_1,t_e]$.} Lemma~\ref{lem:window} applies; put
$t'=t-t_1$. Then $\kappa_1e^{\sigma t'}=\hat a(t)$ and
$\jap{t'}^3e^{-(\sigma-\sigma')t'}=\rhohat(t)$. Since $\varsigma$ is decreasing and
$\kappa_1\ge\kappa$, (L3) gives $C_W\kappa_1\varsigma(\kappa_1)^{N_W}\le\frac12\kappa^{1/2}$.
By Equation~\eqref{eq:A-restart},
$\tilde c(t')e^{\lambda t'}=c(t)e^{\lambda t}+\ellst(v_1(t_1))e^{\lambda t'}$, and
$|\ellst(v_1(t_1))e^{\lambda t'}|\le\|\ellst\|C_{\mathrm{lay}}\kappa^2e^{\sigma t'}\le\|\ellst\|C_{\mathrm{lay}}\kappa\,\hat a(t)$.
Since $\|2\Re(X\phist)\|_{H^s}\le2C|X|$ for $X\in\C$, (L2) gives
\[
\begin{aligned}
\big\|w(t)-2\Re\big(c(t)e^{\lambda t}\phist\big)\big\|_{H^s}
&\le C_W\hat a(t)\big[\rhohat(t)+\hat a(t)+\tfrac12\kappa^{1/2}+\eta_b\big]\\
&\quad+\tfrac12\kappa^{1/2}\hat a(t),
\end{aligned}
\]
which gives the first bound of (b), since $C_W+1\le C_2$. By
Equation~\eqref{eq:A-t1}, the main term $2\Re(c(t)e^{\lambda t}\phist)$ has $H^s$ norm
at most $6C\hat a(t)$. Since $\hat a+\kappa^{1/2}+\eta_b\le3$, we get
\[
\|w(t)\|_{H^s}\le6C\hat a+(C_W+1)\hat a(\rhohat+3)\le C_2\hat a(1+\rhohat).
\]
For the displacement, a trajectory of $\omega_{\mathrm r}$ satisfies, with $z=\flow_{t_1}(y)$,
\[
|\flow_t(y)|-|y|=\big(|\tilde\flow_{t'}(z)|-|z|\big)+\big(|\flow_{t_1}(y)|-|y|\big),
\]
so by Lemma~\ref{lem:window} and Equation~\eqref{eq:A-layer} it moves radially by at
most $C_W\delta_0+\frac12\kappa^{1/2}+C_{\mathrm{lay}}\kappa\le C_2\delta_0+\kappa^{1/2}$.

\emph{Proof of (a).} In both cases the displacement is at most
$C_2\delta_0+\kappa^{1/2}\le\frac a{32}+\frac a{32}=\frac a{16}$, because
$\delta_0\le\dH\le a/(32C_2)$ and $\kappa\le(a/32)^2$. Since
$\supp\omega_{\mathrm r}(0)\subset\NNm$ by (Init), Lemma~\ref{lem:annuli} gives
$\supp(\Phi+w(t))\subset\nobreak\NN$.
\end{proof}

Recall from the expansion after Definition~I.\ref{I-def:multipole} that the outer
moment $\pout_3$ of a function determines the part in the angular modes $\pm3$ of its
field outside its support. So the outer moment of the deviation of a ring
determines the part in these modes of the field that the ring exerts on the rings
outside it. Recall also from Definition~I.\ref{I-def:eigendata} and the discussion
before it that $\pout_*(G)=\pout_3[\phist(G)]$ is the outer moment of the unstable
eigenfunction.
The next corollary gives this moment for an active ring.

\begin{corollary}
\label{cor:active-moment}
Under the assumptions of Theorem~\ref{thm:active}, put
$q(t)=2\pi\,\pout_3[w(t)]$. Then for every $t\in[0,\min(\bar t,\hat T)]$,
\begin{equation}
\label{eq:active-moment}
\big|q(t)-2\pi\pout_*(G)\,c(t)e^{\lambda t}\big|\le2\pi C_3\,\hat a(t)\big[\rhohat(t)+\hat a(t)+\kappa^{1/2}+\eta_b\big].
\end{equation}
\end{corollary}
\begin{proof}
For $X\in\C$, $\pout_3[2\Re(X\phist)]=X\pout_3[\phist]+\bar X\pout_3[\overline{\phist}]$. The
second term vanishes: $\overline{\phist}$ lies in the angular mode $-3$, and
$\int_0^{2\pi}e^{-3i\theta}e^{-3i\theta}\dd\theta=0$. Since
$\pout_*(G)=\pout_3[\phist(G)]$, linearity gives
$\pout_3[w(t)]-\pout_*(G)c(t)e^{\lambda t}=\pout_3\big[w(t)-2\Re(c(t)e^{\lambda t}\phist)\big]$.
The function in brackets is supported in $N$, so Lemma~I.\ref{I-lem:uniform} and
Theorem~\ref{thm:active}\textup{(b)} bound its modulus by
$\Cp C_2\hat a(t)\big[\rhohat(t)+\hat a(t)+\kappa^{1/2}+\eta_b\big]$. Multiplying by $2\pi$
and using $\Cp C_2\le C_3$ proves the corollary.
\end{proof}

\subsection{The Seeded Ring}

The approximate datum of Definition~\ref{def:approx} replaces its innermost ring by
a rearrangement $\Phi^\kappa$ of $\Phi$ whose deviation from $\Phi$ has an unstable
component of size about $\kappa$. We obtain $\Phi^\kappa$ by composing $\Phi$ with
the time-$\kappa$ flow of a divergence-free field $V$. This flow preserves area, so
$\Phi^\kappa$ is equimeasurable with $\Phi$, and $V$ is chosen with
$V\cdot\nabla\Phi=2\Re\phist(0)$, so that $\Phi^\kappa-\Phi$ is $2\kappa\Re\phist(0)$ to
first order in $\kappa$.

\begin{lemma}
\label{lem:seed}
Assume \textup{(N1)}--\textup{(N3)}, let $s\ge3$ be an integer, and let
$M=s+N_G+2$ be the integer of Theorem~\ref{thm:active} for this $s$. There are a
number $\kL>0$ and functions $\Phi^\kappa$, $\kappa\in(0,\kL]$, all depending only
on $\Phi$, and a number $\Cseed<\infty$, depending only on $\Phi$, $\Gone$ and $s$,
such that for every
$\kappa\in(0,\kL]$: $\Phi^\kappa\in C_c^\infty(\R^2)$, $\Phi^\kappa\ge0$,
$\Phi^\kappa$ is $3$-fold rotationally symmetric and equimeasurable with $\Phi$;
$\supp\Phi^\kappa\subset\NNm$;
$\|\Phi^\kappa-\Phi\|_{H^{M-1}}\le\Cseed\kappa$;
and $|\ellst(0)(\Phi^\kappa-\Phi)|\in[\frac\kappa2,2\kappa]$.
\end{lemma}
\begin{proof}
Let $G=0$, $\lambda=\lamst(0)$, and let $\psi_*$ and $n_*=n_*(0)$ be
as in Theorem~I.\ref{I-thm:evans-simple} and Definition~I.\ref{I-def:eigendata}, for the
ring in ring units and $\ell=3$. By these definitions,
$\phist(0)=n_*^{-1}g_*(r)e^{3i\theta}$ with $g_*=\Mres_{\lambda}\kcoef_3\,\psi_*(\log r)$,
where $\kcoef_3(r)=3i\Phi'(r)/r$ and $\Mres_\lambda=(\lambda+3i\Omega(\cdot;0))^{-1}$
(Definition~I.\ref{I-def:sector-setting}). So the mode-$3$ profile of
$\phist(0)$ is
\[
\vphist(r)=n_*^{-1}\,\frac{3i\,\Phi'(r)}{r}\,\frac{\psi_*(\log r)}{\lambda+3i\Omega(r;0)} .
\]
Put $g_{\mathrm s}(r)=n_*^{-1}\psi_*(\log r)/(\lambda+3i\Omega(r;0))$ for $r>0$. It is
smooth: $\psi_*$ solves a linear equation with smooth coefficients
(Lemma~I.\ref{I-lem:sector}\textup{(3)}), $\Omega(\cdot;0)$ is smooth on $(0,\infty)$, and
$|\lambda+3i\Omega|\ge\Re\lambda>0$. Let $\chiseed\in C_c^\infty(\{a/2<|y|<2\})$ be
radial with $\chiseed=1$ on $\{a\le|y|\le1\}$, and put
$\psi_{\mathrm{seed}}=-2\Re\big(\chiseed(r)\,g_{\mathrm s}(r)\,e^{3i\theta}\big)$ and
$V=\nabla^\perp\psi_{\mathrm{seed}}$. Then $\psi_{\mathrm{seed}}$ is real, smooth,
$3$-fold rotationally symmetric and supported in the compact set
$K_V=\supp\chiseed\subset\{a/2<|y|<2\}$, and $V\in C_c^\infty(\R^2;\R^2)$ is
divergence-free and $3$-fold rotationally symmetric. Let $Y_\kappa=\flow_V(\kappa,0,\cdot)$ be the
time-$\kappa$ flow of $V$, with $\flow_V$ as in Lemma~\ref{lem:flow} for the
time-independent bounded field $V$. By Lemma~\ref{lem:flow}\textup{(a)}, \textup{(b)},
\textup{(d)}, $Y$ is defined for all $\kappa\in\R$ and is $C^\infty$ in $(\kappa,y)$,
each $Y_\kappa$ is a $C^\infty$ diffeomorphism of $\R^2$, and $\det DY_\kappa=1$
because $\divop V=0$. Since $V$ does not depend on time, uniqueness for the flow
gives $Y_{\kappa+\kappa'}=Y_\kappa\circ Y_{\kappa'}$ and, since
$V(\Orot y)=\Orot V(y)$, $Y_\kappa\circ\Orot=\Orot\circ Y_\kappa$. For $y\notin K_V$
the constant path $y$ solves $\dot x=V(x)$, since $V(y)=0$; so $Y_\kappa$ is the identity outside
$K_V$.

Put $\Phi^\kappa=\Phi\circ Y_\kappa$. Then $\Phi^\kappa\in C_c^\infty(\R^2)$ is nonnegative and $3$-fold
rotationally symmetric. It is equimeasurable with $\Phi$: for $c>0$,
$\{\Phi^\kappa>c\}=Y_\kappa^{-1}(\{\Phi>c\})$ has the measure of $\{\Phi>c\}$ because
$Y_\kappa$ preserves area, and $\{\Phi^\kappa<-c\}=\{\Phi<-c\}=\emptyset$. Since
$\Phi$ is radial, $V\cdot\nabla\Phi=(V\cdot\er)\Phi'$, and
$V\cdot\er=-r^{-1}\partial_\theta\psi_{\mathrm{seed}}=2\Re\big(\frac{3i}r\chiseed g_{\mathrm s}e^{3i\theta}\big)$.
On $\supp\Phi'\subset\{a\le|y|\le1\}$ we have $\chiseed=1$, so
$(V\cdot\er)\Phi'=2\Re(\vphist e^{3i\theta})=2\Re\phist(0)$; off $\supp\Phi'$ both
sides vanish, because $\vphist$ is a multiple of $\Phi'$. Thus
\begin{equation}
\label{eq:seed-first}
V\cdot\nabla\Phi=2\Re\phist(0).
\end{equation}

We now expand in $\kappa$. Fix $y$ and put $f_y(t)=\Phi(Y_t(y))$. Then
$f_y'(t)=(V\cdot\nabla\Phi)(Y_t(y))$ and $f_y''(t)=(V\cdot\nabla(V\cdot\nabla\Phi))(Y_t(y))$,
so with $\Psi_{\mathrm s}=V\cdot\nabla(V\cdot\nabla\Phi)\in C_c^\infty(\R^2)$,
\[
\Phi^\kappa-\Phi=\int_0^\kappa(V\cdot\nabla\Phi)\circ Y_t\dd t,\qquad
\Phi^\kappa-\Phi-\kappa V\cdot\nabla\Phi=\int_0^\kappa(\kappa-t)\,\Psi_{\mathrm s}\circ Y_t\dd t .
\]
The derivatives of $Y_t(y)-y$ of every order vanish outside $K_V$ and are
continuous on $[-1,1]\times K_V$, hence bounded there. By the chain rule and
$\det DY_t=1$, for every $k$, $\|g\circ Y_t\|_{H^k}\lesssim\|g\|_{H^k}$
with a constant depending only on $\Phi$ and $k$, for
$g\in C_c^\infty(\R^2)$ and $|t|\le1$.
Equation~\eqref{eq:seed-first} and the two identities therefore give constants
$C$, depending only on $\Phi$, and $\Cseed$, depending only on $\Phi$ and $M$,
such that for $0<\kappa\le1$,
$\|\Phi^\kappa-\Phi-2\kappa\Re\phist(0)\|_{L^2}\le C\kappa^2$ and
$\|\Phi^\kappa-\Phi\|_{H^{M-1}}\le\Cseed\kappa$.

Finally choose $\kL\in(0,1]$ sufficiently small that
$\kL\|V\|_{L^\infty}\le a/8$ and
$\kL\|\ellst(0)\|_{(X^0)^*}C\le1/2$,
and let $\kappa\in(0,\kL]$. Since $|Y_\kappa^{-1}(x)-x|\le\kappa\|V\|_{L^\infty}\le a/8$ by
Lemma~\ref{lem:flow}\textup{(b)}, \textup{(c)}, and
$\supp\Phi^\kappa=Y_\kappa^{-1}(\supp\Phi)$ with $\supp\Phi\subset\{a\le|y|\le1\}$, we get
$\supp\Phi^\kappa\subset\{\frac78a\le|y|\le1+\frac a8\}\subset\NNm$. The function
$\Phi^\kappa-\Phi$ is supported in $N$ and $3$-fold rotationally symmetric, so it lies in $X^0$.
By Definition~I.\ref{I-def:eigendata}, $\ellst(0)(2\Re\phist(0))=\ellst(0)(\phist(0))+\ellst(0)(\overline{\phist(0)})=1$.
Hence
$|\ellst(0)(\Phi^\kappa-\Phi)-\kappa|\le\|\ellst(0)\|_{(X^0)^*}C\kappa^2\le\frac\kappa2$, and
$|\ellst(0)(\Phi^\kappa-\Phi)|\in[\frac\kappa2,\frac{3\kappa}2]\subset[\frac\kappa2,2\kappa]$.
\end{proof}

\section{The Datum, the Approximants and the Phase Clock}
\label{sec:datum}

We define the datum $f$ of Theorem~\ref{thm:main} and its approximants
$f^{(n)}$; the solutions $\omega^{(n)}$ with these data form the
$2$-admissible sequence to which Corollary~I.\ref{I-cor:defect} is applied in
the proof of Proposition~\ref{prop:reduction}. We then describe each ring in
its own units and define the phase clock, which fixes activation and
saturation times from $\omega^{(n)}$ itself.

\begin{definition}
\label{def:datum}
Let $A_1>0$, $\mu\in(1,2]$ and $\sbar\in(0,a/16]$, and define $A_k$, $s_k$,
$r_k$ ($k\ge1$) and $f$ by Equation~\eqref{eq:datum}.
\end{definition}

\begin{proposition}
\label{prop:integrability}
In Definition~\ref{def:datum}, the sets $\{ar_k\le|x|\le r_k\}$, $k\ge1$, are
pairwise disjoint.
The function $f$ is radial and nonnegative, $f\in (L^1\cap L^q)(\R^2)$ for every
$q\in[1,\infty)$, and $f\notin L^\infty(\R^2)$. Moreover
$\sum_{k\ge1}A_k^qr_k^2<\infty$ for every $q\in[1,\infty)$.
\end{proposition}
\begin{proof}
Since $\mu>1$ and $\sbar\le a/16$, $s_k\le\sbar/2<a$, so $r_{k+1}<ar_k$ and every
annulus of index $i>k$ lies in the open disk $B_{ar_k}$, because $r_i\le r_{k+1}$;
this disk does not meet the annulus of index $k$, so the annuli are pairwise
disjoint. Moreover
$r_k=\prod_{j<k}s_j=\sbar^{k-1}(2\mu)^{-k(k-1)/2}$.

By Lemma~I.\ref{I-lem:ring-facts}, $\Phi$ is radial and continuous, $0\le\Phi<1$,
$\Phi\not\equiv0$ and $\supp\Phi\subset\{a\le|y|\le1\}$. So the $k$-th term of
the series is radial, nonnegative and supported in the $k$-th annulus; since the
annuli are disjoint, at every point at most one term is nonzero, $f$ is radial and nonnegative, and, by
monotone convergence and the formula for $r_k$, for $q\in[1,\infty)$,
\begin{gather*}
\int f^q\dd x=\sum_{k\ge1}A_k^qr_k^2\|\Phi\|_{L^q}^q,\qquad
\|\Phi\|_{L^q}^q\le|\{a\le|y|\le1\}|\le\pi,\\
A_k^qr_k^2=A_1^q\mu^{q(k-1)}\sbar^{2(k-1)}(2\mu)^{-k(k-1)}=A_1^q\big(\mu^q\sbar^2(2\mu)^{-k}\big)^{k-1}.
\end{gather*}
If $(2\mu)^k\ge2\mu^q$, which holds for all large $k$ because $2\mu>1$, the
bracket is at most $\frac12$, because $\sbar<1$, so
$A_k^qr_k^2\le A_1^q2^{1-k}$ for all large $k$. Hence $\sum_kA_k^qr_k^2<\infty$
and $f\in L^q$. Finally, $f\ge A_k\Phi(\cdot/r_k)$, and, $\Phi$ being continuous,
$A_k\Phi(\cdot/r_k)>\frac12A_k\max\Phi>0$ on a nonempty open set; since $\mu>1$
gives $A_k\to\infty$, $f\notin L^\infty$.
\end{proof}

The approximants keep the first $n$ rings of $f$ and replace ring $n$, the
innermost of them, whose time scale $1/A_n$ is the shortest, by the seeded ring
$\Phi^{\kseed}$ of Lemma~\ref{lem:seed}. This is a nonnegative, $3$-fold
rotationally symmetric function, equimeasurable with $\Phi$, whose coefficient
$\ellst(0)(\Phi^{\kseed}-\Phi)$ along the unstable eigenfunction $\phist(0)$ has
modulus between $\kseed/2$ and $2\kseed$. The relay described in the
Introduction starts from this seed.

\begin{definition}
\label{def:approx}
Assume \textup{(N1)}--\textup{(N3)}, and let $\kL$ and $\Phi^\kappa$ be as in
Lemma~\ref{lem:seed}. For an integer $n\ge2$ and $\kseed\in(0,\kL]$ put
$f^{(n)}=\sum_{k<n}A_k\Phi(\cdot/r_k)+A_n\Phi^{\kseed}(\cdot/r_n)$, let $\omega^{(n)}$ be the solution of Theorem~\ref{thm:smooth} with datum
$f^{(n)}$, let $v^{(n)}=\KBS\ast\omega^{(n)}$, and let $\flow_t$ be the flow of
$v^{(n)}$.
\end{definition}

\begin{lemma}
\label{lem:adm}
In Definition~\ref{def:approx}, $f^{(n)}\in C_c^\infty(\R^2)$,
$f^{(n)}$ and $\omega^{(n)}(t)$ are $3$-fold rotationally symmetric for every
$t$, and, for fixed $\kseed$, $f^{(n)}\to f$ in $(L^1\cap L^q)(\R^2)$ for every
$q\in[1,\infty)$ as $n\to\infty$. Moreover $(f^{(n)},\omega^{(n)})_{n\ge2}$ is a
$2$-admissible sequence for $f$ on $[0,\infty)$ in the sense of
Definition~I.\ref{I-def:adm} (the sequence starts at
$n=2$; the shift $n\mapsto n+1$ gives the indexing of Definition~I.\ref{I-def:adm}).
\end{lemma}
\begin{proof}
By Lemmas~I.\ref{I-lem:ring-facts} and \ref{lem:seed}, $\Phi$ and $\Phi^{\kseed}$ lie
in $C_c^\infty(\R^2)$, are $3$-fold rotationally symmetric ($\Phi$ is radial)
and equimeasurable, so
they have the same $L^q$ norms. Hence $f^{(n)}\in C_c^\infty(\R^2)$ is
$3$-fold rotationally symmetric, and every $\omega^{(n)}(t)$ is $3$-fold
rotationally symmetric, by Theorem~\ref{thm:smooth}\textup{(6)} with the rotation $\Orot$. Let
$q\in[1,\infty)$. In
$f-f^{(n)}=A_n(\Phi-\Phi^{\kseed})(\cdot/r_n)+\sum_{k>n}A_k\Phi(\cdot/r_k)$, the
first term has $L^q$ norm at most $2A_nr_n^{2/q}\|\Phi\|_{L^q}$, by the triangle
inequality and the equality of the $L^q$ norms of $\Phi$ and $\Phi^{\kseed}$; the
terms of the sum have disjoint supports, by Proposition~\ref{prop:integrability}, so
\[
\|f-f^{(n)}\|_{L^q}\le\Big(2A_nr_n^{2/q}
+\Big(\sum_{k>n}A_k^qr_k^2\Big)^{1/q}\Big)\|\Phi\|_{L^q},
\]
and both terms on the right tend to $0$ as $n\to\infty$ because
$\sum_kA_k^qr_k^2<\infty$ (Proposition~\ref{prop:integrability}): $A_n^qr_n^2$ is
a term, and $\sum_{k>n}A_k^qr_k^2$ a tail, of this convergent series. Finally, by
Theorem~\ref{thm:smooth}\textup{(7)}, for $q=2$ and every
$T>0$ the restriction of $\omega^{(n)}$ to $[0,T)$ is a renormalized weak
solution in $X_2(T)$ with datum $f^{(n)}$; so $\omega^{(n)}$ is a renormalized
weak solution on $[0,\infty)$, and with the convergence just proved this is the
$2$-admissibility of Definition~I.\ref{I-def:adm}.
\end{proof}

Each ring is studied in its own spatial and vorticity scale, while the
indexed quantities in the global construction all retain physical time $t$.
Thus ring-$k$ units mean $y=x/r_k$, vorticity divided by $A_k$, and velocity
divided by $A_kr_k$. For a velocity field $V(t,x)$ we write
$V_{[k]}(t,y)=V(t,r_ky)/(A_kr_k)$, suppressing $t$ when it is fixed.

When applying the single-ring results of Section~\ref{sec:single}, we also
rescale time. The standard Euler rescaling, with a physical reference time
$t_*$ and scales $A,r>0$, is
\[
\widetilde\omega(s,y)=A^{-1}\omega^{(n)}(t_*+s/A,ry),\qquad
\widetilde v(s,y)=(Ar)^{-1}v^{(n)}(t_*+s/A,ry).
\]
The homogeneity of $\KBS$ gives $\widetilde v=\KBS\ast\widetilde\omega$, and
\[
\partial_s\widetilde\omega+\widetilde v\cdot\nabla_y\widetilde\omega
=A^{-2}(\partial_t\omega^{(n)}+v^{(n)}\cdot\nabla_x\omega^{(n)})(t_*+s/A,ry)=0.
\]
Hence Theorem~\ref{thm:smooth}\textup{(5)} identifies $\widetilde\omega$ with
the smooth solution whose datum is $A^{-1}\omega^{(n)}(t_*,r\,\cdot)$.
For ring $k$ we take $A=A_k$ and $r=r_k$. A time-dependent ring profile $f_k$
is then written locally as $\widetilde f_k(s)=f_k(t_*+s/A_k)$; this time
substitution is stated explicitly whenever it is used. The spatially
normalized profile $f_k(t)$ itself continues to use physical time.

In ring-$k$ units, the perturbation $w_k$ is the
deviation of ring $k$ from $\Phi$; $2\pi G_k$ is the circulation of the point
vortex by which the linearization $\Lin(G_k)$ of Theorem~I.\ref{I-thm:linear}
represents the rings inside ring $k$; and the external field $h_k$ is the field
of all the other rings on $N'$ minus the field of this point vortex;
its local time rescaling is the field $h$ of Lemma~\ref{lem:reduction}.

The number $\varepsilon_k$ is the small factor in the bounds of
Lemma~\ref{lem:transfer}\textup{(1)}--\textup{(3)} on the field, minus its
monopole, that the material inside ring $k$ induces on ring $k$. The
normalized moment $q_{k+1}$ is the mode-$3$ moment through which this material
drives the unstable coefficient of ring $k$ (Lemma~\ref{lem:coupling}).
For an unseeded quiet ring $k<n$, Theorem~\ref{thm:quiet}\textup{(c)} gives
$\|w_k(t)\|_{H^{s_h}}\le C_Q\maj_k(t)$ under its hypotheses, after rescaling
to ring-$k$ units.

\begin{definition}
\label{def:regions}
Fix $n\ge2$ and the objects of Definition~\ref{def:approx}. For $1\le k\le n$:
\begin{itemize}
\item $f^{(n)}_k$ is the $k$-th term of $f^{(n)}$, that is,
  $f^{(n)}_k=A_k\Phi(\cdot/r_k)$ for $k<n$ and
  $f^{(n)}_n=A_n\Phi^{\kseed}(\cdot/r_n)$;
  \[
  S_k=\supp f^{(n)}_k,\quad S_k(t)=\flow_t(S_k),\quad D_k(t)=\bigcup_{i=k}^nS_i(t),\quad
  D_{n+1}(t)=\emptyset,
  \]
  and $\omega_{\ge k}(t)=\omega^{(n)}(t)\mathbf 1_{D_k(t)}$;
  we set $\omega_{\ge n+1}(t)=0$;
\item $N_k=r_kN$, $N'_k=r_kN'$, $\Nflat_k=r_k\Nflat$, $\NN_k=r_k\NN$,
  $\NNm_k=r_k\NNm$ (the annuli of Section~I.\ref{I-sec:ring-def} and of Section~\ref{sec:II-prelim} dilated by $r_k$);
\item while $S_k(t)\subset N_k$, the \textbf{perturbation} of ring $k$ at the
  time $t$, in ring-$k$ units, is
  $w_k(t,y)=A_k^{-1}\omega^{(n)}(t,r_ky)\mathbf 1_{S_k(t)}(r_ky)-\Phi(y)$;
\item the \textbf{inner circulation} $G_k$ is defined by
  $2\pi G_kA_kr_k^2=\int\omega_{\ge k+1}(0)$; we put
  $\varepsilon_k=\mu s_k^5$, and, if $G_k\le\Gone$, $\sigma_k=\sigma(G_k)$,
  $\lambda_k=\lamst(G_k)$, $x_k=\sigma_k/\mu$ and
  $\beta_k=\beta(G_k)s_k^5/(2\pi)$;
\item for $k<n$, the \textbf{normalized moment} of the inner region is
  \[
  q_{k+1}(t)=\overline{Q_3[\omega_{\ge k+1}(t)]}/(A_{k+1}r_{k+1}^5)
  \]
  (the support of $\omega_{\ge k+1}(t)$ is compact and does not contain $0$, by
  Lemma~\ref{lem:regions-basic}\textup{(2)} below, so $Q_3$ applies);
\item at a time $t$ with $D_{k+1}(t)\subset B_{ar_k/8}$, the \textbf{inner field}
  of ring $k$ at physical time $t$ is
  \[
  \hin_k(t)=\KBS\ast\nu^{\mathrm{in}}_k(t,\cdot)-G_k\frac{y^\perp}{|y|^2}\quad\text{on }N',\qquad
  \nu^{\mathrm{in}}_k(t,y)=A_k^{-1}\omega_{\ge k+1}(t,r_ky);
  \]
  at a time $t$ with $\big(S_1(t)\cup\dots\cup S_{k-1}(t)\big)\cap B_{6r_k}=\emptyset$
  (the rings $i<k$ lie outside $B_{6r_k}$), the \textbf{outer field} is
  $\hout_k(t)=\KBS\ast\nu^{\mathrm{out}}_k(t,\cdot)$ on $N'$, with
  $\nu^{\mathrm{out}}_k(t,y)=A_k^{-1}\omega^{(n)}(t,r_ky)\mathbf 1_{S_1(t)\cup\dots\cup S_{k-1}(t)}(r_ky)$;
  when both are defined, $h_k=\hin_k+\hout_k$ is the \textbf{external field} of
  ring $k$;
\item for an integer $M_h\ge1$ (from Section~\ref{sec:relay} on, $M_h$ is the
  integer of \textup{(SA2)} in Definition~\ref{def:standing}) we put
  $\etin_k(t)=\|\hin_k(t)\|_{C^{M_h}(N')}$,
  $\etout_k(t)=\|\hout_k(t)\|_{C^{M_h}(N')}$ and $\eta_k=\etin_k+\etout_k$,
  as functions of physical time;
\item if $G_k\le\Gone$ and $\eta_k$ is defined on $[0,t]$ (it is then
  continuous there, as we show after Lemma~\ref{lem:regions-basic}), the
  \textbf{Duhamel bound} of ring $k$ is
  \[
  \maj_k(t)=A_k\int_0^t e^{A_k\sigma_k(t-v)}\eta_k(v)\dd v ,
  \]
  and $\maj^{\mathrm{in}}_k$ is defined in the same way from $\etin_k$.
\end{itemize}
\end{definition}

\begin{lemma}
\label{lem:regions-basic}
In Definition~\ref{def:regions} the following hold for every $t\ge0$.
\begin{enumerate}[label=\textup{(\arabic*)},beginpenalty=10000]
\item The sets $S_1,\dots,S_n$ are compact, pairwise disjoint and invariant
  under $\Orot$, and none contains $0$. Moreover
  $S_k\subset\{ar_k\le|x|\le r_k\}$ for $k<n$, $S_n\subset\NNm_n$, and
  $f^{(n)}\mathbf 1_{S_k}=f^{(n)}_k$.
\item $\flow_t\circ\Orot=\Orot\circ\flow_t$ and $\flow_t(0)=0$. The sets
  $S_1(t),\dots,S_n(t)$ are compact, pairwise disjoint and invariant under
  $\Orot$, none contains $0$, and
  $\supp\omega^{(n)}(t)\subset S_1(t)\cup\dots\cup S_n(t)$.
\item For $1\le k\le n$,
  \begin{equation}
  \label{eq:omega-ge}
  \omega^{(n)}(t)\mathbf 1_{S_k(t)}=f^{(n)}_k\circ\flow_t^{-1},\qquad
  \omega_{\ge k}(t)=\Big(\sum_{i=k}^nf^{(n)}_i\Big)\circ\flow_t^{-1}
  =\sum_{i=k}^n\omega^{(n)}(t)\mathbf 1_{S_i(t)},
  \end{equation}
  and $\int\omega^{(n)}(t)\mathbf 1_{S_k(t)}=A_kr_k^2\|\Phi\|_{L^1}$. These
  functions lie in $C_c^\infty(\R^2)$, are nonnegative and are $3$-fold
  rotationally symmetric.
\item If $S_k(t)\subset N_k$, then $w_k(t)\in C_c^\infty(\R^2)$ is real valued,
  $3$-fold rotationally symmetric and supported in $N$; so $w_k(t)\in X^s$ for every integer
  $s\ge0$. Moreover $\int w_k(t)=0$, $\|w_k(t)\|_{L^1}\le2\|\Phi\|_{L^1}$, and
  $\|w_k(t)\|_{L^1}\le|N|^{1/2}\|w_k(t)\|_{L^2}$.
\end{enumerate}
\end{lemma}

\begin{proof}
(1) For $k<n$, $f^{(n)}_k$ is radial and supported in $\{ar_k\le|x|\le r_k\}$
(Lemma~I.\ref{I-lem:ring-facts}); $f^{(n)}_n$ is $3$-fold rotationally symmetric and supported in
$\NNm_n$ (Lemma~\ref{lem:seed}). So each $S_k$ is compact, $\Orot$-invariant and
avoids $0$. Let $i<k$. By these support properties, $S_k\subset\overline B_{5r_k/4}$, and
$\frac54r_k\le\frac54s_ir_i<ar_i$, because $r_k\le r_{i+1}=s_ir_i$ and $s_i\le\sbar/2\le a/32$, while
$S_i\subset\{|x|\ge ar_i\}$; so $S_k\cap S_i=\emptyset$, and
$f^{(n)}\mathbf 1_{S_k}=f^{(n)}_k$.

(2) By Lemma~\ref{lem:adm}, $\omega^{(n)}(t)\circ\Orot=\omega^{(n)}(t)$; with
$\KBS(\Orot x)=\Orot\KBS(x)$ and the change of variables $y=\Orot y'$ this gives
$v^{(n)}(t,\Orot x)=\Orot v^{(n)}(t,x)$. So $t\mapsto\Orot^{-1}\flow_t(\Orot y)$
solves the equation of $\flow_t(y)$ with the same initial value; since
$v^{(n)}$ is $C^1$ and $\nabla v^{(n)}$ is bounded on $[0,t]\times\R^2$
(Theorem~\ref{thm:smooth}\textup{(2)}), Gronwall's inequality gives uniqueness,
and $\flow_t\circ\Orot=\Orot\circ\flow_t$. At $x=0$, $v^{(n)}(t,0)=\Orot
v^{(n)}(t,0)$, so $v^{(n)}(t,0)=0$ and $\flow_t(0)=0$. As $\flow_t$ is a
homeomorphism (Theorem~\ref{thm:smooth}\textup{(4)}), the sets
$S_k(t)=\flow_t(S_k)$ inherit compactness, disjointness and $\Orot$-invariance
from (1), $0\notin S_k(t)$, and
$\supp\omega^{(n)}(t)=\flow_t(\supp f^{(n)})\subset\bigcup_kS_k(t)$, because
$\omega^{(n)}(t)=f^{(n)}\circ\flow_t^{-1}$.

(3) Since $\omega^{(n)}(t)=f^{(n)}\circ\flow_t^{-1}$ and
$\mathbf 1_{S_k(t)}=\mathbf 1_{S_k}\circ\flow_t^{-1}$, (1) gives the first
identity, and (2), by which $D_k(t)$ is the disjoint union of
$S_k(t),\dots,S_n(t)$, gives the other two. Each $f^{(n)}_k$ lies in
$C_c^\infty(\R^2)$, is nonnegative and $3$-fold rotationally symmetric (Lemmas~I.\ref{I-lem:ring-facts}
and \ref{lem:seed}), and $\flow_t^{-1}$ is
a $C^\infty$ diffeomorphism commuting with $\Orot$; this gives the regularity
and symmetry. Since $\det D\flow_t=1$,
$\int\omega^{(n)}(t)\mathbf 1_{S_k(t)}=\int f^{(n)}_k=A_kr_k^2\int\Phi
=A_kr_k^2\|\Phi\|_{L^1}$, as $\Phi\ge0$; for $k=n$ we use that the nonnegative equimeasurable
functions $\Phi^{\kseed}$ and $\Phi$ have the same integral.

(4) The two terms of
$w_k(t,y)=A_k^{-1}(\omega^{(n)}(t)\mathbf 1_{S_k(t)})(r_ky)-\Phi(y)$ are real,
lie in $C_c^\infty(\R^2)$, are $3$-fold rotationally symmetric (by (3) and
Lemma~I.\ref{I-lem:ring-facts}) and are supported in $N$, because
$S_k(t)\subset N_k$ and $\supp\Phi\subset\{a\le|y|\le1\}\subset N$
(Lemma~\ref{lem:annuli}). By (3) and the change of variables $x=r_ky$, both are
nonnegative with integral, hence $L^1$ norm, $\|\Phi\|_{L^1}$; so $\int w_k(t)=0$ and $\|w_k(t)\|_{L^1}\le2\|\Phi\|_{L^1}$.
The last inequality is the Cauchy-Schwarz inequality on $N$.
\end{proof}

The fields of Definition~\ref{def:regions} depend continuously on time: if
$J\subset[0,\infty)$ is a compact interval with $D_{k+1}(t)\subset B_{ar_k/8}$
for $t\in J$, then $t\mapsto\hin_k(t)$ is continuous from $J$ into
$C^M(N')$ for every integer $M\ge0$. Indeed, $\hin_n=-G_ny^\perp/|y|^2$ is
constant, because $D_{n+1}(t)=\emptyset$. Let $k<n$ and $g=\sum_{i>k}f^{(n)}_i$, supported in the compact set
$D_{k+1}(0)$. By Equation~\eqref{eq:omega-ge}, $\omega_{\ge k+1}(t)=g\circ\flow_t^{-1}$,
so the change of variables $y=\flow_t(y')$, with $\det D\flow_t=1$, gives
$(\KBS\ast\omega_{\ge k+1}(t))(x)=\int\KBS(x-\flow_t(y'))g(y')\dd y'$. As
$(t,y')\mapsto\flow_t(y')$ is continuous (Theorem~\ref{thm:smooth}\textup{(4)}),
$\rho=\max\{|\flow_t(y')|:t\in J,\ y'\in D_{k+1}(0)\}$ exists, and $\rho<\frac a8r_k$,
because the maximum is attained at a point of one of the sets
$\flow_t(D_{k+1}(0))=D_{k+1}(t)\subset B_{ar_k/8}$.
Let $U$ be the set of the $x$ at distance less than $\frac a{16}r_k$ from $N'_k$;
as the inner radius of $N'_k$ is $\frac a4r_k$, every $x\in U$ has
$|x|>\frac3{16}ar_k$. So, for $x\in U$, $t\in J$ and $y'\in D_{k+1}(0)$,
$|x-\flow_t(y')|\ge\frac3{16}ar_k-\rho>\frac a{16}r_k$. Since
$|\partial^\beta\KBS(z)|\le C_\beta|z|^{-1-|\beta|}$, the derivatives
$\partial^\beta_x$, $|\beta|\le M$, of the integrand are bounded, for $x\in U$
and $t\in J$, by a constant times $|g(y')|$, which is integrable. So we may
differentiate under the integral sign, and each of these derivatives of the field is
continuous on $J\times U$ by dominated convergence, hence uniformly continuous on
the compact set $J\times N'_k$; since $G_ky^\perp/|y|^2$ does not depend on
$t$, passing to ring-$k$ spatial and velocity units gives the claim. In the same way, if the rings $i<k$ lie outside $B_{6r_k}$
for $t\in J$, then $t\mapsto\hout_k(t)\in C^M(N')$ is continuous on $J$:
replace $g$ by $\sum_{i<k}f^{(n)}_i$ and $U$ by $B_{4r_k}$, which contains $N'_k$,
and use $|x-\flow_t(y')|\ge2r_k$ for
$x\in B_{4r_k}$ and $y'\in S_1\cup\dots\cup S_{k-1}$, which holds because
$\flow_t(y')$ then lies in $S_1(t)\cup\dots\cup S_{k-1}(t)$, outside $B_{6r_k}$. Consequently $\etin_k$, $\etout_k$ and
$\eta_k$ are continuous on every interval on which they are defined, so the
integrals that define $\maj_k$ and $\maj^{\mathrm{in}}_k$ exist.

If on $[0,T]$ both $D_{k+1}(t)\subset B_{ar_k/8}$ and the rings $i<k$ lie outside
$B_{6r_k}$, then, by Lemma~\ref{lem:regions-basic}, the fully rescaled solution
with $A=A_k$, $r=r_k$ and $t_*=0$, split into its parts on $S_k(t)$, on $D_{k+1}(t)$ and on
$S_1(t)\cup\dots\cup S_{k-1}(t)$ (rescaled), satisfies the hypotheses of
Lemma~\ref{lem:reduction} on $[0,A_kT]$ with $G=G_k$ and the explicitly
rescaled field $h(s)=h_k(s/A_k)$ there. Indeed, by Lemma~\ref{lem:regions-basic}\textup{(2)} and \textup{(3)}, the
three parts have pairwise disjoint supports, are $3$-fold rotationally symmetric
and are transported by the flow of the solution in ring-$k$ units. The part on
$D_{k+1}(t)$ is $\nu^{\mathrm{in}}_k$, whose integral at $t=0$ is $2\pi G_k$ by
the definition of $G_k$, so $G=G_k$. The two conditions on $[0,T]$, read in ring-$k$
units, imply the support conditions of Lemma~\ref{lem:reduction}. While
$S_k(t)\subset N_k$, the function $w=\omega_{\mathrm r}-\Phi$ of
Lemma~\ref{lem:reduction}\textup{(3)} is $w(s)=w_k(s/A_k)$. In the same
local coordinate, its field norm and Duhamel bound are
$\eta(s)=\eta_k(s/A_k)$ and $\maj(s)=\maj_k(s/A_k)$: the factor $A_k$ in
Definition~\ref{def:regions} follows from the substitution $s=A_kt$.

Part~I and Section~\ref{sec:single} treat the linearization $\Lin(G)$ for
$G\in[0,\Gone]$. The next lemma shows that every inner circulation $G_k$ lies in
this range when $\sbar$ is small, with a bound independent of $n$.

\begin{lemma}
\label{lem:inner-circ}
In Definition~\ref{def:regions}, $\int\omega_{\ge k+1}(t)$ is independent of
$t$, and $0\le G_k\le\mu s_k^2\le\mu\sbar^2$ for $1\le k\le n$;
moreover $G_n=0$. In particular, if $\mu\sbar^2\le\Gone$, then
$G_k\in[0,\Gone]$ for $1\le k\le n$.
\end{lemma}
\begin{proof}
Let $F\colon\R\to\R$ be continuous with $F(0)=0$. Since $\flow_t$ preserves
area (Theorem~\ref{thm:smooth}\textup{(4)}), Equation~\eqref{eq:omega-ge} gives
$\int F(\omega_{\ge k}(t))\dd x=\int F\big(\sum_{i=k}^nf^{(n)}_i\big)\dd x$. The
$f^{(n)}_i$ are supported in the pairwise disjoint sets $S_i$
(Lemma~\ref{lem:regions-basic}\textup{(1)}), so, as $F(0)=0$, the function
$F\big(\sum_{i=k}^nf^{(n)}_i\big)$ is the sum of the $F(f^{(n)}_i)$. Hence, for
every such $F$ and every $1\le k\le n+1$,
\begin{equation}
\label{eq:conserved}
\int F\big(\omega_{\ge k}(t)\big)\dd x=\sum_{i=k}^n\int F\big(f^{(n)}_i\big)\dd x .
\end{equation}
With $F(s)=s$, $\int\omega_{\ge k+1}(t)$ does not depend on $t$, and, as
$\int f^{(n)}_i=A_ir_i^2\|\Phi\|_{L^1}$ (Lemma~\ref{lem:regions-basic}\textup{(3)}
at $t=0$), $2\pi G_kA_kr_k^2=\|\Phi\|_{L^1}\sum_{i=k+1}^nA_ir_i^2$. So $G_k\ge0$
and $G_n=0$, as the sum is then empty. Since $A_{i+1}r_{i+1}^2=\mu s_i^2\,A_ir_i^2$ and
$\mu s_i^2\le2(a/16)^2<\frac12$, each term of the sum is less than half the
previous one, and comparison with a geometric series shows that the sum is at most
$2A_{k+1}r_{k+1}^2=2\mu s_k^2A_kr_k^2$; so
$G_k\le\|\Phi\|_{L^1}\mu s_k^2/\pi\le\mu s_k^2$, because
$0\le\Phi<1$ and $\supp\Phi\subset\overline B_1$ give
$\|\Phi\|_{L^1}\le\pi$; and $s_k\le\sbar$.
\end{proof}

The next lemma bounds the vorticity and velocity of ring $k$ and the rings
inside it, uniformly in $n$ and $t$. Proposition~\ref{prop:confine} uses
the bound on $\KBS\ast\omega_{\ge k+1}$ to control the expansion of $D_{k+1}$.

\begin{lemma}
\label{lem:speed}
In Definition~\ref{def:regions}, there is an absolute constant $\Cstar$ such that for all $1\le k\le n$, $t\ge0$ and $q\in\{1,4\}$,
\[
\begin{gathered}
\|\omega_{\ge k}(t)\|_{L^q}\le\Cstar^{1/q}A_kr_k^{2/q},\qquad
\|\KBS\ast\omega_{\ge k}(t)\|_{L^\infty}\le\Cstar A_kr_k,\\
\|\omega_{\ge k}(t)\|_{L^1}=\|\Phi\|_{L^1}\sum_{i=k}^nA_ir_i^2\le2\|\Phi\|_{L^1}A_kr_k^2 .
\end{gathered}
\]
\end{lemma}
\begin{proof}
By Equation~\eqref{eq:conserved} with $F(s)=|s|^q$ and the equimeasurability of
$\Phi^{\kseed}$ and $\Phi$ (Lemma~\ref{lem:seed}),
$\|\omega_{\ge k}(t)\|_{L^q}^q=\|\Phi\|_{L^q}^q\sum_{i=k}^nA_i^qr_i^2$. For
$q\in\{1,4\}$ and $i\ge k$ we have
$A_{i+1}^qr_{i+1}^2=\mu^qs_i^2A_i^qr_i^2\le\mu^4\sbar^2A_i^qr_i^2\le\frac18A_i^qr_i^2$,
because $\mu\le2$ and $\sbar\le a/16<\frac1{16}$ give $\mu^4\sbar^2<\frac1{16}$;
so $\sum_{i=k}^nA_i^qr_i^2\le\frac87A_k^qr_k^2$. With $q=1$ this gives the
identity and the bound for $\|\omega_{\ge k}(t)\|_{L^1}$. As $0\le\Phi<1$ and
$\supp\Phi\subset\{a\le|y|\le1\}$ (Lemma~I.\ref{I-lem:ring-facts}), $\|\Phi\|_{L^q}^q\le\pi$, and
$\|\omega_{\ge k}(t)\|_{L^q}^q\le\frac{8\pi}7A_k^qr_k^2\le\Cstar A_k^qr_k^2$.
Finally, let $g=\omega_{\ge k}(t)\in (L^1\cap L^4)(\R^2)$; we may assume $g\ne0$. Let
$x\in\R^2$ and $R>0$; by Lemma~I.\ref{I-lem:BS-basic}(a), $2\pi|\KBS(x-y)|=|x-y|^{-1}$.
We split the integral of $|g(y)||x-y|^{-1}$ over $y$ into the near field and the
far field.
On the near field, H\"older's inequality with exponents $4,4/3$ and
integration in polar coordinates give
\[
\int_{|x-y|<R}|g(y)||x-y|^{-1}\dd y\lesssim\|g\|_{L^4}R^{1/2}.
\]
The far-field integral is at most $\|g\|_{L^1}R^{-1}$. Choosing
$R=(\|g\|_{L^1}/\|g\|_{L^4})^{2/3}$ balances these bounds and yields
\[
\|\KBS\ast g\|_{L^\infty}
\lesssim\|g\|_{L^1}^{1/3}\|g\|_{L^4}^{2/3}
\lesssim A_kr_k.
\]
The implicit constants are absolute, using the numerical $L^1$ and $L^4$
bounds above. Choose $\Cstar$ to dominate this velocity coefficient and $8\pi/7$.
\end{proof}

We call the accumulated contribution of the inner velocity field to a ring's
unstable linear response, with its free growth factored out, the \emph{kick}.
Its amplitude is given precisely by Equation~\eqref{eq:kick-amplitude} below.
The clock below records activation and saturation in physical time $t$,
so that events on different rings can be compared directly. The other time
coordinates are used locally for calculations: $s=A_kt$ puts ring $k$ on its
natural unit time scale, making the estimates uniform in $k$;
$u=A_k(t-\tsh_k)$ additionally puts activation at $u=0$, so that Section~\ref{sec:single}
applies with datum $w_k(\tsh_k)$; and
$\tauK=A_{k+1}(t-T_{k+1})$ measures the kick in the inner ring's units, starting
at its saturation. These are changes of coordinate for the same evolution.
All indexed ring quantities, including $w_k(t)$, $h_k(t)$, $\maj_k(t)$ and
$c_k(t)$, are evaluated at physical time; their rescaled versions are
introduced explicitly inside the calculation that uses them.

The clock turns the relay described in the Introduction into definite
times. It is defined recursively, from the innermost ring outward. Ring $n$
becomes active at time $0$. Ring $k$ becomes active at the time $\tsh_k$, and
$T_k$ is the time at which its unstable coefficient is expected to reach the
size $\delta_0$. Ring $k-1$ then waits a further $\Tpin$ units of ring-$k$
time before being declared active. This \emph{pinning phase} allows its kick
amplitude, with the free linear growth factored out, to settle near the value
acquired around $T_k$. After $T_k$, the mode-$3$ moment of the inner vorticity
has at most polynomial growth, whereas the weight in the kick integral decays
exponentially. Choosing $\Tpin$ sufficiently large therefore makes the
remaining contribution small: the normalized unstable coefficient of ring
$k-1$ stays essentially stationary during its ensuing active phase
(Lemma~\ref{lem:drift}), while its unstable mode grows. For ring $k-1$, the kick
of Lemma~\ref{lem:coupling}\textup{(2)} is evaluated at $\tauK=\Tpin$.

\begin{definition}
\label{def:clock}
Let $\delta_0,\Tpin>0$. In Definition~\ref{def:regions}, assume $G_k\le\Gone$
for $1\le k\le n$. Define recursively, downward in $k$:
\begin{itemize}
\item $\tsh_n=0$;
\item if $\tsh_k$ is defined and $S_k(\tsh_k)\subset N_k$:
  $c_k(\tsh_k)=\ellst(G_k)\big(w_k(\tsh_k)\big)$ and $\ksh_k=|c_k(\tsh_k)|$;
\item if moreover $\ksh_k\in(0,\delta_0)$:
  $\hat T_k=\sigma_k^{-1}\log(\delta_0/\ksh_k)$, $T_k=\tsh_k+\hat T_k/A_k$, and, if
  $k\ge2$, $\tsh_{k-1}=T_k+\Tpin/A_k$.
\end{itemize}
Thus $\hat T_k=A_k(T_k-\tsh_k)$ is the dimensionless active duration.
The clock is \textbf{defined up to level $k$} if $S_i(\tsh_i)\subset N_i$ and
$\ksh_i\in(0,\delta_0)$ for every $i\in[k,n]$, so that $T_i$ is defined for
these $i$. Put $T_{n+1}=0$, $W_k=[T_{k+1},T_k]$, the
\textbf{active phase} $\Phase_k=[\tsh_k,T_k]$ and the \textbf{pinning phase}
$\Pin_k=[T_{k+1},\tsh_k]$ ($k<n$). Ring $k$ is \textbf{quiet} on $[0,\tsh_k]$,
\textbf{active} on $\Phase_k$, and \textbf{debris} after $T_k$. If $h_k$ is
defined on $[\tsh_k,t]$, the running coefficient at physical time $t$ is
\begin{equation}
\label{eq:running-k}
c_k(t)=c_k(\tsh_k)+A_k\int_{\tsh_k}^t
 e^{-\lambda_kA_k(v-\tsh_k)}\,\ellst(G_k)\big(F_k(v)\big)\dd v,
\end{equation}
where $F_k(v)=-h_k(v)\cdot\nabla\Phi$ on $N'$ and $F_k(v)=0$ elsewhere.
The exponent removes the free linear evolution since activation; the factor
$A_k$ comes from the change to the local coordinate $u=A_k(t-\tsh_k)$.
\end{definition}

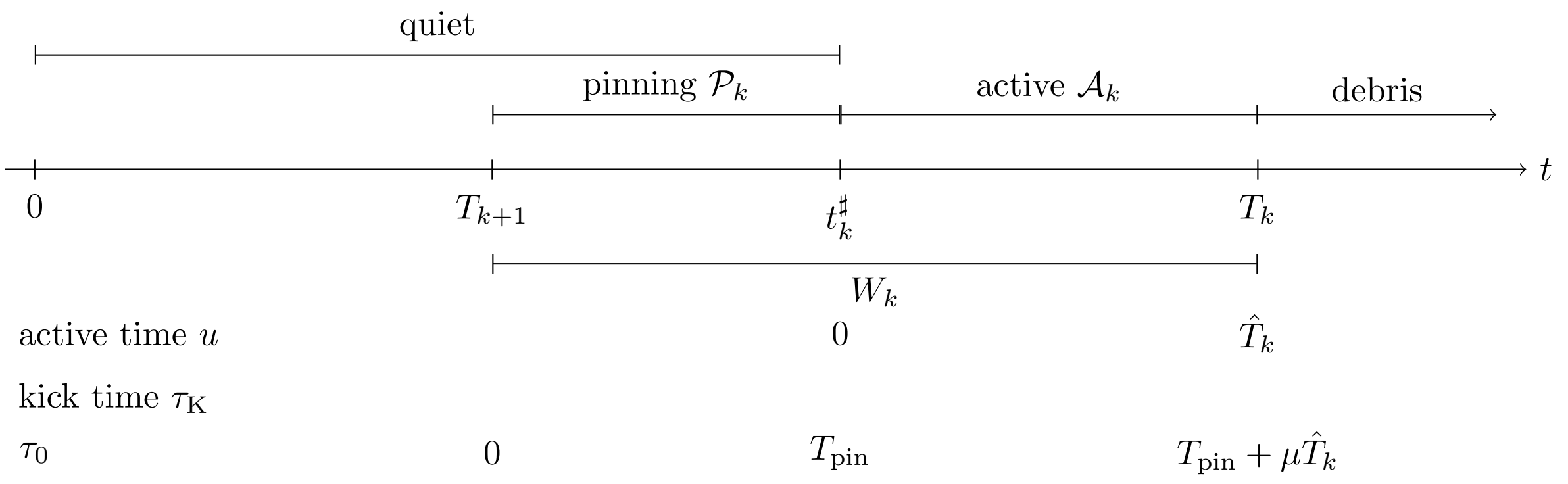
\begin{figure}[H]
\centering
\begin{tikzpicture}[x=1cm,y=1cm,font=\small]
\draw[|-|] (0.5,1.15) -- (8.6,1.15)
  node[midway,above] {quiet};
\draw[|-|] (5.1,0.55) -- (8.6,0.55)
  node[midway,above] {pinning $\Pin_k$};
\draw[|-|] (8.6,0.55) -- (12.8,0.55)
  node[midway,above] {active $\Phase_k$};
\draw[->] (12.8,0.55) -- (15.2,0.55)
  node[midway,above] {debris};
\draw[->] (0.2,0) -- (15.5,0) node[right] {$t$};
\foreach \x in {0.5,5.1,8.6,12.8}
  \draw (\x,0.10) -- (\x,-0.10);
\node[below] at (0.5,-0.12) {$0$};
\node[below] at (5.1,-0.12) {$T_{k+1}$};
\node[below] at (8.6,-0.12) {$\tsh_k$};
\node[below] at (12.8,-0.12) {$T_k$};
\draw[|-|] (5.1,-0.95) -- (12.8,-0.95)
  node[midway,below] {$W_k$};
\node[anchor=west] at (0.2,-1.65) {active time $u$};
\node at (8.6,-1.65) {$0$};
\node at (12.8,-1.65) {$\hat T_k$};
\node[anchor=west] at (0.2,-2.30) {kick time $\tauK$};
\node at (0.5,-2.85) {$\tau_0$};
\node at (5.1,-2.85) {$0$};
\node at (8.6,-2.85) {$\Tpin$};
\node at (12.8,-2.85) {$\Tpin+\mu\hat T_k$};
\end{tikzpicture}
\caption{The intervals of ring $k<n$ in physical time and the local coordinates used in calculations
(positions schematic). The window $W_k$ is the union of the pinning and active
intervals; $\tau_0=-A_{k+1}T_{k+1}$ is the kick-time coordinate of $t=0$.
For the seeded ring $k=n$, $\tsh_n=0$ and $W_n=\Phase_n$, with no
pinning phase.}
\label{fig:ring-clock}
\end{figure}

The integral in Equation~\eqref{eq:running-k} exists, so $c_k(t)$ is well
defined: its integrand is continuous in $v$, since the map
$v\mapsto h_k(v)\in C^0(N')$ is continuous (as shown after
Lemma~\ref{lem:regions-basic}) and
$|\ellst(G_k)(-h\cdot\nabla\Phi)|\le\Cls\|h\|_{C^0(N')}$ for continuous $h$
(Lemma~I.\ref{I-lem:uniform}).

Lemma~\ref{lem:clock} verifies that the recursion gives well-defined,
strictly ordered times whenever its support and amplitude conditions are met.
Knowing the solution through $\tsh_i$ fixes $\ksh_i$, hence $T_i$ and the next
activation time $\tsh_{i-1}$. This lets the continuity argument of
Section~\ref{sec:bootstrap} use the estimates of Section~\ref{sec:relay}
on $[0,t]$ without using later evolution.

To state this dependence precisely, a quantity built from $\omega^{(n)}$ by
Definitions~\ref{def:regions} and \ref{def:clock} is
\textbf{determined by $\omega^{(n)}$ on $[0,\bar t]$} if it keeps its value when
$\omega^{(n)}$ is replaced in these definitions by any
$\tilde\omega\in C^\infty([0,\infty)\times\R^2)$ such that $\tilde\omega=\omega^{(n)}$
on $[0,\bar t]\times\R^2$, every $\tilde\omega(t)$ has compact support, and
$\KBS\ast\tilde\omega$ is smooth with every derivative bounded on $[0,T]\times\R^2$
for every $T$; the flow $\flow$ is then replaced by the flow of
$\KBS\ast\tilde\omega$, and $f^{(n)}$, the sets $S_k$ and the parameters are kept.

\begin{lemma}
\label{lem:clock}
In Definition~\ref{def:clock}, assume that the clock is defined up to level $k$.
\begin{enumerate}[label=\textup{(\arabic*)}]
\item $0=\tsh_n<T_n<\tsh_{n-1}<T_{n-1}<\dots<\tsh_k<T_k$, and $T_k<\tsh_{k-1}$
  if $k\ge2$.
\item For every $i\in[k,n]$, each of $\tsh_i$, $\ksh_i$, $T_i$ and, if $i\ge2$,
  $\tsh_{i-1}$ is determined by $\omega^{(n)}$ on $[0,\tsh_i]$.
\item Hence, for every $t\ge0$, the set $\mathcal J_t$ of the $i\in[k,n]$ with
  $\tsh_i\le t$, and the numbers $\tsh_i$ and $T_i$ for $i\in\mathcal J_t$, are
  determined by $\omega^{(n)}$ on $[0,t]$. In particular, the solution on
  $[0,t]$ determines all activation and saturation times $\tsh_i$, $T_i$
  ($k\le i\le n$) not exceeding $t$, and whether each ring $i\in[k,n]$ is
  quiet, active, or debris at $t$, as defined in Definition~\ref{def:clock}.
\end{enumerate}
\end{lemma}
\begin{proof}
(1) Let $i\in[k,n]$. By Lemma~\ref{lem:inner-circ}, $G_i\ge0$, and $G_i\le\Gone$ by
the assumption of Definition~\ref{def:clock}; so Proposition~I.\ref{I-prop:SU} gives
$\sigma_i=\sigma(G_i)>\frac{211}{2500}$. Moreover $\ksh_i<\delta_0$, as the clock
is defined up to level $k$, so $\hat T_i$ is finite and positive. Thus the
formulas in Definition~\ref{def:clock} uniquely determine $T_i>\tsh_i$ and,
if $i\ge2$, $\tsh_{i-1}$, with $\tsh_{i-1}-T_i=\Tpin/A_i>0$.
These inequalities, for $i=n,n-1,\dots,k$,
give (1).

(2) Let $\bar t\ge\tsh_i$, let $\tilde\omega$ be as above, and mark with a tilde
the objects built from $\tilde\omega$. The velocities $\KBS\ast\tilde\omega$ and
$v^{(n)}$ agree on $[0,\bar t]\times\R^2$ and are Lipschitz in $x$, so, by
uniqueness for the equation of the flow, the two
flows agree on $[0,\bar t]$; and $\tilde G_j=G_j$, because
$\tilde\omega(0)=f^{(n)}$, so the numbers $\sigma_j$ do not change either. We show by downward induction on $j\in[i,n]$ that
$\tilde t^\sharp_j=\tsh_j$, $\tilde\kappa^\sharp_j=\ksh_j$, $\tilde T_j=T_j$, and
$\tilde t^\sharp_{j-1}=\tsh_{j-1}$ if $j\ge2$. First, $\tilde t^\sharp_n=0=\tsh_n$.
If $j\in[i,n]$ and $\tilde t^\sharp_j=\tsh_j$, then $\tsh_j\le\tsh_i\le\bar t$ by (1), so
$\tilde S_j(\tsh_j)=S_j(\tsh_j)\subset N_j$ and $\tilde w_j(\tsh_j)=w_j(\tsh_j)$;
hence $\tilde\kappa^\sharp_j=\ksh_j\in(0,\delta_0)$, and
Definition~\ref{def:clock} gives $\tilde T_j=T_j$ and
$\tilde t^\sharp_{j-1}=\tsh_{j-1}$. This completes the induction and proves (2).

(3) Let $\tilde\omega$ agree with $\omega^{(n)}$ on $[0,t]\times\R^2$, and let
$\tilde{\mathcal J}_t$ be the set of the $i\in[k,n]$ for which
$\tilde t^\sharp_i$ is defined and at most $t$. Let $i_0=\min\mathcal J_t$, which
exists because $\tsh_n=0$. By (1), $\tsh_i$ decreases as $i$ increases, so
$\mathcal J_t=[i_0,n]$, and $\tsh_i,T_i>t$
for $i\notin\mathcal J_t$. By (2) with $\bar t=t$ and $i=i_0$, the numbers
$\tsh_i$, $T_i$, $i\ge i_0$, are the same for $\tilde\omega$, so
$[i_0,n]\subset\tilde{\mathcal J}_t$, and, if $i_0>k$,
$\tilde t^\sharp_{i_0-1}=\tsh_{i_0-1}>t$; the proof of (1) applies to the clock
of $\tilde\omega$ wherever it is defined, so every defined $\tilde t^\sharp_i$ with $i<i_0$ exceeds
$t$. Hence $\tilde{\mathcal J}_t=\mathcal J_t$. These data determine whether
ring $i\in[k,n]$ is quiet, active, or debris at $t$. If $i\notin\mathcal J_t$, then $t<\tsh_i<T_i$ and ring
$i$ is quiet and neither active nor debris, and for $\tilde\omega$ the number
$\tilde t^\sharp_i$ is either undefined or greater than $t$, so ring $i$ is
neither active nor debris and is quiet at $t$ by the convention of
Definition~\ref{def:loc} in Section~\ref{sec:several}, which includes rings with
no assigned activation time among the quiet rings. If $i\in\mathcal J_t$, ring $i$ is
quiet if and only if $t=\tsh_i$, active if and only if $t\le T_i$, and debris if
and only if $t>T_i$.
\end{proof}

\section{Several Rings: Transfer of Moments and Confinement}
\label{sec:several}

We bound the velocity fields that the rings of $\omega^{(n)}$ induce on one
another, and we show that, on a time window after ring $k+1$ saturates, the
material inside ring $k$ stays in the disk $\overline B_{ar_k/8}$.

We assume \textup{(N1)}--\textup{(N3)} of Hypothesis~I.\ref{I-hyp:N}, fix
$n\ge2$, and use the datum and its approximations $\omega^{(n)}$, the ring
quantities and the clock of Definitions~\ref{def:datum}, \ref{def:approx},
\ref{def:regions} and \ref{def:clock}, with $\mu\in(1,2]$ and $\sbar\in(0,a/16]$; we write
$\omega=\omega^{(n)}$, $v=\KBS\ast\omega$, and $\flow_t$ for the flow of $v$.
Since $s_k=\sbar(2\mu)^{-k}\le\sbar/2$ and $r_k/r_i=s_is_{i+1}\cdots s_{k-1}$,
\begin{equation}
\label{eq:sk-small}
s_k\le\frac a{32}\quad(k\ge1),\qquad
\frac{r_k}{r_i}\le s_{k-1}\,\sbar^{\,k-i-1}\quad(1\le i<k).
\end{equation}
Moreover $\sbar<\frac12$ and
$\mu^4\sbar^2\le16\cdot2^{-8}<\frac18$, because $a<1$. The conditions
$\sbar\le a/16$, $\sbar\le\frac12$ and $\mu^4\sbar^2\le\frac18$ appear
among the hypotheses of Lemmas~\ref{lem:transfer} and~\ref{lem:field-bound} and of
Proposition~\ref{prop:confine} so that these can be quoted on
their own; under the standing assumptions of this section they always hold.

In this section $q$, without an index, is the order of a derivative.
For a velocity field $V$, $V_{[j]}$ is $V$ in ring-$j$ units,
that is, $V_{[j]}(y)=V(r_jy)/(A_jr_j)$ (defined before
Definition~\ref{def:regions}). As in
Section~I.\ref{I-sec:prelim}, $\|h\|_{C^k(U)}$ is the sum over $q\le k$ of
$\sup_U|\nabla^qh|$, so the proofs below bound each $\nabla^q$ and then sum over
$q$. We use the \textbf{scaling
rule}: if $c\in\R$, $c'>0$ and $W(y)=cV(c'y)$, then
$\nabla^qW(y)=c\,c'^{\,q}(\nabla^qV)(c'y)$.

In the complex-variable calculations below, we use the spatial variable
$z\in\C$ and the conventions for $V(z)$ and its Cartesian derivatives
from Section~\ref{sec:II-prelim}.
We write $F(z)=V_1(z)-iV_2(z)$ for the complex form of a vector field $V$, and we use
Equation~\eqref{eq:BS-complex} for the complex form of $\KBS\ast g$
(Section~\ref{sec:II-prelim}). The monopole
$\frac1{2\pi}\big(\int g\big)(-\Im z,\Re z)/|z|^2$ of $g$ has complex form $\big(\int g\big)/(2\pi iz)$,
and the rotation $\Orot$ by the angle $2\pi/3$ (Section~I.\ref{I-sec:prelim}) is
multiplication by $e^{2\pi i/3}$. For real $g$, the complex forms
of $\KBS\ast g$ and of $\KBS\ast g$ minus its monopole are holomorphic on
$U=\C\setminus(\supp g\cup\{0\})$, so Equation~\eqref{eq:BS-complex-deriv}
applies to them there.

\subsection{Fields Induced by One Ring on Another}

Off the support of a function $g$, the complex form of the velocity of $g$ is
given by Equation~\eqref{eq:BS-complex}, and expanding $(z-\zeta)^{-1}$ turns it
into a power series: in $z^{-1}$ outside a disk centered at the origin that
contains the support, and in $z$ inside such a disk that does not meet it. For a $3$-fold rotationally
symmetric $g$ only every third coefficient survives. For material inside
ring $j$, apart from the monopole, the series starts at
$z^{-4}$, and on $N'_j$ this gives the factor $\varepsilon_j=\mu s_j^5$ of
Definition~\ref{def:regions} in parts (1)--(3) of the next lemma; for material
outside ring $k$ it starts at $z^2$, which gives the factor $b_{i,k}$ in
part (4).

\begin{lemma}
\label{lem:transfer}
Assume $\sbar\le a/16$, and let $M_h\ge1$ and $s'\ge0$ be integers. There is
$\Cvel\ge1$, depending only on $a$ and $M_h$, such that the following hold at
every time $t\ge0$, for $1\le j<n$. We choose this constant large enough for
Lemmas~\ref{lem:inner-field} and~\ref{lem:field-bound} below as well.
\begin{enumerate}[label=\textup{(\arabic*)}]
\item Let $g\in L^1(\R^2)$ be real valued, $3$-fold rotationally symmetric and supported in
  $\overline B_{ar_{j+1}/8}$, and let $E_g$ be its velocity field minus its
  monopole. Then
  \[
  \|(E_g)_{[j]}\|_{C^{M_h}(N')}
  \le \Cvel\varepsilon_j\sup_{|y|=a/4}|(E_g)_{[j+1]}(y)|.
  \]
\item If $S_{j+1}(t)\subset N_{j+1}$, let
  $V=\KBS\ast\big(A_{j+1}w_{j+1}(t,\cdot/r_{j+1})\big)$. Then
  \[
  \|V_{[j]}\|_{C^{M_h}(N')}
  \le \Cvel\varepsilon_j\|w_{j+1}(t)\|_{H^{s'}}.
  \]
\item If $D_{j+1}(t)\subset\overline B_R$ with $R\le ar_j/8$, let $E$ be the
  velocity field of $\omega_{\ge j+1}(t)$ minus its monopole. Then
  \[
  \|E_{[j]}\|_{C^{M_h}(N')}
  \le \Cvel\varepsilon_j(R/r_{j+1})^3.
  \]
\item For $1\le i<k\le n$ with $S_i(t)\subset N_i$, put
  $V_i=\KBS\ast\big(A_iw_i(t,\cdot/r_i)\big)$ and
  $b_{i,k}=\mu^{i-k}r_k/r_i$. Then
  \[
  \begin{aligned}
  \|(V_i)_{[k]}\|_{C^{M_h}(N')}&\le \Cvel b_{i,k}\|w_i(t)\|_{L^1},\\
  \|(V_i)_{[k]}\|_{C^{M_h}(N')}&\le \Cvel b_{i,k}\|w_i(t)\|_{H^{s'}}.
  \end{aligned}
  \]
  If $S_i(t)\subset N_i$ for every $i<k$, put
  $V_{\mathrm{out}}=\KBS\ast\big(\omega(t)\mathbf 1_{S_1(t)\cup\dots\cup S_{k-1}(t)}\big)$.
  On $N'_k$, $V_{\mathrm{out}}=\sum_{i<k}V_i$, and
  \[
  \|(V_{\mathrm{out}})_{[k]}\|_{C^{M_h}(N')}
  \le \Cvel\sum_{i<k}b_{i,k}\|w_i(t)\|_{L^1}.
  \]
\end{enumerate}
\end{lemma}

\begin{proof}
Throughout this proof, implicit constants depend only on $a$ and $M_h$.
The differentiated geometric series are finite:
\[
\sum_{p\ge1}(3p+q)^q2^{-3(p-1)}<\infty
\qquad(0\le q\le M_h).
\]
We prove the four bounds with implicit constants, which are absorbed into
$\Cvel$.
The $H^{s'}$ bounds then follow from
$\|w\|_{L^1}\le |N|^{1/2}\|w\|_{L^2}\le4\|w\|_{H^{s'}}$,
since each perturbation is supported in $N\subset\overline B_2$
(Lemma~\ref{lem:regions-basic}(4)). Put $\rho_0=a/4$, the inner radius of $N'$.

\emph{Proof of (1).} We bound the coefficients of the series of the field of
$g$ minus its monopole, in ring-$(j+1)$ units, by the values of this field on the circle $|z|=\rho_0$, and then sum
the series at the points of $N'_j$, which in these units lie outside
$B_{\rho_0/s_j}$. Put $V=(E_g)_{[j+1]}$,
$\mathcal M=\sup_{|z|=\rho_0}|V(z)|$ and
$\tilde g(\zeta)=g(r_{j+1}\zeta)/A_{j+1}$. Since $\KBS$ is homogeneous of degree
$-1$, $V$ is the field of $\tilde g$ minus its monopole, and $\tilde g$ is
$3$-fold rotationally symmetric and supported in $\overline B_{a/8}$. Inserting
$(z-\zeta)^{-1}=\sum_{l\ge0}\zeta^lz^{-l-1}$, uniformly convergent for
$|\zeta|\le a/8$ and $z$ in compact subsets of $\{|z|>a/8\}$, in
Equation~\eqref{eq:BS-complex} and integrating term by term, we find that the
complex form of $V$ is
\[
F(z)=\sum_{l\ge1}\alpha_lz^{-l-1},\qquad \alpha_l=\frac1{2\pi i}\int\tilde g(\zeta)\zeta^l\dd A(\zeta)\qquad(|z|>a/8),
\]
without a term $l=0$: in the expansion of the full field that term is the
monopole, which $V$ does not contain. The substitution $\zeta=\Orot\zeta'$ and the $3$-fold
rotational symmetry of $\tilde g$ give
$\alpha_l=e^{2\pi il/3}\alpha_l$, so $\alpha_l=0$ unless $3\mid l$. On
$|z|=\rho_0$ the series converges uniformly, and $\alpha_l\rho_0^{-l-1}$ is the
Fourier coefficient of $\theta\mapsto F(\rho_0e^{i\theta})$ at the frequency
$-(l+1)$; since $V$ is real, $|F|=|V|$, so $|\alpha_l|\le\mathcal M\rho_0^{l+1}$. For
$|z|\ge\rho_0/s_j$, since $(l+1)\cdots(l+q)\le(l+q)^q$ and $s_j\le\frac12$
(Equation~\eqref{eq:sk-small}),
\begin{align*}
|F^{(q)}(z)|&\le\mathcal M\sum_{l\in 3\N}(l+q)^q\rho_0^{l+1}\Big(\frac{s_j}{\rho_0}\Big)^{l+1+q}
=\mathcal M\rho_0^{-q}s_j^{4+q}\sum_{p\ge1}(3p+q)^qs_j^{3(p-1)}\\
&\le\mathcal M\Big(\frac4a\Big)^qC\,s_j^{4+q}.
\end{align*}
In ring-$j$ units
$(E_g)_{[j]}(z)=\mu s_jV(z/s_j)$, since $A_{j+1}r_{j+1}=\mu s_jA_jr_j$, and
$|z/s_j|\ge\rho_0/s_j$ for $z\in N'$. By the scaling rule and
Equation~\eqref{eq:BS-complex-deriv}, on $N'$,
$|\nabla^q(E_g)_{[j]}(z)|=\mu s_j^{1-q}2^{q/2}|F^{(q)}(z/s_j)|$, which is at most
$C\mu s_j^5\mathcal M$. Summing over $q\le M_h$ gives
$\|(E_g)_{[j]}\|_{C^{M_h}(N')}\lesssim\varepsilon_j\mathcal M$,
which gives (1).

\emph{Proof of (2).} Recall that $w_{j+1}(t,\cdot)$ is already expressed in
ring-$(j+1)$ spatial and vorticity units; its first argument $t$ is physical
time. The corresponding perturbation in physical coordinates is
$A_{j+1}w_{j+1}(t,\cdot/r_{j+1})$, whose velocity field is the $V$ of the
statement. By homogeneity,
\[
V_{[j+1]}(z)=\big(\KBS\ast w_{j+1}(t,\cdot)\big)(z).
\]
By Lemma~\ref{lem:regions-basic}(4), $w_{j+1}(t,\cdot)$ is $3$-fold
rotationally symmetric, supported in $N\subset\overline B_2$, and has zero
integral. As in (1), now with $w_{j+1}(t,\cdot)$ in place of $\tilde g$, the
complex form of $V_{[j+1]}$ is
$F(z)=\sum_{l\in 3\N}\alpha_lz^{-l-1}$ for $|z|>2$. Here
$|\alpha_l|\le\frac1{2\pi}2^l\|w_{j+1}(t)\|_{L^1}$, because $|\zeta|\le2$ on $N$,
and the zero integral removes the term $l=0$.
For $|z|\ge\rho_0/s_j$ we have $2/|z|\le8s_j/a\le\frac14$ by
Equation~\eqref{eq:sk-small}, and
\[
\begin{aligned}
|F^{(q)}(z)|&\le\frac{\|w_{j+1}(t)\|_{L^1}}{2\pi}\sum_{l\in 3\N}(l+q)^q2^l\Big(\frac{4s_j}a\Big)^{l+1+q}\\
&\le\frac{\|w_{j+1}(t)\|_{L^1}}{2\pi}\Big(\frac4a\Big)^{1+q}\Big(\frac8a\Big)^3C\,s_j^{4+q}.
\end{aligned}
\]
The second inequality follows, as in (1), by writing $l=3p$ and using
$8s_j/a\le\frac12$. Passing to ring-$j$ units as in (1),
$|\nabla^qV_{[j]}(y)|\lesssim\varepsilon_j\|w_{j+1}(t)\|_{L^1}$
on $N'$, and summing over $q\le M_h$ proves (2).

\emph{Proof of (3).} By Lemma~\ref{lem:regions-basic}(3), $g=\omega_{\ge j+1}(t)$
is $3$-fold rotationally symmetric and supported in $D_{j+1}(t)\subset\overline B_R$, and by
Lemma~\ref{lem:speed}, $\|g\|_{L^1}\le2\|\Phi\|_{L^1}A_{j+1}r_{j+1}^2\le2\pi
A_{j+1}r_{j+1}^2$, since $0\le\Phi<1$ and $\supp\Phi\subset\overline B_1$. Every
$x\in N'_j$ has $|x|\ge ar_j/4\ge2R$. For the field $E$ of $g$ minus its
monopole, Lemma~\ref{lem:ext-multipole} with $\varrho=R$ and the scaling rule
(with $c=1/(A_jr_j)$ and $c'=r_j$) give, for $y\in N'$,
\begin{align*}
|\nabla^qE_{[j]}(y)|&=\frac{r_j^{q}}{A_jr_j}|\nabla^qE(r_jy)|
\le C_{3,q}\frac{2\pi A_{j+1}r_{j+1}^2R^3}{A_jr_j^5}|y|^{-4-q}\\
&\le2\pi C_{3,q}\Big(\frac4a\Big)^{4+q}\mu s_j^2\Big(\frac R{r_j}\Big)^3 .
\end{align*}
The second inequality uses $|y|\ge a/4$ on $N'$ and
$A_{j+1}r_{j+1}^2=\mu s_j^2A_jr_j^2$. Since $\mu s_j^2(R/r_j)^3=\mu s_j^5(R/r_{j+1})^3=\varepsilon_j(R/r_{j+1})^3$,
summing over $q\le M_h$ gives the bound
$\lesssim\varepsilon_j(R/r_{j+1})^3$, proving (3).

\emph{Proof of (4).} Now the field is generated by a ring outside the annulus
on which we evaluate it, so we expand it in nonnegative powers of $z$, in
ring-$i$ units. Recall that $w_i(t,\cdot)$ is already expressed in ring-$i$
spatial and vorticity units, while $t$ is physical time. The perturbation in
physical coordinates is $A_iw_i(t,\cdot/r_i)$, whose velocity field is the
$V_i$ of the statement. Put $s_{i,k}=r_k/r_i\le a/32$
(Equation~\eqref{eq:sk-small}). By homogeneity,
\[
(V_i)_{[i]}(z)=\big(\KBS\ast w_i(t,\cdot)\big)(z).
\]
The function $w_i(t,\cdot)$ is $3$-fold rotationally symmetric and supported in
$N=\{a/2\le|\zeta|\le2\}$. Inserting
$(z-\zeta)^{-1}=-\sum_{l\ge0}z^l\zeta^{-l-1}$, uniformly convergent for
$\zeta\in N$ and $z$ in compact subsets of $B_{a/2}$, in
Equation~\eqref{eq:BS-complex} gives the complex form of $(V_i)_{[i]}$, for $|z|<a/2$:
\[
\begin{gathered}
F(z)=-\sum_{l\ge0}\alpha_lz^l,\qquad
\alpha_l=\frac1{2\pi i}\int w_i(t,\zeta)\zeta^{-l-1}\dd A(\zeta),\\
|\alpha_l|\le\frac{\|w_i(t)\|_{L^1}}{2\pi}\Big(\frac2a\Big)^{l+1},
\end{gathered}
\]
where the bound on $\alpha_l$ uses $|\zeta|\ge a/2$ on $N$, and the
substitution $\zeta=\Orot\zeta'$ and the $3$-fold rotational symmetry of $w_i(t,\cdot)$ show
$\alpha_l=0$ unless $3\mid l+1$.
In ring-$i$ units the annulus $N'_k$ lies in $\overline B_{3s_{i,k}}$, where we
evaluate the series. Let $|z|\le3s_{i,k}$, so that $6s_{i,k}/a\le\frac14$. Writing $l=3p-1$ with
$p\ge1$, and using $l(l-1)\cdots(l-q+1)\le(3p)^q$, we get
\[
\begin{aligned}
|F^{(q)}(z)|&\le\frac{\|w_i(t)\|_{L^1}}{2\pi}\sum_{p\ge1}(3p)^q\Big(\frac2a\Big)^{3p}(3s_{i,k})^{3p-1-q}\\
&\le\frac{\|w_i(t)\|_{L^1}}{2\pi}\Big(\frac2a\Big)^{q+1}\Big(\frac{6s_{i,k}}a\Big)^{2-q}C;
\end{aligned}
\]
in the second inequality, the terms with $3p-1<q$ vanish on the left and are
nonnegative on the right, and the others are bounded by means of the identity
\[
(2/a)^{3p}(3s_{i,k})^{3p-1-q}=(2/a)^{q+1}(6s_{i,k}/a)^{2-q}(6s_{i,k}/a)^{3(p-1)}.
\]
Since $6s_{i,k}/a\le\frac14$, the last factor is at most $2^{-3(p-1)}$, and
$(3p)^q\le(3p+q)^q$; this gives the bound by $C$.
In ring-$k$ units,
$(V_i)_{[k]}(z)=\frac{A_ir_i}{A_kr_k}(V_i)_{[i]}(s_{i,k}z)=\mu^{i-k}s_{i,k}^{-1}(V_i)_{[i]}(s_{i,k}z)$,
and $|s_{i,k}z|\le3s_{i,k}$ for $z\in N'$. Hence, on $N'$, by the scaling rule,
Equation~\eqref{eq:BS-complex-deriv} and $\mu^{i-k}s_{i,k}=b_{i,k}$,
\[
|\nabla^q(V_i)_{[k]}(z)|=\mu^{i-k}s_{i,k}^{q-1}2^{q/2}|F^{(q)}(s_{i,k}z)|
\lesssim b_{i,k}\|w_i(t)\|_{L^1},
\]
and summing over $q\le M_h$ gives (4). For the last
assertion, Lemma~\ref{lem:regions-basic}(2), by which the sets $S_i(t)$ are
pairwise disjoint, and the definition of $w_i$ give
$\omega(t)\mathbf 1_{S_1(t)\cup\dots\cup S_{k-1}(t)}=\sum_{i<k}\big(A_i\Phi(\cdot/r_i)+A_iw_i(t,\cdot/r_i)\big)$.
Each $A_i\Phi(\cdot/r_i)$ is radial and supported in $\{|x|\ge ar_i\}$, so it
equals its angular average, and its field vanishes on $B_{ar_i}$ by
Lemma~\ref{lem:hollow}. Moreover
$N'_k\subset\overline B_{3r_k}\subset B_{ar_i}$, because $3r_k=3s_{i,k}r_i<ar_i$.
Summing the bounds over $i<k$ completes the proof.
\end{proof}

The next lemma gives the estimate of Lemma~\ref{lem:transfer}(3) in the form
used on the active interval in Proposition~\ref{prop:relay-step}, to verify the
inner-field bound required by Theorem~\ref{thm:active}. It makes the $L^1$
input explicit and identifies the estimated field: the velocity induced by
all the inner rings, with vorticity $\omega_{\ge k+1}(t)$, minus their monopole
is precisely $\hin_k$ in ring-$k$ units wherever $\hin_k$ is defined. The
inner vorticity in these units is
$\nu^{\mathrm{in}}_k(t,y)=A_k^{-1}\omega_{\ge k+1}(t,r_ky)$, as defined in
Definition~\ref{def:regions} at physical time $t$. The
same constant $\Cvel$ applies to this bound.

\begin{lemma}
\label{lem:inner-field}
Let $M_h\ge1$ be an integer and use the constant $\Cvel$ of
Lemma~\ref{lem:transfer}. Let $1\le k<n$ and $t\ge0$ be such that
$\omega_{\ge k+1}(t)$ is supported in $\overline B_R$ with $R\le\frac a8r_k$ and
$\|\omega_{\ge k+1}(t)\|_{L^1}\le\Cstar A_{k+1}r_{k+1}^2$. Let $E$ be the
velocity field of $\omega_{\ge k+1}(t)$ minus its monopole. Then
$\|E_{[k]}\|_{C^{M_h}(N')}\le\Cvel\,\mu s_k^2\,(R/r_k)^3$. If moreover
$D_{k+1}(t)\subset B_{ar_k/8}$, then $E_{[k]}=\hin_k(t)$ on $N'$.
\end{lemma}

\begin{proof}
The function $g=\omega_{\ge k+1}(t)$ is $3$-fold
rotationally symmetric by Lemma~\ref{lem:regions-basic}(3), and every $x\in N'_k$ has
$|x|\ge ar_k/4\ge2R$. The computation in the proof of
Lemma~\ref{lem:transfer}(3), with $\|g\|_{L^1}\le\Cstar A_{k+1}r_{k+1}^2$ in
place of $2\pi A_{k+1}r_{k+1}^2$, gives
\[
|\nabla^qE_{[k]}(y)|\lesssim\mu s_k^2(R/r_k)^3\qquad(y\in N',\ 0\le q\le M_h),
\]
with constants depending only on $a$ and $M_h$. Summing over $q$
gives the bound after enlarging the common constant $\Cvel$. Now let
$D_{k+1}(t)\subset B_{ar_k/8}$, so that $\hin_k(t)$ is defined. Since $\KBS$
is homogeneous of degree $-1$, $(\KBS\ast\omega_{\ge k+1}(t))_{[k]}=\KBS\ast\nu^{\mathrm{in}}_k$, and the
monopole of $\omega_{\ge k+1}(t)$ is, in ring-$k$ units,
$\frac1{2\pi}\big(\int\nu^{\mathrm{in}}_k\big)y^\perp/|y|^2$. By
Lemma~\ref{lem:inner-circ}, $\int\omega_{\ge k+1}(t)=\int\omega_{\ge k+1}(0)$, so,
by the definition of $G_k$ (Definition~\ref{def:regions}),
$\int\nu^{\mathrm{in}}_k=A_k^{-1}r_k^{-2}\int\omega_{\ge k+1}(0)=2\pi G_k$. Hence
$E_{[k]}=\KBS\ast\nu^{\mathrm{in}}_k-G_ky^\perp/|y|^2=\hin_k(t)$.
\end{proof}

\subsection{Confinement on a Window}

After ring $k+1$ saturates at the time $T_{k+1}$ of Definition~\ref{def:clock},
it is debris: the bootstrap of Section~\ref{sec:bootstrap} no
longer keeps its support near the annulus $N_{k+1}$. The next proposition shows
that, on a window $[T_{k+1},t_1]$, the material inside ring $k$ nevertheless
stays in the disk $\overline B_{ar_k/8}$, because its speed is at most of
order $A_{k+1}r_{k+1}$; the radius of the disk that contains it grows at most
linearly in the kick time $\tauK$. The bootstrap uses it to keep the inner field
$\hin_k$ defined after $T_{k+1}$ (Lemma~\ref{lem:geometry}\textup{(a)}), and it
bounds the normalized moment $q_{k+1}$, through which the material inside ring
$k$ drives ring $k$ (Section~\ref{sec:kick}). The number $\Lambda r_{k+1}$ in
the proposition is the radius $\varrho(t_1)$ reached at the end of the window.
Condition \textup{(iv)} asks that the disk of this radius lie in
$\overline B_{ar_k/8}$, and that on it the rings $j\le k$, under \textup{(ii)}
and \textup{(iii)}, contribute at most $\Cstar A_{k+1}r_{k+1}$ to the speed.

Here $t_1$ is an endpoint up to which we seek this confinement. In the
bootstrap, it is at most the current time, and the aim is to continue the
control through $T_k$, when ring $k$ saturates
(Lemma~\ref{lem:geometry}\textup{(a)}). The required duration in kick time is
\[
A_{k+1}(T_k-T_{k+1})=\Tpin+\mu\hat T_k.
\]
The later bounds $\ksh_k\asymp s_k^5$ give
$\hat T_k=O(1+\log(1/s_k))$. Thus this interval can be long in kick-time
units; its physical length is divided by $A_{k+1}$. Condition \textup{(iv)}
allows this logarithmic duration: the corresponding $\Lambda$ is
$O(1+\log(1/s_k))$, so both $\Lambda s_k$ and $\Lambda^2s_k$ tend to zero as
$s_k\to0$, with the other parameters fixed.

Assumption \textup{(ii)} keeps ring $k$ and the rings outside it away from
the origin throughout $[T_{k+1},t_1]$, allowing Lemma~\ref{lem:hollow} to
bound their velocity on the inner disk. In the bootstrap these rings are
still quiet or active up to $T_k$, and the stronger annular localization
assumed on the current time interval supplies \textup{(ii)}. The proposition
then supplies confinement of the already saturated inner material on that
same interval. These two controls are combined in the continuity argument
of Section~\ref{sec:bootstrap}.

\begin{proposition}
\label{prop:confine}
Let $1\le k<n$, assume that the clock is defined up to level $k+1$, put
$t_0=T_{k+1}$, and let $t_1\ge t_0$. Put $\tauK=A_{k+1}(t-t_0)$ and
$\Lambda=2+2\Cstar A_{k+1}(t_1-t_0)$. Assume:
\begin{enumerate}[label=\textup{(\roman*)}]
\item $D_{k+1}(t_0)\subset\overline{B}_{2r_{k+1}}$;
\item for every $t\in[t_0,t_1]$ and every $j\le k$,
  $S_j(t)\subset\{|x|\ge\frac a2r_j\}$;
\item $\sbar\le a/16$, $\sbar\le\frac12$ and $\mu^4\sbar^2\le\frac18$;
\item $\Lambda s_k\le a/8$ and
  $\frac83\big(\frac2a\big)^3\Lambda^2s_k\le\mu\Cstar$.
\end{enumerate}
Then for every $t\in[t_0,t_1]$:
\begin{enumerate}[label=\textup{(\alph*)}]
\item $D_{k+1}(t)\subset\overline{B}_{\varrho(t)}$ with
  $\varrho(t)=(2+2\Cstar\tauK)r_{k+1}\le\frac a8r_k$;
\item $|q_{k+1}(t)|\le\Cmom(1+\tauK)^3$, for an absolute constant $\Cmom$.
\end{enumerate}
None of the constants depends on $n$, $k$, $t_0$, $\mu$, $\sbar$, $\delta_0$ or
$\kseed$.
\end{proposition}

\begin{proof}
The constant $\Cstar$ is absolute by Lemma~\ref{lem:speed};
we choose $\Cmom$ below in terms of it. We first bound the
velocity on the disk $\overline B_{\Lambda r_{k+1}}$, and then show by a
continuity argument that $D_{k+1}(t)$ does not leave this disk before $t_1$.

\emph{Speed bound.} Fix $t\in[t_0,t_1]$ and $x\in\overline B_{\Lambda r_{k+1}}$;
by (iv), $|x|\le\Lambda s_kr_k\le\frac a8r_k$. By
Lemma~\ref{lem:regions-basic}(2),(3), $\omega(t)=\omega_{\ge k+1}(t)+\sum_{j\le k}g_j$
with $g_j=\omega(t)\mathbf 1_{S_j(t)}$, so
$v(t,x)=\KBS\ast\omega_{\ge k+1}(t)(x)+\sum_{j\le k}\KBS\ast g_j(x)$, and
$|\KBS\ast\omega_{\ge k+1}(t)(x)|\le\Cstar A_{k+1}r_{k+1}$ by Lemma~\ref{lem:speed}.
It remains to bound the fields of the rings $j\le k$. Let $j\le k$. By Lemma~\ref{lem:regions-basic}(3), $g_j\in C_c^\infty$ is
$3$-fold rotationally symmetric and $\|g_j\|_{L^1}=A_jr_j^2\|\Phi\|_{L^1}\le\pi A_jr_j^2$, since
$0\le\Phi<1$ and $\supp\Phi\subset\overline B_1$. By (ii), $g_j$ is supported in
$\{|y|\ge R_j\}$ with $R_j=\frac a2r_j$, and $|x|\le\frac a8r_k\le R_j/2$.
Lemma~\ref{lem:hollow} with $m_0=3$, $|x|\le\Lambda r_{k+1}$ and
$A_j=\mu^{j-k-1}A_{k+1}$ give
\begin{align*}
|\KBS\ast g_j(x)|&\le\frac4{3\pi}\pi A_jr_j^2|x|^2R_j^{-3}
=\frac43\Big(\frac2a\Big)^3A_jr_j\Big(\frac{|x|}{r_j}\Big)^2\\
&\le\frac43\Big(\frac2a\Big)^3\Lambda^2A_{k+1}r_{k+1}\,\mu^{j-k-1}\frac{r_{k+1}}{r_j}.
\end{align*}
By Equation~\eqref{eq:sk-small}, $r_{k+1}/r_j\le s_k\sbar^{k-j}$, and by (iii),
$\sbar^{k-j}\le2^{j-k}$; with $\mu^{j-k}\le1$,
$\mu^{j-k-1}r_{k+1}/r_j\le\mu^{-1}s_k2^{j-k}$, whose sum over
$j\le k$, a geometric series, is less than $2\mu^{-1}s_k$. By (iv),
$\sum_{j\le k}|\KBS\ast g_j(x)|\le\frac83(\frac2a)^3\Lambda^2s_k\mu^{-1}A_{k+1}r_{k+1}\le\Cstar A_{k+1}r_{k+1}$,
so
\[
|v(t,x)|\le2\Cstar A_{k+1}r_{k+1}\qquad(t\in[t_0,t_1],\ x\in\overline B_{\Lambda r_{k+1}}).\tag{$*$}
\]

\emph{Proof of (a).} By definition $D_{k+1}(t)=\flow_t(D_{k+1}(0))$, and
$D_{k+1}(0)=S_{k+1}\cup\dots\cup S_n$ is compact. Let $t^*$ be the largest
$t\in[t_0,t_1]$ such that $D_{k+1}(t')\subset\overline B_{\Lambda r_{k+1}}$ for
every $t'\in[t_0,t]$; it exists, because the set of such $t$ contains $t_0$, by
(i) and $\Lambda\ge2$, and is closed, because $(t,y)\mapsto\flow_t(y)$ is
continuous (Theorem~\ref{thm:smooth}\textup{(4)}). For $t\in[t_0,t^*]$ and
$y\in D_{k+1}(0)$, the point $\flow_{t'}(y)$ lies in $\overline B_{\Lambda r_{k+1}}$
for $t'\in[t_0,t]$, so \textup{($*$)} gives
$|\partial_{t'}\flow_{t'}(y)|=|v(t',\flow_{t'}(y))|\le2\Cstar A_{k+1}r_{k+1}$.
Integrating from $t_0$ to $t$, and using $|\flow_{t_0}(y)|\le2r_{k+1}$, which
follows from (i), we get
$|\flow_t(y)|\le|\flow_{t_0}(y)|+2\Cstar A_{k+1}r_{k+1}(t-t_0)\le\varrho(t)$. If
$t^*<t_1$, then $\varrho(t^*)<\varrho(t_1)=\Lambda r_{k+1}$, since $\Cstar>0$, and
the uniform continuity of $(t,y)\mapsto|\flow_t(y)|$ on
$[t_0,t_1]\times D_{k+1}(0)$ gives $\delta\in(0,t_1-t^*]$ with
$|\flow_t(y)|<\Lambda r_{k+1}$ for $t\in[t^*,t^*+\delta]$ and $y\in D_{k+1}(0)$,
contradicting the maximality of $t^*$. Hence $t^*=t_1$, and (a) holds, because
$\varrho(t)\le\varrho(t_1)=\Lambda s_kr_k\le\frac a8r_k$ by (iv).

\emph{Proof of (b).} By the definition of $Q_3$
(Definition~I.\ref{I-def:multipole}), (a) and Lemma~\ref{lem:speed},
$|Q_3[\omega_{\ge k+1}(t)]|\le\varrho(t)^3\|\omega_{\ge k+1}(t)\|_{L^1}\le\Cstar A_{k+1}r_{k+1}^2\allowbreak\varrho(t)^3$,
so by the definition of $q_{k+1}$ (Definition~\ref{def:regions}),
$|q_{k+1}(t)|\le\Cstar(\varrho(t)/r_{k+1})^3=\Cstar(2+2\Cstar\tauK)^3$, and
$2+2\Cstar\tauK\lesssim1+\tauK$. Choosing $\Cmom$ large enough proves (b).
\end{proof}

The following definition collects the support conditions under which the
perturbation $w_i$ of every ring that is not debris, and the external field
$h_j$ of every quiet ring, are defined (Definition~\ref{def:regions}), together
with the confinement of
Proposition~\ref{prop:confine}\textup{(a)}.
Lemma~\ref{lem:geometry}\textup{(c)} verifies these conditions in the bootstrap.
Since the clock is constructed recursively, they
must also cover rings whose activation or saturation time has not yet been
assigned.

\begin{definition}
\label{def:loc}
Let $t\ge0$, and consider the times assigned by the recursive clock of
Definition~\ref{def:clock}. A ring $j$ is \textbf{quiet at $t$} if its
activation time $\tsh_j$ has not been assigned, or if it has been assigned
and $t\le\tsh_j$. A ring $i$ is \textbf{debris at $t$} if its saturation time
$T_i$ has been assigned and $t>T_i$. Thus a ring with no assigned saturation
time is included among the rings that are not debris.

The \textbf{localization hypothesis} $(\mathrm{Loc})$ holds at $t$ if both
of the following conditions hold:
\begin{enumerate}[label=\textup{(\arabic*)}]
\item Every ring $i$ that is not debris at $t$ satisfies $S_i(t)\subset N_i$.
\item For every ring $j$ that is quiet at $t$, the outer rings and the inner
  material satisfy
  \[
  S_i(t)\cap B_{6r_j}=\emptyset\quad(i<j),\qquad
  D_{j+1}(t)\subset B_{ar_j/8}.
  \]
  If $T_{j+1}$ has been assigned and $t\ge T_{j+1}$, require in addition
  \[
  D_{j+1}(t)\subset\overline B_{\varrho_j(t)},\qquad
  \varrho_j(t)=\big(2+2\Cstar A_{j+1}(t-T_{j+1})\big)r_{j+1},
  \]
  the confinement radius of Proposition~\ref{prop:confine}\textup{(a)}.
\end{enumerate}
\end{definition}

Under support conditions of this kind the external field $h_j$ of ring $j$ is
small in $\sbar$. The bootstrap uses the next lemma to show that $\eta_j$, the
sum of the norms of $\hin_j$ and $\hout_j$ in $C^{M_h}(N')$, is at most $1$ for every ring $j<n$ that is not debris
(Lemma~\ref{lem:geometry}\textup{(b)}).

\begin{lemma}
\label{lem:field-bound}
Let $M_h\ge1$ be an integer and use the constant $\Cvel$ of
Lemma~\ref{lem:transfer}. Assume
$\sbar\le a/16$ and $\sbar\le\frac12$, let $1\le j<n$ and $t\ge0$, and
assume that $S_i(t)\subset N_i$ for every $i<j$, that $S_i(t)\cap B_{6r_j}=\emptyset$
for every $i<j$, and that $D_{j+1}(t)\subset B_{ar_j/8}$. Then
$\eta_j(t)\le \Cvel\big(\mu\sbar^2+\sbar\big)$.
\end{lemma}

\begin{proof}
All implicit constants in this proof depend only on $a$ and $M_h$, since
$\|\Phi\|_{L^1}\le\pi$. The common constant $\Cvel$ is chosen to dominate
these bounds as well. The hypotheses are those under which
Definition~\ref{def:regions} defines $h_j$ at time $t$, with
$\eta_j(t)=\etin_j(t)+\etout_j(t)$. The inner field is the non-monopole field of
$\omega_{\ge j+1}(t)$ in ring-$j$ units by
Lemma~\ref{lem:inner-field}; the hypotheses of that lemma hold with $R=ar_j/8$, by
Lemma~\ref{lem:speed}. By
Lemma~\ref{lem:transfer}(3) with $R=ar_j/8$,
$\etin_j(t)\lesssim\varepsilon_j(ar_j/(8r_{j+1}))^3=(a/8)^3\mu s_j^2\lesssim\mu\sbar^2$;
here we used $\varepsilon_j=\mu s_j^5$, $r_{j+1}=s_jr_j$ and $s_j\le\sbar$.
If $j=1$, then $\etout_1=0$, as there is no ring outside ring $1$. If $j\ge2$, Lemma~\ref{lem:transfer}(4) with $k=j$ and
Lemma~\ref{lem:regions-basic}(4) give
$\etout_j(t)\lesssim\sum_{i<j}b_{i,j}\|w_i(t)\|_{L^1}\le2\|\Phi\|_{L^1}\sum_{i<j}b_{i,j}$.
By Equation~\eqref{eq:sk-small}, $r_j/r_i\le\sbar^{j-i}$, and $\mu^{i-j}\le1$, so
$b_{i,j}\le\sbar^{j-i}$ and, as $\sbar\le\frac12$,
$\sum_{i<j}b_{i,j}\le\sum_{l\ge1}\sbar^l\le2\sbar$. Hence
$\etout_j(t)\lesssim\sbar$, and adding the two bounds proves
the lemma.
\end{proof}

\section{The Kick}
\label{sec:kick}

We compute the kick that ring $k$ receives from the rings inside it,
for $1\le k<n$. While ring $k$ is quiet, the inner field $\hin_k$ of
Definition~\ref{def:regions} drives its unstable linear response through the
mode-$3$ moment $q_{k+1}$ of the region $D_{k+1}(t)$.
Lemma~\ref{lem:coupling} gives the scalar Equation~\eqref{eq:kick} in the local
kick-time coordinate. Proposition~\ref{prop:kick} bounds the accumulated kick
from above and below. Proposition~\ref{prop:relay-step} uses these bounds
to start ring $k$, with the quantitative comparisons proved in
Lemmas~\ref{lem:handoff} and~\ref{lem:kicktrue}.

\begin{lemma}
\label{lem:coupling}
Assume \textup{(N1)}--\textup{(N3)}, fix $n\ge2$ and the objects of
Definitions~\ref{def:regions} and \ref{def:clock}, and let $1\le k<n$.
\begin{enumerate}[label=\textup{(\arabic*)}]
\item Let $t\ge0$ with $D_{k+1}(t)\subset B_{ar_k/8}$. Then
  $-\hin_k(t)\cdot\nabla\Phi$, extended by $0$ outside $N'$, lies in $X^s$
  for every $s\ge0$, and
  \[
  \ellst(G_k)\big(-\hin_k(t)\cdot\nabla\Phi\big)=\mu\,\beta_k\,q_{k+1}(t).
  \]
  Moreover, with $\nu^{\mathrm{in}}_k$ of Definition~\ref{def:regions} at the
  time $t$, the mode-$3$ profile of $\hin_k(t)\cdot\er$ is
  $\frac i2\pout_3[\nu^{\mathrm{in}}_k]\,r^{-4}$ for $a/4\le r\le3$.
\item Assume that the clock is defined up to level $k+1$ and that
  $D_{k+1}(t)\subset B_{ar_k/8}$ for every $t\in[0,\tsh_k]$. Then
  $t\mapsto-\hin_k(t)\cdot\nabla\Phi$ is continuous from $[0,\tsh_k]$ into
  $X^s$ for every $s$, and the physical-time linear response is
  \[
  \wlin_{\mathrm{in},k}(t)=A_k\int_0^t e^{A_k(t-s)\Lin(G_k)}\big(-\hin_k(s)\cdot\nabla\Phi\big)\dd s
  \quad(0\le t\le\tsh_k).
  \]
  For the scalar calculation, use the local coordinate
  $\tauK=A_{k+1}(t-T_{k+1})\in[\tau_0,\Tpin]$, where
  $\tau_0=-A_{k+1}T_{k+1}$, and define
  \[
  \begin{gathered}
  t(\tauK)=T_{k+1}+\tauK/A_{k+1},\qquad
  \tilde q(\tauK)=q_{k+1}(t(\tauK)),\\
  \cin(\tauK)=\ellst(G_k)\big(\wlin_{\mathrm{in},k}(t(\tauK))\big).
  \end{gathered}
  \]
  Then
  $\cin\in C^1([\tau_0,\Tpin])$ and, exactly,
  \begin{equation}
  \label{eq:kick}
  \frac{d\cin}{d\tauK}=\zeta_0\,\cin+\beta_k\,\tilde q,\qquad \zeta_0=\frac{\lambda_k}{\mu},\qquad
  \cin(\tau_0)=0 .
  \end{equation}
  The \textbf{kick amplitude}, the accumulated inner driving with free growth
  removed, is
  \begin{equation}
  \label{eq:kick-amplitude}
  \Ain(\tauK)=e^{-\zeta_0\tauK}\cin(\tauK)
  =\beta_k\int_{\tau_0}^{\tauK}e^{-\zeta_0s}\tilde q(s)\dd s.
  \end{equation}
  Thus its contribution at activation is
  \begin{equation}
  \label{eq:kick-coefficient}
  \cin_k:=\ellst(G_k)\big(\wlin_{\mathrm{in},k}(\tsh_k)\big)=e^{\zeta_0\Tpin}\Ain(\Tpin).
  \end{equation}
\item If, in addition to the hypotheses of \textup{(2)}, the rings $i<k$ lie
  outside $B_{6r_k}$ for every $t\in[0,\tsh_k]$, then
  $t\mapsto-\hout_k(t)\cdot\nabla\Phi$ is continuous from $[0,\tsh_k]$ into
  $X^s$ for every $s$, and $\wlin_{\mathrm{out},k}$, defined as
  $\wlin_{\mathrm{in},k}$ with $\hout_k$ in place of $\hin_k$, satisfies $\wlin_{\mathrm{in},k}+\wlin_{\mathrm{out},k}=\wlin$ on
  $[0,\tsh_k]$, where $\wlin$ is the physical-time linear response obtained
  from Theorem~\ref{thm:quiet} with $G=G_k$, using the local coordinate
  $\tau=A_kt$ and field $h(\tau)=h_k(\tau/A_k)$.
\end{enumerate}
\end{lemma}
\begin{proof}
\emph{(1).} This part makes precise the computation made before
Definition~I.\ref{I-def:eigendata}, where $b$ is introduced as the mode-$3$
driving term generated by an interior velocity field of unit outer moment, and
$\beta(G)$ as the component of that term along the unstable eigenfunction. By
Lemma~\ref{lem:regions-basic}(2),(3), the function
$\nu^{\mathrm{in}}_k$ of Definition~\ref{def:regions} at time $t$, the vorticity of
$D_{k+1}(t)$ in ring-$k$ units, is real,
$3$-fold rotationally symmetric, in $C_c^\infty(\R^2)$, and its support avoids $0$; by the
hypothesis $D_{k+1}(t)\subset B_{ar_k/8}$ the support is a compact subset of
$B_{a/8}$, so the largest modulus $\rho$ of a point of the support satisfies $\rho<a/8$. By
Lemma~I.\ref{I-lem:BS-basic}(d), $\KBS\ast\nu^{\mathrm{in}}_k=\nabla^\perp\psi$ with
$\psi=\frac1{2\pi}\log|\cdot|\ast\nu^{\mathrm{in}}_k$ smooth on $\{|y|>\rho\}$.
The monopole $G_ky^\perp/|y|^2$ is orthogonal to $\er$, so the polar form
$v\cdot\er=-r^{-1}\partial_\theta\psi$ of the Biot-Savart law
(Section~I.\ref{I-sec:prelim}) gives $\hin_k(t)\cdot\er=-r^{-1}\partial_\theta\psi$
on $N'$. Put $p_3=\pout_3[\nu^{\mathrm{in}}_k]$, the outer mode-$3$ moment of
Definition~I.\ref{I-def:multipole}. By Lemma~\ref{lem:exterior-mode}, at angular order $3$ and with $\varrho=\rho$, the angular mode $3$ of $\psi$ is
$-\frac16p_3r^{-3}$ for $r>\rho$, and by integration by parts in $\theta$ that of
$\partial_\theta\psi$ is $3i$ times that of $\psi$. Since
$N'=\{a/4\le|y|\le3\}\subset\{|y|>\rho\}$, the mode-$3$ profile of
$\hin_k(t)\cdot\er$ is $\frac i2p_3r^{-4}$ for $a/4\le r\le3$, the last claim
of (1). As $\Phi$ is radial and $\Phi'$ vanishes outside $[a,1]$,
$\hin_k(t)\cdot\nabla\Phi=\Phi'(r)\,\hin_k(t)\cdot\er$, so the mode-$3$
profile of $-\hin_k(t)\cdot\nabla\Phi$ is
$-\frac i2p_3r^{-4}\Phi'(r)=p_3\,b(r)$, with $b$ of Definition~I.\ref{I-def:eigendata}.
By Equation~\eqref{eq:BS-complex}, $\hin_k(t)$ is the restriction to
$N'$ of a field smooth on $\{|y|>\rho\}$, and it is $3$-fold rotationally symmetric because
$\nu^{\mathrm{in}}_k$ is; so $-\hin_k(t)\cdot\nabla\Phi$ is real, smooth,
$3$-fold rotationally symmetric and supported in $\{a\le|y|\le1\}\subset N$, hence in $X^s$
for every $s$. By Definition~I.\ref{I-def:eigendata}, $\ellst(G_k)(w)=n_*(G_k)\,\ellst(w_3)$
depends linearly on the mode-$3$ profile $w_3$ of $w$ only, and
$\beta(G_k)=\ellst(G_k)(b\,e^{3i\theta})=n_*(G_k)\ellst(b)$; so
$\ellst(G_k)(-\hin_k(t)\cdot\nabla\Phi)=p_3\,\beta(G_k)$. It remains to
express $p_3$ through the normalized moment $q_{k+1}(t)$ of
Definition~\ref{def:regions}. Since $\omega_{\ge k+1}(t)$ is real-valued,
the change of variables $x=r_ky$ gives
\[
p_3=\frac1{2\pi A_k}\int\omega_{\ge k+1}(t,r_ky)\,\bar y^3\dd y
=\frac{\overline{Q_3[\omega_{\ge k+1}(t)]}}{2\pi A_kr_k^5}
=\frac{A_{k+1}r_{k+1}^5}{2\pi A_kr_k^5}\,q_{k+1}(t)=\frac{\mu s_k^5}{2\pi}\,q_{k+1}(t),
\]
and (1) follows from $\beta_k=\beta(G_k)s_k^5/(2\pi)$
(Definition~\ref{def:regions}).

\emph{(2).} By the hypothesis and the paragraph after
Lemma~\ref{lem:regions-basic}, $t\mapsto\hin_k(t)$ is continuous from
$[0,\tsh_k]$ into $C^M(N')$ for every integer $M\ge0$. Put
$g(t)=-\hin_k(t)\cdot\nabla\Phi$. By (1), $g(t)\in X^s$ for every $s$,
and by Lemma~I.\ref{I-lem:uniform},
$\|g(t)-g(t')\|_{H^M}\le\Cnab(M)\|\hin_k(t)-\hin_k(t')\|_{C^M(N')}$;
so $g$ is continuous into $X^s$ for every $s$. Since $e^{\theta\Lin(G_k)}$ is a
$C_0$-semigroup on $X^s$ (Lemma~I.\ref{I-lem:semigroup-construct}), the integrand
$s\mapsto e^{A_k(t-s)\Lin(G_k)}g(s)$ is continuous into $X^s$, so
$\wlin_{\mathrm{in},k}$ is well defined.

For the scalar calculation we change from physical time to the local kick
coordinate. Since $A_{k+1}=\mu A_k$ by Equation~\eqref{eq:datum} and
$\tsh_k=T_{k+1}+\Tpin/A_{k+1}$ (Definition~\ref{def:clock}), $\tauK$ runs over
$[\tau_0,\Tpin]$ as $t$ runs over $[0,\tsh_k]$. As $\ellst(G_k)$ is bounded on
$X^0$ (Lemma~I.\ref{I-lem:uniform}), it commutes with the integral defining
$\wlin_{\mathrm{in},k}$. Proposition~I.\ref{I-prop:semigroup} turns the semigroup
under this functional into multiplication by $e^{A_k\lambda_k(t-s)}$.
Using (1), followed by the substitution $r=A_{k+1}(s-T_{k+1})$, gives
\begin{align*}
\cin(\tauK)
&=A_k\int_0^{t(\tauK)}e^{A_k\lambda_k(t(\tauK)-s)}\ellst(G_k)\big(g(s)\big)\dd s\\
&=A_k\mu\beta_k\int_0^{t(\tauK)}e^{A_k\lambda_k(t(\tauK)-s)}q_{k+1}(s)\dd s\\
&=\beta_k\int_{\tau_0}^{\tauK}e^{\zeta_0(\tauK-r)}\tilde q(r)\dd r.
\end{align*}
Here the change of variables uses $\dd s=\dd r/A_{k+1}$ and
$A_k\mu/A_{k+1}=1$.
By Equation~\eqref{eq:omega-ge} and area preservation,
$Q_3[\omega_{\ge k+1}(t)]=\int\sum_{i>k}f^{(n)}_i(x')\,\flow_t(x')^3\dd x'$.
The integrand is supported in the compact set $D_{k+1}(0)$, is bounded uniformly
for $t$ in compact intervals, and is continuous in $t$, as
$(t,x')\mapsto\flow_t(x')$ is continuous (Theorem~\ref{thm:smooth}\textup{(4)}).
Thus $q_{k+1}$, and hence its pullback $\tilde q$, are continuous by dominated
convergence. Taking $e^{\zeta_0\tauK}$ out of the last integral and applying the
fundamental theorem of calculus gives $\cin\in C^1([\tau_0,\Tpin])$,
$\cin(\tau_0)=0$, and Equation~\eqref{eq:kick}. At physical time $t=\tsh_k$,
$\tauK=\Tpin$, so
$\cin_k=\cin(\Tpin)=e^{\zeta_0\Tpin}\Ain(\Tpin)$.

\emph{(3).} By the hypotheses and the paragraph after
Lemma~\ref{lem:regions-basic}, $t\mapsto\hout_k(t)$ is continuous from
$[0,\tsh_k]$ into $C^M(N')$ for every $M$. By Lemma~\ref{lem:regions-basic}(3),
$\nu^{\mathrm{out}}_k$ is real and $3$-fold rotationally symmetric, and, as the
rings $i<k$ lie outside $B_{6r_k}$, it vanishes on $B_6$. Hence
$\hout_k(t)$ is real, $3$-fold rotationally symmetric and smooth on $B_6$, by
Equation~\eqref{eq:BS-complex}, and $B_6$ contains $N'$. As in (1),
$-\hout_k(t)\cdot\nabla\Phi\in X^s$ for every $s$, and this function is
continuous in $t$ into $X^s$ by Lemma~I.\ref{I-lem:uniform}, as in (2).
For the application of Section~\ref{sec:single}, use the local coordinate
$\tau=A_kt$ and the explicit pullback $h(\tau)=h_k(\tau/A_k)$.
The discussion after Lemma~\ref{lem:regions-basic} gives the hypotheses of
Lemma~\ref{lem:reduction} on $[0,A_k\tsh_k]$ with $G=G_k$ and this $h$.
Pulling the resulting linear response back to physical time gives
\[
\wlin(t)=A_k\int_0^t e^{A_k(t-s)\Lin(G_k)}
\big(-h_k(s)\cdot\nabla\Phi\big)\dd s.
\]
Since $h_k=\hin_k+\hout_k$, linearity of this integral proves
$\wlin=\wlin_{\mathrm{in},k}+\wlin_{\mathrm{out},k}$.
\end{proof}

Proposition~\ref{prop:kick} works in the following setting, which isolates
Equation~\eqref{eq:kick}. Let $\zeta_0,\beta\in\C$, put $x=\Re\zeta_0$, let
$\tau_0<\tau_1$ be real, put $I=[\tau_0,\tau_1]$, let
$q\colon I\to\C$ be continuous, let $c$ solve $c'=\zeta_0c+\beta q$ on $I$ with
$c(\tau_0)=0$, and put $A(\tau)=e^{-\zeta_0\tau}c(\tau)$. The factor
$e^{-\zeta_0\tau}$ removes the free growth of $c$, so that $A$ is a weighted
integral of the drive $q$, as the next display shows. Since
$A'=e^{-\zeta_0\tau}(c'-\zeta_0c)=\beta e^{-\zeta_0\tau}q$ and $A(\tau_0)=0$,
\begin{equation}
\label{eq:kick-A}
A(\tau)=\beta\int_{\tau_0}^\tau e^{-\zeta_0s}\,q(s)\dd s\qquad(\tau\in I).
\end{equation}
Under the hypotheses of Lemma~\ref{lem:coupling}\textup{(2)}, the functions
$\cin$, $\tilde q$ and $\Ain$, with $\zeta_0=\lambda_k/\mu$, $\beta=\beta_k$ and
$\tau_1=\Tpin$, are an instance of this setting; here
$\tau_0=-A_{k+1}T_{k+1}<0<\Tpin$ by Lemma~\ref{lem:clock}(1) and
Definition~\ref{def:clock}.

Write $\lamref=\lamst(0)$ and $\sigma=\Re\lamref>0$
(Definition~I.\ref{I-def:eigendata}, Proposition~I.\ref{I-prop:U}(a)).
The reference eigenvalue $\lamref$ is a fixed complex number throughout this
section and the subsequent construction. We also use
$I_p(x)=\int_0^\infty e^{-xs}(1+s)^p\dd s$ for $x>0$
(Equation~(I.\ref{I-eq:Ip})); it decreases as $x$ increases.

The next proposition is an elementary estimate for Equation~\eqref{eq:kick-A}.
Its purpose is to give a quantitative lower bound on the kick, which is used
to start the next ring. The drive is close to an exponential on
$[-T_g,0]$. After the receiving ring's free growth is removed, this gives a
main contribution of size $|\beta\qs|/(|\lamref|(\mu-1))$, because the
modulus of the difference between the two growth exponents is comparable to
$|\lamref|(\mu-1)/\mu$.

The assumptions below make this main contribution survive cancellation.
The length condition \textup{(Kick2)} ensures that the exponential interval
captures a definite part of the accumulated kick. The other three bounds
control, respectively, the later polynomially bounded drive
\textup{(Kick1)}, the error in the exponential approximation
\textup{(Kick3)}, and the earlier history \textup{(Kick4)}. Their contributions
are estimated in absolute value and kept small compared with the main term.

\begin{samepage}
\begin{proposition}
\label{prop:kick}
Assume \textup{(N1)}--\textup{(N3)} and the scalar setting above. Let
$1<\mu\le2$, $G,G'\in[0,\Gone]$, and set
\[
\zeta_0=\frac{\lamst(G)}\mu,\qquad
\lambda_{\mathrm{in}}=\lamst(G'),\qquad
w=\lambda_{\mathrm{in}}-\zeta_0,\qquad x=\Re\zeta_0.
\]
Here $\zeta_0$ is the receiving exponent, $\lambda_{\mathrm{in}}$ is the driving
exponent, and $w$ is their difference.
Assume
\[
\max\big\{|\lamst(G)-\lamref|,|\lambda_{\mathrm{in}}-\lamref|\big\}
\le\frac{\sigma(\mu-1)}{8\mu}.
\]
Let $\tau_0<-T_g<0<\tau_1$, $\bar Q>0$, and let $p\ge0$ be an integer.
Suppose there are $\qs\in\C\setminus\{0\}$ and nonnegative measurable functions
$\err$ on $[-T_g,0]$ and $Q_{\mathrm{pre}}$ on $[\tau_0,-T_g]$ such that
\[
\begin{aligned}
|q(\tau)-\qs e^{\lambda_{\mathrm{in}}\tau}|&\le\err(\tau)
&&(-T_g\le\tau\le0),&&\textup{(Q1)}\\
|q(\tau)|&\le Q_{\mathrm{pre}}(\tau)
&&(\tau_0\le\tau\le-T_g),&&\textup{(Q2)}\\
|q(\tau)|&\le\bar Q(1+\tau)^p
&&(0\le\tau\le\tau_1).&&\textup{(Q3)}
\end{aligned}
\]
Define the weighted errors
\[
\Errg=\int_{-T_g}^0e^{-x\tau}\err(\tau)\dd\tau,\qquad
\Errpre=\int_{\tau_0}^{-T_g}e^{-x\tau}Q_{\mathrm{pre}}(\tau)\dd\tau.
\]
Require the following length and error bounds:
\[
\begin{aligned}
\bar Q\,\Ip{p}{\sigma/4}&\le\frac{|\qs|}{10|\lamref|(\mu-1)},
&&\textup{(Kick1)}\\
e^{-\Re w\,T_g}&\le\tfrac12,
&&\textup{(Kick2)}\\
\max\{\Errg,\Errpre\}&\le\frac{|\qs|}{20|\lamref|(\mu-1)}.
&&\textup{(Kick3)--(Kick4)}
\end{aligned}
\]
Then $A(\tau_1)$ satisfies
\[
\frac{|\beta|\,|\qs|}{10|\lamref|(\mu-1)}
\le|A(\tau_1)|\le
\frac{3|\beta|\,|\qs|\,\mu}{|\lamref|(\mu-1)}.
\]
\end{proposition}
\end{samepage}

In Lemma~\ref{lem:kicktrue}, $G'=G_{k+1}$ and $T_g=\hat T_{k+1}$.
For $k+1<n$, the three intervals $[\tau_0,-T_g]$, $[-T_g,0]$ and
$[0,\tau_1]$ describe the quiet and active intervals of ring $k+1$, followed
by the pinning interval of ring $k$, in the local kick coordinate. At the
seeded level $k+1=n$, that lemma uses a slightly shorter exponential
interval.

\begin{proof}
We first compare $w$ with $\lamref$ and bound $x$ from below. Put
$e_\mu=\frac{\mu-1}{8\mu}\le\frac1{16}$ and write $\lambda_{\mathrm{in}}=\lamref+d_1$,
$\lamst(G)=\lamref+d_2$ with $|d_1|,|d_2|\le e_\mu\sigma$. Then
$w=\lamref\frac{\mu-1}\mu+d_1-\frac{d_2}\mu$ with
$|d_1-\frac{d_2}\mu|\le2e_\mu\sigma=\frac{\sigma(\mu-1)}{4\mu}$, so, since
$\Re \lamref=\sigma\le|\lamref|$,
\[
\Re w\ge\frac{\sigma(\mu-1)}\mu-\frac{\sigma(\mu-1)}{4\mu}=\frac34\sigma\frac{\mu-1}\mu>0,\qquad
\frac34\,\frac{|\lamref|(\mu-1)}\mu\le|w|\le\frac54\,\frac{|\lamref|(\mu-1)}\mu .
\]
As $\mu\le2$, $x=\Re\lamst(G)/\mu\ge(1-e_\mu)\sigma/\mu\ge\frac{15}{32}\sigma\ge\frac\sigma4>0$,
and as $\Ip{p}{\cdot}$ is decreasing, $\Ip{p}{x}\le\Ip{p}{\sigma/4}$.

Next we split the integral of Equation~\eqref{eq:kick-A} at $-T_g$ and $0$,
according to the three hypotheses on $q$. Put $\mathsf K=\beta\qs(1-e^{-wT_g})/w$,
the main term, and let $\tau=\tau_1$. Then
\[
A(\tau)=\beta\int_{\tau_0}^{-T_g}e^{-\zeta_0s}q(s)\dd s+\beta\int_{-T_g}^0e^{-\zeta_0s}q(s)\dd s
+\beta\int_0^\tau e^{-\zeta_0s}q(s)\dd s ,
\]
and $|e^{-\zeta_0s}|=e^{-xs}$. By (Q2), the
first term on the right has modulus at most $|\beta|\Errpre$. Since
\[
\int_{-T_g}^0e^{-\zeta_0s}\qs e^{\lambda_{\mathrm{in}}s}\dd s=\qs\int_{-T_g}^0e^{ws}\dd s=\qs(1-e^{-wT_g})/w,
\]
the second term is $\mathsf K$ plus a remainder of modulus at most
$|\beta|\Errg$, by (Q1). By (Q3), $x>0$ and Equation~(I.\ref{I-eq:Ip}), the third
has modulus at most
$|\beta|\bar Q\int_0^\infty e^{-xs}(1+s)^p\dd s=|\beta|\bar Q\,\Ip{p}{x}$.

Finally set $S=|\beta||\qs|/(|\lamref|(\mu-1))$, the natural size of the kick.
By (Kick2) and $|e^{-wT_g}|=e^{-\Re w\,T_g}$, $|1-e^{-wT_g}|\in[\frac12,\frac32]$,
and the bounds on $|w|$ from the first paragraph give
\[
|\mathsf K|\ge\tfrac12\cdot\tfrac45\mu S=\tfrac25\mu S,\qquad
|\mathsf K|\le\tfrac32\cdot\tfrac43\mu S=2\mu S.
\]
By (Kick1) and $\Ip{p}{x}\le\Ip{p}{\sigma/4}$,
$|\beta|\bar Q\,\Ip{p}{x}\le S/10$, and by (Kick3) and (Kick4),
$|\beta|\Errg,|\beta|\Errpre\le S/20$. Thus the three remainders have total
modulus at most $S/5$. The triangle inequality and $\mu>1$ give
\begin{gather*}
|A(\tau)|\ge\big(\tfrac25\mu-\tfrac15\big)S\ge\tfrac15 S\ge\tfrac1{10}S,\\
|A(\tau)|\le\big(2\mu+\tfrac15\big)S\le3\mu S.\qedhere
\end{gather*}
\end{proof}

\section{One Step of the Relay}
\label{sec:relay}

When ring $k+1$ saturates, we have controlled the solution through
$T_{k+1}$. The next step extends these bounds across the pinning interval
$[T_{k+1},\tsh_k]$ and the active interval $[\tsh_k,T_k]$. The moment of ring
$k+1$ supplies the initial amplitude of ring $k$. During the active interval,
we confine the material inside ring $k$ and bound the feedback from the outer
quiet rings. At saturation, the same estimates are available for the next
step outward.

Assume (N1)--(N3), and fix parameters satisfying (P0)--(P6), listed in
Appendix~\ref{sec:param-lists}. Section~\ref{sec:param} shows that they can
be chosen in order and give (SA) and (Det). Fix $n\ge2$ and $A_1>0$, write
$\omega=\omega^{(n)}$, and use the rings and clock of
Definitions~\ref{def:datum}--\ref{def:clock}. All ring quantities use physical
time. The reference growth rate is $\sigma=\sigma(0)$, while ring $i$ grows
at rate $A_i\sigma_i$ with $\sigma_i=\sigma(G_i)$.

With $C_0$ from (P4), let $\rhohat_{C_0,G_i}$ be the transient of
Theorem~\ref{thm:active} in ring-$i$ time. For a ring with a positive
assigned amplitude $\ksh_i$, the envelope and transient in physical time are
\begin{equation}
\label{eq:relay-envelope}
\hat a_i(v)=\ksh_i e^{A_i\sigma_i(v-\tsh_i)},\qquad
\rhohat_i(v)=\rhohat_{C_0,G_i}(A_i(v-\tsh_i)),\qquad v\ge\tsh_i.
\end{equation}
In particular, $\rhohat_i(\tsh_i)=C_0$.

\begin{definition}
\label{def:stop}
If the clock is defined up to level $1$, put $\Tend=T_1$. Otherwise let
$k_{\mathrm s}$ be the level $k\in\{1,\dots,n\}$ such that the clock is defined up
to level $k+1$ but not up to level $k$ (for $k=n$ the first condition is void),
and put $\Tend=\tsh_{k_{\mathrm s}}$. We then say that the clock stops at
level $k_{\mathrm s}$.
\end{definition}

The clock is constructed from $n$ downward. If it stops at $k_{\mathrm s}$,
its activation times are assigned for $i\ge k_{\mathrm s}$ and its saturation
times for $i>k_{\mathrm s}$. In particular, the failed level has an assigned
activation time. We use the quiet, active and debris conventions of
Definition~\ref{def:loc}, including quiet at activation and active at
saturation. Any ring farther out than a quiet or active ring is quiet by
Lemma~\ref{lem:clock}(1).

For $t\in[0,\Tend]$, the bootstrap assumptions (BA)$(t)$ are the following
bounds at every $v\in[0,t]$.
\begin{description}[font=\normalfont]
\item[(BA1)] Rings that have not finished their active interval stay in their
prescribed annuli:
\begin{equation}
\label{eq:BA-support}
S_j(v)\subset\Nflat_j\qquad\text{if ring $j$ is not debris at $v$}.
\end{equation}
\item[(BA2)] The accumulated driving of each quiet ring stays on the scale
$\varepsilon_j$. For every quiet ring $j<n$, $h_j$ is defined on
$[0,v]$ and
\begin{equation}
\label{eq:BA-quiet}
\maj_j(v)\le2\Abar\varepsilon_j.
\end{equation}
\item[(BA3)] After activation, the perturbation grows no faster than its
unstable envelope, with the initial transient allowed. For every ring $i$
active at $v>\tsh_i$,
\begin{equation}
\label{eq:BA-active}
\|w_i(v)\|_{H^s}\le2C_2\hat a_i(v)(1+\rhohat_i(v)).
\end{equation}
\end{description}
Here $\varepsilon_j=\mu s_j^5$ is the scale of the mode-$3$ driving from
the next inner ring, and $\Abar$ is the fixed constant bounding the quiet
response at that scale.
By the \emph{improved bounds through $t$} we mean the following estimates
at every $v\in[0,t]$:
\begin{align}
S_j(v)&\subset\NN_j
&&\text{if ring $j$ is not debris at $v$},\label{eq:improved-support}\\
\maj_j(v)&\le\Abar\varepsilon_j
&&\text{if ring $j<n$ is quiet at $v$},\label{eq:improved-quiet}\\
\|w_i(v)\|_{H^s}&\le C_2\hat a_i(v)(1+\rhohat_i(v))
&&\text{if ring $i$ is active or }v=\tsh_i.\label{eq:improved-active}
\end{align}
In \eqref{eq:improved-quiet}, $h_j$ must still be defined on $[0,v]$.
Equation~\eqref{eq:improved-active} also includes every assigned activation
time in $[0,t]$, with $\ksh_i$ defined and positive. The seed or the hand-off
supplies this amplitude and bound before the clock assigns $T_i$. Since
$\NN\subset\operatorname{int}\Nflat$ and the size bounds are positive,
these strict improvements of (BA1)--(BA3) leave room for continuation in
physical time.

\begin{proposition}[One relay step]
\label{prop:relay-step}
Let $0\le t\le\Tend$ and assume (BA)$(t)$. For $1\le k<n$, suppose the clock is
defined up to level $k+1$, with $T_{k+1}\le t$, and that
\eqref{eq:improved-support}--\eqref{eq:improved-active} hold through $T_{k+1}$.
Assume the conclusions in item~(c) for ring $k+1$ through $T_{k+1}$,
and the following bounds on the amplitudes at the completed activations:
\begin{equation}
\label{eq:relay-amplitudes}
\tfrac12\kseed\le\ksh_n\le2\kseed,\qquad
\Klo s_i^5\le\ksh_i\le\Kup s_i^5\quad(k<i<n).
\end{equation}
Then the following hold.
\begin{enumerate}[label=\textup{(\alph*)}]
\item The improved bounds extend through $\min(t,\tsh_k)$.
\item If $\tsh_k\le t$, the clock is defined up to level $k$ and the
improved bounds hold through $b=\min(t,T_k)$. At activation,
\begin{equation}
\label{eq:relay-startup}
\begin{aligned}
\Klo s_k^5&\le\ksh_k\le\Kup s_k^5,
&\quad \maj_k(\tsh_k)&\le K_1\ksh_k,\\
\|w_k(\tsh_k)\|_{H^{M-1}}&\le K_1 C_Q\ksh_k,
& S_k(\tsh_k)&\subset\NNm_k.
\end{aligned}
\end{equation}
\item In case (b), retain the initial Duhamel bound in
\eqref{eq:relay-startup}. If $b>\tsh_k$, Theorem~\ref{thm:active} and
Corollary~\ref{cor:active-moment} apply with $G=G_k$, $\kappa=\ksh_k$ and
$\eta_b=\eta_{b,k}$ of \eqref{eq:etab} (with $\eta_{b,1}=0$).
More explicitly, for $\tsh_k\le v\le b$ and $u=A_k(v-\tsh_k)$, besides
\eqref{eq:improved-support} and \eqref{eq:improved-active} for ring $k$,
\begin{equation}
\label{eq:relay-active}
\begin{aligned}
&\big\|w_k(v)-2\Re\big(c_k(v)e^{\lambda_k u}\phist(G_k)\big)\big\|_{H^s}\\
&\qquad\le C_2\hat a_k(v)
\big[\rhohat_k(v)+\hat a_k(v)+(\ksh_k)^{1/2}+\eta_{b,k}\big],
\end{aligned}
\end{equation}
\begin{equation}
\label{eq:relay-active-moment}
\begin{aligned}
&\big|2\pi\pout_3[w_k(v)]-2\pi\pout_*(G_k)c_k(v)e^{\lambda_k u}\big|\\
&\qquad\le2\pi C_3\hat a_k(v)
\big[\rhohat_k(v)+\hat a_k(v)+(\ksh_k)^{1/2}+\eta_{b,k}\big].
\end{aligned}
\end{equation}
The inner velocity field satisfies
\begin{equation}
\label{eq:relay-inner-field}
\|\hin_k(v)\|_{C^{M_h}(N')}\le C_0\ksh_k\jap{u}^3
\qquad(\tsh_k\le v\le b).
\end{equation}
All hypotheses of Lemma~\ref{lem:drift} hold with endpoint $b$ and
$\eta=\eta_{b,k}$, and its two smallness conditions are
\begin{equation}
\label{eq:relay-drift-small}
(1+\Tpin)^3e^{-\sigma\Tpin/4}\le\frac1{64},\qquad
\Cls\eta_{b,k}\hat T_k\le\frac1{32}.
\end{equation}
In particular, at $b=T_k$,
\begin{equation}
\label{eq:relay-saturation}
|c_k(T_k)|e^{\sigma_k\hat T_k}
\in[\tfrac{15}{16}\delta_0,\tfrac{17}{16}\delta_0].
\end{equation}
\end{enumerate}
On the portion of the pinning and active intervals reached by time $t$,
the confinement and moment estimates
\eqref{eq:geometry-confinement}--\eqref{eq:geometry-moment} hold for $D_{k+1}$.

For the seeded ring $k=n$, assume instead that the clock is defined up to
$n$, the improved bounds hold at $0$, and
\begin{equation}
\label{eq:relay-seed}
\tfrac12\kseed\le\ksh_n\le2\kseed,\qquad
S_n(0)\subset\NNm_n,\qquad
\|w_n(0)\|_{H^{M-1}}\le C_0\ksh_n.
\end{equation}
Then \eqref{eq:improved-support}--\eqref{eq:improved-active} hold through
$\min(t,T_n)$, and the conclusions in (c) hold if this interval has positive length, with
$\maj_n(0)=0$ as the initial Duhamel bound.
\end{proposition}

\begin{figure}[H]
\centering
\begin{tikzpicture}[x=1cm,y=1cm,font=\small]
\draw[|-|,thick] (0,0.75) -- (4.1,0.75)
  node[midway,above] {bounds already proved};
\draw[|-|,thick] (4.1,0.75) -- (8.8,0.75)
  node[midway,above,align=center] {confinement and\\quiet response};
\draw[|-|,thick] (8.8,0.75) -- (13.4,0.75)
  node[midway,above,align=center] {active growth and\\weak feedback};
\draw[->] (-0.15,0) -- (13.9,0) node[right] {$v$};
\foreach \x in {0,4.1,8.8,13.4}
  \draw (\x,0.10) -- (\x,-0.10);
\node[below] at (0,-0.12) {$0$};
\node[below] at (4.1,-0.12) {$T_{k+1}$};
\node[below,align=center] at (8.8,-0.12) {$\tsh_k$\\hand-off to ring $k$};
\node[below,align=center] at (13.4,-0.12) {$T_k$\\next step outward};
\end{tikzpicture}
\caption{A relay step in physical time. The part already proved supplies the
moment that starts ring $k$. At $T_k$ the estimates have been extended
through the active interval of ring $k$.}
\label{fig:relay-step}
\end{figure}
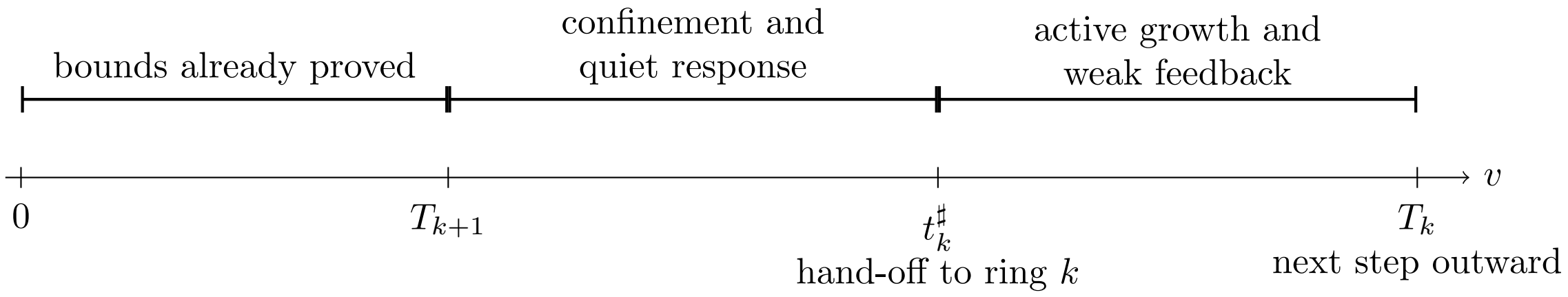

Figure~\ref{fig:relay-step} shows the extension from a completed inner
interval. If $t$ lies within the step, the proof stops at $t$.

\begin{proof}
For the seed with $t=0$, the conclusion is the assumed
\eqref{eq:improved-support}--\eqref{eq:improved-active} at $0$.

\phantomsection\label{par:relay-geometry}
\emph{Confinement.}
At an endpoint $b\le t$ where $\ksh_i\ge\Klo s_i^5$ is known for every
$i<n$ activated strictly before $b$, (BA1) and Lemma~\ref{lem:geometry}
give the support separations
\eqref{eq:geometry-quiet-support}--\eqref{eq:geometry-quiet-inner}, the
field bound \eqref{eq:geometry-field}, and (Loc) in
\eqref{eq:geometry-loc} on $[0,b]$. For $k<n$, the confinement estimate starts
at $T_{k+1}$ with the material inside ring $k$ in
$\overline B_{2r_{k+1}}$. The support of ring $k+1$ is in $\Nflat_{k+1}$ at
saturation, and the material inside it lies in $B_{ar_{k+1}/8}$.
Proposition~\ref{prop:confine} then gives \eqref{eq:geometry-confinement},
with radius $\varrho_k(v)=\big(2+2\Cstar A_{k+1}(v-T_{k+1})\big)r_{k+1}$.
On the pinning interval, $A_{k+1}(v-T_{k+1})\le\Tpin$. During the subsequent
active interval, it is at most $\Tpin+\mu\hat T_k$, where the lower amplitude
bound at activation gives
\begin{equation}
\label{eq:relay-duration}
\hat T_k=\frac{\log(\delta_0/\ksh_k)}{\sigma_k}
\le\frac{\log_+(\delta_0/\Klo)+5\log(1/s_k)}{\sigma'}.
\end{equation}
Thus (B5) keeps this radius below $ar_k/8$ even across the active interval.
Equation~\eqref{eq:geometry-moment} bounds $|q_{k+1}(v)|$ cubically in
$1+A_{k+1}(v-T_{k+1})$. For $k=n$ the inner region is empty.

\phantomsection\label{par:relay-quiet}
\emph{Control of the quiet rings.}
Suppose the support separations and field bound
\eqref{eq:geometry-quiet-support}--\eqref{eq:geometry-field} hold on
$[0,b]$, with $b>0$, and ring $j<n$ is quiet through $b$ with an assigned
activation time. Ring $i<j$ grows more slowly:
Lemma~\ref{lem:detuning} gives
\begin{equation}
\label{eq:relay-growth-gap}
A_j\sigma_j-A_i\sigma_i
=A_j(\sigma_j-\mu^{i-j}\sigma_i)
\ge\frac{A_j\sigb(\mu-1)}{2\mu}.
\end{equation}
In the Duhamel bound for ring $i$, this gap allows us to use the endpoint
value of $\maj_j$. By (M1) of Lemma~\ref{lem:majorant},
\begin{equation}
\label{eq:relay-duhamel-comparison}
A_i\int_0^b e^{A_i\sigma_i(b-v)}\maj_j(v)\dd v
\le \frac{A_i\maj_j(b)}{A_j\sigma_j-A_i\sigma_i}
\le\frac{2\mu}{\sigb(\mu-1)}\maj_j(b).
\end{equation}
Together with the radius gain in Lemma~\ref{lem:transfer-step}, this gives
the product comparison \eqref{eq:O-comparison} of Lemma~\ref{lem:O}.
Its support and field hypotheses are exactly \eqref{eq:O-support}, supplied
by \eqref{eq:geometry-quiet-support}--\eqref{eq:geometry-field}.
The remaining hypotheses \eqref{eq:O-smallness} follow from (BA2), (B7)
and (B9), the last giving $2\Abar\varepsilon_j\le\majst$.
The outer-field estimate \eqref{eq:O-feedback} gives
$\etout_j\le\theta_j\maj_j$, and Lemma~\ref{lem:O} applies Theorem~\ref{thm:quiet} to every
$i\le j$, using the local time $A_iv$. Consequently, for $v\in[0,b]$,
\begin{equation}
\label{eq:boot-quiet}
S_i(v)\subset\NNm_i,\qquad \|w_i(v)\|_{H^{s_h}}\le C_Q\maj_i(v).
\end{equation}
Appendix~\ref{sec:bootstrap-quiet-control} gives the precise rescaling and
checks for \eqref{eq:O-feedback} and \eqref{eq:boot-quiet}.
These estimates also hold through earlier assigned activation times, by
taking each such time as the endpoint $b$.

\phantomsection\label{par:relay-handoff}
\emph{Hand-off.}
Let $k<n$. Apply Lemma~\ref{lem:geometry} through $\min(t,\tsh_k)$,
using the lower bounds in \eqref{eq:relay-amplitudes}, and then
Lemma~\ref{lem:O} with $j=k$ as just verified. This gives
\eqref{eq:O-feedback} and \eqref{eq:boot-quiet} on that interval.
If $\tsh_k\le t$, \eqref{eq:relay-active-moment} for ring $k+1$ identifies
the leading outgoing moment
\begin{equation}
\label{eq:relay-outgoing}
\qs=2\pi\pout_*(G_{k+1})c_{k+1}(T_{k+1})e^{\lambda_{k+1}\hat T_{k+1}}.
\end{equation}
Lemma~\ref{lem:kicktrue} compares the true incoming field to this moment.
In particular $|\qs|\ge(15/16)2\pi|\pout_*(G_{k+1})|\delta_0$, and
Equation~\eqref{eq:boot-kick} bounds the resulting kick. To apply the lemma,
we use \eqref{eq:geometry-moment} during pinning and the conclusions of
item~(c) already assumed for ring $k+1$ through $T_{k+1}$.
Its earlier quiet-ring bounds are \eqref{eq:O-feedback} and
\eqref{eq:boot-quiet}, and its inner support satisfies
\eqref{eq:geometry-confinement} through saturation. For the seed there is no earlier quiet
interval. The remaining smallness conditions are checked in
Appendix~\ref{sec:bootstrap-handoff-checks}.

The factor $|\beta_k|=|\beta(G_k)|s_k^5/(2\pi)$ in the kick turns the outgoing
moment into an incoming coefficient of size $s_k^5$. More precisely,
Lemma~I.\ref{I-lem:uniform} and Equation~\eqref{eq:boot-kick} give
\begin{equation}
\label{eq:relay-incoming-lower}
|\cin_k|=e^{x_k\Tpin}|\Ain(\Tpin)|
\ge \frac{3c_N}{32\mu}\,
\frac{\varepsilon_k\delta_0 e^{x_k\Tpin}}{|\lamref|(\mu-1)}.
\end{equation}
Lemmas~\ref{lem:quiet-time} and~\ref{lem:theta-quiet} give the
bound on accumulated feedback~\eqref{eq:boot-theta} using only amplitudes at levels
$i>k$. With this bound, Proposition~\ref{prop:relay} bounds the entire
response at $\tsh_k$. The seed-size condition ($\kappa$7) bounds the remainder
$\Xi_{k+1}$ by $\delta_0/(\sigb(\mu-1))$. The resulting inequality is
\begin{equation}
\label{eq:relay-startup-duhamel}
\maj_k(\tsh_k)\le 2\Cvel\varepsilon_k e^{x_k\Tpin}
\left[\frac{(4C_2+1)\delta_0}{\sigb(\mu-1)}+C_D(\Tpin)\right]
\le\Abar\varepsilon_k.
\end{equation}
The choice of $K_1$ in Equation~\eqref{eq:K1} therefore gives
$4\maj_k(\tsh_k)\le K_1|\cin_k|$. With these bounds,
Lemma~\ref{lem:handoff} compares the actual coefficient $\ksh_k$
with $|\cin_k|$ within factors $1/2$ and $3/2$. It also gives the support,
initial norm and Duhamel inequalities in \eqref{eq:relay-startup}.
The upper and lower bounds on $\ksh_k$ in \eqref{eq:relay-startup} are
the two sides of \eqref{eq:boot-kick} with these factors included,
as checked in Appendix~\ref{sec:bootstrap-handoff-checks}.

By (B6) and ($\kappa$1), $0<\ksh_k\le\Kup\sbar^5\le\kseed<\delta_0$.
Together with $S_k(\tsh_k)\subset\NNm_k$, this verifies the next clause of
the clock recursion. The active size bound at that time follows directly:
\begin{equation}
\label{eq:boot-switch}
\|w_k(\tsh_k)\|_{H^s}\le K_1C_Q\ksh_k\le C_0\ksh_k
<C_2\hat a_k(\tsh_k)(1+\rhohat_k(\tsh_k)).
\end{equation}

\phantomsection\label{par:relay-active}
\emph{Active growth and the inner vorticity.}
Consider an active interval of positive length. For $k<n$ this means
$\tsh_k<t$ after the hand-off above. For the seed it means $t>0$. Put
$b=\min(t,T_k)$ and apply Lemma~\ref{lem:geometry} through $b$,
using the new lower bound in \eqref{eq:relay-startup} if $k<n$.
This gives \eqref{eq:geometry-confinement}--\eqref{eq:geometry-loc}
on their stated time intervals. The local active time
is $u=A_k(v-\tsh_k)$. The circulation field of $\omega_{\ge k+1}$ is
included in $L(G_k)$. By Lemma~\ref{lem:inner-field}, $\hin_k$ is
bounded by $\Cvel\mu s_k^2$ times the cube of the
confinement radius in ring-$k$ units. This gives the polynomial field
bound
\begin{equation}
\label{eq:relay-inner-calculation}
\|\hin_k(v)\|_{C^{M_h}(N')}
\le \Cvel\mu s_k^5\big(2+2\Cstar(\Tpin+\mu u)\big)^3
\le C_0\ksh_k\jap{u}^3
\end{equation}
by (P4) and the lower bound in \eqref{eq:relay-startup}, as in
Equation~\eqref{eq:boot-inner}.
For $k=n$, $\hin_n=0$.

For $k>1$, apply Lemma~\ref{lem:O} with $j=k-1$ through $b$, obtaining
\eqref{eq:O-feedback} and \eqref{eq:boot-quiet} for the outer quiet rings.
Lemma~\ref{lem:Oprime} then bounds their feedback, using (BA3) for
$\tsh_k<v\le b$ and \eqref{eq:boot-switch} at $\tsh_k$ (or
\eqref{eq:relay-seed} for $k=n$). The term $\|w_k(v)\|_{H^s}$ in the
bound for $\etin_{k-1}(v)$ carries the factor $\varepsilon_{k-1}$;
$\|w_{k-1}(v)\|_{H^s}$ in the bound for $\etout_k(v)$ carries
$b_{k-1,k}=\mu^{-1}s_{k-1}$
(Equations~\eqref{eq:transfer-step-inner}--\eqref{eq:transfer-step-outer}).
Their product is $s_{k-1}^6$, giving
\begin{equation}
\label{eq:relay-outer-field}
\|\hout_k(v)\|_{C^{M_h}(N')}
\le \Crelay\frac{s_{k-1}^6}{\sigb(\mu-1)}\hat a_k(v)
=\eta_{b,k}\hat a_k(v).
\end{equation}
For $k=1$, $\hout_1=0$.

We now bound $c_k(v)-c_k(\tsh_k)$ by splitting
$F_k=-\hin_k\cdot\nabla\Phi-\hout_k\cdot\nabla\Phi$ in
Equation~\eqref{eq:running-k}. For $k<n$, the integral containing
$\hin_k$ starts at $\tauK=\Tpin$ after changing from physical time $s$
to $\tauK=A_{k+1}(s-T_{k+1})$. By \eqref{eq:geometry-moment} and
$x_k\ge\sigma/4$, its modulus is at most $e^{x_k\Tpin}|\beta_k|\Cmom$ times
\begin{equation}
\label{eq:relay-moment-tail}
\int_{\Tpin}^{\infty} e^{-x_k\tau}(1+\tau)^3\dd\tau
\le I_3(\sigma/4)(1+\Tpin)^3e^{-\sigma\Tpin/4}.
\end{equation}
The lower bound in \eqref{eq:boot-kick} and (P2) give
$e^{x_k\Tpin}|\beta_k|\Cmom I_3(\sigma/4)\le |\cin_k|\le2\ksh_k$.
Thus the modulus of the integral containing $\hin_k$ is at most
$2\rT(\Tpin)\ksh_k\le\ksh_k/32$ by (P3).
For the integral containing $\hout_k$, Equation~\eqref{eq:relay-outer-field}
gives the bound $\Cls\eta_{b,k}\hat T_k\ksh_k\le\ksh_k/32$ by (B8).
These give the two terms on the right of \eqref{eq:drift-interval};
for the seed, $\hin_n=0$. The smallness conditions
\eqref{eq:relay-drift-small} give (Coef) and the saturation estimate
\eqref{eq:relay-saturation}. The checks for these applications and
for the rescaling are in Appendix~\ref{sec:bootstrap-active-checks}.

Theorem~\ref{thm:active} and Corollary~\ref{cor:active-moment} now give
\eqref{eq:relay-active}--\eqref{eq:relay-active-moment}, as well as
\eqref{eq:improved-support} and \eqref{eq:improved-active} for ring $k$
through $b$. If the interval consists of the activation time alone,
\eqref{eq:improved-support} and \eqref{eq:improved-active}
come from the hand-off or the seed. No positive-length active estimate is
needed in that case.

\phantomsection\label{par:relay-quiet-size}
\emph{Quiet sizes.}
It remains to improve (BA2) on the intervals just reached. At a positive
time $v$ on the pinning interval it suffices to bound $\maj_j(v)$ for $j=k$.
On the active interval after activation take $j=k-1$ if $k>1$. This ring
has an assigned activation time and is quiet through $v$.
Equations~\eqref{eq:geometry-quiet-support}--\eqref{eq:geometry-field},
\eqref{eq:O-feedback} and \eqref{eq:boot-quiet} hold through $v$.
The conclusions in item~(c) for ring $j+1$ have been proved through
$\min(v,T_{j+1})$; its initial Duhamel bound is
\eqref{eq:relay-startup}, or $\maj_n(0)=0$ for the seed.
The bound on accumulated feedback \eqref{eq:boot-theta}
again follows from Lemmas~\ref{lem:quiet-time} and~\ref{lem:theta-quiet}.

By (M1) of Lemma~\ref{lem:majorant} and \eqref{eq:O-feedback},
\begin{equation}
\label{eq:relay-quiet-absorption}
\maj_j(v)\le\maj_j^{\mathrm{in}}(v)+\theta_j A_jv\maj_j(v)
\le\maj_j^{\mathrm{in}}(v)+\tfrac12\maj_j(v),
\end{equation}
where the last inequality is \eqref{eq:boot-theta}.
Proposition~\ref{prop:relay} and ($\kappa$7) then give
\eqref{eq:improved-quiet} for $j$, as at the hand-off in
\eqref{eq:relay-startup-duhamel}. Appendix~\ref{sec:bootstrap-quiet-checks}
checks the endpoint cases. For every $i<j$, \eqref{eq:O-comparison} gives
\begin{equation}
\label{eq:relay-quiet-outer}
\maj_i(v)\le\tfrac34\Abar\varepsilon_i
\prod_{p=i+1}^{j}(K_O\varepsilon_p)
\le\tfrac38\Abar\varepsilon_i,
\end{equation}
since each factor is at most $1/2$ by (B7). This includes the outer rings
whose activation times have yet to be assigned. When $k=1$ is active
there are no quiet rings after activation.

Equation~\eqref{eq:boot-quiet} for quiet rings and Theorem~\ref{thm:active}(a)
for active rings prove \eqref{eq:improved-support} on the new interval;
Theorem~\ref{thm:active}(b) and
\eqref{eq:relay-quiet-absorption}--\eqref{eq:relay-quiet-outer} prove
\eqref{eq:improved-quiet}--\eqref{eq:improved-active} there.
These three estimates at the initial time of the seed were assumed.
\end{proof}

\section{Choice of Parameters}
\label{sec:param}

Proposition~\ref{prop:relay-step} uses one choice of parameters for every
ring and every depth of the construction. We first set the saturation threshold and the
pinning time, then fix the constants in the activation bounds
\eqref{eq:relay-startup} and the active-ring estimates
\eqref{eq:relay-active}--\eqref{eq:relay-inner-field}.
The seed size and radius ratios are made small after those constants are
fixed. The precise requirements (P0)--(P6) are in
Appendix~\ref{sec:param-lists}.

\begin{lemma}
\label{lem:nonempty}
Assume \textup{(N1)}--\textup{(N3)}. The choices (P0)--(P6) can be made in this
order, with all quantities finite and positive (except $\eta_{b,1}=\theta_1=0$),
and each quantity is determined by quantities fixed at earlier steps or earlier in the same step, as in
Table~\ref{tab:param}. None of the choices depends on $n$, $A_1$ or
$\omega^{(n)}$. Finally, if (P0)--(P6) hold, then for every integer $n\ge2$ and
every $A_1>0$ the standing assumptions (SA) hold with these parameters, and
(Det) holds.
\end{lemma}

\begin{proof}
First fix the ring, its spectral data and the regularity orders. The
interval for $\sigma'$ is nonempty because $\sigb>211/2500>9/200$ and
$3\sigma/4\ge3\sigb/4>9/200$. Choose $\delta_0>0$ small so the quadratic
error is a small fraction of the outgoing moment. Choose $\mu-1>0$
small compared with $\delta_0$ to control
$|c_k(v)-c_k(\tsh_k)|$ relative to $\ksh_k$
(Equation~\eqref{eq:drift-interval}). Choose $\Tpin$ large enough
to make $\rT(\Tpin)$ small.

We now fix the constant $K_1$ in the Duhamel bound
$\maj_k(\tsh_k)\le K_1\ksh_k$ of Equation~\eqref{eq:relay-startup}.
Its formula in Definition~\ref{def:K1} replaces $\Xi_{j+1}$ in
Equation~\eqref{eq:relay} by $\delta_0/(\sigb(\mu-1))$;
condition $(\kappa7)$ will ensure this bound. Thus $K_1$ uses only the data already fixed.
Choose $C_0$ to satisfy Equation~\eqref{eq:param-C0}, which bounds the initial
perturbation and inner field required by Theorem~\ref{thm:active}.
The constants controlling the remainder, including
$\Crelay(C_0,K_1)$, can then be fixed.

The seed conditions either bound $\kappa$ by a fixed positive threshold or
bound expressions of the form
$C\kappa^\alpha\langle\log(\delta_0/\kappa)\rangle^N$ by fixed positive
constants, with $\alpha>0$ and $N\ge0$. The exponents are positive because
$\mu>1$ and $\sigma'<3\sigma/4$.
These expressions tend to zero as $\kappa\to0$. The function $\rhohat$ also
tends to zero at saturation, uniformly in $G\in[0,\Gone]$: its initial layer
has uniformly bounded length, while the saturation time grows like
$\log(1/\kappa)$ and its decay rate is at least $\sigb-\sigma'>0$.
Thus every condition holds throughout $(0,2\kseed]$ once $\kseed$ is small
enough. In particular, $\Xi_{j+1}\le\delta_0/(\sigb(\mu-1))$ at every
hand-off, as required in the choice of $K_1$.

After the seed is fixed, the radius conditions have the same form,
$C\sbar^\alpha\log^N(1/\sbar)$, except for spectral continuity and the
conditions imposed at all levels. Spectral continuity gives (Det) on the
shrinking interval $0\le G\le\mu\sbar^2$. For the conditions at all levels, the
confinement radius grows at most like $\log(1/s_k)$ and the outer feedback
has a factor $s_{k-1}^6$. Positive powers of $s_k$ dominate these logarithms,
uniformly for $s_k\le\sbar$. Hence one sufficiently small $\sbar$ works at
every level. It also puts the new amplitudes, of size $s_k^5$, below the
seed threshold. The exact inequalities and these uniform estimates are checked
in Appendix~\ref{sec:param-verification}; they also give (SA) and (Det).
No choice uses $n$, $A_1$ or the solution $\omega^{(n)}$.
\end{proof}

\begin{table}[H]
\centering
\small
\begin{tabular}{@{}p{0.20\textwidth}p{0.51\textwidth}p{0.21\textwidth}@{}}
\hline
Choice & Purpose & Fixed after \\
\hline
(P0): ring and regularity & Supply the spectral, velocity and single-ring estimates & fixed data \\
(P1): $\delta_0$ & Keep nonlinear errors below the outgoing moment & (P0) \\
(P2): $\mu$ & Control the change in $c_k$ after activation & (P1) \\
(P3): $\Tpin$, $K_1$ & Control the tail of the kick and fix the hand-off constant & (P2) \\
(P4): $C_0$, error constants & Bound the initial perturbation, inner field and remainder & (P3) \\
(P5): $\kseed$ & Control pre-activation and transient errors uniformly & (P4) \\
(P6): $\sbar$ & Ensure separation, confinement and small outer feedback & (P5) \\
\hline
\end{tabular}
\caption{Order and purpose of the parameter choices. The depth $n$, the outer
amplitude $A_1$ and the solution $\omega^{(n)}$ enter none of the choices.
The exact inequalities are in Appendix~\ref{sec:param-lists}.}
\label{tab:param}
\end{table}

\section{Completion of the Bootstrap}
\label{sec:bootstrap}

It remains to iterate the bootstrap through all scales, from the seeded ring
$n$ to ring $1$, using Proposition~\ref{prop:relay-step} at each step.
The time bound and the nonzero moment of ring $1$ will then give
Proposition~\ref{prop:reduction} by the compactness criterion of Part~I.
We use the setting and bootstrap assumptions of Section~\ref{sec:relay}.

\begin{proposition}
\label{prop:bootstrap}
Assume \textup{(N1)}--\textup{(N3)} and \textup{(P0)}--\textup{(P6)}. Then for
every $A_1>0$ and every $n\ge2$:
\begin{enumerate}[label=\textup{(\arabic*)}]
\item $\Tend=T_1$; that is, the clock is defined up to level $1$. Moreover
  $\ksh_1\le\Kup s_1^5$;
\item the conclusions of Theorem~\ref{thm:active} hold for ring $1$ on
  $\Phase_1$, in ring-$1$ units and active time, with $G=G_1$,
  $\kappa=\ksh_1$, the constant $C_0$ and $\eta_b=\eta_{b,1}$; and
  $|c_1(T_1)|e^{\sigma_1\hat T_1}\ge\frac{15}{16}\delta_0$;
\item $D_2(T_1)\subset B_{a/8}$;
\item $T_1\le\Tstar$, with $\Tstar$ of Lemma~\ref{lem:time}.
\end{enumerate}
\end{proposition}

\begin{proof}
\emph{The seed.}
At $t=0$, $\tsh_n=0$, $G_n=0$, and $w_n(0)=\Phi^{\kseed}-\Phi$ in ring-$n$
units. By ($\kappa$6), $\kseed\le\kL$. Lemma~\ref{lem:seed} and (P4) give
\begin{equation}
\label{eq:boot-seed-data}
S_n(0)\subset\NNm_n,\qquad
\tfrac12\kseed\le\ksh_n\le2\kseed,\qquad
\|w_n(0)\|_{H^{M-1}}\le\Cseed\kseed\le C_0\ksh_n.
\end{equation}
By ($\kappa$1), $2\kseed<\delta_0$, so the clock is defined up to level $n$.
By its definition in Theorem~\ref{thm:active}, $t_1(C_0,G)>0$ for every $G$
since $2^{3/2}C_0>1$. Thus $0<t_1(C_0,0)$ and the first formula for
$\rhohat$ applies at $0$, giving $\rhohat_n(0)=C_0$.
With $s\le M-1$ and $C_2\ge1$ the
active size at the initial time satisfies
\begin{equation}
\label{eq:boot-start}
\|w_n(0)\|_{H^s}\le C_0\ksh_n
<C_2\hat a_n(0)(1+\rhohat_n(0)).
\end{equation}
For $j<n$, the initial support $S_j(0)$ is contained in $\NNm_j$. Also
$D_{j+1}(0)\subset\overline B_{5r_{j+1}/4}\subset B_{ar_j/8}$ and the outer
rings lie outside $B_{6r_j}$. These separations define $h_j(0)$ and give
$\maj_j(0)=0$. The same holds for $h_n(0)$ and $\maj_n(0)$, with an empty
inner region. Thus Equations~\eqref{eq:improved-support}--\eqref{eq:improved-active}
hold at $0$.

\emph{Improvement through a fixed time.}
Fix $t\in[0,\Tend]$ with (BA)$(t)$. Apply Proposition~\ref{prop:relay-step}
first to the seeded interval. If $t>T_n$, take $k=n-1$ in the relay step.
After a completed interval ending at $T_{k+1}<t$, we have
Equations~\eqref{eq:improved-support}--\eqref{eq:improved-active}, the amplitudes
in Equation~\eqref{eq:relay-amplitudes}, and the conclusions of
Proposition~\ref{prop:relay-step}(c) for ring $k+1$. These are the hypotheses
for the next level $k$. This is the locomotive induction
shown in Figure~\ref{fig:relay-step}. If the current interval reaches
$t$, we stop there, including the hand-off when $t=\tsh_k$. Otherwise the
clock conditions have been verified
through $k$ and we proceed to $k-1$. There are at most $n$ levels, and if
we reach $T_1$ then $t\le\Tend=T_1$. This proves
Equations~\eqref{eq:improved-support}--\eqref{eq:improved-active} through $t$,
Equation~\eqref{eq:relay-startup} at each nonseed activation reached by $t$,
and Proposition~\ref{prop:relay-step}(c) on each positive-length active
interval reached by $t$. The seed satisfies Equation~\eqref{eq:boot-seed-data}.

\emph{Continuation in time.}
Let $t_*$ be the supremum of the times in $[0,\Tend]$ at which (BA) holds.
The seed gives (BA)$(0)$, and (BA) holds at every time below $t_*$ by its
interval definition. The flow is continuous on $[0,\Tend]\times\R^2$ by
Theorem~\ref{thm:smooth}(2),(4). If a compact support lies in an open set at
one time, its image remains in that set for nearby times. When
$S_k(v)\subset N_k$, $w_k(v)$ is continuous in $H^s$ because
$\omega(v)\mathbf 1_{S_k(v)}=f^{(n)}_k\circ\flow_v^{-1}$
(Equation~\eqref{eq:omega-ge}) is smooth with support in a fixed compact
set on a bounded time interval. When $h_j$ is defined, it and its Duhamel
bound are continuous by Lemma~\ref{lem:reduction}(1).

For a ring that has activated, write
\begin{equation}
\label{eq:boot-continuation-envelope}
R_k(v)=2C_2\hat a_k(v)(1+\rhohat_k(v)),\qquad v\ge\tsh_k,
\end{equation}
for the right side of Equation~\eqref{eq:BA-active}.
It is continuous except at $v=\tsh_k+t_1(C_0,G_k)/A_k$, where the transient
jumps from a left limit of at most $2^{-3/2}$ to $1$, by the definition of
$t_1$ in Theorem~\ref{thm:active}. Thus $R_k$ is right-continuous, and at
every point with a left limit its value is at least that limit.

We first obtain (BA) at $t_*$. If $t_*>0$, a ring that is not debris at
$t_*$ was not debris earlier. Its support is in the closed annulus
$\Nflat_j$ at $t_*$ by continuity, so Equation~\eqref{eq:BA-support} holds.
Every nonseed activation strictly before $t_*$ satisfies
$\ksh_k\ge\Klo s_k^5$ by Equation~\eqref{eq:relay-startup} at earlier times.
Lemma~\ref{lem:geometry}(b) therefore gives the separation of the outer supports
in Equation~\eqref{eq:geometry-quiet-support} and the inner confinement in
Equation~\eqref{eq:geometry-quiet-inner} through $t_*$. These define $h_j$ for all
rings quiet at $t_*$, and continuity gives Equation~\eqref{eq:BA-quiet} there.

For a ring active at $t_*$ with $\tsh_k<t_*$, the left limit bound for
$R_k$ and continuity of $w_k$ give Equation~\eqref{eq:BA-active} at $t_*$.
This proves (BA)$(t_*)$.

Suppose $t_*<\Tend$. The relay step gives
Equations~\eqref{eq:improved-support}--\eqref{eq:improved-active} at $t_*$,
including Equation~\eqref{eq:improved-active} at any activation there.
Choose $\Delta>0$ with
$t_*+\Delta\le\Tend$ and no assigned activation or saturation time in
$(t_*,t_*+\Delta]$. The rings that are not debris on this interval were not
debris at $t_*$. Their supports are compact subsets of
$\NN_j\subset\operatorname{int}\Nflat_j$, so Equation~\eqref{eq:BA-support}
persists for small $\Delta$. Equation~\eqref{eq:relay-startup} gives
$\ksh_k\ge\Klo s_k^5$ at every nonseed activation before $t_*+\Delta$,
including any at $t_*$. Lemma~\ref{lem:geometry}(b), through
Equations~\eqref{eq:geometry-quiet-support}--\eqref{eq:geometry-quiet-inner},
defines the fields for the quiet rings through $t_*+\Delta$.
Equation~\eqref{eq:improved-quiet} gives
$\maj_j(t_*)\le\Abar\varepsilon_j$, so Equation~\eqref{eq:BA-quiet}
persists by continuity. Any ring active on $(t_*,t_*+\Delta]$ has activated
by $t_*$ and satisfies, by Equation~\eqref{eq:improved-active},
\begin{equation}
\label{eq:boot-continuation-margin}
\|w_k(t_*)\|_{H^s}\le C_2\hat a_k(t_*)(1+\rhohat_k(t_*))<R_k(t_*).
\end{equation}
Continuity of $w_k$ and right-continuity of $R_k$ give
Equation~\eqref{eq:BA-active} for small $\Delta$.
Since there are finitely many rings, one $\Delta>0$ works for all
of them, contradicting the definition of $t_*$. Thus $t_*=\Tend$ and
Equations~\eqref{eq:improved-support}--\eqref{eq:improved-active} hold through
$\Tend$.

\emph{The final ring and the elapsed time.}
If the clock stopped at $k_{\mathrm s}<n$, then
$\Tend=\tsh_{k_{\mathrm s}}$ would be reached under (BA). The hand-off in
Proposition~\ref{prop:relay-step}, through Equation~\eqref{eq:relay-startup},
verifies the support and amplitude conditions at that same time. This
contradicts the failure at
$k_{\mathrm s}$ in the fixed clock. Hence $\Tend=T_1$.
For ring $1$, Equation~\eqref{eq:relay-startup} gives $\ksh_1\le\Kup s_1^5$,
Proposition~\ref{prop:relay-step}(c) gives Theorem~\ref{thm:active} through
$T_1$ and the coefficient bound \eqref{eq:relay-saturation}, and
Equation~\eqref{eq:geometry-confinement} gives $D_2(T_1)\subset B_{a/8}$,
since $r_1=1$.

For each $k<n$, Equation~\eqref{eq:relay-startup} gives
\begin{equation}
\label{eq:boot-duration}
T_k-T_{k+1}
=\frac{\Tpin}{A_{k+1}}+\frac{\hat T_k}{A_k}
\le\frac{1}{A_k}\left(\Tpin+
\frac{\log_+(\delta_0/\Klo)+5\log(1/s_k)}{\sigma'}\right).
\end{equation}
The numerator in Equation~\eqref{eq:boot-duration} grows at most linearly
in $k$, while $A_k=A_1\mu^{k-1}$ grows geometrically. The seed satisfies
\begin{equation}
\label{eq:boot-seed-duration}
T_n\le\frac{\log(2\delta_0/\kseed)}{A_n\sigma'}.
\end{equation}
Summing Equations~\eqref{eq:boot-duration}--\eqref{eq:boot-seed-duration}
as in Lemma~\ref{lem:time} gives $T_1\le\Tstar$ uniformly in $n$.
\end{proof}

\begin{proof}[Proof of Proposition~\ref{prop:reduction}]
Assume \textup{(N1)}--\textup{(N3)}. By Lemma~\ref{lem:nonempty}, we may choose
parameters satisfying (P0)--(P6). This fixes $\mu\in(1,\mu_{\max}]$ with
$\mu_{\max}\le\frac{17}{16}$ (by (P2)) and $\sbar\in(0,a/16]$ (by (P6) and
(B1)), and none of (P0)--(P6) involves $A_1$. Let $A_1>0$, let $f$ be the
function of Definition~\ref{def:datum} with these parameters, and let $\Tstar$ be the number of
Lemma~\ref{lem:time}. By Proposition~\ref{prop:integrability}, $f$ is radial
and nonnegative, $f\in (L^1\cap L^q)(\R^2)$ for every $q\in[1,\infty)$, and
$f\notin L^\infty(\R^2)$. The bracket in the definition of $\Tstar$ is a finite
sum of the term $1$ and of nonnegative terms, because $\Tz>0$, $\mu>1$,
$0<\sbar\le a/16<1$ and $\kseed<\delta_0$ by ($\kappa$1) at $\kappa=2\kseed$;
so $\Tstar\in[1/A_1,\infty)$, and in particular $\Tstar\in(0,\infty)$.

Fix $\chiout\in C_c^\infty((a/4,3))$ with $0\le\chiout\le1$ and $\chiout=1$ on
$[a/2,2]$, and put $\Psione(x)=\chiout(|x|)\bar x^3$; it does not depend on
$n$. The cutoff $\chiout$ is chosen so that, at the time $T_1$ below, it equals
$1$ on the support of ring $1$ and vanishes on the material inside ring $1$; so
the moment $\int\omega^{(n)}(T_1)\Psione$ sees only ring $1$.

For $n\ge2$ let $f^{(n)}$ and $\omega^{(n)}$ be as in Definition~\ref{def:approx},
with the seed size $\kseed$ of (P5), which lies in $(0,\kL]$ by ($\kappa$6) at
$\kappa=2\kseed$. By Lemma~\ref{lem:adm}, $(f^{(n)},\omega^{(n)})_{n\ge2}$ is a
$2$-admissible sequence for $f$ on $[0,\infty)$ (indexed from $n=2$), and
$\omega^{(n)}(t)$ is $3$-fold rotationally symmetric for every $n$ and $t$.
For $n\ge2$ put $T^{(n)}=T_1$, which is defined by
Proposition~\ref{prop:bootstrap}(1) and lies in $[0,\Tstar]$ by
Proposition~\ref{prop:bootstrap}(4). (SA) and (Det) hold by
Lemma~\ref{lem:nonempty}, and the clock is defined up to level $1$. By
Proposition~\ref{prop:bootstrap}(2), the conclusions of
Theorem~\ref{thm:active}(a)--(b) hold for ring $1$ on $\Phase_1$ with $G=G_1$,
$\kappa=\ksh_1$, the constant $C_0$ and $\eta_b=\eta_{b,1}=0$
(Definition~\ref{def:P6}). Ring $1$ has
$r_1=1$, so ring-$1$ units are the original lengths; at $t=T_1$ its active time
is $\hat T_1$, so $w(\hat T_1)=w_1(T_1)$ in Theorem~\ref{thm:active}. The
vorticity $\omega^{(n)}(T_1)$ is supported in $D_1(T_1)$, the disjoint union of
$S_1(T_1)$ and $D_2(T_1)$. By Proposition~\ref{prop:bootstrap}(3),
$D_2(T_1)\subset B_{a/8}$, so the part on $D_2(T_1)$ does not
contribute, as $\chiout$ vanishes on $[0,a/4]$. On $S_1(T_1)$,
$\omega^{(n)}(T_1)=A_1(\Phi+w_1(T_1))$ (Definition~\ref{def:regions}), which is
supported in $\NN\subset\{a/2\le|x|\le2\}$ by Theorem~\ref{thm:active}(a), where
$\chiout(|x|)=1$. As $\Phi$ is radial, $\pout_3[\Phi]=0$, with the
multipole moment $\pout_3$ of Definition~I.\ref{I-def:multipole}, so
$\int\omega^{(n)}(T_1)\Psione=A_1\int(\Phi+w_1(T_1))\,\bar y^3\dd y=2\pi A_1\,\pout_3[w_1(T_1)]$.

Write $c=c_1(T_1)$, $\lambda=\lambda_1$ and $\phist=\phist(G_1)$, the
unstable eigenfunction of Definition~I.\ref{I-def:eigendata}. We bound
$|\pout_3[w_1(T_1)]|$ from below by comparing $w_1(T_1)$ with
$2\Re(ce^{\lambda\hat T_1}\phist)$, as in Theorem~\ref{thm:active}\textup{(b)}. As
$\pout_3[\phist]=\pout_*(G_1)$ and $\pout_3[\overline{\phist}]=0$ (its integrand
carries the factor $e^{-6i\theta}$),
$\pout_3[2\Re(ce^{\lambda\hat T_1}\phist)]=ce^{\lambda\hat T_1}\pout_*(G_1)$. The
function $w_1(T_1)-2\Re(ce^{\lambda\hat T_1}\phist)$ is supported in $N$, because
$\Phi$, $\Phi+w_1(T_1)$ and $\phist\in X^0$ are. (P1) gives
$\Cp C_2\delta_0\le\cout/(16\pi)$, and ($\kappa$5) with $G=G_1\in[0,\Gone]$
(Lemma~\ref{lem:inner-circ} and (SA3)) and
$\kappa=\ksh_1$, which lies in $(0,\kseed]$ by
Proposition~\ref{prop:bootstrap}(1) and (B6) (it is positive since the clock is
defined up to level $1$, Definition~\ref{def:clock}), gives
$\Cp C_2(\rhohat_1(T_1)+\ksh_1{}^{1/2})\le\cout/(16\pi)$. As
$\hat a_1(T_1)=\delta_0$, Lemma~I.\ref{I-lem:uniform},
Theorem~\ref{thm:active}(b) with $\eta_b=0$ and these two inequalities give
\begin{equation}
\label{eq:boot-observable-error}
\begin{aligned}
\big|\pout_3[w_1(T_1)]-ce^{\lambda\hat T_1}\pout_*(G_1)\big|
&\le\Cp\big\|w_1(T_1)-2\Re(ce^{\lambda\hat T_1}\phist)\big\|_{H^s}\\
&\le\Cp C_2\delta_0\big[\rhohat_1(T_1)+\delta_0+\ksh_1{}^{1/2}\big]
\le\frac{\cout\delta_0}{8\pi}.
\end{aligned}
\end{equation}
As $G_1\in[0,\Gone]$, $|\pout_*(G_1)|\ge\cout/(2\pi)$ (Lemma~I.\ref{I-lem:uniform}), and
$|ce^{\lambda\hat T_1}|=|c_1(T_1)|e^{\sigma_1\hat T_1}\ge\frac{15}{16}\delta_0$
by Proposition~\ref{prop:bootstrap}(2). Hence Equation~\eqref{eq:boot-observable-error} gives
$|\pout_3[w_1(T_1)]|\ge\frac{\cout}{2\pi}\cdot\frac{15}{16}\delta_0-\frac{\cout\delta_0}{8\pi}
=\frac{11}{16}\cdot\frac{\cout\delta_0}{2\pi}\ge\frac{\cout\delta_0}{4\pi}$, and
multiplying by $2\pi A_1$ gives
\begin{equation}
\label{eq:boot-observable-defect}
\left|\int\omega^{(n)}(T^{(n)})\Psione\right|
\ge c_0:=\frac12A_1\cout\delta_0>0.
\end{equation}
The number $c_0$
involves only $A_1$, $\cout$ and $\delta_0$, and $\Tstar$ only $A_1$, $\mu$,
$\sbar$, $\sigma'$, $\Tz$, $\delta_0$ and $\kseed$; none of them depends on $n$.
Moreover $\Psione$ has the form required in
Corollary~I.\ref{I-cor:defect}, because $\chiout\in C_c^\infty((a/4,3))$ is real.
Corollary~I.\ref{I-cor:defect}, using Equation~\eqref{eq:boot-observable-defect}
and $\bar T=\Tstar$, gives a renormalized weak
solution $\omega_\infty$ on $[0,\infty)$ with datum $f$, in the sense of
Definition~I.\ref{I-def:weak} with $q=2$, such that
$\omega_\infty\in C([0,\infty);L^p(\R^2))$ for every $p\in[1,\infty)$ and, for
every $t\ge0$, $\omega_\infty(t)$ is equimeasurable with $f$ and $3$-fold
rotationally symmetric, and a time $t_0\in[0,\Tstar)$ such that $\omega_\infty$ and $f$ are
distinct elements of $X_q(T')$ for every $T'>t_0$ and every $q\in(4/3,\infty)$.

Let $T>t_0$ and $q\in(4/3,\infty)$; then $T>0$. By
Lemma~I.\ref{I-lem:radial-stationary}, the stationary function $t\mapsto f$ is a
renormalized weak solution in $X_q(T)$ with datum $f$. Also $\omega_\infty\in
X_q(T)$, since, by equimeasurability, $\|\omega_\infty(t)\|_{L^r}=\|f\|_{L^r}$ for every
$r\in[1,\infty)$ and $\omega_\infty\in C([0,T];L^1\cap L^q)$. By
Definition~I.\ref{I-def:weak}, the restriction of $\omega_\infty$ to $[0,T)$
satisfies Equations~(I.\ref{I-eq:weak}) and (I.\ref{I-eq:renorm}) for every test function
$\varphi\in\DT$ and every renormalization $\Ren$; these identities involve only
$\omega_\infty$, its velocity $\KBS\ast\omega_\infty$ and the datum $f$, and not
the exponent of the space. So $\omega_\infty$ is a renormalized weak solution in
$X_q(T)$ with datum $f$, and, as $T>t_0$, the two solutions are distinct elements
of $X_q(T)$. This holds in particular for every $T>\Tstar$, since $t_0<\Tstar$.
So $f$, $\Tstar$, the stationary solution, $\omega_\infty$ and $t_0$ have every
property stated in Theorem~\ref{thm:main}, for every $A_1>0$ with the same $\mu$
and $\sbar$, which were fixed before $A_1$.
\end{proof}

\appendix
\section{Verification of the Numerical Hypotheses}
\label{app:verify}

This appendix proves Theorem~\ref{thm:N-holds}, that is, the statements
\textup{(N1)}--\textup{(N3)} of Hypothesis~I.\ref{I-hyp:N}. The inequality
\textup{(N3)} is between explicit numbers, and we prove it by hand
(Proposition~\ref{prop:N3}). The statements \textup{(N1)} and \textup{(N2)} concern
the zeros of $\Evans_3(\cdot;G)$ for every $G$ in an interval, and we prove them with
the argument principle. For this we compute winding numbers of $\Evans_3(\cdot;G)$
along the boundaries of finitely many rectangles, from enclosures of $\Evans_3$ in
ball arithmetic.

Throughout, $a_0,d,\xiL,\xia,\xib,\xiR$, $H$, $\chi$, $\prof$, $\Omega_G$,
$\psi_\ell$, $\Evans_\ell$, $\Thl{\ell}$ and $\epsst$ are as in
Definition~I.\ref{I-def:ring} and Hypothesis~I.\ref{I-hyp:N}, and $G\ge0$. Only the
properties of $\chi$ listed in Definition~I.\ref{I-def:ring} are used. Recall that, for
$\Re\lambda>0$, $\psi_\ell(\cdot;\lambda,G)$ is the solution of
$\psi''=(\ell^2+i\ell\prof'/(\lambda+i\ell\Omega_G))\psi$ on $[\xiL,\xiR]$ with
$\psi(\xiL)=1$ and $\psi'(\xiL)=\ell$, and that
$\Evans_\ell(\lambda;G)=e^{-\ell(\xiR-\xiL)}(\psi'+\ell\psi)(\xiR)$.
We identify $\C$ with $\R^2$, write $|Y|_w=\max(\ell|u|,|v|)$ for
$\ell\ge1$ and $Y=(u,v)\in\C^2$, and put
\[
x_a=\frac{\xia-a_0}{d},\quad y_b=-\frac{\xib}{d},\quad \Gmax=\frac1{50},\quad \Gmax'=d,\quad
\Rect=\Big[\frac9{200},3\Big]\times[-5,3].
\]
The weighted norm $|Y|_w$ is used for $Y=(\psi,\psi')$: in it, a solution of
$\psi''=\ell^2\psi$ grows at most by the factor $e^{\ell t}$ over an interval of
length $t$ (proof of Lemma~\ref{lem:end-regions}). The numbers $x_a$ and $y_b$ are
the arguments of the two factors $H$ of $\prof$ at $\xia$ and at $\xib$, so that
$\epsst=\max(H(x_a),H(y_b))$.
Definition~I.\ref{I-def:ring} gives exactly $\xia-a_0=-\frac12+2^{-54}$ and
$d=\frac1{50}(1+6\cdot2^{-58})$. So $d>\frac1{50}$, whence the range
$[0,\Gmax']$ of the computations contains the range $[0,\Gmax]$ of $G$ in
\textup{(N1)} and \textup{(N2)}, and $50(1-6\cdot2^{-58})<1/d<50$, so
$x_a=-(\frac12-2^{-54})/d$ and $y_b=-1/(2d)$ satisfy
\begin{equation}
\label{eq:xayb}
-25<x_a<-25+4\cdot10^{-15},\qquad -25<y_b<-25+10^{-15}.
\end{equation}
For the argument principle we use winding numbers in the following form.
For a continuous closed curve $c\colon[0,1]\to\C\setminus\{0\}$,
$\wind(c,0)=(\theta(1)-\theta(0))/2\pi$, where $\theta$ is any continuous
function with $c=|c|e^{i\theta}$. If $\Gamma$ is the boundary of a closed
rectangle in $\{\Re\lambda>0\}$ and $\Evans_3(\cdot;G)\ne0$ on $\Gamma$,
$\wind(\Evans_3(\Gamma;G),0)$ is the winding number of
$t\mapsto\Evans_3(\gamma(t);G)$, where $\gamma$ traverses $\Gamma$ once
counterclockwise.

\subsection{Ball Arithmetic and the Trusted Base}
\label{app:arith}

The computations use \textbf{ball arithmetic}: a real number is enclosed in a
ball $[m\pm r]=[m-r,m+r]$ whose midpoint $m$ and radius $r\ge0$ are binary
floating-point numbers, and a complex number in a box $X+iY$, the product of two
real balls. We use the implementation of Johansson \cite{Johansson2017}, part of
FLINT \cite{FLINT} from version~3 on, in FLINT~3.6.0 as bundled with
version~0.9.0 of its Python interface python-flint \cite{pythonflint}, at a
working precision of $128$ bits. For a real ball $X$, $\mathrm{lo}(X)$ and
$\mathrm{up}(X)$ are outward dyadic lower and upper bounds; for a box $Z$,
$|Z|^{\rm up}$ is a dyadic upper bound of $|z|$ over $z\in Z$; for a ball $\rho$,
$\mathrm{box}(\rho)$ is a computed box containing
$[-\mathrm{up}(\rho),\mathrm{up}(\rho)]+i[-\mathrm{up}(\rho),\mathrm{up}(\rho)]$,
and hence the closed disk of
radius $\rho'$ about $0$ for every $\rho'\in\rho$ with $\rho'\ge0$;
$\mathrm{hull}(X,X')$ is a ball computed from $\mathrm{lo}(X)$ and
$\mathrm{up}(X')$; and $X\cup X'$ is the computed union of two real balls, again
a ball. Every enclosure rests on the following \textbf{containment property}.
\begin{enumerate}[label=\textup{(CP)}, leftmargin=*]
\item For every arithmetic operation $f$ that we use (sums, differences, products
  and quotients of real balls and of boxes, integer powers of a real ball, the
  maximum of two real balls, and the exponential and the square root of a real
  ball), and for all input balls $X_1,\dots,X_k$, the output ball contains
  $f(x_1,\dots,x_k)$ for every choice of $x_j\in X_j$ at which $f$ is defined; a
  ball with infinite radius, or with an undefined midpoint, is regarded as the
  whole real line. A ball built from an integer contains that integer, a ball
  built from an exact midpoint $m$ and an exact radius $r$ contains $[m-r,m+r]$, and a
  box built from two real balls is their product. Moreover
  $\mathrm{lo}(X)\le x\le\mathrm{up}(X)$ for every $x\in X$, $|z|\le|Z|^{\rm up}$ for
  every $z\in Z$, $\mathrm{hull}(X,X')$ contains every real number between
  $\mathrm{lo}(X)$ and $\mathrm{up}(X')$, $X\cup X'$ contains $X$, $X'$ and every
  real number between a point of $X$ and a point of $X'$, the computed
  intersection of a ball $X$ with $[0,\infty)$ contains $X\cap[0,\infty)$, and a
  comparison such as $X>0$ or $X\le c$ returns true only if every point of $X$
  satisfies it.
\end{enumerate}
The data of Definition~I.\ref{I-def:ring}, the contour points, and the endpoints of
the intervals of $G$ and the rectangles of Section~\ref{app:runs} are rational.
They are held exactly, as pairs of integers, in the arbitrary-precision integer
arithmetic of CPython, and every comparison between them is made there. A
rational $n/k$ in lowest terms, $k\ge1$, enters ball arithmetic as follows. If
$k=2^e$ and $|n|<2^{125}$, it enters as the product of the ball built from $n$ and
the power $-e$ of the ball built from $2$; the program checks that this product
has radius $0$, so by (CP) it is the point $n/k$. Otherwise it enters as the
quotient of the balls built from $n$ and $k$, which contains $n/k$ by (CP). No
hardware floating-point number of the programming language enters any enclosure:
every ball is built from integers or from other balls.

The proof of Proposition~\ref{prop:runs} is a computation. It is correct provided
the following \textbf{trusted base} holds.
\begin{enumerate}[label=\textup{(TB\arabic*)}, leftmargin=*]
\item The ball arithmetic of FLINT~3.6.0, as bundled with python-flint~0.9.0,
  satisfies (CP) for the operations listed there.
\item CPython integer arithmetic, and the rational arithmetic of its standard
  module \texttt{fractions}, are correct.
\item The program runs as written, with its assertions enabled, and its routines
  carry out the enclosure procedure of Section~\ref{app:procedure} and the
  certifications and runs of Section~\ref{app:runs}.
\end{enumerate}
The function \texttt{evans} of the module \texttt{ringcert.evans} follows steps
(i)--(iv) of Section~\ref{app:procedure} line by line, and the part of (TB3) that
concerns it is established by reading it against them. A certification (Section~\ref{app:runs})
is carried out by the function \texttt{certify\_cell} of the module
\texttt{ringcert.winding}, which follows the description of the certification in
Section~\ref{app:runs} line by line, and a run
by the function \texttt{certify\_cell\_adaptive} of the same module. Everything
else in this paper has a written proof that uses no machine computation.

\subsection{Two Elementary Lemmas}

The first estimate controls the end regions in the enclosure procedure;
the second confines the relevant zeros to $\Rect$.

\begin{lemma}
\label{lem:end-regions}
The following two estimates hold.
\begin{enumerate}[label=\textup{(\arabic*)}]
\item We have $\int_{-\infty}^{\xia}|\prof'|\dd\xi\le\frac94e^{2x_a}$ and
  $\int_{\xib}^{\infty}|\prof'|\dd\xi\le\frac94e^{2y_b}$, and both bounds are less
  than $3e^{-50}$.
\item Let $\ell\ge1$, $G\ge0$ and $\Re\lambda>0$, let $\psi=\psi_\ell(\cdot;\lambda,G)$,
  and let $[\alpha,\beta]$ be one of the intervals $[\xiL,\xia]$ and $[\xib,\xiR]$.
  Let $\psi_0$ solve $\psi_0''=\ell^2\psi_0$ with $\psi_0(\alpha)=\psi(\alpha)$ and
  $\psi_0'(\alpha)=\psi'(\alpha)$, and put $Y=(\psi,\psi')$, $Y_0=(\psi_0,\psi_0')$
  and $\epsend=3e^{-50}/\Re\lambda$. Then, for every $\xi\in[\alpha,\beta]$,
  \[
  |Y(\xi)-Y_0(\xi)|_w\le e^{\ell(\xi-\alpha)}\,|Y(\alpha)|_w\,\big(e^{\epsend}-1\big).
  \]
\end{enumerate}
\end{lemma}
\begin{proof}
(1) Write $A(\xi)=H(\frac{\xi-a_0}d)$, $B(\xi)=H(-\frac\xi d)$, so that
$\prof=\chi AB$, and $x=\frac{\xi-a_0}d$, $y=-\frac\xi d$. Since
$H'=\frac12\operatorname{sech}^2$, $A'=\frac1{2d}\operatorname{sech}^2x\ge0$ and
$B'=-\frac1{2d}\operatorname{sech}^2y\le0$; we use $H(x)\le e^{2x}$,
$\operatorname{sech}^2x\le4e^{2x}$, $0\le H\le1$ and $\operatorname{sech}^2\le1$. On
$(-\infty,\xia]$, $\prof'=\chi'AB+\chi(A'B+AB')$, and we bound the integrals of the
moduli of the two terms separately. For the first, $\chi$ is nondecreasing from
$\chi(\xiL)=0$ to $\chi(\xia)=1$, so $\int_{-\infty}^{\xia}|\chi'|=1$, and $AB\le A\le A(\xia)=H(x_a)$ since $A$ is
increasing; so $\int_{-\infty}^{\xia}|\chi'|AB\le H(x_a)\le e^{2x_a}$. For the second,
$0\le\chi\le1$, and the bounds just listed give $|A'B+AB'|\le\frac1{2d}(\operatorname{sech}^2x+H(x)\operatorname{sech}^2y)
\le\frac5{2d}e^{2x}$, and $\int_{-\infty}^{\xia}e^{2(\xi-a_0)/d}\dd\xi=\frac d2e^{2x_a}$,
so the second term contributes at most $\frac54e^{2x_a}$;
adding gives $\frac94e^{2x_a}$. On $[\xib,\infty)$ exchange the roles of $x$ and
$y$: $\chi$ is nonincreasing from $1$ to $0$, $AB\le B\le B(\xib)=H(y_b)\le e^{2y_b}$,
$|A'B+AB'|\le\frac1{2d}(\operatorname{sech}^2x\,H(y)+\operatorname{sech}^2y)\le\frac5{2d}e^{2y}$
and $\int_{\xib}^\infty e^{-2\xi/d}\dd\xi=\frac d2e^{2y_b}$, and the same addition gives
$\frac94e^{2y_b}$. By Equation~\eqref{eq:xayb}, both
bounds are at most $\frac94e^{-50}e^{8\cdot10^{-15}}<3e^{-50}$, as
$e^{8\cdot10^{-15}}<\frac43$.

(2) By the equation for $\psi$ in Hypothesis~I.\ref{I-hyp:N}, the vector
$Z=(\ell\psi,\psi')$ satisfies $Z'=A_0Z+\mathcal QZ$ with
$A_0=\begin{psmallmatrix}0&\ell\\ \ell&0\end{psmallmatrix}$ and
$\mathcal Q=\begin{psmallmatrix}0&0\\ \delta/\ell&0\end{psmallmatrix}$, where
$\delta=i\ell\prof'/(\lambda+i\ell\Omega_G)$, and $Z_0=(\ell\psi_0,\psi_0')$ satisfies
$Z_0'=A_0Z_0$; and $|Y|_w=\|Z\|_\infty$, the maximum of the moduli of the
components. Since $\Omega_G$ is real, $|\lambda+i\ell\Omega_G|\ge\Re\lambda$, so
the $\|\cdot\|_\infty$ operator norm of $\mathcal Q$ is
$b=|\prof'|/|\lambda+i\ell\Omega_G|\le|\prof'|/\Re\lambda$, and
$\int_\alpha^\beta b\le\epsend$ by (1). For $t\ge0$,
$e^{tA_0}=\begin{psmallmatrix}\cosh\ell t&\sinh\ell t\\ \sinh\ell t&\cosh\ell t\end{psmallmatrix}$
has norm $e^{\ell t}$. By variation of constants,
$Z(\xi)-Z_0(\xi)=\int_\alpha^\xi e^{(\xi-\eta)A_0}\mathcal Q(\eta)Z(\eta)\dd\eta$. Put
$u(\xi)=e^{-\ell(\xi-\alpha)}\|Z(\xi)-Z_0(\xi)\|_\infty$ and $c=\|Z(\alpha)\|_\infty$;
then $\|Z_0(\eta)\|_\infty\le e^{\ell(\eta-\alpha)}c$. Taking norms in the
variation of constants formula, with
$\|Z(\eta)\|_\infty\le\|Z_0(\eta)\|_\infty+\|Z(\eta)-Z_0(\eta)\|_\infty$, and multiplying
by $e^{-\ell(\xi-\alpha)}$ gives
$u(\xi)\le\int_\alpha^\xi b(\eta)(c+u(\eta))\dd\eta$, so by Gronwall's inequality
$c+u(\xi)\le c\exp(\int_\alpha^\xi b)$ and, as $\int_\alpha^\xi b\le\epsend$,
$u(\xi)\le c\,(e^{\epsend}-1)$. Since $c=|Y(\alpha)|_w$ and
$e^{\ell(\xi-\alpha)}u(\xi)=|Y(\xi)-Y_0(\xi)|_w$, this is the claim.
\end{proof}

\begin{lemma}
\label{lem:apriori}
Let $0\le G\le\Gmax'$. Every zero $\lambda$ of $\Evans_3(\cdot;G)$ with
$\Re\lambda\ge\frac9{200}$ lies in $\Rect$.
\end{lemma}
\begin{proof}
The proof compares $\psi_\ell$, $\ell=3$, with the solution $e^{\ell(\xi-\xiL)}$ of
$\psi''=\ell^2\psi$ on the whole interval $[\xiL,\xiR]$. If $|\lambda+i\ell\Omega_G|$
is large on $[\xiL,\xiR]$ compared with the total variation of $\prof$, the
coefficient of the equation for $\psi_\ell$ differs from $\ell^2$ by a term whose
modulus has small integral over $[\xiL,\xiR]$; then
$|\Evans_\ell(\lambda;G)-2\ell|<2\ell$, so $\Evans_\ell(\lambda;G)\ne0$.

Let $\ell=3$, let $\lambda$ be a zero of $\Evans_\ell(\cdot;G)$ with
$\Re\lambda\ge\frac9{200}$, so $\Re\lambda>0$, and put $TV=\int_{\R}|\prof'|\dd\xi$ and
$\Omega_{\max}=\frac12+\Gmax'e^{-2\xiL}$. We first bound $TV$. On $[\xia,\xib]$, $\chi=1$, so, with $A,B$ as in the proof of
Lemma~\ref{lem:end-regions}, $\prof'=A'B+AB'$ with $A'\ge0\ge B'$ and $0\le A,B\le1$.
Hence $\int_{\xia}^{\xib}|\prof'|\le\int_{\xia}^{\xib}(A'-B')
=A(\xib)-A(\xia)+B(\xia)-B(\xib)<2$, and Lemma~\ref{lem:end-regions}(1) gives
$TV<2+6e^{-50}$. Since $e>2.7$ and $2.7^7>1000$, $e^{50}>2.7^{49}=(2.7^7)^7>10^{21}$, so
$6e^{-50}<10^{-20}$; and
$\log2=\sum_{k\ge1}\frac1{k2^k}>\sum_{k=1}^6\frac1{k2^k}>0.69$. So
$TV/\log2<(2+10^{-20})/0.69<2.9$.

We next bound $\Omega_G$ on $[\xiL,\infty)$. For $\xi\in\R$, $0\le\prof\le1$ gives
$0\le\int_{-\infty}^\xi\prof(\xi')e^{2\xi'}\dd\xi'\le\frac12e^{2\xi}$, so, as
$G\le\Gmax'$,
$0\le\Omega_G(\xi)\le\frac12+Ge^{-2\xi}\le\Omega_{\max}$ for $\xi\ge\xiL$. By
Definition~I.\ref{I-def:ring}, $-2\xiL=4528212286209965\cdot2^{-51}<2.011$; $e<2.7183$ gives $e^2<7.3892$, and
$e^x\le1+x+x^2$ ($0\le x\le1$) gives $e^{0.011}<1.01113$, so
$e^{-2\xiL}<7.3892\cdot1.01113<7.4715$. As $\Gmax'=d$ and $d<0.0200001$,
$\Omega_{\max}<\frac12+0.0200001\cdot7.4715<0.6495$.

Finally we compare $\psi$ with $e^{\ell(\xi-\xiL)}$. Let $\psi=\psi_\ell(\cdot;\lambda,G)$ and
$\mu=\min_{\xi\in[\xiL,\xiR]}|\lambda+i\ell\Omega_G(\xi)|$, positive as $\Re\lambda>0$. The proof of
Lemma~\ref{lem:end-regions}(2), run on $[\xiL,\xiR]$ with
$\psi_0(\xi)=e^{\ell(\xi-\xiL)}$ (which has the same data at $\xiL$),
$\int_{\xiL}^{\xiR}b\le TV/\mu=:\varepsilon'$ and $|Y(\xiL)|_w=\ell$, gives
$|Y(\xiR)-Y_0(\xiR)|_w\le e^{\ell(\xiR-\xiL)}\ell\,(e^{\varepsilon'}-1)$. Since
$\psi_0'+\ell\psi_0=2\ell e^{\ell(\xi-\xiL)}$ and
$|(\psi'-\psi_0')+\ell(\psi-\psi_0)|\le2|Y-Y_0|_w$, we get
$|\Evans_\ell(\lambda;G)-2\ell|\le2\ell(e^{\varepsilon'}-1)$. If $\mu>TV/\log2$, then
$\varepsilon'<\log2$, the right side is less than $2\ell$, and
$\Evans_\ell(\lambda;G)\ne0$. As $\lambda$ is a zero, $\mu\le TV/\log2$, and $\mu$ is attained at some
$\xi\in[\xiL,\xiR]$, so $|\lambda+i\ell\,\Omega_G(\xi)|\le TV/\log2$. At that $\xi$,
$\Re\lambda\le TV/\log2<2.9$ and $|\Im\lambda+\ell\Omega_G(\xi)|\le TV/\log2$, so
$0\le\Omega_G(\xi)\le\Omega_{\max}$ gives $-\ell\,\Omega_{\max}-2.9<\Im\lambda<2.9$. As
$\ell=3$, $-3\Omega_{\max}-2.9>-3\cdot0.6495-2.9>-4.85$. Hence
$\lambda\in[\frac9{200},2.9)\times(-4.85,2.9)\subset\Rect$.
\end{proof}

\subsection{The Enclosure Procedure}
\label{app:procedure}

The procedure encloses $\Evans_\ell(\lambda;G)$ for all $\lambda$ in a box and
all $G$ in an interval. It crosses the end regions by
Lemma~\ref{lem:end-regions}, and $[\xia,\xib]$, where $\chi=1$, by a Taylor
method with coefficients and remainders enclosed in ball arithmetic. The Taylor
method needs the derivatives of $\prof$, which are explicit on $[\xia,\xib]$, where
$\prof$ is the product of the two factors $H$, but not elsewhere, as $\chi$ is known
only through the properties listed in Definition~I.\ref{I-def:ring};
Lemma~\ref{lem:end-regions}(2) needs only the bounds of
Lemma~\ref{lem:end-regions}(1) on $\prof'$. Fix the
Taylor order $\NT=16$, written $N$ in this subsection, and the points $s_0<s_1<\dots<s_n$, $n=360$, given by
$s_0=\xia$, $s_n=\xib$ and $s_j=\RN(\xia+\frac jn(\xib-\xia))$ for $0<j<n$, where
$\RN(q)$ is the binary64 number nearest to $q$ (ties to even). Only three facts are
used: the $s_j$ are dyadic rationals, they are strictly increasing, and the steps
$h_j=s_{j+1}-s_j$ are exact dyadic rationals.

The input is an integer $\ell\ge1$, a box $\Lambda^{\rm b}=\mathsf X+i\mathsf Y$ and a
real ball $\mathcal G^{\rm b}$. Each operation below is a ball operation; the
procedure \textbf{fails} if a comparison that it requires is not certified as
true, or if its output is not a box of finite radius.
\begin{enumerate}[label=\textup{(\roman*)}]
\item \emph{Coefficient data, independent of $\lambda$ and $G$.} For each $j<n$,
  and for $S=\{s_j\}$ (the \emph{point data}) and for $S=[s_j,s_{j+1}]$ (the
  \emph{step data}), compute
  \[
  \begin{gathered}
  \mathsf A_0=\mathrm{hull}\big(H(\tfrac{\min S-a_0}d),H(\tfrac{\max S-a_0}d)\big),\quad
  \bar{\mathsf A}_0=\mathrm{hull}\big(H(\tfrac{a_0-\max S}d),H(\tfrac{a_0-\min S}d)\big),\\
  \mathsf B_0=\mathrm{hull}\big(H(\tfrac{-\max S}d),H(\tfrac{-\min S}d)\big),\quad
  \bar{\mathsf B}_0=\mathrm{hull}\big(H(\tfrac{\min S}d),H(\tfrac{\max S}d)\big),
  \end{gathered}
  \]
  with $H(x)$ evaluated as $1/(1+\exp(-2x))$; then, for $0\le k\le N$,
  \[
  \begin{gathered}
  \mathsf A_{k+1}=\frac{2}{d(k+1)}\sum_{r=0}^k\mathsf A_r\bar{\mathsf A}_{k-r},\qquad
  \bar{\mathsf A}_{k+1}=-\mathsf A_{k+1},\\
  \mathsf B_{k+1}=-\frac{2}{d(k+1)}\sum_{r=0}^k\mathsf B_r\bar{\mathsf B}_{k-r},\qquad
  \bar{\mathsf B}_{k+1}=-\mathsf B_{k+1},
  \end{gathered}
  \]
  and $\mathsf w_k=\sum_{r\le k}\mathsf A_r\mathsf B_{k-r}$ for $k\le N+1$,
  $\mathsf f_k=(k+1)\mathsf w_{k+1}$, $\mathsf E^\pm_k=(\pm2)^ke^{\pm2S}/k!$ (with $e^{\pm2S}$ a
  ball containing $e^{\pm2\xi}$ for $\xi\in S$) and
  $\mathsf q_k=\sum_{r\le k}\mathsf w_r\mathsf E^+_{k-r}$ for $k\le N$. Then compute the
  balls
  \[
  \gamma_0=\mathrm{hull}\Big(0,\frac{e^{2x_a+2\xia}}{2/d+2}\Big),\qquad
  \gamma_{j+1}=\gamma_j+\Big(\sum_{k<N}\mathsf q_k^{\rm pt}\frac{h_j^{k+1}}{k+1}
  +\mathsf q_N^{\rm st}\frac{h_j^{N+1}}{N+1}\Big)\cap[0,\infty),
  \]
  where the superscripts ${\rm pt}$ and ${\rm st}$ refer to the point and the step
  data of step $j$, and, with $\mathsf g_0=\gamma_j$ for the point data,
  $\mathsf g_0=\mathrm{hull}(\gamma_j,\gamma_{j+1})$ for the step data, and
  $\mathsf g_r=\mathsf q_{r-1}/r$ for $1\le r<N$, put
  $\mathsf o_k=\sum_{r\le k}\mathsf g_r\mathsf E^-_{k-r}$ for $k<N$. The point data
  serve for the Taylor polynomials at $s_j$, and the step data, which enclose the
  coefficients at every point of the step, for the a priori enclosures on the step
  and the remainders.
\item \emph{Left end.} Require $\mathrm{lo}(\mathsf X)>0$. Let $\varepsilon$ be an upper
  bound of $3e^{-50}/\mathrm{lo}(\mathsf X)$, and require $\varepsilon\le1$. With
  $L_1=\xia-\xiL$, put
  $P=e^{\ell L_1}+\mathrm{box}(2\varepsilon\,\mathrm{up}(e^{\ell L_1}))$ and $P'=\ell P$.
  The condition $\varepsilon\le1$ gives $e^{\varepsilon}-1\le2\varepsilon$, which
  bounds the factor $e^{\epsend}-1$ of Lemma~\ref{lem:end-regions}(2) (proof of
  Proposition~\ref{prop:enclosure}); this is the origin of $2\varepsilon$ here and
  of $4\varepsilon$ in step (iv).
\item \emph{Middle steps.} For $j=0,1,\dots,n-1$ in turn, do the following, and then
  replace $(P,P')$ by the new values. For the point data and for the step data of
  step $j$, compute
  \[
  \begin{gathered}
  \mathsf D_0=\Lambda^{\rm b}+i\ell(\mathsf o_0+\mathcal G^{\rm b}\mathsf E^-_0),\qquad
  \mathsf D_k=i\ell(\mathsf o_k+\mathcal G^{\rm b}\mathsf E^-_k)\quad(1\le k<N),\\
  \mathsf R_0=1/\mathsf D_0,\qquad \mathsf R_k=-\mathsf R_0\sum_{r=1}^k\mathsf D_r\mathsf R_{k-r}\quad(1\le k<N),\\
  \mathsf V_k=i\ell\sum_{r\le k}\mathsf f_r\mathsf R_{k-r}+\ell^2[k=0]\quad(0\le k<N),
  \end{gathered}
  \]
  where $[k=0]$ is $1$ if $k=0$ and $0$ otherwise. Let $q_*=|\mathsf V^{\rm st}_0|^{\rm up}$,
  require that $q_*$ is finite, let $\kst$ be an upper bound of
  $(\max(q_*,\ell^2))^{1/2}$, require $\kst h_j\le700$, choose a dyadic number
  $M\ge\max(\kst|P|^{\rm up},|P'|^{\rm up})\,e^{\kst h_j}$, put
  $P^{\rm st}=P+\mathrm{box}(h_jM)$ and $P'^{\rm st}=P'+\mathrm{box}(h_jq_*M/\kst)$,
  and compute $\Psi_0=P$, $\Psi_1=P'$ and
  $\Psi_{k+2}=\frac1{(k+1)(k+2)}\sum_{r\le k}\mathsf V^{\rm pt}_r\Psi_{k-r}$, and in the same
  way $\Psi^{\rm st}_k$ from $P^{\rm st}$, $P'^{\rm st}$ and $\mathsf V^{\rm st}$, for
  $k\le N+1$. The new values are
  \[
  P_{\rm new}=\sum_{k<N}\Psi_kh_j^k+\Psi^{\rm st}_Nh_j^N,\qquad
  P'_{\rm new}=\sum_{k=1}^Nk\Psi_kh_j^{k-1}+(N+1)\Psi^{\rm st}_{N+1}h_j^N .
  \]
  Here $P^{\rm st}$ and $P'^{\rm st}$ are a priori enclosures of $\psi_\ell$ and
  $\psi_\ell'$ on the whole step, and the terms with $\Psi^{\rm st}$ are the Taylor
  remainders.
\item \emph{Right end and output.} With $L_2=\xiR-\xib$ and
  $W=\max(\ell|P|^{\rm up},|P'|^{\rm up})$, return
  \[
  \Evans^{\rm b}=e^{-\ell(\xiR-\xiL)}\Big(e^{\ell L_2}(P'+\ell P)
  +\mathrm{box}\big(4\varepsilon W\,\mathrm{up}(e^{\ell L_2})\big)\Big).
  \]
\end{enumerate}

\begin{proposition}
\label{prop:enclosure}
Let $\ell\ge1$ be an integer, let $\Lambda\subset\{\Re\lambda>0\}$ be a compact
axis-parallel box, and let $\Cell\subset[0,\infty)$ be a compact interval. Assume
that every ball operation satisfies \textup{(CP)}. Run the procedure of
Section~\ref{app:procedure} with $\ell$, a box $\Lambda^{\rm b}\supset\Lambda$ and a ball
$\mathcal G^{\rm b}\supset\Cell$. If it does not fail, it returns a box
$\Evans^{\rm b}$ of finite radius such that $\Evans_\ell(\lambda;G)\in\Evans^{\rm b}$ for every
$\lambda\in\Lambda$ and every $G\in\Cell$.
\end{proposition}
\begin{proof}
Fix $\lambda\in\Lambda$ and $G\in\Cell$, and let $\psi=\psi_\ell(\cdot;\lambda,G)$. We
show that every computed ball contains a corresponding exact quantity, defined
from $\lambda$, $G$ and $\chi$; the last one is $\Evans_\ell(\lambda;G)$. By (CP) it
suffices to check, at each step, that the exact quantity is the same function of
the earlier exact quantities, or else to justify the enclosure separately. For a
smooth function $F$ we write $F_k(\xi)=F^{(k)}(\xi)/k!$.

\emph{Step \textup{(i)}.} On $[\xia,\xib]$, $\chi=1$ and, $\chi$ being smooth, all its
derivatives vanish there, so $\prof^{(k)}=(AB)^{(k)}$, where
$A(\xi)=H(\frac{\xi-a_0}d)$, $B(\xi)=H(-\frac\xi d)$, $\bar A=1-A=H(\frac{a_0-\xi}d)$ and
$\bar B=1-B=H(\frac\xi d)$, as $1-H(x)=H(-x)$. Since $H'=2H(1-H)$, $A'=\frac2dA\bar A$, $\bar A'=-A'$,
$B'=-\frac2dB\bar B$ and $\bar B'=-B'$. The Leibniz rule gives
$A_{k+1}=\frac2{d(k+1)}\sum_{r\le k}A_r\bar A_{k-r}$, likewise for $B$, and
$\prof_k=\sum_{r\le k}A_rB_{k-r}$, $(\prof')_k=(k+1)\prof_{k+1}$,
$(e^{\pm2\xi})_k=(\pm2)^ke^{\pm2\xi}/k!$; the function $q=\prof e^{2\xi}$ has
$q_k=\sum_{r\le k}\prof_r(e^{2\xi})_{k-r}$. For $\xi\in S$, since $H$ is increasing,
$A(\xi)\in\mathsf A_0$, $\bar A(\xi)\in\bar{\mathsf A}_0$, $B(\xi)\in\mathsf B_0$ and
$\bar B(\xi)\in\bar{\mathsf B}_0$; by induction on $k$, each of
$\mathsf A_k,\bar{\mathsf A}_k,\mathsf B_k,\bar{\mathsf B}_k,\mathsf w_k,
\mathsf f_k,\mathsf E^\pm_k,\mathsf q_k$ contains the Taylor coefficient at $\xi$ of
$A,\bar A,B,\bar B,\prof,\prof',e^{\pm2\xi},q$ respectively. So the point data enclose
the Taylor coefficients at $s_j$, and the step data at every point of
$[s_j,s_{j+1}]$.

The function $g(\xi)=\int_{-\infty}^\xi\prof(\xi')e^{2\xi'}\dd\xi'$, for which
$\Omega_G=e^{-2\xi}(G+g)$, has $g'=q\ge0$, so it is
nondecreasing. On $(-\infty,\xia]$, $0\le\prof\le A\le e^{2(\xi-a_0)/d}$, so
$0\le g(\xia)\le\int_{-\infty}^{\xia}e^{2(\xi-a_0)/d+2\xi}\dd\xi=e^{2x_a+2\xia}/(2/d+2)$,
that is, $g(s_0)\in\gamma_0$. Taylor's theorem with the Lagrange remainder,
applied to $t\mapsto\int_{s_j}^{s_j+t}q$, gives
\[
g(s_{j+1})-g(s_j)=\sum_{k<N}q_k(s_j)\frac{h_j^{k+1}}{k+1}+q_N(\eta_j)\frac{h_j^{N+1}}{N+1}
\]
for some $\eta_j\in[s_j,s_{j+1}]$. Since $\mathsf q_k^{\rm pt}\ni q_k(s_j)$ and
$\mathsf q_N^{\rm st}\ni q_N(\eta_j)$, (CP) puts this increment in the ball in
parentheses defining $\gamma_{j+1}$; as it is nonnegative ($q\ge0$), it lies in the
intersection of that ball with $[0,\infty)$, so $g(s_j)\in\gamma_j$ for every $j$ by
induction, and $g(\xi)\in\mathrm{hull}(\gamma_j,\gamma_{j+1})$ for
$\xi\in[s_j,s_{j+1}]$ since $g$ is nondecreasing. As $g_r=q_{r-1}/r$ for $r\ge1$, the
balls $\mathsf o_k$ contain the Taylor coefficients of $ge^{-2\xi}$, at $s_j$ for the
point data and at every point of the step for the step data.

\emph{Step \textup{(ii)}.} Since $\lambda\in\Lambda^{\rm b}$,
$\Re\lambda\ge\mathrm{lo}(\mathsf X)>0$, so $\varepsilon\ge\epsend$ (as in
Lemma~\ref{lem:end-regions}) and, as $\varepsilon\le1$,
$e^{\epsend}-1\le e^\varepsilon-1\le(e-1)\varepsilon\le2\varepsilon$.
Lemma~\ref{lem:end-regions}(2) on $[\xiL,\xia]$, where $\psi_0(\xi)=e^{\ell(\xi-\xiL)}$
and $|Y(\xiL)|_w=\ell$, gives $|\psi(\xia)-e^{\ell L_1}|\le2\varepsilon e^{\ell L_1}$ and
$|\psi'(\xia)-\ell e^{\ell L_1}|\le2\ell\varepsilon e^{\ell L_1}$. Hence
$\psi(\xia)\in P$ and $\psi'(\xia)\in\ell P=P'$.

\emph{Step \textup{(iii)}.} Assume $\psi(s_j)\in P$ and $\psi'(s_j)\in P'$. On
$[\xia,\xib]$ the coefficient $V=\ell^2+i\ell\prof'/(\lambda+i\ell\Omega_G)$ is smooth, so
$\psi$ is smooth there. Put $D=\lambda+i\ell\Omega_G=\lambda+i\ell(g+G)e^{-2\xi}$ and
$R=1/D$, where $D\ne0$ as $\Re D=\Re\lambda>0$. The identity $DR=1$ gives
$R_0=1/D_0$ and $R_k=-R_0\sum_{r=1}^kD_rR_{k-r}$ for $k\ge1$, and
$V=\ell^2+i\ell\prof'R$ gives $V_k=i\ell\sum_{r\le k}(\prof')_rR_{k-r}+\ell^2[k=0]$. By
step (i) and (CP), $\mathsf V^{\rm pt}_k\ni V_k(s_j)$ and
$\mathsf V^{\rm st}_k\ni V_k(\xi)$ for every $\xi\in[s_j,s_{j+1}]$; in particular
$|V|\le q_*\le\kst^2$ on the step. We next justify the a priori enclosures
$P^{\rm st}$ and $P'^{\rm st}$. Put
$u(t)=\max(\kst|\psi(s_j+t)|,|\psi'(s_j+t)|)$ for $t\in[0,h_j]$. Then
$\kst|\psi(s_j+t)|\le\kst|\psi(s_j)|+\int_0^t\kst|\psi'|$ and
$|\psi'(s_j+t)|\le|\psi'(s_j)|+\int_0^t\kst\cdot\kst|\psi|$, so
$u(t)\le u(0)+\kst\int_0^tu$, and Gronwall's inequality gives
$u(t)\le u(0)e^{\kst t}\le M$; the last inequality holds because $\psi(s_j)\in P$
and $\psi'(s_j)\in P'$ give $u(0)\le\max(\kst|P|^{\rm up},|P'|^{\rm up})$, and
$t\le h_j$. Hence, for $\xi\in[s_j,s_{j+1}]$,
$|\psi(\xi)-\psi(s_j)|\le h_jM$ and, since $|\psi|\le M/\kst$ on the step,
$|\psi'(\xi)-\psi'(s_j)|\le\int|V||\psi|\le h_jq_*M/\kst$; so $\psi(\xi)\in P^{\rm st}$
and $\psi'(\xi)\in P'^{\rm st}$. The equation $\psi''=V\psi$ gives
$(k+1)(k+2)\psi_{k+2}=\sum_{r\le k}V_r\psi_{k-r}$, so, by induction on $k$ with
(CP), starting from $\Psi_0=P\ni\psi(s_j)$ and $\Psi_1=P'\ni\psi'(s_j)$, and for the
step data from $P^{\rm st}$ and $P'^{\rm st}$, $\Psi_k\ni\psi_k(s_j)$ and
$\Psi^{\rm st}_k\ni\psi_k(\xi)$ for every $\xi\in[s_j,s_{j+1}]$ and $k\le N+1$.
Taylor's theorem with the integral remainder gives
\[
\psi(s_{j+1})=\sum_{k<N}\psi_k(s_j)h_j^k+h_j^N\int_0^1N(1-\tau)^{N-1}\psi_N(s_j+\tau h_j)\dd\tau .
\]
The weight $N(1-\tau)^{N-1}$ is a probability density on $[0,1]$, so the integral
lies in the closed convex hull of $\{\psi_N(\xi):\xi\in[s_j,s_{j+1}]\}$, which lies in
the closed convex box $\Psi^{\rm st}_N$. Hence $\psi(s_{j+1})\in P_{\rm new}$; the same
argument for $\psi'$, with $(\psi')_r=(r+1)\psi_{r+1}$, gives
$\psi'(s_{j+1})\in P'_{\rm new}$. Step (ii) gives the assumption for $j=0$, as
$s_0=\xia$, so by induction on $j$, $\psi(\xib)\in P$ and
$\psi'(\xib)\in P'$ at the end of step (iii).

\emph{Step \textup{(iv)}.} Apply Lemma~\ref{lem:end-regions}(2) on $[\xib,\xiR]$, where
$\psi_0$ has the data of $\psi$ at $\xib$; by step (iii),
$|Y(\xib)|_w\le W$, and, as in step (ii), $e^{\epsend}-1\le2\varepsilon$.
Since $(\psi_0'+\ell\psi_0)'=\ell(\psi_0'+\ell\psi_0)$,
$(\psi_0'+\ell\psi_0)(\xiR)=e^{\ell L_2}(\psi'+\ell\psi)(\xib)$, and
\[
|(\psi'+\ell\psi)(\xiR)-e^{\ell L_2}(\psi'+\ell\psi)(\xib)|\le2|Y(\xiR)-Y_0(\xiR)|_w
\le2e^{\ell L_2}|Y(\xib)|_w\cdot2\varepsilon\le4\varepsilon We^{\ell L_2}.
\]
Multiplying by $e^{-\ell(\xiR-\xiL)}$ gives $\Evans_\ell(\lambda;G)\in\Evans^{\rm b}$.
\end{proof}

\subsection{Certified Winding Numbers}

We read the winding number of $\Evans_3(\cdot;G)$ along the boundary of a rectangle
off enclosures of its values on the segments of the boundary, simultaneously for all
$G$ in an interval. A set $B\subset\C$ is of \textbf{type P} if it does not meet $(-\infty,0]$, and
of \textbf{type N} if it does not meet $[0,\infty)$. On a set of type P the principal
argument is continuous, and on a set of type N a continuous argument with values in
$(0,2\pi)$ exists.

\begin{proposition}
\label{prop:winding}
Let $\Gamma$ be the boundary of a closed rectangle $\Rect'$ in $\{\Re\lambda>0\}$,
traversed counterclockwise and cut into segments $[p_i,p_{i+1}]$,
$i=0,\dots,I-1$, with $p_I=p_0$. Let $\Cell$ be a set, let $E(\cdot;G)$ be continuous
on $\Gamma$ for each $G\in\Cell$, and let $\mathcal B_i\subset\C$ be sets with
$E(\lambda;G)\in\mathcal B_i$ for all $\lambda\in[p_i,p_{i+1}]$ and $G\in\Cell$. Assume that
each $\mathcal B_i$ is of type P or of type N, and assign to each $i$ one type for which this
holds. Then $E(\cdot;G)\ne0$ on $\Gamma$ and, for every $G\in\Cell$,
\[
\wind\big(E(\Gamma;G),0\big)=\sum_{i\ \text{of type N}}\Big(\mathbf 1\big[\Im E(p_{i+1};G)<0\big]-\mathbf 1\big[\Im E(p_i;G)<0\big]\Big).
\]
In this sum the two terms at a point shared by two segments of type N cancel.
\end{proposition}
\begin{proof}
Fix $G\in\Cell$. Every $\mathcal B_i$ misses $0$, so $E(\cdot;G)\ne0$ on $\Gamma$. The
principal argument $\mathrm{Arg}\colon\C\setminus\{0\}\to(-\pi,\pi]$ is continuous on
$\C\setminus(-\infty,0]$, and
$\theta_{\rm N}(z)=\mathrm{Arg}(z)+2\pi\,\mathbf 1[\Im z<0]$ is continuous on
$\C\setminus[0,\infty)$, with values in $(0,2\pi)$. Since $E(\cdot;G)$ maps the
$i$-th segment into $\mathcal B_i$, $\mathrm{Arg}(E(\cdot;G))$ is a continuous
argument on a segment of type P, and $\theta_{\rm N}(E(\cdot;G))$ on one of type N. Each differs by a constant from a continuous argument along the
whole curve, so the winding number is $\frac1{2\pi}$ times the sum over the
segments of their increments. On a segment of type N, by the definition of
$\theta_{\rm N}$, the increment is that of $\mathrm{Arg}$ plus $2\pi$ times the
difference of the indicators at the two ends. The increments of $\mathrm{Arg}$ telescope around the
closed curve to $0$, and what remains is the indicator sum. At a point shared by
two segments of type N the indicator enters once with each sign.
\end{proof}

\subsection{The Two Computations}
\label{app:runs}

We reduce \textup{(N1)} and \textup{(N2)} to the statements \textup{(N1$'$)} and
\textup{(N2$'$)} below, about winding numbers along the boundaries of finitely many
rectangles, and verify these by two computations. The line
$\Re\lambda=\frac{211}{2500}$ in \textup{(N2$'$)} is the right edge of the strip of
\textup{(N2)}.

\begin{definition}
\label{def:Nprime}
We consider the following two statements.
\begin{enumerate}[label=\textup{(N\arabic*$'$)}]
\item For every $G\in[0,\Gmax']$, $\Evans_3(\cdot;G)\ne0$ on $\partial\Rect$ and
  $\wind(\Evans_3(\partial\Rect;G),0)=1$.
\item There are $236$ closed intervals $\Cell_i$, $1\le i\le236$, with
  union $[0,\Gmax']$, and closed rectangles
  $\Rect_1,\dots,\Rect_{236}\subset\Rect$, the left edge of
  each of which lies on the line $\Re\lambda=\frac{211}{2500}$, such that for
  every $i$ and every $G\in\Cell_i$, $\Evans_3(\cdot;G)\ne0$ on $\partial\Rect_i$
  and $\wind(\Evans_3(\partial\Rect_i;G),0)=1$.
\end{enumerate}
\end{definition}

\begin{lemma}
\label{lem:Nprime-N}
\textup{(N1$'$)} implies \textup{(N1)}. \textup{(N1)} and \textup{(N2$'$)} imply
\textup{(N2)}.
\end{lemma}
\begin{proof}
Let $G\in[0,\Gmax]\subset[0,\Gmax']$. By
Lemma~I.\ref{I-lem:evans-basic}\textup{(1)} and its last assertion, $\Evans_3(\cdot;G)$ is holomorphic on
$\{\Re\lambda>0\}$, which contains $\Rect$ and every $\Rect_i$.

Assume \textup{(N1$'$)}. By the argument principle, $\Evans_3(\cdot;G)$ has exactly
one zero, counted with multiplicity, in the interior of $\Rect$. By
Lemma~\ref{lem:apriori}, every zero with $\Re\lambda\ge\frac9{200}$
lies in $\Rect$, and by \textup{(N1$'$)} none lies on $\partial\Rect$; so the zeros in
$\{\Re\lambda\ge\frac9{200}\}$ are exactly those in the interior of $\Rect$. This is
\textup{(N1)}.

Assume \textup{(N1)} and \textup{(N2$'$)}, let $\lamst$ be the zero of
\textup{(N1)}, and choose $i$ with $G\in\Cell_i$. By \textup{(N2$'$)} and the argument
principle, $\Evans_3(\cdot;G)$ has a zero $\lambda'$ in the interior of $\Rect_i$, whose left edge
lies on $\{\Re\lambda=\frac{211}{2500}\}$, so $\Re\lambda'>\frac{211}{2500}>\frac9{200}$, and $\lambda'=\lamst$ by \textup{(N1)}.
So the only zero of $\Evans_3(\cdot;G)$ in $\{\Re\lambda\ge\frac9{200}\}$ has real
part greater than $\frac{211}{2500}$, and there is none in
$\{\frac9{200}\le\Re\lambda\le\frac{211}{2500}\}$. This is \textup{(N2)}.
\end{proof}

Both computations use the procedure of Section~\ref{app:procedure} with $\ell=3$.
For a closed interval $\Cell$ with rational endpoints, a closed rectangle
$\Rect'=[x_0,x_1]\times[y_0,y_1]\subset\{\Re\lambda>0\}$ with rational vertices, a
number $\nu>0$ and integers $J,K\ge1$, the \textbf{certification} of
$(\Cell,\Rect')$ with density $\nu$, at most $J$ rounds and at most $K$ segments
is the following computation. The boundary of $\Rect'$ is traversed
counterclockwise from $(x_0,y_0)$, and each edge, of length $L$, is cut into
$\max(2,\lceil\nu L\rceil)$ segments of equal length; all cut points are rational.
For a segment $[p,q]$, the procedure is run with
$\Lambda^{\rm b}=(X_p\cup X_q)+i(Y_p\cup Y_q)$, where $X_p,X_q,Y_p,Y_q$ are the balls
that contain $\Re p$, $\Re q$, $\Im p$, $\Im q$, and with $\mathcal G^{\rm b}$ the
union of the balls that contain the two endpoints of $\Cell$; by (CP),
$\Lambda^{\rm b}\supset[p,q]$ and $\mathcal G^{\rm b}\supset\Cell$. The resulting box
$\Evans^{\rm b}$ is accepted as of type P if it is certified that
$\Im\Evans^{\rm b}>0$, or $\Im\Evans^{\rm b}<0$, or $\Re\Evans^{\rm b}>0$; otherwise it
is accepted as of type N if it is certified that $\Re\Evans^{\rm b}<0$. The first
round treats every segment; each later round treats the two halves, cut at the
midpoint, of every segment not accepted in the previous round, including those for
which the procedure failed. A segment is cut only if there are at most $K$
segments at that moment, and at most $J$ rounds are made; otherwise the
certification does not succeed. If every segment has been accepted, then at every
point $p$ where a segment of type N meets a segment of type P the procedure is run
with $\Lambda^{\rm b}$ the box of the balls that contain $\Re p$ and $\Im p$ and
with the same $\mathcal G^{\rm b}$, and the
indicator $\mathbf 1[\Im\Evans_3(p;G)<0]$ is taken to be $1$ if the resulting box
is certified to have negative imaginary part and $0$ if it is certified to have
positive imaginary part. By Proposition~\ref{prop:winding}, only endpoints of
segments of type N enter the indicator sum, and at a point shared by two such
segments the terms cancel; so the sum needs the sign of $\Im\Evans_3$ only at
endpoints of segments of type N that are not shared with another segment of
type N. The certification \textbf{succeeds} if every segment is
accepted and every such indicator is decided; its result is then the indicator sum
of Proposition~\ref{prop:winding}.

A \textbf{run} with at most $D$ halvings certifies a list of pairs
$(\Cell,\Rect')$: if the certification of $(\Cell,\Rect')$ does not succeed, $\Cell$
is halved at its midpoint and each half is certified with the same $\Rect'$, at
most $D$ times in succession; the certified pieces take the place of $\Cell$, and a
piece still not certified after $D$ halvings is reported as failed.

The run for \textup{(N1$'$)} uses the rectangle $\Rect$ and the $25$ intervals
\[
\Big[\frac{j}{20}d,\frac{j+1}{20}d\Big]\quad(0\le j\le14),\qquad
\Big[\frac{30+j}{40}d,\frac{31+j}{40}d\Big]\quad(0\le j\le9),
\]
whose union is $[0,\Gmax']$, with density $\nu=20$, at most $16$ rounds, at most
$6000$ segments and at most $4$ halvings. The run for \textup{(N2$'$)} uses $236$
intervals $\Cell_1,\dots,\Cell_{236}$, numbered from left to right, and rectangles
$\Rect_i=[\frac{211}{2500},x_i]\times[y_i,y_i']$, with density
$\nu=3/(x_i-\frac{211}{2500})$, at most $10$ rounds, at most $800$ segments and at
most $5$ halvings. These $236$ pairs, with every number an exact fraction, form the file
\texttt{sector3\_gap\_cells.json} of the reproduction package
(Section~\ref{app:reproduce}), and the file \texttt{sector3\_count\_cells.json}
lists the $25$ intervals of the first run, also numbered from left to right. In each of these lists the endpoints of
the intervals are rational, the first interval starts at $0$, the last ends at
$d=\Gmax'$, and consecutive intervals share an endpoint; every $\Rect_i$ has
rational vertices, satisfies $x_i>\frac{211}{2500}$ and $y_i<y_i'$, and lies in
$\Rect$. The re-check of Section~\ref{app:reproduce} confirms these facts, and that
every rectangle of the first run is $\Rect$, for the output files of the runs, in
exact rational arithmetic.

\begin{proposition}
\label{prop:runs}
\textup{(N1$'$)} and \textup{(N2$'$)} hold.
\end{proposition}
\begin{proof}
The proof is the computation just described, carried out under the trusted base of
Section~\ref{app:arith}. Each of the $25$ certifications of the first run and each
of the $236$ of the second succeeded at once, with result $1$, so no interval was
halved and the output files list exactly the intervals and rectangles described
above.

It remains to show that a successful certification proves the corresponding
winding statement. Consider a certification of $(\Cell,\Rect')$ that succeeds with result $1$. By
Lemma~I.\ref{I-lem:evans-basic}\textup{(1)} and its last assertion, $\Evans_3(\cdot;G)$ is continuous on
$\{\Re\lambda>0\}\supset\partial\Rect'$ for every $G\in\Cell$. By
Proposition~\ref{prop:enclosure}, the box of each accepted segment $[p,q]$
contains $\Evans_3(\lambda;G)$ for all $\lambda\in[p,q]$ and $G\in\Cell$, and each decided
indicator equals $\mathbf 1[\Im\Evans_3(p;G)<0]$ for every $G\in\Cell$. By (CP) for
the comparisons, a box certified to have $\Im\Evans^{\rm b}>0$ or
$\Im\Evans^{\rm b}<0$ misses the real line, and one certified to have
$\Re\Evans^{\rm b}>0$ misses $(-\infty,0]$, so a box accepted as of type P misses
$(-\infty,0]$; one accepted as of type N, certified to have $\Re\Evans^{\rm b}<0$,
misses $[0,\infty)$. The accepted segments are consecutive and cover
$\partial\Rect'$. By Proposition~\ref{prop:winding}, $\Evans_3(\cdot;G)\ne0$ on
$\partial\Rect'$ and $\wind(\Evans_3(\partial\Rect';G),0)=1$ for every $G\in\Cell$.
Applied to the first run, whose intervals cover $[0,\Gmax']$, this gives
\textup{(N1$'$)}; applied to the second, with the properties of the $\Cell_i$ and
$\Rect_i$ stated above, it gives \textup{(N2$'$)}.
\end{proof}

\subsection{The Closed-Form Inequality}

The proof of \textup{(N3)} is a computation by hand from the exact values of
$a_0$ and $d$ in Definition~I.\ref{I-def:ring} and from
Equation~\eqref{eq:xayb}.

\begin{proposition}
\label{prop:N3}
\textup{(N3)} holds.
\end{proposition}
\begin{proof}
We first bound $\epsst$. Since $H(x)\le e^{2x}$, Equation~\eqref{eq:xayb} gives
$\epsst\le e^{2\max(x_a,y_b)}<e^{-49}<10^{-21}$, as $e^{49}>2.7^{49}>10^{21}$
(proof of Lemma~\ref{lem:apriori}).

We next bound $e^{a_0}$. From Definition~I.\ref{I-def:ring}, $a_0=-7304210039150669\cdot2^{-54}<-0.4054$. The partial
sum $\sum_{k\le5}0.4054^k/k!$ of the exponential series is a lower bound for
$e^{0.4054}$, and it exceeds $1.49989$, so $e^{a_0}<e^{-0.4054}<1/1.49989<0.66672$.
Squaring and multiplying, $0.66672^2<0.444516$, $0.66672^3<0.296368$,
$0.66672^6<0.087834$ and $0.66672^9<0.026032$.

We now bound $\Thl{\ell}$ for $\ell\in\{6,9\}$. Put $u=\pi\ell d$. Since $\pi<3.1416$ and $0<d<0.0200001$, we have $0<u<0.377$ for $\ell=6$ and
$0<u<0.5655$ for $\ell=9$; in particular $0<u<\pi$. For $u\ge0$,
$\sin u\ge u-u^3/6$, since the difference vanishes at $0$ and has derivative
$\cos u-1+u^2/2$, which is nonnegative (it vanishes at $0$ and its derivative $u-\sin u$ is
nonnegative). So $\Thl{\ell}=u/\sin u\le(1-u^2/6)^{-1}$ as
$u^2<6$. Hence $\Thl{6}\le(1-0.377^2/6)^{-1}<(1-0.023689)^{-1}<1.024264$
and $\Thl{9}<(1-0.0533)^{-1}<1.0564$. Therefore
\[
\begin{gathered}
e^{6a_0}\Thl{6}+\epsst<0.087834\cdot1.024264+10^{-21}<0.08997,\\
e^{9a_0}\Thl{9}+\epsst<0.026032\cdot1.0564+10^{-21}<0.0276,
\end{gathered}
\]
and both are less than $\frac9{100}$, which is \textup{(N3)}.
\end{proof}

\begin{proof}[Proof of Theorem~\ref{thm:N-holds}]
By Proposition~\ref{prop:runs}, \textup{(N1$'$)} and \textup{(N2$'$)} of
Definition~\ref{def:Nprime} hold; by Lemma~\ref{lem:Nprime-N}, \textup{(N1)} holds, and
then \textup{(N2)}; by Proposition~\ref{prop:N3}, \textup{(N3)} holds.
\end{proof}

\subsection{Reproduction}
\label{app:reproduce}

The computations of Section~\ref{app:runs} can be repeated with version~1.0.0 of
the package \texttt{euler-ring-certificate}, which is available under the Apache
License 2.0 in the public repository \url{https://github.com/mcoiculescu/NonUniqueEuler}
and contains the code, the two
certificate files of Section~\ref{app:runs} and the reference outputs. It requires
only python-flint~0.9.0, which bundles FLINT~3.6.0, under CPython~3.10 or later.
The program refuses to run with assertions disabled, and, unless explicitly
allowed, with another FLINT version. From the top directory of the package, the
command \mbox{\texttt{./reproduce.sh 8}} runs, with $8$ worker processes, among other
items: an environment check; the unit tests; the certification of the $25$
intervals of \textup{(N1$'$)} and of the $236$ intervals and rectangles of
\textup{(N2$'$)} (commands \texttt{ringcert sector3-count} and
\texttt{ringcert sector3-gap}); a re-check of the outputs in exact rational
arithmetic (\texttt{ringcert check}); and a comparison with the reference
outputs. A successful run reports all items as passed.

The outputs record, for every certified interval, its endpoints, its rectangle,
its result, and the numbers of segments and indicators. From the output files
alone, the re-check confirms that the intervals of each run cover $[0,\Gmax']$,
with consecutive intervals sharing an endpoint, that every result is $1$, that no
indicator is undecided and no interval failed, that every rectangle of the first
run is $\Rect$, and that every rectangle of the second run is nondegenerate, has
left edge on the line $\Re\lambda=\frac{211}{2500}$ and lies in $\Rect$; it also
checks that these output files are the ones recorded in the summary of the run.

The output files \texttt{sector3\_count.json} and \texttt{sector3\_gap.json}
contain no timing or environment data; on an Apple arm64 machine with
CPython~3.13.5 they are byte-identical to the reference outputs of the package.
The SHA-256 digests of the certificate files and of these outputs are the following.
\begin{center}
\scriptsize
\begin{tabular}{ll}
\texttt{sector3\_count\_cells.json} &
\texttt{518a0ca4f05fd24bae273aebea20862fad1c1beb46dbdee9f8b5fb76c8729dc1}\\
\texttt{sector3\_gap\_cells.json} &
\texttt{fb95a3f72eee1f017ddc5c26e7dfca52841759603e34b0131a440c8992086c91}\\
\texttt{sector3\_count.json} &
\texttt{2aa43b18d7289b7288f4bc660b1fdbec547458949ee675ef4b23f65cbcf1a1c9}\\
\texttt{sector3\_gap.json} &
\texttt{f06feec05cc6633b390bb9ee6191ce59c3bd5231a01ffa0f48263c8c9173a786}
\end{tabular}
\end{center}
The reference outputs were produced on that machine only. With another FLINT
version or on other hardware the numbers of segments may differ, and an interval
may be halved. The proof needs only that the re-check passes, that is, that the output files
have the properties listed above. If some intervals
were halved, the argument of Proposition~\ref{prop:runs} gives \textup{(N1$'$)},
and \textup{(N2$'$)} with the certified pieces, each paired with the rectangle of
its interval, in place of the $236$ intervals; the proof of
Lemma~\ref{lem:Nprime-N} does not use their number.

\section{Global Smooth Solutions}
\label{app:smooth}

We prove Theorem~\ref{thm:smooth}: global existence and uniqueness for smooth,
compactly supported data, and properties (4)--(7) of the solution. The argument
is a Picard iteration through linear transport. The iterates obey uniform
bounds on every bounded time interval and converge with all their spatial
derivatives. Lemma~\ref{lem:grad-bound} bounds the velocity gradient by a
logarithm of the vorticity gradient. This gives a double-exponential bound on
the latter, finite at every time. We use $L^2$-Sobolev estimates and, besides
elementary real analysis, only
Section~I.\ref{I-sec:prelim}, Section~\ref{sec:II-prelim} and the facts on ordinary differential equations
collected in Lemma~\ref{lem:ode}. For
the classical theory of smooth solutions see \cite{Wolibner}, \cite{Holder} and
\cite{Kato1967}.

\subsection{Ordinary Differential Equations}

We use the following ODE facts from \cite{Teschl} to construct the smooth flow.

\begin{lemma}[Theorems 2.2, 2.10, 2.13 and 2.17 in \cite{Teschl}]
\label{lem:ode}
Let $n\ge1$ and $k\ge1$ be integers, let $U\subset\R\times\R^n$ be open, and let
$F\in C^k(U;\R^n)$.
\begin{enumerate}[label=\textup{(\roman*)}]
\item For every $(s,y)\in U$ the problem $\dot x=F(t,x)$, $x(s)=y$, has a unique
  maximal solution, defined on an open interval containing $s$. We write
  $\flow_F(t,s,y)$ for its value at time $t$.
\item Around each point of $U$ there is an open set $I\times\mathcal U\subset U$,
  with $I$ an interval, such that $\flow_F$ and $\partial_t\flow_F$ are defined
  and of class $C^k$ on $I\times I\times\mathcal U$, and there $D_y\flow_F$
  solves the variational equation
  \[
  \partial_tD_y\flow_F(t,s,y)=D_xF\big(t,\flow_F(t,s,y)\big)\,D_y\flow_F(t,s,y).
  \]
\item If $U=\R\times\R^n$ and for every $T>0$ there are numbers $M(T)$, $L(T)$
  with $|F(t,x)|\le M(T)+L(T)|x|$ on $[-T,T]\times\R^n$, then every maximal
  solution is defined on all of $\R$.
\end{enumerate}
\end{lemma}

Item (i) is Theorem 2.2 of \cite{Teschl} (the Picard-Lindel\"of theorem; a
$C^1$ function is locally Lipschitz) together with Theorem 2.13 there, which
produces the maximal solution. Item (ii) is Theorem 2.10 of \cite{Teschl}, and
item (iii) is Theorem 2.17 there. We use these facts through the following
lemma, which we apply to the velocity of each iterate and of the limit;
Section~\ref{sec:single} also uses it directly.

\begin{lemma}
\label{lem:flow}
Let $u\in C^\infty(\R\times\R^2;\R^2)$ be bounded on $[-T,T]\times\R^2$ for every
$T>0$, and write $\flow_u$ for the solution map of $\dot x=u(t,x)$ from
Lemma~\ref{lem:ode}(i). Then:
\begin{enumerate}[label=\textup{(\alph*)}]
\item $\flow_u(t,s,y)$ is defined for all $(t,s,y)\in\R\times\R\times\R^2$, and
  $\flow_u\in C^\infty(\R\times\R\times\R^2)$;
\item $\flow_u(t,\tau,\flow_u(\tau,s,y))=\flow_u(t,s,y)$ for all $t,\tau,s,y$, and
  each map $\flow_u(t,s,\cdot)$ is a $C^\infty$ diffeomorphism of $\R^2$ with
  inverse $\flow_u(s,t,\cdot)$;
\item $|\flow_u(t,s,y)-y|\le|t-s|\sup_{[-T,T]\times\R^2}|u|$ whenever
  $s,t\in[-T,T]$;
\item $Y(t)=D_y\flow_u(t,s,y)$ solves $\partial_tY=\nabla u(t,\flow_u(t,s,y))\,Y$
  with $Y(s)=I$, and if $\divop u=0$, then $\det D_y\flow_u\equiv1$.
\end{enumerate}
\end{lemma}

\begin{proof}
Lemma~\ref{lem:ode}(iii), with $M(T)=\sup_{[-T,T]\times\R^2}|u|$ and $L(T)=0$,
shows that $\flow_u$ is defined on $\R\times\R\times\R^2$. Both sides of the
identity in (b) solve $\dot x=u(t,x)$ as functions of $t$ and agree at $t=\tau$,
so they agree by the uniqueness in Lemma~\ref{lem:ode}(i). With $t=s$ the
identity gives $\flow_u(s,\tau,\flow_u(\tau,s,y))=y$, so $\flow_u(\tau,s,\cdot)$ is
a bijection with inverse $\flow_u(s,\tau,\cdot)$. Part (c) follows from the
integral form
$\flow_u(t,s,y)-y=\int_s^tu(\sigma,\flow_u(\sigma,s,y))\dd\sigma$ of the
equation, since $\sigma$ lies in $[-T,T]$ when $s$ and $t$ do.

We prove the smoothness in (a). Fix $k\ge1$. Lemma~\ref{lem:ode}(ii) gives smoothness only
near each point, for times in a small interval. We therefore cover a trajectory
by finitely many such neighborhoods and write $\flow_u$ near a given point, by
(b), as a composition of finitely many maps, each of which is $C^k$.
Fix $(t_*,s_*,y_*)$, let $J$ be the
closed interval with endpoints $s_*$ and $t_*$, and put
$\gamma(\tau)=\flow_u(\tau,s_*,y_*)$. For every $\tau\in J$,
Lemma~\ref{lem:ode}(ii) gives an open set
$I_\tau\times\mathcal U_\tau\ni(\tau,\gamma(\tau))$ on which $\flow_u$ is $C^k$ in
the sense stated there. The sets
$\{\sigma\in J\cap I_\tau:\gamma(\sigma)\in\mathcal U_\tau\}$, $\tau\in J$, are
relatively open in $J$ and cover $J$. By the Lebesgue number lemma there are
points $s_*=\tau_0,\tau_1,\dots,\tau_N=t_*$, monotone in the index, and points
$\tau'_i\in J$ such that $\tau_i,\tau_{i+1}\in I_{\tau'_i}$ and
$\gamma(\tau_i)\in\mathcal U_{\tau'_i}$ for $0\le i<N$. By (b), for $(t,s,y)$ near $(t_*,s_*,y_*)$,
\[
\flow_u(t,s,y)=\flow_u\Big(t,\tau_{N-1},\flow_u\big(\tau_{N-1},\tau_{N-2},\cdots
\flow_u(\tau_2,\tau_1,\flow_u(\tau_1,s,y))\cdots\big)\Big).
\]
The innermost map $(s,y)\mapsto\flow_u(\tau_1,s,y)$ is $C^k$ near
$(\tau_0,\gamma(\tau_0))$; each map $y\mapsto\flow_u(\tau_{i+1},\tau_i,y)$ is $C^k$
on $\mathcal U_{\tau'_i}$, which is a neighborhood of $\gamma(\tau_i)$, and it
maps $\gamma(\tau_i)$ to $\gamma(\tau_{i+1})$; and
$(t,x)\mapsto\flow_u(t,\tau_{N-1},x)$ is $C^k$ for $t\in I_{\tau'_{N-1}}\ni t_*$
and $x\in\mathcal U_{\tau'_{N-1}}$. By
continuity of these maps, for $(t,s,y)$ near $(t_*,s_*,y_*)$ each of them is
evaluated in the set on which it is $C^k$, so, by the chain rule, $\flow_u$ is $C^k$ near $(t_*,s_*,y_*)$. Since $k$
is arbitrary, $\flow_u\in C^\infty$, and by (b) the maps $\flow_u(t,s,\cdot)$ are
$C^\infty$ diffeomorphisms.

We prove (d). By (a), $\flow_u$ is $C^\infty$. Differentiating
$\partial_t\flow_u(t,s,y)=u(t,\flow_u(t,s,y))$ in $y$ and exchanging the order
of the derivatives gives
\[
\partial_tD_y\flow_u(t,s,y)=\nabla u\big(t,\flow_u(t,s,y)\big)\,D_y\flow_u(t,s,y),
\]
and $Y(s)=I$ because $\flow_u(s,s,y)=y$. The matrix
$Y(t)$ is invertible by (b), since $\flow_u(t,s,\cdot)$ is a $C^\infty$
diffeomorphism, and Jacobi's formula gives
$\partial_t\det Y=\det Y\operatorname{tr}(Y^{-1}\partial_tY)=\det
Y\operatorname{tr}\nabla u=(\divop u)\det Y$. If $\divop u=0$, then $\det Y$ is
constant in $t$ and equal to $\det Y(s)=1$.
\end{proof}

\subsection{Two Estimates on the Plane}

Lemma~\ref{lem:grad-bound} is a pointwise
bound for the velocity gradient, which Lemma~\ref{lem:iter-lip} uses to bound
the gradients of the iterates for all times. Lemma~\ref{lem:comm-R2} is a
commutator estimate, which the energy estimate of Lemma~\ref{lem:iter-Hs} uses.

\begin{lemma}
\label{lem:grad-bound}
Let $g\in C^1_c(\R^2)$ be real with $g\not\equiv0$, and let $u=\KBS\ast g$. Then
for all $x\in\R^2$ and $i,l\in\{1,2\}$,
\[
|\partial_lu_i(x)|\le\|g\|_{L^\infty}\Big(2+\log_+\frac{\|\nabla g\|_{L^\infty}}{\|g\|_{L^\infty}}\Big)
+\frac{\|g\|_{L^1}}{2\pi},
\]
where $\|\nabla g\|_{L^\infty}=\sup_x|\nabla g(x)|$.
\end{lemma}

\begin{proof}
Since $\KBS$ is locally integrable and $g\in C^1_c$, differentiation under the
integral sign gives $\partial_lu_i(x)=\int(\KBS)_i(z)\,\partial_lg(x-z)\dd z$.
Fix $\delta\in(0,1]$, and split the integral into the near field
$\{|z|<\delta\}$ and the far field $\{|z|>\delta\}$. By
Lemma~I.\ref{I-lem:BS-basic}(a), $\int_{B_\delta}|\KBS|\dd z=\delta$, so the part of the integral
over $\{|z|<\delta\}$ is at most $\delta\|\nabla g\|_{L^\infty}$. On
$\{|z|>\delta\}$ we write $\partial_lg(x-z)=-\partial_{z_l}[g(x-z)]$ and apply
the divergence theorem on $\{\delta<|z|<\rho'\}$, with $\rho'$ so large that
$g(x-z)=0$ for $|z|\ge\rho'$. The outward normal on $\{|z|=\delta\}$ is
$-z/|z|$, so
\begin{align*}
\int_{|z|>\delta}(\KBS)_i(z)\,\partial_lg(x-z)\dd z
&=\int_{|z|>\delta}\partial_l(\KBS)_i(z)\,g(x-z)\dd z\\
&\quad+\int_{|z|=\delta}(\KBS)_i(z)\,g(x-z)\frac{z_l}{|z|}\dd S .
\end{align*}
By Lemma~I.\ref{I-lem:BS-basic}(a), $|\KBS(z)|=(2\pi|z|)^{-1}$ and
$|\partial_l(\KBS)_i(z)|\le(2\pi|z|^2)^{-1}$. Hence the boundary term, an integral
over a circle of length $2\pi\delta$, is at most
$(2\pi\delta)^{-1}\|g\|_{L^\infty}2\pi\delta=\|g\|_{L^\infty}$; integration in
polar coordinates shows that the volume integral
over $\{\delta<|z|<1\}$ is at most
$\|g\|_{L^\infty}\int_\delta^1r^{-1}\dd r=\|g\|_{L^\infty}\log\frac1\delta$; and the
volume integral over $\{|z|\ge1\}$, where $|\partial_l(\KBS)_i|\le(2\pi)^{-1}$, is at most $\|g\|_{L^1}/(2\pi)$. Hence
$|\partial_lu_i(x)|\le\delta\|\nabla g\|_{L^\infty}+\|g\|_{L^\infty}\big(1+\log\frac1\delta\big)+\|g\|_{L^1}/(2\pi)$.
The choice $\delta=\min\{1,\|g\|_{L^\infty}/\|\nabla g\|_{L^\infty}\}$ gives the
claim: then $\delta\|\nabla g\|_{L^\infty}\le\|g\|_{L^\infty}$ and
$\log\frac1\delta=\log_+(\|\nabla g\|_{L^\infty}/\|g\|_{L^\infty})$.
\end{proof}

\begin{lemma}
\label{lem:comm-R2}
Let $s\ge2$ be an integer, $g\in C_c^\infty(\R^2)$, and let
$u\in C^\infty(\R^2;\R^2)$ satisfy
$|\nabla u|_{\max}:=\max_{i,l}\sup_x|\partial_lu_i(x)|<\infty$ and
$\partial^\beta u\in L^2$ for $1\le|\beta|\le s$. Then for every multi-index
$\alpha$ with $|\alpha|=s$,
\[
\|[\partial^\alpha,u\cdot\nabla]g\|_{L^2}\le2^{s+1}\Big[\big(|\nabla u|_{\max}+\|\nabla u\|_{H^{s-1}}\big)\|g\|_{H^s}
+\|\nabla u\|_{H^{s-1}}\|\nabla g\|_{L^\infty}\Big].
\]
\end{lemma}

\begin{proof}
By the Leibniz formula,
$[\partial^\alpha,u\cdot\nabla]g=\sum_{l=1}^2\sum_{0<\beta\le\alpha}\binom\alpha\beta\partial^\beta u_l\,\partial_l\partial^{\alpha-\beta}g$,
and $\sum_l\sum_{\beta\le\alpha}\binom\alpha\beta=2^{s+1}$. It suffices to bound
the $L^2$ norm of each product by the bracket on the right. We distinguish three
cases according to the order $|\beta|$ of the derivative that falls on $u$. If $|\beta|=1$, the product is at
most $|\nabla u|_{\max}\|g\|_{H^s}$, since
$\partial_l\partial^{\alpha-\beta}g$ is a derivative of order $s$ of $g$. If $\beta=\alpha$, it is at most
$\|\partial^\alpha u_l\|_{L^2}\|\partial_lg\|_{L^\infty}\le\|\nabla u\|_{H^{s-1}}\|\nabla g\|_{L^\infty}$,
since $\partial^\alpha u_l$ is a derivative of order $s-1$ of an entry of
$\nabla u$. If $2\le|\beta|\le s-1$ (a case that occurs only for $s\ge3$),
H\"older's inequality, which bounds the $L^2$ norm of the product by the
product of the $L^4$ norms of the factors, and (E1) of
Lemma~\ref{lem:embed} give the bound
$\|\partial^\beta u_l\|_{H^1}\|\partial_l\partial^{\alpha-\beta}g\|_{H^1}
\le\|\nabla u\|_{H^{|\beta|}}\|g\|_{H^{s-|\beta|+2}}\le\|\nabla u\|_{H^{s-1}}\|g\|_{H^s}$,
because $|\beta|\le s-1$ and $|\beta|\ge2$.
\end{proof}

\subsection{The Iteration}

Throughout the rest of this appendix, $\omega_0$ and $R_0$ are as in
Theorem~\ref{thm:smooth}, and we assume $\omega_0\not\equiv0$. We put
\begin{gather*}
A=\|\omega_0\|_{L^\infty},\qquad
N_0=\sup_x|\nabla\omega_0(x)|,\\
U_*=\Big(\frac{2A\|\omega_0\|_{L^1}}{\pi}\Big)^{1/2},\qquad
c_*=2+\frac{\|\omega_0\|_{L^1}}{2\pi A},
\end{gather*}
and, for $t\in\R$,
\begin{gather*}
\mathcal L(t)=\Big(\log_+\frac{N_0}{A}+c_*\Big)e^{2A|t|}-c_*,\qquad
\Psi(t)=Ae^{\mathcal L(t)},\\
V(t)=2A\big(\mathcal L(t)+c_*\big),\qquad R(t)=R_0+U_*|t| .
\end{gather*}
The functions $\mathcal L$, $\Psi$ and $V$ are even and nondecreasing in
$|t|$; moreover $\mathcal L\ge\mathcal L(0)=\log_+(N_0/A)\ge0$, and
$\mathcal L'=V$ on $(0,\infty)$. These are the bounds that the iterates satisfy
uniformly in the iteration index: $U_*$ bounds their velocities
(Lemma~\ref{lem:iterates}(iii)), $\Psi(t)$ the gradients of their vorticities
and $V(t)$ the operator norms of their velocity gradients
(Lemma~\ref{lem:iter-lip}), and the disk $B_{R(t)}$ contains their supports. The
identity $\mathcal L'=V$ is what closes the induction: by
Lemma~\ref{lem:grad-bound}, a vorticity with the $L^1$ and $L^\infty$ norms of
$\omega_0$ whose gradient is at most $\Psi(t)$ has a velocity gradient of operator norm at most
$V(t)$, and Gronwall's inequality along the flow of that velocity returns the
bound $\Psi(t)$ for the next iterate. Since $\supp\omega_0$ is a compact subset of
$B_{R_0}$, there is
$R_0'<R_0$ with $\omega_0=0$ outside $\overline B_{R_0'}$. We use $R_0'$ so
that the supports lie in the open disk $B_{R(t)}$ required in
Theorem~\ref{thm:smooth}(1).

\begin{definition}
\label{def:iterates}
Put $\omega^{[0]}(t,x)=\omega_0(x)$ for $(t,x)\in\R\times\R^2$. If
$\omega^{[j]}\in C^\infty(\R\times\R^2)$ is defined and
$v^{[j]}=\KBS\ast\omega^{[j]}$ satisfies the hypotheses of Lemma~\ref{lem:flow},
let $\flow^{[j]}=\flow_{v^{[j]}}$ and put
$\omega^{[j+1]}(t,x)=\omega_0\big(\flow^{[j]}(0,t,x)\big)$.
\end{definition}

Each iterate is transported by the preceding velocity field. The next lemma
shows that Lemma~\ref{lem:flow} applies at every step, so the construction is
well-defined.

\begin{lemma}
\label{lem:iterates}
The iterates $\omega^{[j]}$ are defined for every $j\ge0$, and:
\begin{enumerate}[label=\textup{(\roman*)}]
\item $\omega^{[j]}\in C^\infty(\R\times\R^2)$ is real, and $\omega^{[j]}(t,x)=0$
  for $|x|>R_0'+U_*|t|$;
\item $\|\omega^{[j]}(t)\|_{L^p}=\|\omega_0\|_{L^p}$ for every $t\in\R$ and every
  $p\in[1,\infty]$;
\item $v^{[j]}\in C^\infty(\R\times\R^2;\R^2)$, every derivative of $v^{[j]}$ is
  bounded on $[-T,T]\times\R^2$ for every $T>0$, $\divop v^{[j]}=0$, and
  $|v^{[j]}|\le U_*$;
\item $\flow^{[j]}$ is defined and $C^\infty$ on $\R\times\R\times\R^2$,
  $\det D_y\flow^{[j]}\equiv1$, $\omega^{[j+1]}(0)=\omega_0$, and
  $\partial_t\omega^{[j+1]}+v^{[j]}\cdot\nabla\omega^{[j+1]}=0$ on $\R\times\R^2$.
\end{enumerate}
\end{lemma}

\begin{proof}
We argue by induction on $j$. Items (i) and (ii) hold for $j=0$, since
$\omega^{[0]}(t)=\omega_0$ for every $t$. Assume (i) and
(ii) for some $j$. From them we derive (iii) and (iv) for $j$, and then (i) and
(ii) for $j+1$.

We prove (iii) for $j$. Fix $T>0$. By (i), and since $R_0'<R_0$, every derivative
$\partial_t^i\partial_x^\gamma\omega^{[j]}$ is continuous and vanishes outside
the compact set $[-T,T]\times\overline B_{R(T)}$ on $[-T,T]\times\R^2$, so it is
bounded there. When $(t,x)$ ranges over a bounded subset of
$(-T,T)\times\R^2$, the functions
$z\mapsto\KBS(z)\,\partial_t^i\partial_x^\gamma\omega^{[j]}(t,x-z)$ are dominated
by a constant multiple of the integrable function
$|\KBS(z)|\mathbf 1_{\{|z|\le\rho\}}$, for a suitable $\rho$: one so large
that $|x-z|>R(T)$ for all $|z|>\rho$ and all $x$ in that subset. Differentiation under the integral sign
and dominated convergence give
$\partial_t^i\partial_x^\gamma v^{[j]}=\KBS\ast\partial_t^i\partial_x^\gamma\omega^{[j]}$,
continuous, and bounded on $[-T,T]\times\R^2$ by Lemma~I.\ref{I-lem:BS-basic}(b).
Lemma~I.\ref{I-lem:BS-basic}(d) gives $\divop v^{[j]}=0$, and
Lemma~I.\ref{I-lem:BS-basic}(b) with (ii) gives
\[
|v^{[j]}(t,x)|\le\big(2\|\omega^{[j]}(t)\|_{L^1}\|\omega^{[j]}(t)\|_{L^\infty}/\pi\big)^{1/2}=U_* .
\]

We prove (iv) for $j$. By (iii), Lemma~\ref{lem:flow} applies to $v^{[j]}$ and
gives the properties of $\flow^{[j]}$. Hence $\omega^{[j+1]}$ is defined, real
and $C^\infty$, and $\omega^{[j+1]}(0,x)=\omega_0(x)$. By Lemma~\ref{lem:flow}(b),
$\omega^{[j+1]}(t,\flow^{[j]}(t,0,a))=\omega_0(a)$ for all $t$ and $a$.
Differentiating in $t$, with the chain rule and
$\partial_t\flow^{[j]}(t,0,a)=v^{[j]}(t,\flow^{[j]}(t,0,a))$, gives
$(\partial_t\omega^{[j+1]}+v^{[j]}\cdot\nabla\omega^{[j+1]})(t,\flow^{[j]}(t,0,a))=0$,
and $\flow^{[j]}(t,0,\cdot)$ is onto, so the transport equation holds at every
point of $\R\times\R^2$.

We prove (i) and (ii) for $j+1$. If $\omega^{[j+1]}(t,x)\ne0$, then
$a=\flow^{[j]}(0,t,x)$ satisfies $|a|\le R_0'$, and by Lemma~\ref{lem:flow}(c),
$|x-a|\le U_*|t|$; so $|x|\le R_0'+U_*|t|$. The substitution
$x=\flow^{[j]}(t,0,a)$, whose Jacobian determinant is $1$ by (iv), gives
$\int|\omega^{[j+1]}(t,x)|^p\dd x=\int|\omega_0(a)|^p\dd a$ for $p<\infty$, and the
case $p=\infty$ follows because $\flow^{[j]}(0,t,\cdot)$ is a bijection.
\end{proof}

\begin{lemma}
\label{lem:iter-lip}
For every $j\ge0$, $t\in\R$ and $x\in\R^2$,
\[
|\nabla\omega^{[j]}(t,x)|\le\Psi(t),\qquad
\max_{i,l}|\partial_lv^{[j]}_i(t,x)|\le A\big(\mathcal L(t)+c_*\big),\qquad
\|\nabla v^{[j]}(t,x)\|_{\mathrm{op}}\le V(t).
\]
\end{lemma}

\begin{proof}
The third bound follows from the second, because the operator norm of a
$2\times2$ matrix is at most twice the largest modulus of its entries. We show,
by induction on $j$, that the first bound holds for all $t$ and that it implies
the second for the same $j$.

For $j=0$ we have
$|\nabla\omega_0|\le N_0\le\max\{A,N_0\}=Ae^{\mathcal L(0)}\le\Psi(t)$. Suppose
that $|\nabla\omega^{[j]}(\tau,x)|\le\Psi(\tau)$ for all $\tau$ and $x$. Since
$\mathcal L\ge0$, we have $\log_+(\Psi(\tau)/A)=\mathcal L(\tau)$. By
Lemma~\ref{lem:iterates}(ii), $\|\omega^{[j]}(\tau)\|_{L^\infty}=A$ and
$\|\omega^{[j]}(\tau)\|_{L^1}=\|\omega_0\|_{L^1}$. So $\omega^{[j]}(\tau)$ is real,
lies in $C^1_c(\R^2)$ by Lemma~\ref{lem:iterates}(i), and is not identically
zero, since $A>0$; and, as $\log_+$ is nondecreasing,
$\log_+(\|\nabla\omega^{[j]}(\tau)\|_{L^\infty}/A)\le\mathcal L(\tau)$. Thus Lemma~\ref{lem:grad-bound},
applied to $g=\omega^{[j]}(\tau)$, gives
\[
|\partial_lv^{[j]}_i(\tau,x)|\le A\big(2+\mathcal L(\tau)\big)+\frac{\|\omega_0\|_{L^1}}{2\pi}=A\big(\mathcal L(\tau)+c_*\big).
\]
Hence $\|\nabla v^{[j]}(\tau,x)\|_{\mathrm{op}}\le V(\tau)$.

Next we deduce the first bound for $j+1$. By the chain rule,
$\nabla\omega^{[j+1]}(t,x)=D_y\flow^{[j]}(0,t,x)^{\mathsf T}\,\nabla\omega_0\big(\flow^{[j]}(0,t,x)\big)$.
By Lemma~\ref{lem:flow}(d), $Y(\tau)=D_y\flow^{[j]}(\tau,t,x)$ solves
$\partial_\tau Y=\nabla v^{[j]}(\tau,\flow^{[j]}(\tau,t,x))Y$ with $Y(t)=I$, so
$\|Y(\tau)\|_{\mathrm{op}}\le1+\big|\int_t^\tau V(\sigma)\|Y(\sigma)\|_{\mathrm{op}}\dd\sigma\big|$,
and Gronwall's inequality gives
$\|Y(0)\|_{\mathrm{op}}\le\exp\big|\int_0^tV(\sigma)\dd\sigma\big|=e^{\mathcal L(t)-\mathcal L(0)}$,
because $\mathcal L'=V$ on $(0,\infty)$ and both functions are even. Hence,
as a matrix and its transpose have the same operator norm,
$|\nabla\omega^{[j+1]}(t,x)|\le N_0e^{\mathcal L(t)-\mathcal L(0)}\le Ae^{\mathcal L(t)}=\Psi(t)$,
since $N_0e^{-\mathcal L(0)}=\min\{N_0,A\}$.
\end{proof}

The bounds of Lemma~\ref{lem:iter-lip} control the commutator in an energy
estimate for the derivatives of every order, which bounds the iterates in every
$H^s$, uniformly in $j$.

\begin{lemma}
\label{lem:iter-Hs}
For every integer $s\ge1$ and every $T>0$, there is a constant $C$, depending
only on $s$, $T$, $R_0$, $A$, $N_0$ and the norms $\|\omega_0\|_{L^1}$ and
$\|\omega_0\|_{H^{s'}}$, $0\le s'\le s$, such that
\[
\|\omega^{[j]}(t)\|_{H^s}\le C\qquad(j\ge0,\ |t|\le T).
\]
\end{lemma}

\begin{proof}
We argue by induction on $s$. For $s=1$, Lemma~\ref{lem:iterates}(i)--(ii) and
Lemma~\ref{lem:iter-lip} give $\|\omega^{[j]}(t)\|_{L^2}=\|\omega_0\|_{L^2}$ and
$\|\nabla\omega^{[j]}(t)\|_{L^2}^2\le|B_{R(T)}|\Psi(T)^2$, since $\nabla\omega^{[j]}(t)$
vanishes outside $B_{R(T)}$ and $\Psi(t)\le\Psi(T)$. Adding the two bounds gives
the claim for $s=1$, as $|B_{R(T)}|=\pi R(T)^2$.

Let $s\ge2$ and assume the bound at level $s-1$ for all iterates. Throughout
this step, $C$ has only the dependencies stated in the lemma and is uniform
in $j$ and $|t|\le T$. For $j=0$ the bound follows from the datum. Let $j\ge0$; we
prove the bound for $\omega^{[j+1]}$. Let
$|\alpha|=s$. Differentiating the equation of Lemma~\ref{lem:iterates}(iv) gives
\[
\partial_t\partial^\alpha\omega^{[j+1]}+v^{[j]}\cdot\nabla\partial^\alpha\omega^{[j+1]}
=-[\partial^\alpha,v^{[j]}\cdot\nabla]\omega^{[j+1]} .
\]
We multiply by $\partial^\alpha\omega^{[j+1]}$ and integrate; all functions are
smooth and real, and $\omega^{[j+1]}(t)$ vanishes outside $B_{R(T)}$ for
$|t|\le T$. Since $\divop v^{[j]}=0$, the transport term gives, after an
integration by parts in $x$,
$-\frac12\int(\divop v^{[j]})(\partial^\alpha\omega^{[j+1]})^2\dd x=0$. Put
$\mathcal H(t)=\sum_{|\alpha|=s}\|\partial^\alpha\omega^{[j+1]}(t)\|_{L^2}^2$ and
$Z(t)=\|\omega^{[j+1]}(t)\|_{H^s}$; there are $s+1$ multi-indices with
$|\alpha|=s$. We apply Lemma~\ref{lem:comm-R2} with
$u=v^{[j]}(t)$ and $g=\omega^{[j+1]}(t)$, using $|\nabla v^{[j]}|_{\max}\le V(T)$
and $|\nabla\omega^{[j+1]}|\le\Psi(T)$ (Lemma~\ref{lem:iter-lip}), and
\[
\|\nabla v^{[j]}(t)\|_{H^{s-1}}=\|\omega^{[j]}(t)\|_{H^{s-1}}\le C
\]
(Lemma~I.\ref{I-lem:BS-basic}(d) and the induction hypothesis). By the
Cauchy-Schwarz inequality, the absolute value of the time derivative of
$\|\partial^\alpha\omega^{[j+1]}(t)\|_{L^2}^2$ is at most $2Z$ times the
$L^2$ norm of the commutator; summing over $\alpha$ gives, for
$|t|\le T$,
\begin{align*}
|\mathcal H'(t)|&\le(s+1)2^{s+2}\big[(V(T)+C)Z+C\Psi(T)\big]Z\\
&\le C(Z^2+1)\le C(\mathcal H+1),
\end{align*}
because $Z\le Z^2+1$ and $Z^2\le\mathcal H+\|\omega^{[j+1]}(t)\|_{H^{s-1}}^2\le\mathcal H+C$.
Since $\mathcal H(0)\le\|\omega_0\|_{H^s}^2$, as $\omega^{[j+1]}(0)=\omega_0$,
Gronwall's inequality, applied to $\mathcal H+1$, gives
$\mathcal H(t)\le(\|\omega_0\|_{H^s}^2+1)e^{C|t|}$. Hence $Z^2\le\mathcal H+C$
is bounded uniformly in $j$ and $|t|\le T$, proving the induction step.
\end{proof}

\subsection{Convergence}

An energy estimate for differences of consecutive iterates gives convergence
in $L^2$, uniformly on bounded time intervals (Lemma~\ref{lem:iter-cauchy}).
The higher Sobolev bounds of Lemma~\ref{lem:iter-Hs} and interpolation then give
convergence of all spatial derivatives.

\begin{lemma}
\label{lem:iter-cauchy}
Let $T>0$ and put $e_j(t)=\|\omega^{[j]}(t)-\omega^{[j-1]}(t)\|_{L^2}$ for $j\ge1$.
Then for every $j\ge1$ and every $t\in[-T,T]$,
\[
e_{j}(t)\le2\|\omega_0\|_{L^2}\frac{\big(2R(T)\Psi(T)|t|\big)^{j-1}}{(j-1)!} .
\]
\end{lemma}

\begin{proof}
Put $R=R(T)$ and $\Psi=\Psi(T)$. For $j=1$, $e_1(t)\le2\|\omega_0\|_{L^2}$ by
Lemma~\ref{lem:iterates}(ii) and the triangle inequality. Let $j\ge1$. Subtracting the equations of
Lemma~\ref{lem:iterates}(iv) for $\omega^{[j+1]}$ and $\omega^{[j]}$ shows that
$\delta=\omega^{[j+1]}-\omega^{[j]}$ satisfies
\[
\partial_t\delta+v^{[j]}\cdot\nabla\delta=-\big(\KBS\ast(\omega^{[j]}-\omega^{[j-1]})\big)\cdot\nabla\omega^{[j]} .
\]
We multiply by $\delta$ and integrate; the transport term vanishes, because
$\divop v^{[j]}=0$ and $\delta(t)$ has compact support. For $|t|\le T$, $\nabla\omega^{[j]}(t)$ vanishes outside $B_R$ and is
bounded by $\Psi$, by Lemma~\ref{lem:iterates}(i) and Lemma~\ref{lem:iter-lip};
and $\omega^{[j]}(t)-\omega^{[j-1]}(t)$ vanishes outside $B_R$, so
Lemma~I.\ref{I-lem:BS-basic}(c), with $\rho=R$, gives
$\|\KBS\ast(\omega^{[j]}-\omega^{[j-1]})(t)\|_{L^2(B_R)}\le2Re_j(t)$. Hence, by the
Cauchy-Schwarz inequality,
$|\frac{\dd}{\dd t}e_{j+1}^2|\le4R\Psi e_je_{j+1}$. To divide by $e_{j+1}$, which
may vanish, we regularize: for all $\epsilon>0$,
$\big|\frac{\dd}{\dd t}(e_{j+1}^2+\epsilon)^{1/2}\big|\le2R\Psi e_j$, because
$e_{j+1}\le(e_{j+1}^2+\epsilon)^{1/2}$. Since
$e_{j+1}(0)=0$, letting $\epsilon\to0$ gives
$e_{j+1}(t)\le2R\Psi\big|\int_0^te_j(\tau)\dd\tau\big|$ for $|t|\le T$, and the
claim follows by induction on $j$: inserting the bound for $e_j$ and integrating
gives the bound for $e_{j+1}$.
\end{proof}

\begin{proposition}
\label{prop:iter-limit}
There is a real $\omega\in C^\infty(\R\times\R^2)$ such that, with
$v=\KBS\ast\omega$, for every $T>0$ and every multi-index $\gamma$,
$\partial_x^\gamma\omega^{[j]}\to\partial_x^\gamma\omega$ and
$\partial_x^\gamma v^{[j]}\to\partial_x^\gamma v$ uniformly on $[-T,T]\times\R^2$
as $j\to\infty$. The function $\omega$ has properties \textup{(1)}--\textup{(3)}
of Theorem~\ref{thm:smooth}, and for all $t\in\R$, $x\in\R^2$ and
$p\in[1,\infty]$,
\begin{equation}
\label{eq:limit-bounds}
\begin{gathered}
\omega(t,x)=0\ \text{ for }|x|>R_0'+U_*|t|,\qquad
\|\omega(t)\|_{L^p}=\|\omega_0\|_{L^p},\qquad |\nabla\omega(t,x)|\le\Psi(t),\\
|v(t,x)|\le U_*,\qquad \max_{i,l}|\partial_lv_i(t,x)|\le A\big(\mathcal L(t)+c_*\big).
\end{gathered}
\end{equation}
\end{proposition}

\begin{proof}
We first prove the convergence. Fix $T>0$. By Lemma~\ref{lem:iter-cauchy} the
series $\sum_j\sup_{|t|\le T}e_j(t)$ is at most
$2\|\omega_0\|_{L^2}\exp(2R(T)\Psi(T)T)$. Thus $(\omega^{[j]})$ is a Cauchy sequence
for the norm $\sup_{|t|\le T}\|\cdot\|_{L^2}$, since $\omega^{[j]}-\omega^{[i]}$ is a
sum of differences of consecutive iterates. For an integer
$k\ge1$, (E3) of Lemma~\ref{lem:embed}, and Lemma~\ref{lem:iter-Hs} and the triangle
inequality, which bound the $H^{2k}$ norm of $\omega^{[j]}(t)-\omega^{[i]}(t)$
uniformly in $i,j$ and $|t|\le T$, give
\[
\sup_{|t|\le T}\|\omega^{[j]}(t)-\omega^{[i]}(t)\|_{H^k}^2\le C\sup_{|t|\le T}\|\omega^{[j]}(t)-\omega^{[i]}(t)\|_{L^2},
\]
where $C$ has the dependencies of Lemma~\ref{lem:iter-Hs} with $s=2k$.
The right side tends to $0$ as $i,j\to\infty$. By (E2), applied to
$\partial_x^\gamma(\omega^{[j]}-\omega^{[i]})(t)$ and combined with this bound for
$k=|\gamma|+2$, for all $\gamma$ the functions
$\partial_x^\gamma\omega^{[j]}$, which are continuous on $\R\times\R^2$, form a
uniform Cauchy sequence on $[-T,T]\times\R^2$. Their uniform limits are
continuous. The limit of $\partial_x^\gamma\omega^{[j]}$ is $\partial_x^\gamma$ of
the limit $\omega$ of $\omega^{[j]}$, as one sees by passing to the limit in
$g(x+he_l)-g(x)=\int_0^h\partial_lg(x+\sigma e_l)\dd\sigma$, with
$g=\partial_x^\gamma\omega^{[j]}(t)$. The limits for
different $T$ agree. The first line of Equation~\eqref{eq:limit-bounds} follows
from Lemma~\ref{lem:iterates}(i)--(ii) and Lemma~\ref{lem:iter-lip}, because the
convergence is uniform and the supports lie in a fixed ball for $|t|\le T$. In
particular $\omega(t)\in C_c^\infty(\R^2)$ for all $t$, and $\omega$ has
property (1) of Theorem~\ref{thm:smooth}, since $R_0'+U_*|t|<R(t)$.

We now treat the velocity. For a bounded measurable $g$ vanishing outside
$B_\rho$, Lemma~I.\ref{I-lem:BS-basic}(b) and $\|g\|_{L^1}\le\pi\rho^2\|g\|_{L^\infty}$
give $\|\KBS\ast g\|_{L^\infty}\le\sqrt2\,\rho\|g\|_{L^\infty}$. Applied to
$g=\partial_x^\gamma\omega^{[j]}(t)-\partial_x^\gamma\omega(t)$ with $\rho=R(T)$,
this shows that $\partial_x^\gamma v^{[j]}=\KBS\ast\partial_x^\gamma\omega^{[j]}$
converges to $\KBS\ast\partial_x^\gamma\omega$ uniformly on $[-T,T]\times\R^2$.
By Lemma~I.\ref{I-lem:BS-basic}(d) applied to $\omega(t)$, $\KBS\ast\partial_x^\gamma\omega(t)=\partial_x^\gamma v(t)$
and $\divop v=0$. The bounds on $v$ in Equation~\eqref{eq:limit-bounds} follow
from Lemma~\ref{lem:iterates}(iii) and Lemma~\ref{lem:iter-lip}, by passing to the
limit in the corresponding bounds for $v^{[j]}$.

We now derive the equation. For every $\gamma$, Lemma~\ref{lem:iterates}(iv),
differentiated in $x$ and integrated in time from $0$, using
$\omega^{[j+1]}(0)=\omega_0$,
gives
\[
\partial_x^\gamma\omega^{[j+1]}(t,x)=\partial^\gamma\omega_0(x)
-\int_0^t\partial_x^\gamma\big(v^{[j]}\cdot\nabla\omega^{[j+1]}\big)(\tau,x)\dd\tau .
\]
By the Leibniz formula the integrand is a finite sum of products of
$x$-derivatives of $v^{[j]}$ and of $\omega^{[j+1]}$, each converging uniformly
and bounded on $[-T,T]\times\R^2$. Passing to the limit $j\to\infty$,
\[
\partial_x^\gamma\omega(t,x)=\partial^\gamma\omega_0(x)-\int_0^t\partial_x^\gamma(v\cdot\nabla\omega)(\tau,x)\dd\tau,
\]
with an integrand that is continuous and bounded on $[-T,T]\times\R^2$. Hence,
by the fundamental theorem of calculus,
$\partial_t\partial_x^\gamma\omega$ exists, is continuous and bounded, and
equals $-\partial_x^\gamma(v\cdot\nabla\omega)$. With $\gamma=0$ this is
property (3).

Finally we prove that $\omega$ and $v$ are $C^\infty$ with bounded derivatives.
So far we control the derivatives taken first in $x$ and then once in $t$. We
obtain the higher time derivatives one at a time from the identity
$\partial_t\partial_x^\gamma\omega=-\partial_x^\gamma(v\cdot\nabla\omega)$, and
then use Schwarz's theorem to show that the order of differentiation does not
matter. To keep track of the order, we use the following classes.
For $T>0$ and an integer $l\ge0$, let $\mathcal C^l(T)$ be the set of functions
$g$ on $(-T,T)\times\R^2$ such that $\partial_t^i\partial_x^\gamma g$ exists, is
continuous and is bounded, for all $i\le l$ and all multi-indices $\gamma$
(including $i=0$ and $\gamma=0$); here the
$x$-derivatives are taken first. We have shown $\omega\in\mathcal C^1(T)$. Suppose
$\omega\in\mathcal C^l(T)$. Every $\partial_t^i\partial_x^\gamma\omega$ with
$i\le l$ vanishes for $|x|\ge R(T)$, so differentiation under the integral sign,
as in the proof of Lemma~\ref{lem:iterates}(iii), gives
$\partial_t^i\partial_x^\gamma v=\KBS\ast\partial_t^i\partial_x^\gamma\omega$,
continuous and bounded by $\sqrt2R(T)\sup|\partial_t^i\partial_x^\gamma\omega|$;
so $v\in\mathcal C^l(T)$. By the Leibniz formula $v\cdot\nabla\omega\in\mathcal C^l(T)$,
and the identity $\partial_t\partial_x^\gamma\omega=-\partial_x^\gamma(v\cdot\nabla\omega)$
shows that $\partial_t^{i+1}\partial_x^\gamma\omega$ exists, is continuous and
is bounded for $i\le l$. Hence $\omega\in\mathcal C^{l+1}(T)$. By induction,
$\omega$ and $v$ lie in $\mathcal C^l(T)$ for all $l$ and $T$. By Schwarz's
theorem on mixed partial derivatives (if $\partial_ag$, $\partial_bg$ and
$\partial_b\partial_ag$ exist near a point and $\partial_b\partial_ag$ is
continuous at it, then $\partial_a\partial_bg$ exists at that point and equals
$\partial_b\partial_ag$ there), an induction on the order shows that every
partial derivative of $\omega$ or $v$, taken in any order, exists, is continuous,
and equals the corresponding $\partial_t^i\partial_x^\gamma$ derivative. Thus
$\omega,v\in C^\infty(\R\times\R^2)$, $v$ and every derivative of $v$ are bounded
on $[-T,T]\times\R^2$ for all $T$, and $\omega$ has property (2).
\end{proof}

\subsection{Proof of Theorem~\ref{thm:smooth}}

The proof of Theorem~\ref{thm:smooth} also gives the following explicit bounds,
in terms of the constants defined at the start of the iteration above.

\begin{lemma}
\label{lem:smooth-bounds}
In Theorem~\ref{thm:smooth}, assume $\omega_0\not\equiv0$ and put
$A=\|\omega_0\|_{L^\infty}$, $N_0=\sup_x|\nabla\omega_0(x)|$,
$U_*=\sqrt{2A\|\omega_0\|_{L^1}/\pi}$, $c_*=2+\|\omega_0\|_{L^1}/(2\pi A)$,
$\mathcal L(t)=(\log_+(N_0/A)+c_*)e^{2A|t|}-c_*$ and $\Psi(t)=Ae^{\mathcal L(t)}$.
Then for every $t\in\R$ and $x\in\R^2$,
\begin{gather*}
|\nabla\omega(t,x)|\le\Psi(t),\qquad|v(t,x)|\le U_*,\\
\max_{i,l}|\partial_lv_i(t,x)|\le A\big(\mathcal L(t)+c_*\big),\qquad
\|\nabla v(t,x)\|_{\mathrm{op}}\le2A\big(\mathcal L(t)+c_*\big).
\end{gather*}
\end{lemma}

\begin{proof}[Proof of Theorem~\ref{thm:smooth} and Lemma~\ref{lem:smooth-bounds}]
We first prove (5), for any $\omega$ with properties (1)--(3); it gives the
uniqueness in Theorem~\ref{thm:smooth}, and the proof of (6) uses it. Let $\tilde\omega$
be as in (5), and fix a compact interval $J\subset I$ containing $0$. Let
$\rho>0$ be such that $\omega(t)$ and $\tilde\omega(t)$ vanish outside $B_\rho$ for
$t\in J$; this is possible by (1) and the hypothesis on $\tilde\omega$. Put
$\delta=\omega-\tilde\omega\in C^1(J\times\R^2)$. Subtracting the two
vorticity equations and using $v-\KBS\ast\tilde\omega=\KBS\ast\delta$ gives
\[
\partial_t\delta+v\cdot\nabla\delta=-(\KBS\ast\delta)\cdot\nabla\tilde\omega
\qquad\text{on }J\times\R^2 .
\]
We multiply by $\delta$ and integrate in $x$. The transport term vanishes,
because $\divop v=0$ and $\delta(t)$ has compact support. The function
$\nabla\tilde\omega$ is continuous on $J\times\R^2$ and vanishes there for
$|x|\ge\rho$, so it is bounded by some $L<\infty$, and
Lemma~I.\ref{I-lem:BS-basic}(c) gives $\|\KBS\ast\delta(t)\|_{L^2(B_\rho)}\le2\rho\|\delta(t)\|_{L^2}$.
Hence, by the Cauchy-Schwarz inequality, as $\nabla\tilde\omega(t)$ vanishes
outside $B_\rho$,
$|\frac{\dd}{\dd t}\|\delta(t)\|_{L^2}^2|\le4\rho L\|\delta(t)\|_{L^2}^2$ on
$J$. Since $\delta(0)=0$, Gronwall's inequality gives $\delta=0$ on $J$. Every
point of $I$ lies in such an interval $J$, which proves (5). In particular at
most one $\omega\in C^\infty(\R\times\R^2)$ has properties (1)--(3), since such
an $\omega$ satisfies the hypotheses on $\tilde\omega$ in (5) with $I=\R$.

If $\omega_0\equiv0$, the function $\omega\equiv0$ has properties (1)--(4), (6)
and (7). From now on $\omega_0\not\equiv0$, and $\omega$ is the function of
Proposition~\ref{prop:iter-limit}, which has properties (1)--(3). The bounds of
Lemma~\ref{lem:smooth-bounds} are contained in Equation~\eqref{eq:limit-bounds},
together with $\|M\|_{\mathrm{op}}\le2\max_{i,l}|M_{il}|$ for $2\times2$ matrices
$M$.

We prove (4). By (2), Lemma~\ref{lem:flow} applies to $v$; let
$\flow(t,y)=\flow_v(t,0,y)$ and $\flow_t=\flow(t,\cdot)$. By Lemma~\ref{lem:flow},
$\flow$ is defined for all $t$; by (a) and (b) there, $\flow$ and
$(t,x)\mapsto\flow_t^{-1}(x)=\flow_v(0,t,x)$ lie in $C^\infty(\R\times\R^2)$ and
each $\flow_t$ is a $C^\infty$ diffeomorphism; and by (d),
$\det D\flow_t\equiv1$ because $\divop v=0$. By (3),
$\frac{\dd}{\dd t}\omega(t,\flow(t,a))=(\partial_t\omega+v\cdot\nabla\omega)(t,\flow(t,a))=0$,
so $\omega(t)=\omega_0\circ\flow_t^{-1}$. For a continuous $F$ with $F(0)=0$,
the function $F(\omega(t))$ is continuous with compact support, and the change
of variables $x=\flow_t(a)$, whose Jacobian determinant is $1$, gives
$\int F(\omega(t))\dd x=\int F(\omega_0)\dd a$, since
$\omega(t)\circ\flow_t=\omega_0$.

We prove (6). Let $\Orot_\theta\in SO(2)$. Since $(\Orot_\theta z)^\perp=\Orot_\theta z^\perp$,
we have $\KBS(\Orot_\theta z)=\Orot_\theta\KBS(z)$. Put
$\omega_\theta(t,x)=\omega(t,\Orot_\theta^{-1}x)$. The substitution
$y=\Orot_\theta y'$ gives
$\KBS\ast\omega_\theta(t,x)=\Orot_\theta v(t,\Orot_\theta^{-1}x)$, and
$\nabla\omega_\theta(t,x)=\Orot_\theta\nabla\omega(t,\Orot_\theta^{-1}x)$. Hence
$(\KBS\ast\omega_\theta)\cdot\nabla\omega_\theta=(v\cdot\nabla\omega)(t,\Orot_\theta^{-1}x)$,
and $\omega_\theta$ is a $C^\infty$ solution of the equation in (3) with datum
$\omega_0\circ\Orot_\theta^{-1}$, vanishing outside $B_{R(t)}$ at time $t$. If
$\omega_0\circ\Orot_\theta=\omega_0$, then $\omega_0\circ\Orot_\theta^{-1}=\omega_0$,
and (5) with $I=\R$ gives $\omega_\theta=\omega$, that is,
$\omega(t)\circ\Orot_\theta^{-1}=\omega(t)$ for all $t$, which is equivalent to
$\omega(t)\circ\Orot_\theta=\omega(t)$.

We prove (7): we verify the conditions of
Definition~I.\ref{I-def:weak}, and the membership in $L^\infty_t(L^1\cap L^\infty)_x$. Let $T>0$, $\varphi\in\DT$, and let $\Ren\in C^1(\R)$ with
$\Ren'$ bounded and $\Ren(0)=0$. The function $\Ren(\omega)$ lies in
$C^1(\R\times\R^2)$, vanishes outside $B_{R(t)}$ at time $t$, and by (3)
satisfies $\partial_t\Ren(\omega)+v\cdot\nabla\Ren(\omega)=\Ren'(\omega)(\partial_t\omega+v\cdot\nabla\omega)=0$.
We multiply by $\varphi$ and integrate over $(0,T)\times\R^2$. Since
$\varphi(t,\cdot)=0$ for $t$ near $T$, an integration by parts in $t$ gives
$\int_0^T\!\int\partial_t\Ren(\omega)\,\varphi=-\int\Ren(\omega_0)\varphi(0,\cdot)-\int_0^T\!\int\Ren(\omega)\,\partial_t\varphi$,
and since $\divop v=0$, an integration by parts in $x$ gives
$\int(v\cdot\nabla\Ren(\omega))\varphi\dd x=-\int\Ren(\omega)\,v\cdot\nabla\varphi\dd x$.
Hence Equation~(I.\ref{I-eq:renorm}) holds, and with $\Ren(s)=s$ so does
Equation~(I.\ref{I-eq:weak}). Finally, $\|\omega(t)\|_{L^1}=\|\omega_0\|_{L^1}$ by
(4) with $F(s)=|s|$, and $\|\omega(t)\|_{L^\infty}=\|\omega_0\|_{L^\infty}$
because $\omega(t)=\omega_0\circ\flow_t^{-1}$. Moreover $\omega$ is continuous and
$\omega(t)$ vanishes outside $B_{R(T)}$ for $|t|\le T$, so $\omega$ is uniformly
continuous on $[0,T]\times\R^2$, and $t\mapsto\omega(t)$ is continuous from
$[0,T]$ into $(L^1\cap L^q)(\R^2)$ for all $q\in[1,\infty]$. Hence the
restriction of $\omega$ to $[0,T)$ lies in $L^\infty_t(L^1\cap L^\infty)_x$, and
in the space $X_q(T)$ of Section~I.\ref{I-sec:prelim} for every
$q\in(4/3,\infty)$. This proves (7) and completes the
proof.
\end{proof}

\section{Quantitative Details of the Relay Construction}
\label{sec:construction-details}

This appendix gives the parameter conditions and estimates used in the
relay step of Section~\ref{sec:relay}. After fixing the data and recording the
parameter choices, we prove the velocity, coefficient and time estimates.
The last subsections verify their hypotheses in the common bootstrap setting.

\subsection{Fixed data for the estimates}
\label{sec:relay-data}

The parameter choices in Section~\ref{sec:param} use the following data.
The estimates below are stated with their individual hypotheses, so their
constants can be chosen in the order required there.
Here $\sigma=\sigma(0)=\Re\lamref$ is the reference growth rate, and
$\sigma_i=\sigma(G_i)=\Re\lambda_i$ is the rate in ring-$i$ units.
Its physical growth rate is $A_i\sigma_i$.

\begin{definition}
\label{def:standing}
We say that the \textbf{standing assumptions (SA)} hold if \textup{(N1)}--\textup{(N3)}
hold and the following objects are fixed:
\begin{enumerate}[label=\textup{(SA\arabic*)},leftmargin=*]
\item an integer $s\ge3$ and $\sigma'\in\big(\frac9{200},\min(\sigb,\frac34\sigma)\big)$;
  the corresponding regularity level $M$, constants and threshold functions
  supplied by Theorem~\ref{thm:active};
\item the integers $s_h=M-1$, $N_Q=\max(2,\lceil\Lam(s_h)/\sigb\rceil)$ and
  $M_h=s_h+N_Q+1$, with $\Lam$ of Lemma~\ref{lem:energy}; the associated
  quiet-ring constants of Theorem~\ref{thm:quiet} and the velocity-estimate
  constant $\Cvel$ at order $M_h$ from Lemma~\ref{lem:transfer};
\item numbers $\delta_0\in(0,\dH]$, $\mu\in(1,\frac{17}{16}]$ and $\Tpin>0$; a
  number $\kseed\in(0,\kL]$, with $\kL$ of Lemma~\ref{lem:seed}; a number
  $\sbar\in(0,a/16]$ with $\mu\sbar^2\le\Gone$
  (Lemma~\ref{lem:inner-circ}); an integer $n\ge2$; and the objects of
  Definitions~\ref{def:datum}--\ref{def:clock} with these parameters.
\end{enumerate}
Under (SA), put
\begin{equation}
\label{eq:relay-tail-weight}
\rT(T)=(1+T)^3e^{-\sigma T/4}\qquad(T\ge0).
\end{equation}
We use the envelopes $\hat a_k$ and $\rhohat_k$ of Section~\ref{sec:relay},
with the constant $C_0\ge1$ named in each estimate when it is allowed to vary.
We say that the
\textbf{detuning hypothesis (Det)} holds if
\begin{equation}
\label{eq:detuning-hypothesis}
\sup_{0\le G\le \mu\sbar^2}|\lamst(G)-\lamref|
\le\frac{\sigma(\mu-1)}{8\mu}.
\end{equation}
\end{definition}

For applications of Section~\ref{sec:single} beginning at physical time $t_*$,
we set $\theta=A_k(t-t_*)$ and pull back $w_k,h_k,\eta_k$ by
$t=t_*+\theta/A_k$. For a quiet ring $t_*=0$, and its Duhamel bound is
$\theta\mapsto\maj_k(\theta/A_k)$. Thus Theorem~\ref{thm:quiet}(c) gives
$\|w_k(t)\|_{H^{s_h}}\le C_Q\maj_k(t)$ whenever its hypotheses hold.
The fields and their norms, and hence their Duhamel bounds, are continuous
wherever defined, by the paragraph after Lemma~\ref{lem:regions-basic}.
Lemma~\ref{lem:inner-circ} and (SA3) give
$0\le G_k\le\mu\sbar^2\le\Gone$ at every level. All the eigendata are
therefore defined, and $\sigma_k\ge\sigb$.

\subsection{The parameter inequalities}
\label{sec:param-lists}

The common setting in Section~\ref{sec:relay} uses the choices (P0)--(P6)
below. They are made in the order of Table~\ref{tab:param}. The inequalities
ensure the hypotheses of the estimates used in the relay.

\begin{definition}
\label{def:P0-P4}
\emergencystretch=3em
The \textbf{structural list} consists of the following ordered choices.
The auxiliary constants are defined at the cited locations.
\begin{description}[font=\normalfont]
\item[\textup{(P0.1)}] $\Phi$, $a$, $\Gone$ of Definition~I.\ref{I-def:ring};
  $\sigma$, $\sigb$ of Definition~I.\ref{I-def:eigendata}; and
  $\lamref=\lamst(0)$ of Section~\ref{sec:kick}.
\item[\textup{(P0.2)}] An integer $s$ and a number $\sigma'$ as in (SA1); then
  $c_\gamma$ of Appendix~\ref{sec:relay-constants}.
\item[\textup{(P0.3)}] The integers $N_G$, $M$, the constants
  $\dH,C_2,C_3$ and $\sigma_{\max}$, and the functions
  $\kappa_0(\cdot)$, $\bar t_1(\cdot)$, $C(\cdot)$ of Theorem~\ref{thm:active};
  $\kL$ and $\Cseed$ of Lemma~\ref{lem:seed}; $\|\ellst\|$, $\Cls$ and $\Cp$ of
  Lemma~I.\ref{I-lem:uniform}.
\item[\textup{(P0.4)}] The integers $s_h$, $N_Q$, $M_h$ and the constants
  $\majst,C_Q$ of (SA2).
\item[\textup{(P0.5)}] The common velocity constant $\Cvel$ of
  Lemmas~\ref{lem:transfer}, \ref{lem:inner-field} and~\ref{lem:field-bound},
  at order $M_h$; the absolute constant $\Cstar$ of Lemma~\ref{lem:speed}; and
  $\Cmom$ of Proposition~\ref{prop:confine}.
\item[\textup{(P0.6)}] $c_N$, $C_N$, $\cout$ of Lemma~I.\ref{I-lem:uniform};
  the constants of Lemma~\ref{lem:kicktrue} are chosen at (P4).
\item[\textup{(P1)}] Choose $\delta_0\in(0,\dH]$ sufficiently small that
  $\frac{16}{15}\cdot\frac{2\pi C_3\delta_0^2}{\sigma}\le\frac14\cdot\frac{15}{16}\cdot\frac{\cout\delta_0}{20|\lamref|}$
  and $\Cp C_2\delta_0\le\frac{\cout}{16\pi}$.
\item[\textup{(P2)}] $\mu\in(1,\mu_{\max}]$, with
  $\mu_{\max}=1+\min\big(c'\delta_0,\frac1{16}\big)$ and
  a fixed $c'>0$ sufficiently small that
  $c'\le\frac{15}{16}\cdot\frac{\cout}{10|\lamref|\,\Cmom\,\Ip{3}{\sigma/4}}$; then
  $\vartheta_\mu$, $\gamma_\mu$, $\gamma'_\mu$ of Appendix~\ref{sec:relay-constants},
  and $\CDuh$, $\KO$ of Lemma~\ref{lem:O}.
\item[\textup{(P3)}] $\Tpin>0$ with $\rT(\Tpin)\le\frac1{64}$, with the function
  $\rT$ of Definition~\ref{def:standing}. Then $C_D(\Tpin)$
  (Proposition~\ref{prop:relay}); $K_1$, $\Abar$, $\Klo$, $\Kup$
  (Definition~\ref{def:K1});
  $\Tz$ (Lemma~\ref{lem:time}) and
  choose $C_\tau\ge\max\{3,\frac{15}{\sigma'}\}$.
\item[\textup{(P4)}] Choose $C_0$ so that
  \begin{equation}
  \label{eq:param-C0}
  \begin{gathered}
  C_0\ge\max\{1,2\Cseed,K_1C_Q\},\\
  C_0\ge\Cvel \,\mu\,(2+2\Cstar\Tpin)^3(1+2\Cstar)^3/\Klo.
  \end{gathered}
  \end{equation}
  Then fix $\kappa_0(C_0)$, $\bar t_1(C_0)$ and $C(C_0)$ of
  Theorem~\ref{thm:active}, the common constant $\Crelay(C_0,K_1)$ of
  Proposition~\ref{prop:relay} and Lemma~\ref{lem:Oprime}, and $\Cerr$ of
  Lemma~\ref{lem:kicktrue}.
\end{description}
\end{definition}

The constant $C_0$ serves (Init) and (Field) in Theorem~\ref{thm:active}.
The seed estimate bounds the initial perturbation by $2\Cseed\ksh_n$.
Equation~\eqref{eq:relay-startup}
bounds it by $K_1C_Q\ksh_k$ at ring $k$'s activation.
Equation~\eqref{eq:boot-inner} bounds the inner field using the second
inequality of (P4). These give (Init) and the inner-field bound in (Field) in the local coordinate
$u=A_k(t-\tsh_k)$.
The coefficients that $C_0$ must dominate are fixed at (P0)--(P3).
The constants $\Crelay$ and $\Cerr$ are fixed after $C_0$ at (P4), before
choosing the seed and radius ratios.

\begin{definition}
\label{def:P5}
The \textbf{$\kappa$-list (P5)}: choose $\kseed>0$ so small that, for every
$\kappa\in(0,2\kseed]$, the following hold.
\begin{description}[font=\normalfont]
\item[\textup{($\kappa$1)}] $\kappa\le\delta_04^{-\mu/(\mu-1)}$.
\item[\textup{($\kappa$2)}] \begin{multline*}
  \frac{2\pi C_3\kappa^{1/2}\mu}{\frac34\sigma(\mu-1)}+\Cerr\Big(\frac\kappa{\delta_0}\Big)^{\gamma_\mu}\log^4\frac{\delta_0}\kappa\\
  +\Cerr\big[1+\rT(\Tpin)\Ip{3}{\tfrac\sigma8}\big]\Big(\frac\kappa{\delta_0}\Big)^{\vartheta_\mu}\log^3\frac{\delta_0}\kappa
  \le\frac12\cdot\frac{15}{16}\cdot\frac{\cout}{20|\lamref|(\mu-1)}.
  \end{multline*}
\item[\textup{($\kappa$3)}] $\Cerr(\kappa/\delta_0)^{\vartheta_\mu}\le\frac{15}{16}\cdot\frac{\cout}{20|\lamref|(\mu-1)}$.
\item[\textup{($\kappa$4)}] $\kappa\le\kappa_0(C_0)$.
\item[\textup{($\kappa$5)}] For every $G\in[0,\Gone]$, with $\hat T=\sigma(G)^{-1}\log(\delta_0/\kappa)$
  and $\rhohat$ the transient of Theorem~\ref{thm:active} for $(C_0,G)$:
  $\Cp C_2\big(\rhohat(\hat T)+\kappa^{1/2}\big)\le\frac{\cout}{16\pi}$.
\item[\textup{($\kappa$6)}] $\kappa\le\kL$.
\item[\textup{($\kappa$7)}] $\Crelay(C_0,K_1)\big(1+\frac1{\sigb(\mu-1)}\big)(\kappa/\delta_0)^{\gamma'_\mu}\langle\log(\delta_0/\kappa)\rangle^4\le\frac1{\sigb(\mu-1)}$.
\end{description}
\end{definition}

In (P5), $\kappa$ stands for an amplitude $\ksh_k$ with which a ring starts its
active phase and, in ($\kappa$1) and ($\kappa$6), also for $2\kseed$.
Theorem~\ref{thm:active} is applied to ring $k$ with $\kappa=\ksh_k$, and
Lemma~\ref{lem:kicktrue}, at level $k$, with $\kappa=\ksh_{k+1}$.
Conditions ($\kappa$1)--($\kappa$3) are hypotheses (c)--(e) of
Lemma~\ref{lem:kicktrue}, ($\kappa$4) is hypothesis (Small) of
Theorem~\ref{thm:active}, ($\kappa$5) is used in the proof of
Proposition~\ref{prop:reduction}, and ($\kappa$6) makes the seed of
Definition~\ref{def:approx} admissible. The seed estimate
\eqref{eq:boot-start} and the amplitude bounds in Equation~\eqref{eq:relay-startup}
put every starting amplitude
in $(0,2\kseed]$, so the conditions are imposed on that whole interval.

The formula for $K_1$ in Definition~\ref{def:K1} uses quantities fixed at
(P0)--(P3). The $+1$ in $4C_2+1$ allows the remainder $\Xi_{j+1}$ in
Proposition~\ref{prop:relay} to be as large as $\delta_0/(\sigb(\mu-1))$.
After $C_0$ and $\Crelay(C_0,K_1)$ are fixed, ($\kappa$7) gives
$\Xi_{j+1}\le\delta_0/(\sigb(\mu-1))$ for every
$\kappa=\ksh_{j+1}\in(0,2\kseed]$.
The comparison of exponents is checked in
Appendix~\ref{sec:bootstrap-handoff-checks}. The same $K_1$ and $\Abar$
therefore work at every hand-off.

The number $\sbar$, which fixes the ratios $s_k$ of consecutive radii, is chosen
last, because the conditions (B1)--(B11) on it involve constants fixed at all the
earlier steps; each of them holds for every
sufficiently small $\sbar$ (proof of Lemma~\ref{lem:nonempty}).

\begin{definition}
\label{def:P6}
The \textbf{$\sbar$-list (P6)}: choose $\sbar>0$ such that the following hold,
where $s_k=\sbar(2\mu)^{-k}$, $\eta_{b,1}=0$ and $\eta_{b,k}$, $\eta_b$ are as in
Equation~\eqref{eq:etab}.
\begin{description}[font=\normalfont]
\item[\textup{(B1)}] $\sbar\le a/16$ and $\sbar\le\frac12$.
\item[\textup{(B2)}] $\mu^4\sbar^2\le\frac18$.
\item[\textup{(B3)}] $\mu\sbar^2\le\Gone$.
\item[\textup{(B4)}] $\sup_{0\le G\le \mu\sbar^2}|\lamst(G)-\lamref|\le\sigma(\mu-1)/(8\mu)$
  (this is (Det)).
\item[\textup{(B5)}] For every $k\ge1$: $\Lconf_ks_k\le a/8$ and
  $\frac83(\frac2a)^3(\Lconf_k)^2s_k\le\mu\Cstar$, where
  \[
  \Lconf_k=2+2\Cstar\mu\Big[\Tpin+1+\frac1{\sigma'}\Big(5\log\frac1{s_k}+\log_+\frac{\delta_0}{\Klo}\Big)\Big].
  \]
\item[\textup{(B6)}] $\Kup\sbar^5\le\kseed$.
\item[\textup{(B7)}] $16\Cvel C_Q\CDuh\KO\sbar^6\le1$, $\KO\sbar^6\le\frac14$,
  $\KO\mu\sbar^5\le\frac12$, and
  \begin{equation}
  \label{eq:E22}
  \KO\sbar^6\Big[\Tz+2\log\frac1{\sbar}+\frac1{\mu-1}\Big]\le\frac{\mu-1}{4\Cvel C_QC_\tau}\min\Big\{\frac12,\frac1{4\Cls K_1}\Big\}.
  \end{equation}
\item[\textup{(B8)}] $\eta_b\le\frac12$; for every $k\ge2$,
  \[
  \frac{\Cls\eta_{b,k}}{\sigma'}\Big[5\log\frac1{s_k}+\log_+\frac{\delta_0}{\Klo}\Big]\le\frac1{32};
  \]
  and $\Cls\eta_b\sigma'^{-1}\log(2\delta_0/\kseed)\le\frac1{32}$.
\item[\textup{(B9)}] $2\Abar\mu\sbar^5\le\majst$ and $\Cvel (\mu\sbar^2+\sbar)\le1$.
\item[\textup{(B10)}] $4\Cvel \big(\frac{2C_Q}{\sigb(\mu-1)}+1\big)\mu\sbar^5\le\frac12$.
\item[\textup{(B11)}] $\|\ellst\|C_QK_1\Abar\mu\sbar^5\le\frac12$, and
  \[
  \Cerr\eta_b\delta_0\frac{\mu^2}{\sigma^2(\mu-1)^2}+2\pi C_3\eta_b\frac{4\mu\delta_0}{3\sigma(\mu-1)}
  \le\frac14\cdot\frac{15}{16}\cdot\frac{\cout\delta_0}{20|\lamref|(\mu-1)} .
  \]
\end{description}
\end{definition}

The second part of (B1), and (B2), hold automatically (see the remark after
Equation~\eqref{eq:sk-small}); they are listed because Proposition~\ref{prop:confine} and
Lemma~\ref{lem:field-bound} state them as hypotheses.

\subsection{The accumulated outer field}
\label{sec:param-feedback}

The outer-field bound \eqref{eq:O-feedback} has the small factor $\theta_j$,
but its effect accumulates throughout the quiet interval. The next lemma
uses (B7), Equation~\eqref{eq:E22}, to bound $\theta_jA_jt$ for
$A_jt\le\bar\tau^\sharp_j$. The two thresholds are those required by
Proposition~\ref{prop:relay}\textup{(H2)} and Lemma~\ref{lem:Oprime}, and by
Lemma~\ref{lem:handoff}\textup{(ii)}, respectively.

\begin{lemma}
\label{lem:theta-quiet}
Assume (SA) and \textup{(B7)}, and let $K_1$ be as in Definition~\ref{def:K1},
$\Tz$ and $C_\tau$ as in \textup{(P3)}, $\theta_j$ as in Lemma~\ref{lem:O} and
$\bar\tau^\sharp_j$ as in Lemma~\ref{lem:quiet-time}. Then
\begin{equation}
\label{eq:theta-quiet}
\theta_j\bar\tau^\sharp_j\le\min\Big\{\frac12,\frac1{4\Cls K_1}\Big\}
\qquad(1\le j<n).
\end{equation}
\end{lemma}

\begin{proof}
For $j=1$, $\theta_1=0$ by Lemma~\ref{lem:O}, and the claim holds. Let
$2\le j<n$. By Lemma~\ref{lem:O}, $\theta_j=4\Cvel C_Q\KO s_{j-1}^6$.

We first bound $\bar\tau^\sharp_j$ by a multiple of $C_\tau$. By (SA3),
$1<\mu\le\frac{17}{16}$, so $\mu\log(2\mu)\le\frac{17}{16}\log\frac{17}8<1$.
The three terms of $\bar\tau^\sharp_j$ (Lemma~\ref{lem:quiet-time}) are those in
$\Tz$, in $\log(1/s_j)$ and in $(\mu-1)^{-2}$; we compare them with
those of $C_\tau[\Tz+\log(1/s_j)+(\mu-1)^{-1}]/(\mu-1)$. By the definition of
$C_\tau$ in \textup{(P3)} and the bound on $\mu\log(2\mu)$, the coefficients
satisfy $3\le C_\tau$,
$11/\sigma'\le C_\tau$ and $5\mu\log(2\mu)/\sigma'\le5/\sigma'\le C_\tau$;
this gives $\bar\tau^\sharp_j\le C_\tau\big[\Tz+\log(1/s_j)+(\mu-1)^{-1}\big]/(\mu-1)$.
Next we pass from $s_j$ to $s_{j-1}$. Since
$s_{j-1}\le\sbar\le a/16<\frac1{2\mu}$, we have
$\log\frac1{s_j}=\log\frac1{s_{j-1}}+\log(2\mu)\le2\log\frac1{s_{j-1}}$, so
\[
\theta_j\bar\tau^\sharp_j\le\frac{4\Cvel C_QC_\tau}{\mu-1}\,\KO\,g(s_{j-1}),\qquad
g(x)=x^6\Big[\Tz+2\log\frac1x+\frac1{\mu-1}\Big].
\]
It remains to replace $s_{j-1}$ by $\sbar$, since Equation~\eqref{eq:E22}
bounds $\KO g(\sbar)$; for this we show that $g$ is nondecreasing near $0$. For $0<x<1$,
$g'(x)=x^5\big[6\big(\Tz+2\log\frac1x+\frac1{\mu-1}\big)-2\big]$, which is
nonnegative when $\log(1/x)\ge1/6$, because $\Tz>0$ and $\mu>1$. So $g$ is
nondecreasing on $(0,e^{-1/6}]$. Since $s_{j-1}\le\sbar\le a/16<\frac1{16}<e^{-1/6}$,
$g(s_{j-1})\le g(\sbar)$, and Equation~\eqref{eq:E22} gives
\[
\theta_j\bar\tau^\sharp_j\le\frac{4\Cvel C_QC_\tau}{\mu-1}\KO g(\sbar)\le\min\{\frac12,\frac1{4\Cls K_1}\}.\qedhere
\]
\end{proof}

\subsection{Verification of the parameter choices}
\label{sec:param-verification}

\begin{proof}[Detailed verification of Lemma~\ref{lem:nonempty}]
We take the steps in order. At each step the quantities of the earlier steps are
fixed, and we show that the conditions of the step can be met.

\emph{(P0)--(P4).} The constants of (P0) exist and are finite by the results
cited there; in particular, by the formulas of (SA2) and $M\ge s+4$
(Theorem~\ref{thm:active}), $s_h=M-1\ge s+3\ge3$ and $N_Q\sigb\ge\Lam(s_h)$, as
Theorem~\ref{thm:quiet} requires. The interval $\big(\frac9{200},\min(\sigb,\frac34\sigma)\big)$ for $\sigma'$
in (SA1) is nonempty: by Proposition~I.\ref{I-prop:U}\textup{(a)}, $\sigb$ is attained and
$\sigb>\frac{211}{2500}>\frac9{200}$; and
$\frac34\sigma\ge\frac34\sigb>\frac{633}{10000}>\frac9{200}$. Moreover
$c_\gamma>0$, because $\sigma'<\frac34\sigma<\frac{15}{16}\sigma$. The conditions in (P1) hold for sufficiently small $\delta_0>0$ because
$\sigma\ge\sigb>0$, $\cout>0$ and the other constants have already been fixed.
In (P2), $c'>0$ because $\Cmom>0$ and $0<\Ip{3}{\sigma/4}<\infty$, so
$\mu_{\max}>1$. In (P3), $\Tpin$ exists because $\rT(T)=(1+T)^3e^{-\sigma T/4}$
(Definition~\ref{def:standing}) tends to $0$ as $T\to\infty$, as $\sigma>0$. The formulas of (P3) and (P4) have positive finite denominators:
$c_N>0$, $\delta_0>0$, $\sigb>0$, $\mu-1>0$ and $\Klo>0$. So all the
constants of (P0)--(P4) are finite and positive, with $C_0\ge1$ as
Theorem~\ref{thm:active} requires.
Each can be chosen in terms of the quantities in its row of Table~\ref{tab:param}, which are fixed at
earlier rows or earlier in the same step.

\emph{(P5).} Since (P5) is imposed for every $\kappa\in(0,2\kseed]$, it suffices
to find $x_0>0$ such that each condition holds for every $\kappa\in(0,x_0]$; then every $\kseed\in(0,x_0/2]$ is admissible. For ($\kappa$1), ($\kappa$4), ($\kappa$6) this holds
because $\delta_04^{-\mu/(\mu-1)}$, $\kappa_0(C_0)$, $\kL$ are positive. For
($\kappa$2), ($\kappa$3), ($\kappa$7) it suffices that each left side tends to
$0$ as $\kappa\to0$. The left sides are finite sums of terms $C\kappa^{1/2}$ and
$C(\kappa/\delta_0)^\gamma\langle\log(\delta_0/\kappa)\rangle^N$ with $C$ fixed at
(P0)--(P4), integers $N\ge0$, and
$\gamma\in\{\gamma_\mu,\vartheta_\mu,\gamma'_\mu\}$, positive because $\mu>1$ and
$c_\gamma>0$; each term tends to $0$. For ($\kappa$5) we need a bound uniform in
$G$. Let $\kappa<\delta_0$ and $G\in[0,\Gone]$, and put $c=\sigb-\sigma'>0$, where
$\sigma'<\sigb$ by (SA1). By Theorem~\ref{thm:active},
$t_1(C_0,G)\le\bar t_1(C_0)$, and $\sigma(G)\in[\sigb,\sigma_{\max}]$, so
$\hat T\ge\sigma_{\max}^{-1}\log(\delta_0/\kappa)$. If moreover
$\kappa\le\delta_0e^{-\sigma_{\max}\bar t_1(C_0)}$, then $\hat T\ge\bar t_1(C_0)$, so
$y=\hat T-t_1(C_0,G)\ge\hat T-\bar t_1(C_0)\ge0$. Hence $\hat T$ lies in the range
$t\ge t_1$ of the second formula for the transient in Theorem~\ref{thm:active},
and, as $\sigma(G)-\sigma'\ge c$,
\[
\rhohat(\hat T)=\jap{y}^3e^{-(\sigma(G)-\sigma')y}\le\jap{y}^3e^{-cy}\le C(c)e^{-cy/2}
\le C(c)e^{c\bar t_1(C_0)/2}\Big(\frac\kappa{\delta_0}\Big)^{c/(2\sigma_{\max})},
\]
where $C(c)$ can be chosen because $\sup_{y\ge0}\jap{y}^3e^{-cy/2}<\infty$; the last step uses
$y\ge\hat T-\bar t_1(C_0)$ and the lower bound on $\hat T$. The right side does not depend on
$G$ and tends to $0$ as $\kappa\to0$, and so does $\kappa^{1/2}$; so ($\kappa$5)
holds for all small $\kappa$, uniformly in $G\in[0,\Gone]$. We take for $x_0$
the least of the thresholds found for ($\kappa$1)--($\kappa$7).

\emph{(P6)}, with (P0)--(P5) fixed. It suffices that each condition holds for
every $\sbar$ in some interval $(0,\sbar_i]$, $\sbar_i>0$; then
every $\sbar\in(0,\min_i\sbar_i]$ is admissible. The conditions (B1)--(B3), (B6), (B7), (B9),
(B10), (B11), the first part of (B8) and its last part (with $\kseed$ fixed)
have left sides that are finite sums of terms $C\sbar^\alpha\log^N(1/\sbar)$ with
$\alpha>0$, $N\in\{0,1\}$ and $C$ fixed at (P0)--(P5); here we use
$\eta_b=\Crelay\sbar^6/(\sigb(\mu-1))$. Each such term tends to $0$ as
$\sbar\to0$, while the right sides are fixed positive numbers. For (B4), the
supremum tends to $0$ as $\sbar\to0$, because $\lamst$ is continuous at $G=0$
(Proposition~I.\ref{I-prop:U}\textup{(a)}) and $\lamst(0)=\lamref$; the right side
$\sigma(\mu-1)/(8\mu)$ is positive. The conditions (B5) and (B8) are imposed at
every level $k$, and we reduce each of them to one condition on $\sbar$. For
(B5), put
$\alpha_0=2+2\Cstar\mu\big[\Tpin+1+\sigma'^{-1}\log_+(\delta_0/\Klo)\big]$ and
$\beta_0=10\Cstar\mu/\sigma'$, so that $\Lconf_k=\alpha_0+\beta_0\log(1/s_k)$.
The functions $x\mapsto x(\alpha_0+\beta_0\log(1/x))$ and
$x\mapsto x(\alpha_0+\beta_0\log(1/x))^2$
tend to $0$ as $x\to0$, so there is $\sbar_{\mathrm{B5}}>0$ such that both inequalities of
(B5) hold at level $k$ whenever $s_k\le\sbar_{\mathrm{B5}}$. Since $s_k\le\sbar$ for every
$k$, (B5) holds for every $\sbar\le\sbar_{\mathrm{B5}}$. For (B8) and $k\ge2$, we have
$s_k=s_{k-1}/(2\mu)$, and, by the formula for $\eta_{b,k}$ in
Equation~\eqref{eq:etab}, the $k$-th inequality reads $F(s_{k-1})\le\frac1{32}$
with
\[
F(x)=\frac{\Cls\Crelay}{\sigma'\sigb(\mu-1)}\,x^6\Big[5\Big(\log\frac1x+\log(2\mu)\Big)+\log_+\frac{\delta_0}{\Klo}\Big].
\]
Since $F(x)\to0$ as $x\to0$, there is $\sbar_{\mathrm{B8}}>0$ with $F\le\frac1{32}$ on
$(0,\sbar_{\mathrm{B8}}]$, and, since $s_{k-1}\le\sbar$ for $k\ge2$, the
inequalities hold for all $k\ge2$ once $\sbar\le\sbar_{\mathrm{B8}}$.

No formula or condition of (P0)--(P6) contains $n$, $A_1$ or $\omega^{(n)}$; the
levels in (B5) and (B8) range over all integers $k\ge1$ and $k\ge2$, not only
up to $n$.

\emph{(SA) and (Det).} Assume (P0)--(P6) and let $n\ge2$ and $A_1>0$. (SA1) holds by (P0.2) and (P0.3), and (SA2) by (P0.4) and (P0.5). For
(SA3): $\mu\in(1,\frac{17}{16}]$ because
$\mu_{\max}\le1+\frac1{16}$; $\sbar\in(0,a/16]$ and $\mu\sbar^2\le\Gone$ by (B1) and (B3);
$\delta_0\in(0,\dH]$ by (P1); $\Tpin>0$ by (P3); and $\kseed\in(0,\kL]$ by
($\kappa$6) at $\kappa=2\kseed$. The objects of
Definitions~\ref{def:datum}--\ref{def:clock} are then defined with these
parameters: Definition~\ref{def:datum} needs $\mu\in(1,2]$ and $\sbar\in(0,a/16]$,
Definition~\ref{def:approx} needs $\kseed\in(0,\kL]$, and
Definition~\ref{def:clock} needs $\delta_0,\Tpin>0$ and $G_k\le\Gone$ for
$1\le k\le n$, which holds by Lemma~\ref{lem:inner-circ} and (B3). Hypothesis
(Det) is (B4).
\end{proof}

\subsection{Growth rates and transmission of the velocity field}
\label{sec:relay-transmission}

The growth gap makes the response of an outer quiet ring smaller than the
response of the inner ring driving it. Recall that $x_i=\sigma_i/\mu$.

\begin{lemma}
\label{lem:detuning}
Assume (SA) and (Det). Put $e_\mu=\frac{\mu-1}{8\mu}$. Then for all levels $i,l$,
\begin{equation}
\label{eq:detuning-rates}
|\sigma_i-\sigma|\le e_\mu\sigma,\qquad e_\mu\le\frac1{136},\qquad
\sigma_i\ge\frac{15}{16}\sigma>\sigma',\qquad
\frac\sigma4\le x_i\le\frac{17}{16}\sigma,
\end{equation}
and
\begin{equation}
\label{eq:detuning-gap}
\mu\sigma_i-\sigma_l\ge\sigma\big((\mu-1)-e_\mu(\mu+1)\big)\ge\tfrac34\sigma(\mu-1)\ge\tfrac12\sigb(\mu-1);
\end{equation}
in particular,
\begin{equation}
\label{eq:detuning-ordered}
\sigma_i\mu^{i-j}\le\sigma_j-\frac{\sigb(\mu-1)}{2\mu}\qquad(i<j).
\end{equation}
\end{lemma}
\begin{proof}
By Lemma~\ref{lem:inner-circ} and (SA3), $0\le G_i\le \mu\sbar^2$, so (Det)
gives $|\sigma_i-\sigma|=|\Re(\lamst(G_i)-\lamref)|\le e_\mu\sigma$. As
$1<\mu\le\frac{17}{16}$, $e_\mu=\frac18(1-\frac1\mu)\le\frac1{136}$, so
$\sigma_i\ge(1-e_\mu)\sigma\ge\frac{15}{16}\sigma>\frac34\sigma>\sigma'$ by
(SA1) and $\sigma=\sigma(0)\ge\sigb>0$ (Proposition~I.\ref{I-prop:U}(a)); and
$x_i=\sigma_i/\mu$ lies between $\frac{15}{16}\cdot\frac{16}{17}\sigma\ge\frac\sigma4$
and $(1+e_\mu)\sigma\le\frac{17}{16}\sigma$. For levels $i,l$,
$\mu\sigma_i-\sigma_l\ge\mu(1-e_\mu)\sigma-(1+e_\mu)\sigma=\sigma((\mu-1)-e_\mu(\mu+1))$,
and $e_\mu(\mu+1)=\frac{(\mu-1)(\mu+1)}{8\mu}\le\frac14(\mu-1)$ because
$\mu+1\le2\mu$; with $\sigma\ge\sigb$ this gives Equation~\eqref{eq:detuning-gap}. For $i<j$,
$\mu^{i-j}\le\mu^{-1}$, and Equation~\eqref{eq:detuning-gap} with $(j,i)$ in place of $(i,l)$ gives
$\sigma_i\mu^{i-j}\le\sigma_i/\mu\le\sigma_j-\sigb(\mu-1)/(2\mu)$.
\end{proof}

Recall that $\etin_i$ and $\etout_i$ are the $C^{M_h}(N')$ norms of
$\hin_i$ and $\hout_i$ in ring-$i$ units. The inner field has its monopole
removed. The next lemma bounds $\etin_i$ using $w_{i+1}$ and $\etin_{i+1}$,
and bounds $\etout_i$ by summing over $w_l$ for $l<i$.

\begin{lemma}
\label{lem:transfer-step}
Assume (SA), and let $t\ge0$.
\begin{enumerate}[label=\textup{(\arabic*)}]
\item Let $1\le i<n$. Assume that $S_{i+1}(t)\subset N_{i+1}$ and that the fields
  $h_i$ and $h_{i+1}$ are defined at $t$ (Definition~\ref{def:regions}). Then for
  every integer $s'\ge0$,
  \begin{equation}
  \label{eq:transfer-step-inner}
  \etin_i(t)\le \Cvel \varepsilon_i\big(\|w_{i+1}(t)\|_{H^{s'}}+\etin_{i+1}(t)\big).
  \end{equation}
\item Let $1\le k\le n$. Assume that $h_k$ is defined at $t$ and that
  $S_i(t)\subset N_i$ for every $i<k$. Then for every integer $s'\ge0$,
  \begin{equation}
  \label{eq:transfer-step-outer}
  \etout_k(t)\le \Cvel \sum_{i<k}b_{i,k}\|w_i(t)\|_{H^{s'}},\qquad b_{i,k}=\mu^{i-k}r_k/r_i;
  \end{equation}
  in particular $\etout_1=0$.
\end{enumerate}
\end{lemma}
\begin{proof}
We use (SA2) for the order $M_h$ and constant $\Cvel$ of
Lemma~\ref{lem:transfer}, and (SA3) for $\sbar\le a/16$ and the ring
construction whose properties we recall next.
We use three facts. (a) $\omega_{\ge l}(t)$ is real and $3$-fold rotationally symmetric
(Lemma~\ref{lem:regions-basic}(3)). (b) $\int w_l(t)=0$ if $S_l(t)\subset N_l$
(Lemma~\ref{lem:regions-basic}(4)). (c) For real radial $g\in L^1(\R^2)$
supported in $\overline B_\rho$, the radial Biot--Savart formula of
Lemma~I.\ref{I-lem:radial-stationary} gives
$\KBS\ast g(x)=\frac1{2\pi}(\int g)\frac{x^\perp}{|x|^2}$ for $|x|\ge2\rho$.
For bounded $g$ this follows directly, with pointwise equality outside its
support by continuity. Truncation of $g$ and dominated convergence give the
general $L^1$ case, since the kernel is bounded on the support for these $x$.

(1) We split the vorticity inside ring $i$ into ring $i+1$ and the region inside
ring $i+1$; in physical units,
$\omega_{\ge i+1}(t)=A_{i+1}\Phi(\cdot/r_{i+1})+A_{i+1}w_{i+1}(t,\cdot/r_{i+1})+\omega_{\ge i+2}(t)$.
The first term is radial and supported in $\overline B_{r_{i+1}}$, and
$N'_i\subset\{|x|\ge\frac a4r_i\}$ with $2r_{i+1}=2s_ir_i\le\frac a8r_i$, so by (c)
its velocity on $N'_i$ is a monopole; the second term has zero integral by (b).
Hence on $N'_i$ the velocity of $\omega_{\ge i+1}(t)$ minus its monopole is the
velocity of the second term, bounded in $C^{M_h}(N'_i)$, in ring-$i$ units, by
$\Cvel \varepsilon_i\|w_{i+1}(t)\|_{H^{s'}}$ (Lemma~\ref{lem:transfer}(2)), plus
$E_g$, the velocity of $g=\omega_{\ge i+2}(t)$ minus its monopole. The function
$g$ is $3$-fold rotationally symmetric by (a) and supported in $D_{i+2}(t)\subset B_{ar_{i+1}/8}$,
as $h_{i+1}$ is defined at $t$; so by Lemma~\ref{lem:transfer}(1), $E_g$ is
bounded in $C^{M_h}(N'_i)$, in ring-$i$ units, by $\Cvel \varepsilon_i$ times
$\sup|E_g|$ on $\{|x|=ar_{i+1}/4\}$ in ring-$(i+1)$ units. In these units $E_g$ is
$\hin_{i+1}$ (its monopole is $G_{i+1}x^\perp/|x|^2$ by the definition of
$G_{i+1}$ and Lemma~\ref{lem:inner-circ}), and $\{|y|=a/4\}\subset N'$, so the
supremum is at most $\etin_{i+1}(t)$.

(2) By definition, $\hout_k$ is the velocity of
$\sum_{i<k}\omega^{(n)}(t)\mathbf 1_{S_i(t)}$ on $N'_k$, in ring-$k$ units. By the
last assertion of Lemma~\ref{lem:transfer}(4), this is the sum over $i<k$ of the
fields induced on $N'_k$ by $A_iw_i(t,\cdot/r_i)$, and by its first assertion each
of them is bounded in $C^{M_h}(N'_k)$, in ring-$k$ units, by
$\Cvel b_{i,k}\|w_i(t)\|_{H^{s'}}$. Sum over $i<k$; for $k=1$ the sum is empty.
\end{proof}

Under (SA) and (Det), for levels $i<j$ the mode of ring $j$ grows faster in physical time.
This lets its Duhamel bound control those of the outer rings. On an interval
$[0,\bar t]$ on which the fields are defined, the comparison uses just two
integral estimates.

Recall that
$\eta_i=\etin_i+\etout_i\ge0$ is the sum of the inner and outer field norms.
Thus the Duhamel bound in physical time is
\begin{equation}
\label{eq:O-duh}
\maj_i(t)=A_i\int_0^t e^{A_i\sigma_i(t-v)}\eta_i(v)\dd v.
\end{equation}
For $i<j$, first factor out $A_j$ and use $A_i/A_j=\mu^{i-j}$.
Lemma~\ref{lem:detuning} then gives
\[
A_j\sigma_j-A_i\sigma_i
=A_j\big(\sigma_j-\mu^{i-j}\sigma_i\big)
\ge\frac{A_j\sigb(\mu-1)}{2\mu}>0.
\]
We first record a generic integral estimate. Let $F$ be any nonnegative
integrable function on $[0,\bar t]$ satisfying
$e^{A_j\sigma_j(t-v)}F(v)\le F(t)$ for $0\le v\le t\le\bar t$.
Then
\begin{equation}
\label{eq:O-fast}
\begin{aligned}
A_i\int_0^t e^{A_i\sigma_i(t-v)}F(v)\dd v
&\le A_iF(t)\int_0^t e^{-(A_j\sigma_j-A_i\sigma_i)(t-v)}\dd v\\
&\le\frac{A_i}{A_j\sigma_j-A_i\sigma_i}F(t)\le\frac{2\mu}{\sigb(\mu-1)}F(t),
\end{aligned}
\end{equation}
since $A_i\le A_j$. Separately, $\etin_j\le\eta_j$ by the definition of
$\eta_j$, and $A_i\sigma_i\le A_j\sigma_j$ by Equation~\eqref{eq:detuning-ordered}.
Using these inequalities and then $A_i\le A_j$ gives
\begin{equation}
\label{eq:O-same}
\begin{aligned}
A_i\int_0^t e^{A_i\sigma_i(t-v)}\etin_j(v)\dd v
&\le A_i\int_0^t e^{A_j\sigma_j(t-v)}\eta_j(v)\dd v\\
&\le A_j\int_0^t e^{A_j\sigma_j(t-v)}\eta_j(v)\dd v=\maj_j(t).
\end{aligned}
\end{equation}

\subsection{Comparison of quiet rings}
\label{sec:proof-O}

The outer quiet rings respond more weakly than the innermost quiet ring.
Combining the growth gap with spatial separation gives a product of small
factors at each outward step and makes their combined feedback small.

\begin{lemma}
\label{lem:O}
Assume (SA) and (Det). Put $\CDuh=2\mu/(\sigb(\mu-1))$ and
\[
\KO=8\Cvel (C_Q\CDuh+1)+32\Cvel C_Q\CDuh .
\]
Let $1\le j<n$, assume that the clock is
defined up to level $j+1$, and let $\bar t\in(0,\tsh_j]$. Assume
that for every $t\in[0,\bar t]$ and every $i\le j$,
\begin{equation}
\label{eq:O-support}
S_l(t)\cap B_{6r_i}=\emptyset\ (l<i),\qquad
D_{i+1}(t)\subset B_{ar_i/8},\qquad \eta_i(t)\le1.
\end{equation}
Assume moreover
\begin{equation}
\label{eq:O-smallness}
\maj_j\le\majst\text{ on }[0,\bar t],\qquad
\KO\sbar^6\le\frac14,\qquad
16\Cvel C_Q\CDuh\KO\sbar^6\le1,\qquad \KO\mu\sbar^5\le\frac12.
\end{equation}
Then for every
$t\in[0,\bar t]$ and every $i<j$:
\begin{equation}
\label{eq:O-comparison}
\maj_i(t)\le\tfrac34\Ptil_{i,j}\,\maj_j(t),\qquad \Ptil_{i,j}=\prod_{p=i}^{j-1}\KO\varepsilon_p ,
\end{equation}
and
\begin{equation}
\label{eq:O-feedback}
\etout_j(t)\le\theta_j\maj_j(t),\qquad
\theta_j=4\Cvel C_Q\KO s_{j-1}^6\ (j\ge2),\qquad \theta_1=0.
\end{equation}
Moreover, for
every $t\in(0,\bar t]$ and every $i\le j$, the hypotheses of
Theorem~\ref{thm:quiet} hold for ring $i$ on $[0,A_it]$, in ring-$i$ units, with
the pulled-back functions of Section~\ref{sec:relay-data} and $G=G_i$; in particular its conclusions
\textup{(a)}--\textup{(c)} hold there.
\end{lemma}

\begin{proof}[Quantitative proof of Lemma~\ref{lem:O}]
The support conditions in Equation~\eqref{eq:O-support} define $h_i$,
$i\le j$, on $[0,\bar t]$ by Definition~\ref{def:regions}, so the
integral estimates~\eqref{eq:O-duh}--\eqref{eq:O-same} apply.
Put $\Pprod_{i,l}=\prod_{p=i}^{l-1}\Cvel \varepsilon_p$ for $i<l$, so that
$\Pprod_{i,l}=(\Cvel /\KO)^{l-i}\Ptil_{i,l}$ and
$\Ptil_{i,l}\Ptil_{l,j}=\Ptil_{i,j}$.
We prove Equation~\eqref{eq:O-comparison} by a continuity argument in $t$.
Equation~\eqref{eq:O-BH} below is this bound with the factor $1$ in place of
$\frac34$ and with a small term $\epsilon e^{A_j\sigma_jt}$ added to $\maj_j$, so
that its right side is positive at $t=0$, and we improve the factor to $\frac14$.

Fix $\epsilon\in(0,\majst e^{-A_j\sigma_j\bar t}]$ (unrelated to the numbers
$\varepsilon_p$) and put $\maj^\epsilon_j(t)=\maj_j(t)+\epsilon e^{A_j\sigma_jt}$.
In applying the generic estimate~\eqref{eq:O-fast}, we take $F=\maj_j$ or
$F=\maj^\epsilon_j$. The required inequality follows from (M1) of
Lemma~\ref{lem:majorant}, rescaled to physical time; the added term
$\epsilon e^{A_j\sigma_jt}$ satisfies it with equality.
Let $t_\epsilon$ be the supremum of
the $t_1\in[0,\bar t]$ such that
\begin{equation}
\label{eq:O-BH}
\maj_i(t)\le\Ptil_{i,j}\,\maj^\epsilon_j(t)\qquad\text{for every }t\in[0,t_1]\text{ and every }i<j .
\end{equation}
Both sides are continuous in $t$, since $h_i$, $i\le j$, is defined on
$[0,\bar t]$; at $t=0$ the left side vanishes and the right side is
$\Ptil_{i,j}\epsilon>0$. So Equation~\eqref{eq:O-BH} holds on $[0,t_\epsilon]$, and
$t_\epsilon>0$.

We first verify the hypotheses of Theorem~\ref{thm:quiet} on $[0,t_\epsilon]$.
Every factor of $\Ptil_{i,j}$ satisfies
$\KO\varepsilon_p=\KO\mu s_p^5\le\KO\mu\sbar^5\le\frac12$, so
$\Ptil_{i,j}\le\frac12$, and Equation~\eqref{eq:O-BH} with $\maj_j\le\majst$ gives
$\maj_i(t)\le\frac12(\majst+\epsilon e^{A_j\sigma_j\bar t})\le\majst$ for
$t\in[0,t_\epsilon]$ and $i<j$. Ring $i<n$ starts at $\Phi$ in ring-$i$ units,
$G_i\in[0,\Gone]$, the support conditions in Equation~\eqref{eq:O-support}
give the separation required by Lemma~\ref{lem:reduction} for ring $i$, and
$\|h_i\|_{C^{M_h}(N')}\le\etin_i+\etout_i=\eta_i\le1$ with $\eta_i$ continuous. So
Theorem~\ref{thm:quiet} applies to ring $i<j$ on $[0,A_it_\epsilon]$ and to ring
$j$ on $[0,A_j\bar t]$ by the hypothesis $\maj_j\le\majst$; in particular, for $t\in[0,t_\epsilon]$ and $i\le j$,
$S_i(t)\subset\NNm_i\subset N_i$ and
\begin{equation}
\label{eq:O-Q}
\|w_i(t)\|_{H^{s_h}}\le C_Q\maj_i(t).
\end{equation}

Let $t\in[0,t_\epsilon]$ and $i<j$. Lemma~\ref{lem:transfer-step}(1)
at the levels $i,\dots,j-1$ with $s'=s_h$ gives by induction
$\etin_i(t)\le\sum_{l=i+1}^{j}\Pprod_{i,l}\|w_l(t)\|_{H^{s_h}}+\Pprod_{i,j}\etin_j(t)$,
and Lemma~\ref{lem:transfer-step}(2) gives
$\etout_i(t)\le \Cvel \sum_{l<i}b_{l,i}\|w_l(t)\|_{H^{s_h}}$. We insert
Equations~\eqref{eq:O-Q} and \eqref{eq:O-BH} into these two bounds and group the
terms of $\eta_i(t)=\etin_i(t)+\etout_i(t)$ by the ring they come from.
Ring $j$ and the region inside it (the terms with $\|w_j(t)\|_{H^{s_h}}$ and
$\etin_j(t)$) give
at most $\frac{\Cvel}{\KO}\Ptil_{i,j}(C_Q\maj^\epsilon_j(t)+\etin_j(t))$, and the
rings $l\in(i,j)$ give
$C_Q\Ptil_{i,j}\maj^\epsilon_j(t)\sum_{l=i+1}^{j-1}(\Cvel /\KO)^{l-i}\le\frac87C_Q\frac{\Cvel}{\KO}\Ptil_{i,j}\maj^\epsilon_j(t)$,
as $\Cvel /\KO\le\frac18$. For the rings $l<i$ we have
$b_{l,i}\Ptil_{l,j}=b_{l,i}\Ptil_{l,i}\Ptil_{i,j}$ and
\[
b_{l,i}\Ptil_{l,i}=\mu^{l-i}\frac{r_i}{r_l}\KO^{i-l}\mu^{i-l}\prod_{p=l}^{i-1}s_p^5
=\KO^{i-l}\Big(\frac{r_i}{r_l}\Big)^6\le\big(\KO\sbar^6\big)^{i-l}\le4^{l-i}
\]
because $s_l\cdots s_{i-1}=r_i/r_l\le\sbar^{i-l}$ and $\KO\sbar^6\le\frac14$;
summing the geometric series, these rings give at most $\frac43\Cvel C_Q\KO\sbar^6$ times $\Ptil_{i,j}\maj^\epsilon_j(t)$. We insert
these bounds into Equation~\eqref{eq:O-duh} and use Equation~\eqref{eq:O-fast}
with $F=\maj^\epsilon_j$ and Equation~\eqref{eq:O-same}:
\[
\maj_i(t)\le\Big[\frac{\Cvel}{\KO}\big(C_Q\CDuh+1\big)+\frac87\cdot\frac{\Cvel C_Q\CDuh}{\KO}+\frac43\Cvel C_Q\CDuh\KO\sbar^6\Big]\Ptil_{i,j}\,\maj^\epsilon_j(t).
\]
As $\KO\ge8\Cvel (C_Q\CDuh+1)$, $\KO\ge32\Cvel C_Q\CDuh$ and
$16\Cvel C_Q\CDuh\KO\sbar^6\le1$, the bracket is at most
$\frac18+\frac1{28}+\frac1{12}<\frac14$, so
\begin{equation}
\label{eq:O-improved}
\maj_i(t)\le\tfrac14\Ptil_{i,j}\,\maj^\epsilon_j(t)\qquad(t\in[0,t_\epsilon],\ i<j).
\end{equation}
The right side is positive, so Equation~\eqref{eq:O-BH} is strict at
$t=t_\epsilon$ and, by continuity, holds beyond $t_\epsilon$ unless
$t_\epsilon=\bar t$. Hence $t_\epsilon=\bar t$ for every admissible $\epsilon$.

So, for every $t\in(0,\bar t]$, the hypotheses of Theorem~\ref{thm:quiet} hold
for ring $i\le j$ on $[0,A_it]$, since they hold on $[0,A_i\bar t]$; and
$\epsilon\to0$ in Equation~\eqref{eq:O-improved} gives
$\maj_i(t)\le\frac14\Ptil_{i,j}\maj_j(t)\le\frac34\Ptil_{i,j}\maj_j(t)$.
The identity and estimate above, now with indices $(i,j)$, read
\begin{equation}
\label{eq:O-product}
b_{i,j}\Ptil_{i,j}
=\KO^{j-i}\Big(\frac{r_j}{r_i}\Big)^6
\le(\KO\sbar^6)^{j-i}\le4^{-(j-i)}\qquad(i<j),
\end{equation}
Finally, Equation~\eqref{eq:transfer-step-outer} at level $j$,
Equation~\eqref{eq:O-Q}, and the limit of Equation~\eqref{eq:O-improved} as
$\epsilon\to0$ give
\[
\etout_j(t)\le\tfrac14\Cvel C_Q\sum_{i<j}\KO^{j-i}\Big(\frac{r_j}{r_i}\Big)^6\maj_j(t).
\]
The term $i=j-1$ is $\KO s_{j-1}^6$, and the term $i<j-1$ is
\[
\KO s_{j-1}^6\KO^{j-1-i}(r_{j-1}/r_i)^6\le\KO s_{j-1}^6(\KO\sbar^6)^{j-1-i}\le\KO s_{j-1}^6\,4^{-(j-1-i)} .
\]
So the sum is at most $\frac43\KO s_{j-1}^6$, and
$\etout_j(t)\le\frac13\Cvel C_Q\KO s_{j-1}^6\maj_j(t)\le\theta_j\maj_j(t)$.
For $j=1$ there is no outer ring and $\etout_1=0$.
\end{proof}

\subsection{The field transmitted by the next inner ring}
\label{sec:proof-relay}

To bound $\maj_j^{\mathrm{in}}$, we integrate $\etin_j$ separately before
$\tsh_{j+1}$, on $[\tsh_{j+1},T_{j+1}]$, and after $T_{j+1}$.
The next proposition estimates these integrals and uses
$\etout_j\le\theta_j\maj_j$ to obtain $\maj_j\le2\maj_j^{\mathrm{in}}$.

\begin{proposition}
\label{prop:relay}
Assume (SA) and (Det), and let $C_0\ge1$ and $K_1>0$.
There is a constant $\Crelay=\Crelay(C_0,K_1)$ for the estimates below. It
depends only on $C_0$, $K_1$, $\sigma$, $\sigb$ and the constants of (SA1)--(SA2), and not
on $\mu$, $\delta_0$,
$\Tpin$, $\sbar$, $\kseed$, $n$ or the level. Let $1\le j<n$, assume that the
clock is defined up to level $j+1$, let $\bar t\in[0,\tsh_j]$, and assume:
\begin{description}[font=\normalfont]
\item[\textup{(H1)}] \textup{(Loc)} holds at every $t\in[0,\bar t]$;
\item[\textup{(H2)}] $\etout_j\le\theta_j\maj_j$ on $[0,\bar t]$, with $\theta_j$ of
  Lemma~\ref{lem:O}, and $\theta_jA_j\bar t\le\frac12$;
\item[\textup{(H3)}] if $j+1<n$, the conclusions of Theorem~\ref{thm:quiet},
  applied by the rescaling of Section~\ref{sec:relay-data}, hold for ring $j+1$ through physical time
  $\min(\bar t,\tsh_{j+1})$;
\item[\textup{(H4)}] if $\bar t>\tsh_{j+1}$, then, for
  $t\in[\tsh_{j+1},\min(\bar t,T_{j+1})]$,
  \begin{equation}
  \label{eq:transmitted-active-data}
  \begin{gathered}
  \maj_{j+1}(\tsh_{j+1})\le K_1\ksh_{j+1},\qquad
  D_{j+2}(t)\subset B_{ar_{j+1}/8},\\
  \|w_{j+1}(t)\|_{H^s}\le2C_2\hat a_{j+1}(t)(1+\rhohat_{j+1}(t)),\\
  \|\hin_{j+1}(t)\|_{C^{M_h}(N')}\le C_0\ksh_{j+1}\jap{A_{j+1}(t-\tsh_{j+1})}^3.
  \end{gathered}
  \end{equation}
\end{description}
Recall that the Duhamel bound of the inner field, with its monopole removed, is
\begin{equation}
\label{eq:relay-inner-duhamel}
\maj^{\mathrm{in}}_j(t)=A_j\int_0^t
e^{A_j\sigma_j(t-v)}\etin_j(v)\dd v,
\qquad \etin_j=\|\hin_j\|_{C^{M_h}(N')}.
\end{equation}
Then
\begin{equation}
\label{eq:relay-absorb}
\maj_j(t)\le2\maj^{\mathrm{in}}_j(t)\qquad(0\le t\le\bar t),
\end{equation}
and:
\begin{enumerate}[label=\textup{(\roman*)}]
\item[\textup{(o)}] if $j+1<n$, then for every $t\in[0,\min(\bar t,\tsh_{j+1})]$,
  \begin{equation}
  \label{eq:relay-before-activation}
  \maj^{\mathrm{in}}_j(t)\le \Cvel \Big(\frac{2C_Q}{\sigb(\mu-1)}+1\Big)\varepsilon_j\,\maj_{j+1}(t);
  \end{equation}
\item if $\bar t\le\tsh_{j+1}$, then
  $\maj_j(\bar t)\le2\Cvel \big(\frac{2C_Q}{\sigb(\mu-1)}+1\big)\varepsilon_j\,\maj_{j+1}(\bar t)$;
\item if $\bar t=\tsh_j$, then
  \begin{equation}
  \label{eq:relay}
  \maj_j(\tsh_j)\le2\Cvel \varepsilon_je^{x_j\Tpin}\Big[\frac{4C_2\delta_0}{\sigb(\mu-1)}+C_D(\Tpin)+\Xi_{j+1}\Big],
  \end{equation}
  where $C_D(T)=(2+2\Cstar T)^3T$,
  $\gamma_j=\max(\sigma_j,\mu\sigma')/(\mu\sigma_{j+1})$, and
  \begin{equation}
  \label{eq:relay-remainder}
  \Xi_{j+1}=\Crelay\Big(1+\frac1{\sigb(\mu-1)}\Big)\delta_0\Big(\frac{\ksh_{j+1}}{\delta_0}\Big)^{1-\gamma_j}\Big\langle\log\frac{\delta_0}{\ksh_{j+1}}\Big\rangle^4;
  \end{equation}
\item if $\tsh_{j+1}<\bar t<\tsh_j$, then $\maj_j(\bar t)$ is at most the right side
  of Equation~\eqref{eq:relay}.
\end{enumerate}
All times used are at most $\bar t$ or are determined by $\omega^{(n)}$ on
$[0,\bar t]$.
\end{proposition}

\begin{proof}
For this calculation set $\tau=A_j\bar t$, and write $\tau^{\mathrm{act}}=A_j\tsh_{j+1}$ and $\tau^{\mathrm{sat}}=A_jT_{j+1}$
(the ring-$j$ times at which ring $j+1$ becomes active and saturates),
$\kappa=\ksh_{j+1}$ and $\hat T=\hat T_{j+1}=\sigma_{j+1}^{-1}\log(\delta_0/\kappa)$.
As the clock is defined up to level $j+1$, $\kappa\in(0,\delta_0)$,
$\tau^{\mathrm{sat}}=\tau^{\mathrm{act}}+\hat T/\mu$ and
$A_j\tsh_j=\tau^{\mathrm{sat}}+\Tpin/\mu$; at a ring-$j$ time $v\ge\tau^{\mathrm{act}}$
the active time of ring $j+1$ is $\mu(v-\tau^{\mathrm{act}})$. Every ring-$j$ time
$v\le\tau$ corresponds to $v/A_j\le\tsh_j$, where ring $j$ is quiet; so by
(Loc)(2) the rings $i<j$ lie outside $B_{6r_j}$, $D_{j+1}(v/A_j)\subset B_{ar_j/8}$,
and $h_j$ is defined there.

We first reduce the claims to bounds on $\maj^{\mathrm{in}}_j$, the Duhamel bound of
the inner field. Let $\tau'\in[0,\tau]$. As $\eta_j=\etin_j+\etout_j$, (H2) and (M1) of
Lemma~\ref{lem:majorant} for $\maj_j$ give
$\maj_j(\frac{\tau'}{A_j})\le\maj^{\mathrm{in}}_j(\frac{\tau'}{A_j})+\theta_j\int_0^{\tau'}e^{\sigma_j(\tau'-v)}\maj_j(\frac{v}{A_j})\dd v
\le\maj^{\mathrm{in}}_j(\frac{\tau'}{A_j})+\theta_j\tau'\maj_j(\frac{\tau'}{A_j})$,
and $\theta_j\tau'\le\frac12$, so $\maj_j\le2\maj^{\mathrm{in}}_j$. It remains to
bound $\maj^{\mathrm{in}}_j$ by half the right sides of (i)--(iii), and to prove (o).

\emph{Before activation of ring $j+1$.} If $j+1=n$, then $\tau^{\mathrm{act}}=0$, (o) is
void and (i) concerns only $\tau=0$, where $\maj_j(0)=0$. Let $j+1<n$,
$\tau'\le\min(\tau,\tau^{\mathrm{act}})$, $v\in[0,\tau']$ and $t=v/A_j$. Then
$t\le\tsh_{j+1}<T_{j+1}$, so ring $j+1$ is quiet at $t$: (Loc)(1) gives
$S_{j+1}(t)\subset N_{j+1}$, and (Loc)(2) for ring $j+1$ gives
$D_{j+2}(t)\subset B_{ar_{j+1}/8}$ and $S_i(t)\cap B_{6r_{j+1}}=\emptyset$ for $i\le j$;
so $h_{j+1}$ is defined at $t$. Lemma~\ref{lem:transfer-step}(1) with $s'=s_h$
and Theorem~\ref{thm:quiet}(c) for ring $j+1$, which holds by (H3), give
$\etin_j(\frac{v}{A_j})\le \Cvel \varepsilon_j(C_Q\maj_{j+1}(v/A_j)+\eta_{j+1}(v/A_j))$. By (M1),
$\maj_{j+1}(v/A_j)\le e^{-\mu\sigma_{j+1}(\tau'-v)}\maj_{j+1}(\tau'/A_j)$, and
$\mu\sigma_{j+1}-\sigma_j\ge\frac12\sigb(\mu-1)$ (Lemma~\ref{lem:detuning}), so
$\int_0^{\tau'}e^{\sigma_j(\tau'-v)}\maj_{j+1}(v/A_j)\dd v\le2\maj_{j+1}(\tau'/A_j)/(\sigb(\mu-1))$.
As $\sigma_j/\mu\le\sigma_{j+1}$, the substitution $v'=\mu v$ gives
\[
\int_0^{\tau'}e^{\sigma_j(\tau'-v)}\eta_{j+1}(v/A_j)\dd v=\frac1\mu\int_0^{\mu\tau'}e^{\frac{\sigma_j}\mu(\mu\tau'-v')}\eta_{j+1}(\frac{v'}{A_{j+1}})\dd v'\le\maj_{j+1}(\tau'/A_j) .
\]
Together,
\begin{equation}
\label{eq:relay-quiet}
\maj^{\mathrm{in}}_j(\frac{\tau'}{A_j})\le \Cvel \Big(\frac{2C_Q}{\sigb(\mu-1)}+1\Big)\varepsilon_j\,\maj_{j+1}(\tau'/A_j)\qquad\big(0\le\tau'\le\min(\tau,\tau^{\mathrm{act}})\big),
\end{equation}
which is (o). If $\tau\le\tau^{\mathrm{act}}$, (o) at $\tau'=\tau$ and $\maj_j\le2\maj^{\mathrm{in}}_j$
give (i).

Let now $\tau^{\mathrm{act}}<\tau\le A_j\tsh_j$, and split
$\maj^{\mathrm{in}}_j(\frac{\tau}{A_j})=\int_0^\tau e^{\sigma_j(\tau-v)}\etin_j(\frac{v}{A_j})\dd v$ over
$\Iint_1=[0,\tau^{\mathrm{act}}]$, $\Iint_2=[\tau^{\mathrm{act}},\tau^{\mathrm{sat}}]$ and
$\Iint_3=[\tau^{\mathrm{sat}},A_j\tsh_j]$ (intersected with $[0,\tau]$), on which ring
$j+1$ is quiet, active and debris, respectively. In each part we use
$e^{\sigma_j(\tau-v)}\le e^{\sigma_j(A_j\tsh_j-v)}$ and
\begin{equation}
\label{eq:relay-id}
e^{\sigma_j(A_j\tsh_j-\tau^{\mathrm{sat}})}=e^{x_j\Tpin},\qquad
\kappa\,e^{\sigma_j(\tau^{\mathrm{sat}}-\tau^{\mathrm{act}})}=\kappa\Big(\frac{\delta_0}\kappa\Big)^{\frac{\sigma_j}{\mu\sigma_{j+1}}}\le\delta_0\Big(\frac\kappa{\delta_0}\Big)^{1-\gamma_j},
\end{equation}
where the inequality holds as $\sigma_j/(\mu\sigma_{j+1})\le\gamma_j$ and $\kappa<\delta_0$.

\emph{The part over $\Iint_1$} vanishes if $j+1=n$. Otherwise it equals
$e^{\sigma_j(\tau-\tau^{\mathrm{act}})}\maj^{\mathrm{in}}_j(\frac{\tau^{\mathrm{act}}}{A_j})$. We
bound the exponential factor by its value at $\tau=A_j\tsh_j$ and split its exponent
at $\tau^{\mathrm{sat}}$, so that both parts of Equation~\eqref{eq:relay-id} apply. With
Equation~\eqref{eq:relay-quiet} at $\tau'=\tau^{\mathrm{act}}$, the bound
$\maj_{j+1}(\tsh_{j+1})\le K_1\kappa$ of (H4) and
Equation~\eqref{eq:relay-id}, this bounds the part over $\Iint_1$ by
$\Cvel \varepsilon_je^{x_j\Tpin}K_1\big(\frac{2C_Q}{\sigb(\mu-1)}+1\big)\delta_0(\kappa/\delta_0)^{1-\gamma_j}$.

\emph{The part over $\Iint_2$.} Let $v\in[\tau^{\mathrm{act}},\min(\tau,\tau^{\mathrm{sat}})]$,
$t=v/A_j$ and $u=\mu(v-\tau^{\mathrm{act}})\in[0,\hat T]$. Then $t\le T_{j+1}$, so
$S_{j+1}(t)\subset N_{j+1}$ by (Loc)(1). The field $h_{j+1}$ is defined at $t$:
$D_{j+2}(t)\subset B_{ar_{j+1}/8}$ by (H4), the rings $i<j$ lie outside
$B_{6r_j}\supset B_{6r_{j+1}}$, and ring $j$ lies in $N_j\subset\{|x|\ge\frac a2r_j\}$
by (Loc)(1), which misses $B_{6r_{j+1}}$ as $6s_j\le\frac{6a}{16}<\frac a2$.
Lemma~\ref{lem:transfer-step}(1) with $s'=s$ and (H4) give
$\etin_j(\frac{v}{A_j})\le \Cvel \varepsilon_j(2C_2\hat a_{j+1}(t)+2C_2\hat a_{j+1}(t)\rhohat_{j+1}(t)+C_0\kappa\jap{u}^3)$.
As $\hat a_{j+1}(t)=\kappa e^{\sigma_{j+1}u}=\delta_0e^{-\mu\sigma_{j+1}(\tau^{\mathrm{sat}}-v)}$,
Lemma~\ref{lem:detuning} bounds the contribution of the first term by
$2C_2\Cvel \varepsilon_je^{x_j\Tpin}\delta_0\int_{-\infty}^{\tau^{\mathrm{sat}}}e^{-(\mu\sigma_{j+1}-\sigma_j)(\tau^{\mathrm{sat}}-v)}\dd v
\le \Cvel \varepsilon_je^{x_j\Tpin}\cdot4C_2\delta_0/(\sigb(\mu-1))$.
For the second term, Equation~\eqref{eq:active-transient} (with
$G=G_{j+1}$) gives $\hat a_{j+1}(t)\rhohat_{j+1}(t)\le C(C_0)\kappa\jap{u}^3e^{\sigma'u}$;
moreover $e^{\sigma_j(A_j\tsh_j-v)}=e^{x_j\Tpin}e^{\frac{\sigma_j}\mu(\hat T-u)}$ and,
for $u\in[0,\hat T]$,
$\frac{\sigma_j}\mu(\hat T-u)+\sigma'u\le\max(\frac{\sigma_j}\mu,\sigma')\hat T=\gamma_j\sigma_{j+1}\hat T=\gamma_j\log\frac{\delta_0}\kappa$.
With $\dd v=\dd u/\mu\le\dd u$, the second term contributes at most $2C_2C(C_0)$
times the left side of Equation~\eqref{eq:relay-transient-integral}, and the third (in the same way,
with $e^{\sigma'u}$ replaced by $1$) at most $C_0$ times it:
\begin{equation}
\label{eq:relay-transient-integral}
\Cvel \varepsilon_je^{x_j\Tpin}\kappa\Big(\frac{\delta_0}\kappa\Big)^{\gamma_j}\int_0^{\hat T}\jap{u}^3\dd u
\le \Cvel \varepsilon_je^{x_j\Tpin}\delta_0\Big(\frac\kappa{\delta_0}\Big)^{1-\gamma_j}\hat T\jap{\hat T}^3 .
\end{equation}
Finally $\hat T\le\sigb^{-1}\log(\delta_0/\kappa)$, as $\sigma_{j+1}\ge\sigb$, and
$\jap{cx}\le\max(1,c)\jap{x}$ for $c>0$, so
$\hat T\jap{\hat T}^3\le\langle\log(\delta_0/\kappa)\rangle^4/(\sigb\min(1,\sigb)^3)$.

\emph{The part over $\Iint_3$} matters only if $\tau>\tau^{\mathrm{sat}}$. Let
$v\in[\tau^{\mathrm{sat}},\tau]$ and $t=v/A_j$, so that $T_{j+1}\le t\le\tsh_j$ and
$\tauK=A_{j+1}(t-T_{j+1})=\mu(v-\tau^{\mathrm{sat}})\le\mu(A_j\tsh_j-\tau^{\mathrm{sat}})=\Tpin$. Ring $j$ is quiet at
$t$, so by (Loc)(2) $D_{j+1}(t)\subset B_{ar_j/8}\cap\overline B_{\varrho_j(t)}$ with
$\varrho_j(t)=(2+2\Cstar\tauK)r_{j+1}$. Hence $R_t=\max\{|x|:x\in D_{j+1}(t)\}$
satisfies $R_t<ar_j/8$ and $R_t\le(2+2\Cstar\Tpin)r_{j+1}$, and
Lemma~\ref{lem:transfer}(3) with $R=R_t$ gives
$\etin_j(\frac{v}{A_j})\le \Cvel \varepsilon_j(R_t/r_{j+1})^3\le \Cvel \varepsilon_j(2+2\Cstar\Tpin)^3$. As
$e^{\sigma_j(\tau-v)}\le e^{\sigma_j(A_j\tsh_j-\tau^{\mathrm{sat}})}=e^{x_j\Tpin}$ (as
$v\ge\tau^{\mathrm{sat}}$, $\tau\le A_j\tsh_j$, and by Equation~\eqref{eq:relay-id}) and the
interval has length at most $\Tpin/\mu\le\Tpin$, this part is at most
$\Cvel \varepsilon_je^{x_j\Tpin}C_D(\Tpin)$.

Adding the three parts, with $\langle\log(\delta_0/\kappa)\rangle\ge1$ for the
part over $\Iint_1$, $\maj^{\mathrm{in}}_j(\frac{\tau}{A_j})$ is at most
$\Cvel \varepsilon_je^{x_j\Tpin}[\frac{4C_2\delta_0}{\sigb(\mu-1)}+C_D(\Tpin)
+C(C_0,K_1)\big(1+(\sigb(\mu-1))^{-1}\big)\delta_0(\frac\kappa{\delta_0})^{1-\gamma_j}\langle\log\frac{\delta_0}\kappa\rangle^4]$
where $C(C_0,K_1)$ depends only on the stated fixed data. Taking
$\Crelay(C_0,K_1)$ large enough absorbs this coefficient.
So $\maj^{\mathrm{in}}_j(\frac{\tau}{A_j})$ is at most half the right side of
Equation~\eqref{eq:relay}, and $\maj_j\le2\maj^{\mathrm{in}}_j$ gives (ii) and (iii).

Finally, if $\tau\le\tau^{\mathrm{act}}$, the proof used the solution on
$[0,\tau/A_j]$ and the fact that ring $j+1$ is quiet there, which is determined by
$\omega^{(n)}$ on $[0,\tau/A_j]$ (Lemma~\ref{lem:clock}); otherwise
$\tsh_{j+1}<\tau/A_j$, and $\ksh_{j+1}$, $T_{j+1}$ and $\tsh_j$ are determined by
$\omega^{(n)}$ on $[0,\tsh_{j+1}]$ (Lemma~\ref{lem:clock}).
\end{proof}

\subsection{The feedback of the outer rings}
\label{sec:proof-Oprime}

The active ring drives the outer quiet rings, which in turn act back on it.
The next lemma shows that these two transmissions produce the small factor
$s_{k-1}^6$, so the feedback stays small compared with the active growth envelope.

\begin{lemma}
\label{lem:Oprime}
Assume (SA) and (Det), and let $C_0\ge1$ and $K_1>0$.
After enlarging the common constant $\Crelay(C_0,K_1)$ of
Proposition~\ref{prop:relay} if necessary, the estimate below holds. It
depends only on $C_0$, $K_1$, $\sigma$, $\sigb$ and the constants of
(SA1)--(SA2), and not on $\mu$,
$\delta_0$, $\Tpin$, $\sbar$, $\kseed$, $n$ or the level. Let $1\le k\le n$, assume
that the clock is defined up to level $k$, and let $t\in\Phase_k$ with $t>0$. Assume on
$[0,t]$: \textup{(Loc)}; if $k\ge2$, the hypotheses and conclusions of
Lemma~\ref{lem:O} for $j=k-1$ and $\bar t=t$ (ring $k-1$ is quiet on $[0,t]$,
since $t\le T_k<\tsh_{k-1}$), and $\theta_{k-1}A_{k-1}t\le\frac12$; the active
bound $\|w_k(t')\|_{H^s}\le2C_2\hat a_k(t')(1+\rhohat_k(t'))$ and the inner bound
$\|\hin_k(t')\|_{C^{M_h}(N')}\le C_0\ksh_k\jap{A_k(t'-\tsh_k)}^3$ for
$t'\in[\tsh_k,t]$; and, if $k<n$: the conclusions of
Theorem~\ref{thm:quiet} for ring $k$, after rescaling, through physical time $\tsh_k$,
$\maj_k(\tsh_k)\le K_1\ksh_k$, and $D_{k+1}(t')\subset B_{ar_k/8}$ for
$t'\in[\tsh_k,t]$. Then for every $t'\in[\tsh_k,t]$,
\begin{equation}
\label{eq:etab}
\|\hout_k(t')\|_{C^{M_h}(N')}\le\eta_{b,k}\,\hat a_k(t'),\qquad
\eta_{b,k}=\Crelay\frac{s_{k-1}^6}{\sigb(\mu-1)}\le\eta_b=\Crelay\frac{\sbar^6}{\sigb(\mu-1)}
\end{equation}
for $k\ge2$, and $\hout_1=0$.
\end{lemma}

\begin{proof}
For $k=1$, $\hout_1=0$ by Lemma~\ref{lem:transfer-step}(2). Let $k\ge2$. We use
ring-$(k-1)$ time $\tau'=A_{k-1}t'$ and put $\tau^{\mathrm{act}}=A_{k-1}\tsh_k$; at
$\tau'\ge\tau^{\mathrm{act}}$ the active time of ring $k$ is
$u=\mu(\tau'-\tau^{\mathrm{act}})$, and $\mu\tau^{\mathrm{act}}=A_k\tsh_k$.

\emph{Reduction to $\maj_{k-1}$.} Let $t'\in[\tsh_k,t]$. By (Loc) and the
hypotheses of Lemma~\ref{lem:O}, $h_k$ is defined at $t'$:
$D_{k+1}(t')\subset B_{ar_k/8}$ by hypothesis (trivially if $k=n$), the rings
$i<k-1$ lie outside $B_{6r_{k-1}}\supset B_{6r_k}$, and ring $k-1$ lies in
$N_{k-1}$, which misses $B_{6r_k}$. Lemma~\ref{lem:transfer-step}(2) with
$s'=s_h$, the bound $\|w_i(t')\|_{H^{s_h}}\le C_Q\maj_i(t')$ of
Theorem~\ref{thm:quiet}(c) for $i\le k-1$ (its hypotheses hold by
Lemma~\ref{lem:O} if $t'>0$, and $w_i(0)=0$ as these rings start at $\Phi$), and
$\maj_i\le\frac34\Ptil_{i,k-1}\maj_{k-1}$ (Lemma~\ref{lem:O}) give
$\etout_k(t')\le \Cvel C_Q[b_{k-1,k}+\frac34\sum_{i<k-1}b_{i,k}\Ptil_{i,k-1}]\maj_{k-1}(t')$.
As $b_{i,k}=b_{k-1,k}\,b_{i,k-1}$ and Equation~\eqref{eq:O-product} gives
$b_{i,k-1}\Ptil_{i,k-1}\le4^{-(k-1-i)}$, the bracket is at most
$b_{k-1,k}(1+\frac34\cdot\frac13)=\frac54b_{k-1,k}$, and
\begin{equation}
\label{eq:Op-red}
\etout_k(t')\le2\Cvel C_Q\,b_{k-1,k}\,\maj_{k-1}(t') .
\end{equation}

\emph{Bound on $\maj_{k-1}$.} Proposition~\ref{prop:relay} applies with $j=k-1$,
$\bar t=t$ and the present $C_0$, $K_1$: the clock is defined up to level
$k$, and $A_{k-1}t\le A_{k-1}\tsh_{k-1}$ as $t\le T_k<\tsh_{k-1}$; (H1) is assumed; (H2)
follows from $\etout_{k-1}\le\theta_{k-1}\maj_{k-1}$ on $[0,t]$
(Lemma~\ref{lem:O}) and $\theta_{k-1}A_{k-1}t\le\frac12$; (H3) is assumed, through physical time $\min(t,\tsh_k)=\tsh_k$; for (H4), let $t>\tsh_k$: the bound
$\maj_k(\tsh_k)\le K_1\ksh_k$ is assumed if $k<n$ and holds as
$\maj_n(\tsh_n)=\maj_n(0)=0$ if $k=n$, and on $[\tsh_k,\min(t,T_k)]=[\tsh_k,t]$
the inclusion $D_{k+1}\subset B_{ar_k/8}$ and the two norm inequalities of
Equation~\eqref{eq:transmitted-active-data} are assumed. So Equation~\eqref{eq:relay-absorb} gives
$\maj_{k-1}\le2\maj^{\mathrm{in}}_{k-1}$ on
$[0,t]$, and, if $k<n$, (o) holds at physical time $\tsh_k$.

We split $\maj^{\mathrm{in}}_{k-1}$ at $\tau^{\mathrm{act}}$. For
$\tau'\in[\tau^{\mathrm{act}},A_{k-1}t]$,
$\maj^{\mathrm{in}}_{k-1}(\frac{\tau'}{A_{k-1}})=e^{\sigma_{k-1}(\tau'-\tau^{\mathrm{act}})}\maj^{\mathrm{in}}_{k-1}(\tsh_k)
+\int_{\tau^{\mathrm{act}}}^{\tau'}e^{\sigma_{k-1}(\tau'-v)}\etin_{k-1}(\frac{v}{A_{k-1}})\dd v$, and $e^{\sigma_{k-1}(\tau'-\tau^{\mathrm{act}})}\le e^{\sigma_ku}$, as
$\sigma_{k-1}\le\mu\sigma_k$ (Lemma~\ref{lem:detuning}); recall
$\ksh_ke^{\sigma_ku}=\hat a_k(t')$. The first term vanishes if $k=n$, as then
$\tau^{\mathrm{act}}=0$; if $k<n$, (o) at physical time $\tsh_k$,
$\mu\tau^{\mathrm{act}}=A_k\tsh_k$ and $\maj_k(\tsh_k)\le K_1\ksh_k$ bound it by
$\Cvel \big(\frac{2C_Q}{\sigb(\mu-1)}+1\big)K_1\varepsilon_{k-1}\hat a_k(t')$.

For the second term, let $v\in[\tau^{\mathrm{act}},\tau']$ and
$u_v=\mu(v-\tau^{\mathrm{act}})$. At physical time $v/A_{k-1}\in[\tsh_k,t]$, ring $k$ lies in
$N_k$ by (Loc)(1), $h_{k-1}$ is defined by (Loc)(2) for the quiet ring $k-1$, and
$h_k$ is defined, as shown above; so Lemma~\ref{lem:transfer-step}(1) with $s'=s$
and the hypotheses on ring $k$ give
$\etin_{k-1}(\frac{v}{A_{k-1}})\le \Cvel \varepsilon_{k-1}(2C_2\hat a_k(v/A_{k-1})+2C_2\hat a_k(v/A_{k-1})\rhohat_k(v/A_{k-1})+C_0\ksh_k\jap{u_v}^3)$.
As $\hat a_k(v/A_{k-1})=\ksh_ke^{\mu\sigma_k(v-\tau^{\mathrm{act}})}$ and
$\mu\sigma_k-\sigma_{k-1}\ge\frac12\sigb(\mu-1)$, the first part contributes at most
$4C_2\Cvel \varepsilon_{k-1}\hat a_k(t')/(\sigb(\mu-1))$. For the other two we write
$e^{\sigma_{k-1}(\tau'-v)}=e^{\sigma_{k-1}(\tau'-\tau^{\mathrm{act}})}e^{-\sigma_{k-1}v'}$
with $v'=v-\tau^{\mathrm{act}}$, so $u_v=\mu v'$, and use
$\jap{\mu v'}\le\mu(1+v')$ and $\mu^3\le2$. By
Equation~\eqref{eq:active-transient},
$\hat a_k(v/A_{k-1})\rhohat_k(v/A_{k-1})\le C(C_0)\ksh_k\jap{\mu v'}^3e^{\mu\sigma'v'}$, and
$\sigma_{k-1}-\mu\sigma'>\frac{15}{16}\sigma-\frac{17}{16}\cdot\frac34\sigma=\frac{9\sigma}{64}$
by Lemma~\ref{lem:detuning} and (SA1). Hence, by Equation~(I.\ref{I-eq:Ip}), the
transient part contributes at most
$2C_2C(C_0)\Cvel \varepsilon_{k-1}\ksh_ke^{\sigma_ku}\mu^3\int_0^\infty(1+v')^3e^{-\frac{9\sigma}{64}v'}\dd v'
\le4C_2C(C_0)\Ip{3}{\frac{9\sigma}{64}}\Cvel \varepsilon_{k-1}\hat a_k(t')$,
and, as $\sigma_{k-1}\ge\frac{9\sigma}{64}$, $\Ip{3}{\cdot}$ is decreasing and the inner part contributes at most
$2C_0\Ip{3}{\frac{9\sigma}{64}}\Cvel \varepsilon_{k-1}\hat a_k(t')$.

Using $\sigb(\mu-1)\le\sigb$ and choosing $\Crelay(C_0,K_1)$ sufficiently large, the coefficients of
$\Cvel \varepsilon_{k-1}\hat a_k(t')$ found above add up to at most
$\Crelay/(4\Cvel^2 C_Q\sigb(\mu-1))$; so, for $\tau'\in[\tau^{\mathrm{act}},A_{k-1}t]$,
\[
\maj_{k-1}(t')\le2\maj^{\mathrm{in}}_{k-1}(\frac{\tau'}{A_{k-1}})\le\frac{\Crelay}{2\Cvel C_Q}\cdot\frac{\varepsilon_{k-1}\hat a_k(t')}{\sigb(\mu-1)} .
\]
Finally $b_{k-1,k}\varepsilon_{k-1}=\mu^{-1}s_{k-1}\cdot\mu s_{k-1}^5=s_{k-1}^6$,
so Equation~\eqref{eq:Op-red} gives
\[
\etout_k(t')\le\Crelay s_{k-1}^6\hat a_k(t')/(\sigb(\mu-1))=\eta_{b,k}\hat a_k(t')
\]
for $u=A_k(t'-\tsh_k)$; this is Equation~\eqref{eq:etab}, since
$\|\hout_k(t')\|_{C^{M_h}(N')}=\etout_k(t')$.
Also $\eta_{b,k}\le\eta_b$ as $s_{k-1}\le\sbar$. Every bound used at $t'$ involves
the solution on $[0,t']$ only.
\end{proof}

\subsection{The coefficient after activation}
\label{sec:relay-coefficient}

The coefficient $c_k$ has its free linear growth removed. Its exponential
weight suppresses the polynomial growth of $q_{k+1}$ after $T_{k+1}$.
The next lemma bounds $c_k(t')-c_k(t)$ using $q_{k+1}$ and $\hout_k$.
In item (ii), $q'$ is the amplitude supplied by Lemma~\ref{lem:kicktrue}.

\begin{lemma}
\label{lem:drift}
Assume (SA) and (Det). Let $1\le k\le n$, assume that the clock is defined up to
level $k$, and let $\bar t\in[\tsh_k,T_k]$ and $\eta\ge0$. Assume:
\begin{enumerate}[label=\textup{(\roman*)}]
\item if $k<n$, with $\Cmom$ of Proposition~\ref{prop:confine},
  \begin{equation}
  \label{eq:drift-moment}
  |q_{k+1}(t)|\le\Cmom\big(1+A_{k+1}(t-T_{k+1})\big)^3
  \qquad(\tsh_k\le t\le\bar t);
  \end{equation}
\item if $k<n$: the hypotheses of Lemma~\ref{lem:coupling}\textup{(2)} hold,
  and there is $q'\in\C\setminus\{0\}$ such that
  \begin{equation}
  \label{eq:drift-startup}
  \begin{gathered}
  \ksh_k\ge\tfrac12e^{x_k\Tpin}|\Ain(\Tpin)|,\qquad
  |\Ain(\Tpin)|\ge\frac{|\beta_k||q'|}{10|\lamref|(\mu-1)},\\
  \mu-1\le\frac{|q'|}{10|\lamref|\Cmom\Ip{3}{\sigma/4}};
  \end{gathered}
  \end{equation}
\item
  \begin{equation}
  \label{eq:drift-outer}
  \|\hout_k(t)\|_{C^0(N')}\le\eta\,\hat a_k(t)
  \qquad(\tsh_k\le t\le\bar t);
  \end{equation}
\item for every $t\in[\tsh_k,\bar t]$,
  \begin{equation}
  \label{eq:drift-support}
  D_{k+1}(t)\subset B_{ar_k/8},\qquad
  S_i(t)\cap B_{6r_k}=\emptyset\quad(i<k),
  \end{equation}
  so that $h_k$ and $c_k$ are defined there.
\end{enumerate}
Then for $\tsh_k\le t\le t'\le\bar t$,
\begin{equation}
\label{eq:drift-interval}
|c_k(t')-c_k(t)|\le2\rT\big(\Tpin+A_{k+1}(t-\tsh_k)\big)\,\ksh_k
+\Cls\,\eta\,A_k(t'-t)\,\ksh_k,
\end{equation}
without the first term when $k=n$. If
\begin{equation}
\label{eq:drift-smallness}
\rT(\Tpin)\le\frac1{64},\qquad \Cls\eta\hat T_k\le\frac1{32},
\end{equation}
then hypothesis \textup{(Coef)} of
Theorem~\ref{thm:active} holds on $[0,A_k(\bar t-\tsh_k)]$ for the pulled-back
functions with $G=G_k$ and $\kappa=\ksh_k$. If $\bar t=T_k$, then
\begin{equation}
\label{eq:drift-saturation}
|c_k(T_k)|e^{\sigma_k\hat T_k}\in[\tfrac{15}{16}\delta_0,\tfrac{17}{16}\delta_0].
\end{equation}
\end{lemma}

\begin{proof}
Write $\kappa=\ksh_k$, $\lambda=\lambda_k$, $t(v)=\tsh_k+v/A_k$,
$u=A_k(t-\tsh_k)$ and $u'=A_k(t'-\tsh_k)$. By (iv), the fields are defined
on the required interval. Changing variables in Equation~\eqref{eq:running-k}
shows that $c_k(t')-c_k(t)$ is the sum of
\[
\int_u^{u'}e^{-\lambda v}\ellst(G_k)
\big(-h^\bullet_k(t(v))\cdot\nabla\Phi\big)\dd v,
\qquad \bullet\in\{\mathrm{in},\mathrm{out}\}.
\]
By Lemma~I.\ref{I-lem:uniform} and Equation~\eqref{eq:drift-outer},
$|\ellst(G_k)(-\hout_k(t(v))\cdot\nabla\Phi)|
\le\Cls\eta\kappa e^{\sigma_kv}$.
The modulus of the integral with $\bullet=\mathrm{out}$ is therefore at most
$\Cls\eta\kappa(u'-u)=\Cls\eta\kappa A_k(t'-t)$.

If $k=n$, then $D_{n+1}=\emptyset$ and $G_n=0$, so $\hin_n=0$.
For $k<n$, Lemma~\ref{lem:coupling}\textup{(1)} gives
$\ellst(G_k)(-\hin_k(t(v))\cdot\nabla\Phi)=\mu\beta_kq_{k+1}(t(v))$.
With $\zeta_0=\lambda/\mu$ and $\tau_a=\Tpin+\mu u$, the substitution
$\tauK=\Tpin+\mu v=A_{k+1}(t(v)-T_{k+1})$ gives the integral with
$\bullet=\mathrm{in}$ as
\[
e^{\zeta_0\Tpin}\beta_k\int_{\tau_a}^{\Tpin+\mu u'}
e^{-\zeta_0\tauK}q_{k+1}\big(T_{k+1}+\tauK/A_{k+1}\big)\dd\tauK.
\]
By Equation~\eqref{eq:drift-moment}, its modulus is at most
\begin{align*}
e^{x_k\Tpin}|\beta_k|\Cmom\int_{\tau_a}^\infty e^{-x_k\tauK}(1+\tauK)^3\dd\tauK
&\le e^{x_k\Tpin}|\beta_k|\Cmom(1+\tau_a)^3e^{-x_k\tau_a}\Ip{3}{x_k}\\
&\le e^{x_k\Tpin}|\beta_k|\Cmom\Ip{3}{\sigma/4}\,\rT(\tau_a).
\end{align*}
Here the first inequality uses $1+\tau_a+y\le(1+\tau_a)(1+y)$, and the second
uses Equation~\eqref{eq:detuning-rates} and the monotonicity of $I_3$.
By Equation~\eqref{eq:drift-startup},
\[
|\beta_k|\Cmom\Ip{3}{\sigma/4}\le\frac{|\beta_k||q'|}{10|\lamref|(\mu-1)}
\le|\Ain(\Tpin)|\le2\kappa e^{-x_k\Tpin}.
\]
The modulus of the integral with $\bullet=\mathrm{in}$ is at most $2\kappa\rT(\tau_a)$, proving
Equation~\eqref{eq:drift-interval}.

Apply Equation~\eqref{eq:drift-interval} to the pair $(\tsh_k,t)$, with $t\le\bar t\le T_k$.
Equation~\eqref{eq:drift-smallness} gives
$|c_k(t)-c_k(\tsh_k)|\le(2/64+1/32)\kappa=\kappa/16$.
Since $|c_k(\tsh_k)|=\kappa$, we obtain $|c_k(t)|\le17\kappa/16\le2\kappa$,
which gives (Coef) for $u\mapsto c_k(t(u))$.
Finally, at $t=T_k$ this gives
$|c_k(T_k)|\in[15\kappa/16,17\kappa/16]$, and
$\kappa e^{\sigma_k\hat T_k}=\delta_0$ by Definition~\ref{def:clock}.
\end{proof}

\subsection{The actual coefficient at the hand-off}
\label{sec:relay-handoff}

The inner linear response $\cin_k$ determines the activation amplitude once
the outer response and the quadratic quiet-ring error are small relative to
it. The next lemma also supplies the support and Sobolev bounds needed to
apply Theorem~\ref{thm:active}.

\begin{lemma}
\label{lem:handoff}
Assume (SA) and (Det), let $K_1>0$, let $1\le k<n$, and assume that the clock is
defined up to level $k+1$. Assume:
\begin{enumerate}[label=\textup{(\roman*)}]
\item for every $t\in[0,\tsh_k]$, $D_{k+1}(t)\subset B_{ar_k/8}$ and the rings
  $i<k$ lie outside $B_{6r_k}$, so that $\eta_k$ is defined on $[0,\tsh_k]$
  (Definition~\ref{def:regions}); and the hypotheses of Theorem~\ref{thm:quiet}
  hold for ring $k$ through physical time $\tsh_k$, using the pulled-back functions;
\item for some $\theta_k\ge0$, $\etout_k\le\theta_k\maj_k$ on $[0,\tsh_k]$ and
  $\Cls\theta_kA_k\tsh_kK_1\le\frac14$;
\item $\maj_k(\tsh_k)\le\frac12K_1|\cin_k|$, with $\cin_k$ of Lemma~\ref{lem:coupling};
\item $\|\ellst\|C_QK_1\maj_k(\tsh_k)\le\frac12$.
\end{enumerate}
Then, in ring-$k$ units,
\begin{equation}
\label{eq:handoff}
\begin{gathered}
\tfrac12|\cin_k|\le\ksh_k\le\tfrac32|\cin_k|,\qquad \maj_k(\tsh_k)\le K_1\ksh_k,\\
\|w_k(\tsh_k)\|_{H^{s_h}}\le K_1C_Q\ksh_k,\qquad \supp(\Phi+w_k(\tsh_k))\subset\NNm .
\end{gathered}
\end{equation}
\end{lemma}

\begin{proof}
Write $\maj=\maj_k$ and $c=\cin_k$. By (i) and
Theorem~\ref{thm:quiet}(a), for $t\in[0,\tsh_k]$, $D_{k+1}(t)\subset B_{ar_k/8}$,
the rings $i<k$ lie outside $B_{6r_k}$ and $S_k(t)\subset N_k$. So the hypotheses
of Lemma~\ref{lem:coupling}\textup{(2)} and \textup{(3)} hold, $c$ is defined,
and by Lemma~\ref{lem:coupling}\textup{(3)} the physical-time linear response corresponding to
Theorem~\ref{thm:quiet} for ring $k$ is
$\wlin=\wlin_{\mathrm{in},k}+\wlin_{\mathrm{out},k}$; at physical time $\tsh_k$,
\[
\ellst(G_k)\big(w_k(\tsh_k)\big)=c+\ellst(G_k)\big(\wlin_{\mathrm{out},k}(\tsh_k)\big)+\ellst(G_k)\big(w_k(\tsh_k)-\wlin(\tsh_k)\big),
\]
and $\ksh_k$ is the modulus of the left side (Definition~\ref{def:clock}). As
$G_k\in[0,\Gone]$ (Lemma~\ref{lem:inner-circ} and (SA3)), $\ellst(G_k)$ is bounded on $X^0$ (Lemma~I.\ref{I-lem:uniform}), so
it commutes with the integral defining $\wlin_{\mathrm{out},k}$, and
Proposition~I.\ref{I-prop:semigroup} gives
$\ellst(G_k)(\wlin_{\mathrm{out},k}(\tsh_k))=A_k\int_0^{\tsh_k}e^{A_k\lambda_k(\tsh_k-v)}\ellst(G_k)(-\hout_k(v)\cdot\nabla\Phi)\dd v$.
By Lemma~I.\ref{I-lem:uniform}, $\|\hout_k\|_{C^0(N')}\le\etout_k$, (ii), (M1) of
Lemma~\ref{lem:majorant} and (iii),
\[
\big|\ellst(G_k)\big(\wlin_{\mathrm{out},k}(\tsh_k)\big)\big|
\le\Cls\theta_kA_k\int_0^{\tsh_k}e^{A_k\sigma_k(\tsh_k-v)}\maj(v)\dd v
\le\Cls\theta_kA_k\tsh_k\maj(\tsh_k)\le\frac1{4K_1}\cdot\frac{K_1|c|}2=\frac{|c|}8 .
\]
Since $\ellst(G_k)$ has norm at most $\|\ellst\|$ on $X^0$
(Lemma~I.\ref{I-lem:uniform}) and $\|\cdot\|_{L^2}\le\|\cdot\|_{H^{s_h}}$,
Theorem~\ref{thm:quiet}(b), (iv) and (iii) give
\[
\big|\ellst(G_k)\big(w_k(\tsh_k)-\wlin(\tsh_k)\big)\big|\le\|\ellst\|C_Q\maj(\tsh_k)^2
\le\frac{\maj(\tsh_k)}{2K_1}\le\frac{|c|}4 .
\]
Hence $|\ksh_k-|c||\le\frac38|c|$, which gives the first claim. Then (iii) and
$|c|\le2\ksh_k$ give $\maj(\tsh_k)\le K_1\ksh_k$, and
Theorem~\ref{thm:quiet}(c) and (a) (the ring part is $\Phi+w_k(\tsh_k)$) give
$\|w_k(\tsh_k)\|_{H^{s_h}}\le C_Q\maj(\tsh_k)\le K_1C_Q\ksh_k$ and
$\supp(\Phi+w_k(\tsh_k))\subset\NNm$.
\end{proof}

\subsection{Constants and error exponents}
\label{sec:relay-constants}

The four error exponents below are used only to express the smallness
conditions in Lemma~\ref{lem:kicktrue} and Section~\ref{sec:param}. Under
(SA), put
\[
\begin{gathered}
c_\gamma=1-\frac{16\sigma'}{15\sigma},\qquad
\vartheta_\mu=\frac{\mu-1}3,\qquad
\gamma_\mu=\min(\vartheta_\mu,c_\gamma),\\
\gamma'_\mu=\min\Big(\frac{13}{20}\cdot\frac{\mu-1}\mu,c_\gamma\Big).
\end{gathered}
\]
Under (Det), Lemma~\ref{lem:detuning} shows that $c_\gamma$ bounds
$1-\sigma'/\sigma_k$ from below.
The assumptions on $\sigma'$ and $\mu$ make all four numbers positive.

\begin{definition}
\label{def:K1}
Under (SA), with $c_N,C_N$ of Lemma~I.\ref{I-lem:uniform}, put
\begin{equation}
\label{eq:K1}
\begin{gathered}
K_1=\frac{256\mu \Cvel |\lamref|}{3c_N}\Big[\frac{4C_2+1}{\sigb}+\frac{C_D(\Tpin)(\mu-1)}{\delta_0}\Big],\\
\Abar=2\Cvel e^{17\sigma\Tpin/16}\Big[\frac{(4C_2+1)\delta_0}{\sigb(\mu-1)}+C_D(\Tpin)\Big],
\end{gathered}
\end{equation}
with $C_D$ of Proposition~\ref{prop:relay}, and
\[
\Klo=\frac{c_N\delta_0}{22|\lamref|(\mu-1)},\qquad \Kup=\frac{5C_N\delta_0\mu}{|\lamref|(\mu-1)}e^{17\sigma\Tpin/16}.
\]
\end{definition}

\subsection{Comparison of the true moment with the scalar kick}
\label{sec:proof-kicktrue}

We compare the actual incoming moment with the exponentially growing moment
in Proposition~\ref{prop:kick}, separating the earlier quiet history, the
active interval and the polynomial growth after saturation.

\begin{lemma}
\label{lem:kicktrue}
Assume (SA) and (Det), and let $C_0\ge1$ and $K_1>0$.
Choose one constant $\Cerr\ge1$, large enough for all the error bounds below.
It may depend on $C_0$, $K_1$, $\mu$, $\Tpin$, $\Klo$ and the fixed ring data,
but is independent of $\kappa$, $\sbar$, $n$ and the level. In particular, it
also absorbs the fixed prefactor in the estimate of the earlier history.
Let $1\le k<n$, assume that
the clock is defined up to level $k+1$, put $\kappa=\ksh_{k+1}$,
$T_g=\hat T_{k+1}$ and
\begin{equation}
\label{eq:kicktrue-moment}
\qs=2\pi\pout_*(G_{k+1})\,c_{k+1}(T_{k+1})\,e^{\lambda_{k+1}\hat T_{k+1}} .
\end{equation}
Assume the following.
\begin{itemize}
\item $D_{k+1}(t)\subset B_{ar_k/8}$ for every $t\in[0,\tsh_k]$ (the hypotheses
  of Lemma~\ref{lem:coupling}\textup{(2)}).
\item The conclusions of Theorem~\ref{thm:active} and
  Corollary~\ref{cor:active-moment} hold for ring $k+1$ on $\Phase_{k+1}$ with
  $G=G_{k+1}$, $\kappa=\ksh_{k+1}$, the constant $C_0$ and $\eta_b$ of
  Equation~\eqref{eq:etab}.
\item The hypotheses \textup{(i)}--\textup{(iv)} of Lemma~\ref{lem:drift} hold for
  ring $k+1$ with $\bar t=T_{k+1}$ and some $\eta\in[0,\eta_b]$, and
  $\rT(\Tpin)\le\frac1{64}$ and $\Cls\eta\hat T_{k+1}\le\frac1{32}$.
\item If $k+1<n$: for every $t\in[0,\tsh_{k+1}]$, $D_{k+2}(t)\subset B_{ar_{k+1}/8}$
  and the rings $i\le k$ lie outside $B_{6r_{k+1}}$, so that $\eta_{k+1}$ is
  defined on $[0,\tsh_{k+1}]$ (Definition~\ref{def:regions}); the hypotheses of
  Theorem~\ref{thm:quiet} hold for ring $k+1$ through physical time $\tsh_{k+1}$, using the pulled-back functions; $\maj_{k+1}(\tsh_{k+1})\le K_1\kappa$;
  the hypotheses of Proposition~\ref{prop:confine} hold with $k+1$ in place of $k$
  and $t_1=T_{k+1}$; and $\kappa\ge\Klo s_{k+1}^5$.
\item On the pinning interval, including the case $k+1=n$,
  \begin{equation}
  \label{eq:kicktrue-tail}
  |q_{k+1}(t)|\le\Cmom\big(1+A_{k+1}(t-T_{k+1})\big)^3
  \qquad(T_{k+1}\le t\le\tsh_k).
  \end{equation}
\item The following inequalities hold:
\begin{enumerate}[label=\textup{(\alph*)}]
\item $\frac{16}{15}\cdot\frac{2\pi C_3\delta_0^2}{\sigma}\le\frac14\cdot\frac{15}{16}\cdot\frac{\cout\delta_0}{20|\lamref|}$;
\item $\mu-1\le\frac{15}{16}\cdot\cout\delta_0/\big(10|\lamref|\Cmom\Ip{3}{\sigma/4}\big)$;
\item $\kappa\le\delta_04^{-\mu/(\mu-1)}$;
\item \begin{multline*}
  \frac{2\pi C_3\kappa^{1/2}\mu}{\frac34\sigma(\mu-1)}+\Cerr\Big(\frac\kappa{\delta_0}\Big)^{\gamma_\mu}\log^4\frac{\delta_0}\kappa\\
  +\Cerr\big[1+\rT(\Tpin)\Ip{3}{\tfrac\sigma8}\big]\Big(\frac\kappa{\delta_0}\Big)^{\vartheta_\mu}\log^3\frac{\delta_0}\kappa
  \le\frac12\cdot\frac{15}{16}\cdot\frac{\cout}{20|\lamref|(\mu-1)};
  \end{multline*}
\item $\Cerr\big(\frac\kappa{\delta_0}\big)^{\vartheta_\mu}\le\frac{15}{16}\cdot\frac{\cout}{20|\lamref|(\mu-1)}$;
\item $\Cerr\eta_b\delta_0\frac{\mu^2}{\sigma^2(\mu-1)^2}+2\pi C_3\eta_b\frac{4\mu\delta_0}{3\sigma(\mu-1)}
  \le\frac14\cdot\frac{15}{16}\cdot\frac{\cout\delta_0}{20|\lamref|(\mu-1)}$.
\end{enumerate}
\end{itemize}
Then
\begin{equation}
\label{eq:kicktrue-moment-size}
|\qs|\in[\tfrac{15}{16},\tfrac{17}{16}]\,2\pi|\pout_*(G_{k+1})|\delta_0,
\end{equation}
and the
kick amplitude $\Ain$ of Lemma~\ref{lem:coupling}\textup{(2)} satisfies
\begin{equation}
\label{eq:kicktrue-amplitude}
\frac{|\beta_k|\,|\qs|}{10|\lamref|(\mu-1)}\le|\Ain(\Tpin)|\le\frac{3|\beta_k|\,|\qs|\,\mu}{|\lamref|(\mu-1)} .
\end{equation}
\end{lemma}

\begin{proof}[Quantitative proof of Lemma~\ref{lem:kicktrue}]
Recall the local pullback $\tilde q(\tau)=q_{k+1}(T_{k+1}+\tau/A_{k+1})$
from Lemma~\ref{lem:coupling}. We apply Proposition~\ref{prop:kick} to the instance described after
Lemma~\ref{lem:coupling} (an instance of its setting by
Lemma~\ref{lem:coupling}\textup{(2)}): $G=G_k$, $G'=G_{k+1}$, $\beta=\beta_k$,
$q=\tilde q$, $c=\cin$, $A=\Ain$, $\tau_0=-A_{k+1}T_{k+1}$, $\tau_1=\Tpin$ and
$(\bar Q,p)=(\Cmom,3)$. We verify its hypotheses (Q1)--(Q3) and (Kick1)--(Kick4),
and then the remaining ones; (a)--(f) refer to the inequalities in the statement
above. Write $\sigma_+=\sigma_{k+1}$, $x=x_k=\sigma_k/\mu$,
$L=\log(\delta_0/\kappa)$, $\lambda_{\mathrm{in}}=\lambda_{k+1}$, $\zeta_0=\lambda_k/\mu$ and
$w=\lambda_{\mathrm{in}}-\zeta_0$, so that $\Re w=\sigma_+-x$ and $T_g=L/\sigma_+$, and let
$\mathsf B=\frac{15}{16}\cout\delta_0/(20|\lamref|(\mu-1))$.

By Lemma~\ref{lem:detuning} and its proof,
$\frac{15}{16}\sigma\le\sigma_+\le\frac{17}{16}\sigma$, $x\ge\frac\sigma4$,
$\mu\sigma_+-\sigma_k\ge\frac34\sigma(\mu-1)$, $e_\mu(\mu+1)\le\frac14(\mu-1)$ and
$|\sigma_i-\sigma|\le e_\mu\sigma$ for $i\in\{k,k+1\}$, so
$\mu\sigma_+-\sigma_k\le\mu(1+e_\mu)\sigma-(1-e_\mu)\sigma=\sigma((\mu-1)+e_\mu(\mu+1))\le\frac54\sigma(\mu-1)$; as
$\mu-1\le\frac1{16}$,
\begin{equation}
\label{eq:Rew}
\frac34\sigma\,\frac{\mu-1}\mu\le\Re w=\frac{\mu\sigma_+-\sigma_k}\mu\le\frac54\sigma(\mu-1)\le\frac\sigma8 .
\end{equation}
By Equations~\eqref{eq:detuning-gap} and~\eqref{eq:detuning-rates}, and $\mu\le\frac{36}{17}$,
$1-x/\sigma_+=\Re w/\sigma_+=(\mu\sigma_+-\sigma_k)/(\mu\sigma_+)\ge\frac{12}{17}\cdot\frac{\mu-1}\mu\ge\frac{\mu-1}3=\vartheta_\mu$;
also $1-\sigma'/\sigma_+\ge1-\frac{16\sigma'}{15\sigma}=c_\gamma$, and
$\gamma_\mu=\min(\vartheta_\mu,c_\gamma)$. By (c), $L\ge\frac\mu{\mu-1}\log4\ge1$;
and $T_g\le\frac{16}{15\sigma}L$. As $e^{\sigma_+T_g}=\delta_0/\kappa$,
$\kappa e^{yT_g}=\delta_0(\kappa/\delta_0)^{1-y/\sigma_+}$ for every real $y$.

\emph{The amplitude $\qs$ and (Kick1).} The assumed hypotheses
\textup{(i)}--\textup{(iv)} and smallness conditions of Lemma~\ref{lem:drift}
give Equation~\eqref{eq:drift-saturation} for ring $k+1$:
$|c_{k+1}(T_{k+1})|e^{\sigma_+\hat T_{k+1}}\in[\frac{15}{16}\delta_0,\frac{17}{16}\delta_0]$;
as $|e^{\lambda_{k+1}\hat T_{k+1}}|=e^{\sigma_+\hat T_{k+1}}$, this is the claimed
range of $|\qs|$. As
$G_{k+1}\in[0,\Gone]$ (Lemma~\ref{lem:inner-circ} and (SA3)), $2\pi|\pout_*(G_{k+1})|\ge\cout$ (Lemma~I.\ref{I-lem:uniform}),
so $|\qs|\ge\frac{15}{16}\cout\delta_0>0$, hence $\qs\ne0$, and $\mathsf B\le|\qs|/(20|\lamref|(\mu-1))$. Then
(b) gives $\mu-1\le|\qs|/(10|\lamref|\Cmom\Ip{3}{\sigma/4})$, which is (Kick1).

\emph{(Kick2) and (Q3).} By Equation~\eqref{eq:Rew}, $\sigma_+\le\frac{17}{16}\sigma$
and (c),
$\Re w\,T_g\ge\frac34\cdot\frac{16}{17}\cdot\frac{\mu-1}\mu\,L\ge\frac{12}{17}\log4\ge\log2$,
so $e^{-\Re w\,T_g}\le\frac12$. (Q3) on $[0,\Tpin]$ is
Equation~\eqref{eq:kicktrue-tail}, expressed in the local coordinate $\tauK$.

\emph{(Q1).} Let $\tau\in[-T_g,0]$ and $t=T_{k+1}+\tau/A_{k+1}$. By
Definition~\ref{def:clock}, $T_{k+1}-\tsh_{k+1}=T_g/A_{k+1}$, so
$t\in\Phase_{k+1}$, and $u=\tau+T_g=A_{k+1}(t-\tsh_{k+1})\in[0,T_g]$ is the
active time of ring $k+1$.
We split $q_{k+1}$ into the moment of ring $k+1$ and that of the region inside it.
As $D_{k+1}(t)$ is the disjoint union of $S_{k+1}(t)$ and $D_{k+2}(t)$ and $Q_3$ is
linear,
$Q_3[\omega_{\ge k+1}(t)]=Q_3[\omega^{(n)}(t)\mathbf 1_{S_{k+1}(t)}]+Q_3[\omega_{\ge k+2}(t)]$;
the second term vanishes if $k+1=n$, as $D_{n+1}=\emptyset$. Otherwise, dividing by
$A_{k+1}r_{k+1}^5$ and using
$A_{k+2}r_{k+2}^5=\mu s_{k+1}^5A_{k+1}r_{k+1}^5$, the second term
becomes $\mu s_{k+1}^5q_{k+2}(t)$. As $T_{k+1}-T_{k+2}=\Tpin/A_{k+2}+T_g/A_{k+1}$
(Definition~\ref{def:clock}) and $A_{k+2}=\mu A_{k+1}$, the kick time of ring $k+1$ is
$\tau'=A_{k+2}(t-T_{k+2})=\mu A_{k+1}(t-T_{k+1})+A_{k+2}(T_{k+1}-T_{k+2})=\mu\tau+\Tpin+\mu T_g$.
So $\tilde q(\tau)=q^{\mathrm{ring}}(\tau)+q^{\mathrm{deep}}(\tau)$ with
\[
q^{\mathrm{ring}}(\tau)=\frac{\overline{Q_3[\omega^{(n)}(t)\mathbf 1_{S_{k+1}(t)}]}}{A_{k+1}r_{k+1}^5},\qquad
q^{\mathrm{deep}}(\tau)=\mu s_{k+1}^5\,q_{k+2}(t),
\]
where $t=T_{k+1}+\tau/A_{k+1}$, and
$q^{\mathrm{deep}}\equiv0$ if $k+1=n$.
By the assumed conclusions of Theorem~\ref{thm:active}, $w_{k+1}(t)$ is defined and
$\supp(\Phi+w_{k+1}(t))\subset\NN$, and in ring-$(k+1)$ units
$\omega^{(n)}(t)\mathbf 1_{S_{k+1}(t)}=A_{k+1}(\Phi+w_{k+1}(t))$
(Definition~\ref{def:regions}). As $\Phi$ is radial and $w_{k+1}(t)$ is real,
$q^{\mathrm{ring}}(\tau)=\overline{Q_3[w_{k+1}(t)]}=2\pi\pout_3[w_{k+1}(t)]$, which
is the function $q(u)$ of Corollary~\ref{cor:active-moment} for ring $k+1$, as
$w(u)=w_{k+1}(t)$. As $\qs e^{\lambda_{\mathrm{in}}\tau}=2\pi\pout_*(G_{k+1})c_{k+1}(T_{k+1})e^{\lambda_{\mathrm{in}}u}$,
\begin{align*}
\tilde q(\tau)-\qs e^{\lambda_{\mathrm{in}}\tau}&=\big[q^{\mathrm{ring}}(\tau)-2\pi\pout_*(G_{k+1})c_{k+1}(t)e^{\lambda_{\mathrm{in}}u}\big]\\
&\quad+2\pi\pout_*(G_{k+1})\big(c_{k+1}(t)-c_{k+1}(T_{k+1})\big)e^{\lambda_{\mathrm{in}}u}+q^{\mathrm{deep}}(\tau).
\end{align*}
Equation~\eqref{eq:active-moment} bounds the bracket, and
Equation~\eqref{eq:drift-interval} for ring $k+1$, with the physical-time pair $(t,T_{k+1})$ and
$\eta=\eta_b$, bounds $|c_{k+1}(t)-c_{k+1}(T_{k+1})|$
(hypothesis (iii) of Lemma~\ref{lem:drift} with
$\eta\le\eta_b$ implies (iii) with $\eta_b$, and
Equation~\eqref{eq:drift-interval} does not use the bound on
$\Cls\eta\hat T_{k+1}$). With $\hat a(u)=\kappa e^{\sigma_+u}$,
$\rhohat(u)=\rhohat_{k+1}(\tsh_{k+1}+u/A_{k+1})$ and $P=\max_{G\in[0,\Gone]}|\pout_*(G)|$, (Q1) holds with
\begin{equation}
\label{eq:kicktrue-error}
\begin{aligned}
\err(\tau)&=2\pi C_3\hat a(u)\big[\rhohat(u)+\hat a(u)+\kappa^{1/2}+\eta_b\big]\\
&\quad+2\pi P\big(2\kappa\,\rT(\Tpin+\mu u)\,\mathbf 1_{k+1<n}+\Cls\eta_b\kappa(T_g-u)\big)e^{\sigma_+u}
+|q^{\mathrm{deep}}(\tau)| .
\end{aligned}
\end{equation}

\emph{(Kick3).} We bound
$\Errg=\int_{-T_g}^0e^{-x\tau}\err(\tau)\dd\tau=\int_0^{T_g}e^{x(T_g-u)}\err(u-T_g)\dd u$
term by term in Equation~\eqref{eq:kicktrue-error}.
As $x<\sigma_+$, $2\sigma_+-x\ge\sigma_+\ge\frac{15}{16}\sigma$, so by
(a) and $\mu-1\le1$ the weighted integral of $2\pi C_3\hat a(u)^2$ is
$2\pi C_3\kappa^2e^{xT_g}\int_0^{T_g}e^{(2\sigma_+-x)u}\dd u\le\frac{2\pi C_3\kappa^2e^{2\sigma_+T_g}}{2\sigma_+-x}
\le\frac{16}{15}\cdot\frac{2\pi C_3\delta_0^2}{\sigma}\le\frac{\mathsf B}4$.
In the same way, by Equation~\eqref{eq:Rew}, the terms containing
$\kappa^{1/2}$ and $\eta_b$ in the first line of Equation~\eqref{eq:kicktrue-error}
contribute at most
$2\pi C_3\kappa^{1/2}\kappa e^{\sigma_+T_g}/\Re w\le\delta_0\cdot2\pi C_3\kappa^{1/2}\mu/(\frac34\sigma(\mu-1))$
and $2\pi C_3\eta_b\delta_0/\Re w\le2\pi C_3\eta_b\cdot4\mu\delta_0/(3\sigma(\mu-1))$.
For $2\pi C_3\hat a(u)\rhohat(u)$, Equation~\eqref{eq:active-transient} gives
$\hat a(u)\rhohat(u)\le C(C_0)\kappa\jap{u}^3e^{\sigma'u}$; as
$x(T_g-u)+\sigma'u\le\max(x,\sigma')T_g$ on $[0,T_g]$,
$\jap{T_g}\le(1+\frac{16}{15\sigma})L$, $1-\max(x,\sigma')/\sigma_+\ge\gamma_\mu$ and
$\kappa<\delta_0$, it contributes at most
\[
2\pi C_3C(C_0)\,T_g\jap{T_g}^3\kappa e^{\max(x,\sigma')T_g}
\le \Cerr\,\delta_0
\Big(\frac\kappa{\delta_0}\Big)^{\gamma_\mu}L^4,
\]
where the coefficient is bounded by $\Cerr$ by its choice.
The term containing $\rT(\Tpin+\mu u)$ vanishes if $k+1=n$. If $k+1<n$, we use
$\rT(\Tpin+\mu u)\le\rT(\Tpin)(1+\mu u)^3e^{-\sigma\mu u/4}$,
$(1+\mu u)^3\le\mu^3(1+u)^3$, $\mu\le\frac{17}{16}$, $1-x/\sigma_+\ge\vartheta_\mu$ and
$\sigma_+-x-\frac{\mu\sigma}4\le\frac\sigma8-\frac\sigma4=-\frac\sigma8$
(Equation~\eqref{eq:Rew}); it contributes at most
\[
4\pi P\kappa e^{xT_g}\rT(\Tpin)\mu^3\int_0^{\infty}(1+u)^3e^{-\sigma u/8}\dd u
\le \Cerr\,\rT(\Tpin)\,\Ip{3}{\tfrac\sigma8}\,\delta_0\Big(\frac\kappa{\delta_0}\Big)^{\vartheta_\mu}.
\]
For $2\pi P\Cls\eta_b\kappa(T_g-u)e^{\sigma_+u}$, the substitution $u=T_g-v$
and Equation~\eqref{eq:Rew} bound its weighted integral by
$2\pi P\Cls\eta_b\kappa e^{\sigma_+T_g}\allowbreak\int_0^{T_g}v\,e^{-\Re w\,v}\dd v
\le\frac{2\pi P\Cls\eta_b\delta_0}{(\Re w)^2}
\le \Cerr\eta_b\delta_0\frac{\mu^2}{\sigma^2(\mu-1)^2}$.
We have $q^{\mathrm{deep}}\equiv0$ if $k+1=n$. Let $k+1<n$. As
$\tsh_{k+1}=T_{k+2}+\Tpin/A_{k+2}\ge T_{k+2}$, $\Phase_{k+1}$ lies in
$[T_{k+2},T_{k+1}]$, and Proposition~\ref{prop:confine}(b) at level $k+1$,
with $t_0=T_{k+2}$ and $t_1=T_{k+1}$, gives $|q_{k+2}(t)|\le\Cmom(1+\tau')^3$
there; with $\tau'=\Tpin+\mu(\tau+T_g)$,
$|q^{\mathrm{deep}}(\tau)|\le\mu s_{k+1}^5\Cmom\big(1+\Tpin+\mu(\tau+T_g)\big)^3$
on $[-T_g,0]$. As $x\ge\frac\sigma4$,
$\int_{-T_g}^0e^{-x\tau}\dd\tau\le4e^{xT_g}/\sigma$. On $[-T_g,0]$,
$1+\Tpin+\mu(\tau+T_g)\le1+\Tpin+2T_g\le(1+\Tpin+\frac{32}{15\sigma})L$, as
$\mu\le2$, $T_g\le\frac{16}{15\sigma}L$ and $L\ge1$. With
$s_{k+1}^5\le\kappa/\Klo$, $\kappa e^{xT_g}=\delta_0(\kappa/\delta_0)^{1-x/\sigma_+}$,
$1-x/\sigma_+\ge\vartheta_\mu$ and $\kappa<\delta_0$, the weighted integral of
$|q^{\mathrm{deep}}|$ is at most
\[
\frac{4\mu\Cmom}{\sigma\Klo}\Big(1+\Tpin+\frac{32}{15\sigma}\Big)^3L^3\,\kappa e^{xT_g}
\le\Cerr\,\delta_0\Big(\frac\kappa{\delta_0}\Big)^{1-x/\sigma_+}L^3
\le\Cerr\,\delta_0\Big(\frac\kappa{\delta_0}\Big)^{\vartheta_\mu}L^3 .
\]
As $L\ge1$, the terms containing $\kappa^{1/2}$, $\rhohat$, $\rT$ and
$|q^{\mathrm{deep}}|$ in Equation~\eqref{eq:kicktrue-error} contribute at most
$\delta_0$ times the left side of (d), hence at most $\frac{\mathsf B}2$.
The two terms containing $\eta_b$ contribute at most $\frac{\mathsf B}4$ by (f).
Together with the estimate for $2\pi C_3\hat a^2$,
$\Errg\le\frac{\mathsf B}4+\frac{\mathsf B}2+\frac{\mathsf B}4=\mathsf B\le|\qs|/(20|\lamref|(\mu-1))$, which is (Kick3).

\emph{(Q2) and (Kick4) when $k+1<n$.} By Definition~\ref{def:clock},
$\tau_0=-A_{k+1}T_{k+1}=-A_{k+1}\tsh_{k+1}-T_g<-T_g$, as $A_{k+1}\tsh_{k+1}>0$
(Lemma~\ref{lem:clock}(1)). On $[\tau_0,-T_g]$ we take $Q_{\mathrm{pre}}=|\tilde q|$,
so (Q2) holds. In the ring-$(k+1)$ time $v=\tau-\tau_0\in[0,A_{k+1}\tsh_{k+1}]$,
$-x\tau=xT_g+x(A_{k+1}\tsh_{k+1}-v)$, and the time of $\omega^{(n)}$ is $v/A_{k+1}$, the
argument of $q_{k+1}$ in Definition~\ref{def:regions}.

To bound $\Errpre$ we bound $q_{k+1}$ on the quiet interval of ring $k+1$ by the
perturbation and the inner field of that ring.
Put $j=k+1$, let $t\in[0,\tsh_j]$ and work in ring-$j$ units. As
$D_{j+1}(t)\subset B_{ar_j/8}$, $\hin_j(t)$ is defined, and
Theorem~\ref{thm:quiet}(a) gives $S_j(t)\subset N_j$, so $w_j(t)$ is defined. As
$D_j(t)$ is the disjoint union of $S_j(t)$ and $D_{j+1}(t)$,
$Q_3[\omega_{\ge j}(t)]=A_jr_j^5(Q_3[\Phi+w_j(t)]+Q_3[g])$ with
$g(y)=A_j^{-1}\omega_{\ge j+1}(t,r_jy)$. Here $Q_3[\Phi]=0$ as $\Phi$ is radial,
and $\Phi$ and $\Phi+w_j(t)$ are supported in $\NNm\subset N\subset\overline B_2$
(Lemmas~I.\ref{I-lem:ring-facts} and~\ref{lem:annuli} and Theorem~\ref{thm:quiet}(a)), so
$|Q_3[w_j(t)]|\le8\|w_j(t)\|_{L^1}$. The real function $g$ is
$\nu^{\mathrm{in}}_j$ at physical time $t$ (Definition~\ref{def:regions}); put
$p_3=\pout_3[g]=\frac1{2\pi}\int g(y')\,\bar y'^3\dd y'$
(Definition~I.\ref{I-def:multipole}). As $D_{j+1}(t)\subset B_{ar_j/8}$,
Lemma~\ref{lem:coupling}\textup{(1)} with $j$ in place of $k$ shows that the
mode-$3$ profile of $\hin_j(t)\cdot\er$ is $\frac i2p_3r^{-4}$ for
$a/4\le r\le3$; at $r=1$,
$\frac1{2\pi}\int_0^{2\pi}\hin_j(t)(e^{i\theta})\cdot\er\,e^{-3i\theta}\dd\theta=\frac i2p_3$,
so $|p_3|\le2\max_{|y|=1}|\hin_j(t)(y)|\le2\etin_j(t)$, the unit circle
lying in $N'$. As $g$ is real, $\overline{Q_3[g]}=\int g(y)\,\bar y^3\dd y=2\pi p_3$,
so $|Q_3[g]|\le4\pi\etin_j(t)$. Dividing by $A_jr_j^5$ gives
\begin{equation}
\label{eq:kicktrue-quiet-moment}
|q_j(t)|\le C\big(\|w_j(t)\|_{L^1}+\etin_j(t)\big),
\end{equation}
with $C$ depending only on the fixed ring data.

For $v\in[0,A_j\tsh_j]$, $w_j(v/A_j)$ is supported in $N$, so the
Cauchy-Schwarz inequality, Theorem~\ref{thm:quiet}(c) and (M1) of
Lemma~\ref{lem:majorant} give, using $|N|^{1/2}\le4$,
\begin{equation}
\label{eq:kicktrue-quiet-L1}
\begin{aligned}
\|w_j(v/A_j)\|_{L^1}
&\le 4\|w_j(v/A_j)\|_{H^{s_h}}\le 4C_Q\maj_j(\tfrac{v}{A_j})\\
&\le 4C_Qe^{-\sigma_+(A_j\tsh_j-v)}\maj_j(\tsh_j).
\end{aligned}
\end{equation}
Since $\int_0^{A_j\tsh_j}e^{-(\sigma_+-x)(A_j\tsh_j-v)}\dd v\le(\sigma_+-x)^{-1}$, the
part of $\int_0^{A_j\tsh_j}e^{x(A_j\tsh_j-v)}|q_j(v/A_j)|\dd v$ coming from
$\|w_j\|_{L^1}$ is at most $4CC_Q\maj_j(\tsh_j)/(\sigma_+-x)$.
As $x<\sigma_+$ and $\etin_j\le\eta_j$,
\begin{equation}
\label{eq:kicktrue-quiet-field}
\int_0^{A_j\tsh_j}e^{x(A_j\tsh_j-v)}\etin_j(\frac{v}{A_j})\dd v\le\maj_j(\tsh_j).
\end{equation}
Equations~\eqref{eq:kicktrue-quiet-moment}--\eqref{eq:kicktrue-quiet-field},
the assumed inequality $\maj_{k+1}(\tsh_{k+1})\le K_1\kappa$,
Equation~\eqref{eq:Rew} for $\sigma_+-x=\Re w$, and $1-x/\sigma_+\ge\vartheta_\mu$ give
\[
\begin{aligned}
\Errpre&=e^{xT_g}\int_0^{A_{k+1}\tsh_{k+1}}e^{x(A_{k+1}\tsh_{k+1}-v)}|q_{k+1}(v/A_{k+1})|\dd v\\
&\le C\Big(\frac{16\mu C_Q}{3\sigma(\mu-1)}+1\Big)K_1\,\delta_0\Big(\frac\kappa{\delta_0}\Big)^{\vartheta_\mu}
\le\Cerr\,\delta_0\Big(\frac\kappa{\delta_0}\Big)^{\vartheta_\mu},
\end{aligned}
\]
which is at most $\mathsf B\le|\qs|/(20|\lamref|(\mu-1))$ by (e); this is (Kick4). The
remaining hypotheses of Proposition~\ref{prop:kick} hold: $G_k,G_{k+1}\in[0,\Gone]$
(Lemma~\ref{lem:inner-circ} and (SA3)),
$|\lamst(G_i)-\lamref|\le\sigma(\mu-1)/(8\mu)$ for $i\in\{k,k+1\}$ by (Det),
$1<\mu\le2$, $T_g>0$ and $\tau_1=\Tpin>0$. So Proposition~\ref{prop:kick} applies,
and its conclusion at $\tau_1$ is the claimed bound on $|\Ain|$.

\emph{The seeded level $k+1=n$.} Here $\tsh_n=0$, so $\tau_0=-T_g$, while
Proposition~\ref{prop:kick} asks for $\tau_0<-T_g$. We apply it with $T_g$
replaced by some $T'_g\in[\frac{17}{24}T_g,T_g)$ and the same $\qs$, $\lambda_{\mathrm{in}}$ and
$\err$. (Q1) and $\Errg\le\mathsf B$ hold on the shorter interval $[-T'_g,0]$, so
(Kick3) holds; (Kick1) and (Q3) do not involve $T_g$; and
$\Re w\,T'_g\ge\frac{17}{24}\Re w\,T_g\ge\log2$, because the computation for (Kick2)
gives $\Re w\,T_g\ge\frac{12}{17}\log4=\frac{24}{17}\log2$.
On $[\tau_0,-T'_g]=[-T_g,-T'_g]$ we take $Q_{\mathrm{pre}}=|\tilde q|$; then (Q2) holds
and $\Errpre\le(T_g-T'_g)e^{xT_g}\max_{[-T_g,0]}|\tilde q|$, as $\tilde q$ is continuous.
This bound tends to $0$ as $T'_g\to T_g$, so $T'_g$ close enough to $T_g$ gives
(Kick4). The remaining hypotheses hold as in
the case $k+1<n$, and Proposition~\ref{prop:kick} gives the claimed bound.
\end{proof}

\subsection{Summability of the relay times}
\label{sec:proof-time}

The growth time increases only logarithmically as the activation amplitude
decreases, while conversion to physical time contributes the geometric
factor $A_k^{-1}$. Their sum bounds the total relay duration uniformly in
the number of rings.

\begin{lemma}
\label{lem:time}
Assume (SA) and (Det). Put $\Tz=\Tpin+1+\sigma'^{-1}\log_+(\delta_0/\Klo)$ and
\[
\Tstar=\frac1{A_1}\Big[\frac{\mu}{\mu-1}\Big(\Tz+\frac{5}{\sigma'}\log\frac1{\sbar}\Big)
+\frac{5\log(2\mu)}{\sigma'}\cdot\frac{\mu^2}{(\mu-1)^2}+1+\frac{1}{\sigma'}\log\frac{2\delta_0}{\kseed}\Big].
\]
If the clock is defined up to level $1$, $\ksh_k\ge\Klo s_k^5$ for $k<n$,
and $\ksh_n\ge\frac12\kseed$, then
\begin{equation}
\label{eq:relay-total-time}
T_1\le\Tstar.
\end{equation}
\end{lemma}

\begin{proof}
By Lemma~\ref{lem:detuning} and (SA1), $\sigma_k\ge\frac{15}{16}\sigma>\frac34\sigma>\sigma'$
for every $k$. Since the clock is defined up to level $1$, $\ksh_k\in(0,\delta_0)$,
so $\log(\delta_0/\ksh_k)>0$ and $\hat T_k\le\sigma'^{-1}\log(\delta_0/\ksh_k)$. For
$k<n$, the hypothesis gives
$\hat T_k\le\sigma'^{-1}\big[\log_+(\delta_0/\Klo)+5\log(1/s_k)\big]$, and by
Definition~\ref{def:clock}, with $A_{k+1}=\mu A_k\ge A_k$,
\[
T_k-T_{k+1}=\frac{\Tpin}{A_{k+1}}+\frac{\hat T_k}{A_k}\le\frac1{A_k}\Big[\Tz+\frac{5}{\sigma'}\log\frac1{s_k}\Big].
\]
For $k=n$, $\tsh_n=0$ and $\ksh_n\ge\frac12\kseed$ give
$T_n=\hat T_n/A_n\le\sigma'^{-1}\log(2\delta_0/\kseed)/A_1$. We sum over $k$, using
$T_1=T_n+\sum_{k=1}^{n-1}(T_k-T_{k+1})$, $\log(1/s_k)=\log(1/\sbar)+k\log(2\mu)$,
$A_k=A_1\mu^{k-1}$, and
\[
\sum_{k\ge1}\mu^{1-k}=\frac\mu{\mu-1},\qquad\sum_{k\ge1}k\mu^{1-k}=\frac{\mu^2}{(\mu-1)^2}.
\]
All terms are nonnegative, so $T_1$ is at most $\Tstar$ minus $1/A_1$.
\end{proof}

Before a ring activates, its elapsed quiet time is bounded by the relay
intervals already reached. The next estimate uses only amplitudes assigned
by the observation time, so it remains available during the induction.

\begin{lemma}
\label{lem:quiet-time}
Assume (SA), (Det) and $\Kup\sbar^5\le\kseed$, and let $\Tz$ be as
in Lemma~\ref{lem:time}. Let $1\le j<n$ and put
\[
\bar\tau^\sharp_j=\frac{3\Tz+11\sigma'^{-1}\log(1/s_j)}{\mu-1}+\frac{5\mu\log(2\mu)}{\sigma'(\mu-1)^2}.
\]
Let $t\ge0$, and assume: for every $i\in(j,n]$ such that $\tsh_i$ is defined and
$\tsh_i\le t$, the clock is defined up to level $i$, and moreover
$\ksh_i\ge\Klo s_i^5$ if $i<n$; $t\le\tsh_j$ whenever $\tsh_j$ is defined;
and $\ksh_n\ge\frac12\kseed$. Then
\begin{equation}
\label{eq:relay-quiet-time}
A_jt\le\bar\tau^\sharp_j.
\end{equation}
\end{lemma}

\begin{proof}
Recall (Definition~\ref{def:clock}) that $\tsh_i$ is defined once the clock is
defined up to level $i+1$, and then $\tsh_{i+1}\le T_{i+1}\le\tsh_i$, as
$\hat T_{i+1}\ge0$. Let $i_*$ be the least $i\in(j,n]$ such that $\tsh_i$ is
defined and $\tsh_i\le t$; it exists as $\tsh_n=0$. For $i\in[i_*,n]$, $\tsh_i$ is
defined and $\tsh_i\le\tsh_{i_*}\le t$, so by hypothesis the clock is defined up to level $i$,
and $\ksh_i\ge\Klo s_i^5$ if $i<n$; in particular $T_{i_*}$ is defined. If
$i_*-1>j$, then $\tsh_{i_*-1}=T_{i_*}+\Tpin/A_{i_*}$ is defined and, by the
minimality of $i_*$, larger than $t$; if $i_*-1=j$, then
$t\le\tsh_j=T_{j+1}+\Tpin/A_{j+1}$ by hypothesis. In both cases, by
Definition~\ref{def:clock},
\[
t\le T_{i_*}+\frac{\Tpin}{A_{i_*}}=\frac{\hat T_n}{A_n}+\sum_{i=i_*}^{n-1}\Big(\frac{\Tpin}{A_{i+1}}+\frac{\hat T_i}{A_i}\Big)+\frac{\Tpin}{A_{i_*}}.
\]
We multiply by $A_j$ and use $A_j/A_i=\mu^{j-i}$. The pinning terms contribute at
most $\Tpin\sum_{i>j}\mu^{j-i}=\Tpin/(\mu-1)$. For $j<i<n$ with $i\ge i_*$, as in
the proof of Lemma~\ref{lem:time},
$\hat T_i\le\sigma'^{-1}[5\log(1/s_i)+\log_+(\delta_0/\Klo)]$ with
$\log(1/s_i)=\log(1/s_j)+(i-j)\log(2\mu)$; with $\sum_{i>j}\mu^{j-i}=\frac1{\mu-1}$
and $\sum_{i>j}(i-j)\mu^{j-i}=\frac\mu{(\mu-1)^2}$, these active intervals
contribute at most
\[
\frac{\sigma'^{-1}\log_+(\delta_0/\Klo)+5\sigma'^{-1}\log(1/s_j)}{\mu-1}+\frac{5\mu\log(2\mu)}{\sigma'(\mu-1)^2}.
\]
The seeded interval contributes $\mu^{j-n}\hat T_n\le\hat T_n$. As $c_N\le C_N$ and
$\frac1{22}\le5\mu e^{17\sigma\Tpin/16}$, $\Klo\le\Kup$, so
$\ksh_n\ge\frac12\kseed\ge\frac12\Kup\sbar^5\ge\frac12\Klo\sbar^5$ and
\[
\hat T_n\le\frac1{\sigma'}\Big[\log2+\log_+\frac{\delta_0}{\Klo}+5\log\frac1{\sbar}\Big]
\le\Tz+\frac{6}{\sigma'}\log\frac1{s_j},
\]
because $\log2\le\log(1/\sbar)\le\log(1/s_j)$, by $\sbar\le a/16<\frac12$ and
$s_j\le\sbar$. Adding the three groups and using $\mu-1\le1$,
\[
A_jt\le\frac{2\Tz+11\sigma'^{-1}\log(1/s_j)}{\mu-1}+\frac{5\mu\log(2\mu)}{\sigma'(\mu-1)^2}\le\bar\tau^\sharp_j .\qedhere
\]
\end{proof}

\subsection{Separation and confinement of the supports}
\label{sec:bootstrap-geometry}

We use the common parameter choices and bootstrap bounds of
Section~\ref{sec:relay}. A lower bound on each activation amplitude bounds
its active duration. The inner material therefore has only a limited time
to spread, and the separation of radii keeps it inside the next ring.
Only amplitudes assigned strictly before the current endpoint are needed.

\begin{lemma}
\label{lem:geometry}
Let $t\in[0,\Tend]$. Assume that \textup{(BA1)} holds at every $t'\in[0,t]$, and
that $\ksh_k\ge\Klo s_k^5$ for every level $k<n$ with $\tsh_k<t$. For a level
$k<n$ such that $T_{k+1}$ is defined, put $e_k=T_k$ if $T_k$ is defined and
$e_k=\tsh_k$ otherwise. Then:
\begin{enumerate}[label=\textup{(\alph*)}]
\item for every level $k<n$ such that $T_{k+1}$ is defined and $T_{k+1}\le t$, the
  hypotheses of Proposition~\ref{prop:confine} hold with $t_0=T_{k+1}$ and
  $t_1=\min(t,e_k)$; and for every $t'\in[T_{k+1},\min(t,e_k)]$,
  \begin{equation}
  \label{eq:geometry-confinement}
  \begin{gathered}
  D_{k+1}(t')\subset\overline B_{\varrho_k(t')}\subset B_{ar_k/8},\\
  \varrho_k(t')=\big(2+2\Cstar A_{k+1}(t'-T_{k+1})\big)r_{k+1},
  \end{gathered}
  \end{equation}
  \begin{equation}
  \label{eq:geometry-moment}
  |q_{k+1}(t')|\le\Cmom\big(1+A_{k+1}(t'-T_{k+1})\big)^3;
  \end{equation}
\item for every $t'\in[0,t]$ and every level $j$ that is not debris at $t'$, every
  ring $i<j$ is quiet at $t'$, and
  \begin{equation}
  \label{eq:geometry-quiet-support}
  S_i(t')\subset\Nflat_i\subset\R^2\setminus B_{6r_j}\qquad(i<j),
  \end{equation}
  \begin{equation}
  \label{eq:geometry-quiet-inner}
  D_{j+1}(t')\subset B_{ar_j/8}.
  \end{equation}
  Hence $h_j$ is defined at $t'$ and the support hypotheses of
  Lemma~\ref{lem:reduction} hold in ring-$j$ units. If $j<n$, then
  \begin{equation}
  \label{eq:geometry-field}
  \eta_j(t')=\|\hin_j(t')\|_{C^{M_h}(N')}+\|\hout_j(t')\|_{C^{M_h}(N')}\le1;
  \end{equation}
  the two terms are the norms of the inner velocity field minus its monopole
  and of the outer velocity field, in ring-$j$ units;
\item the localization conditions of Definition~\ref{def:loc} hold throughout:
  \begin{equation}
  \label{eq:geometry-loc}
  \textup{(Loc)}\text{ holds at }t'\qquad(t'\in[0,t]).
  \end{equation}
\end{enumerate}
\end{lemma}

\begin{proof}
Recall that $\Nflat_i=r_i\Nflat$, with the annulus $\Nflat$ of
Section~\ref{sec:II-prelim}. By (BA1), a ring $i$ that is not debris at $t'\in[0,t]$ has
$S_i(t')\subset\Nflat_i\subset\{\frac58ar_i\le|x|\le\frac74r_i\}$; and
$s_k\le a/32$ for every $k$ by Equation~\eqref{eq:sk-small}.

\emph{Outer rings in \textup{(b)}.} Let $j$ be not debris at $t'\in[0,t]$. Every
ring $i<j$ is quiet at $t'$ (Lemma~\ref{lem:clock}(1)), so
$S_i(t')\subset\{|x|\ge\frac58ar_i\}$, and
$\frac58ar_i\ge\frac58ar_{j-1}=\frac{5a}{8s_{j-1}}r_j\ge20r_j$; hence
$S_i(t')\cap B_{6r_j}=\emptyset$.

The material inside ring $j$ is controlled from the inside out: the confinement
of $D_{j+1}$ after the time $T_{j+1}$, at which the active interval of ring $j+1$
ends, starts from the confinement of $D_{j+2}$ at $T_{j+1}$. So we argue by
downward induction on $j=n,\dots,1$ and show that $D_{j+1}(t')\subset
B_{ar_j/8}$ at every $t'\in[0,t]$ at which ring $j$ is not debris, and, if
$j<n$, that (a) holds for $k=j$. For $j=n$, $D_{n+1}=\emptyset$. Let $j<n$ and
assume the claim for $j+1$.

\emph{\textup{(a)} for $k=j$.} Let $T_{j+1}\le t$ be defined and
$t_1=\min(t,e_j)$; the clock is defined up to level $j+1$. We check the
hypotheses of Proposition~\ref{prop:confine}.
(i): ring $j+1$ is not debris at $T_{j+1}$, so
$S_{j+1}(T_{j+1})\subset\overline B_{7r_{j+1}/4}$, and
$D_{j+2}(T_{j+1})\subset B_{ar_{j+1}/8}$ by induction; so
$D_{j+1}(T_{j+1})\subset\overline B_{2r_{j+1}}$. (ii): for
$t''\in[T_{j+1},t_1]$ and $i\le j$, ring $i$ is not debris at $t''$, since
$t''\le e_j$ and $T_i\ge T_j$ whenever $T_i$ is defined
(Lemma~\ref{lem:clock}(1)); so $S_i(t'')\subset\{|x|\ge\frac a2r_i\}$. (iii) is
(B1) and (B2). (iv): we bound the kick time $A_{j+1}(t_1-T_{j+1})$ of ring $j$
at $t_1$. If $t_1\le\tsh_j$, then $A_{j+1}(t_1-T_{j+1})\le\Tpin$, as
$\tsh_j=T_{j+1}+\Tpin/A_{j+1}$ (Definition~\ref{def:clock}). If $t_1>\tsh_j$,
then $e_j=T_j$, so the clock is defined up to level $j$, and $\tsh_j<t$. So
$\ksh_j\ge\Klo s_j^5$ by hypothesis, and $\sigma_j\ge\sigma'$ by
Lemma~\ref{lem:detuning}. As $t_1\le T_j=\tsh_j+\hat T_j/A_j$ and
$A_{j+1}=\mu A_j$, $A_{j+1}(t_1-T_{j+1})\le\Tpin+\mu\hat T_j$. It remains to
bound $\hat T_j$, the length of the active interval of ring $j$ in ring-$j$ time: as
$\ksh_j<\delta_0$ (the clock is defined up to level $j$), $\Klo>0$ (Definition~\ref{def:K1}; $c_N>0$ by Lemma~I.\ref{I-lem:uniform}) and
$\sigma_j\ge\sigma'>0$ (SA1),
$\hat T_j=\sigma_j^{-1}\log(\delta_0/\ksh_j)\le\sigma'^{-1}\log\big(\delta_0/(\Klo s_j^5)\big)\le\sigma'^{-1}\big(5\log\frac1{s_j}+\log_+\frac{\delta_0}{\Klo}\big)$.
In both cases, as $\Tpin<\mu(\Tpin+1)$,
\[
A_{j+1}(t_1-T_{j+1})<\mu\Big[\Tpin+1+\frac1{\sigma'}\Big(5\log\frac1{s_j}+\log_+\frac{\delta_0}{\Klo}\Big)\Big],
\]
so $\Lambda=2+2\Cstar A_{j+1}(t_1-T_{j+1})<\Lconf_j$, with $\Lconf_j$ of (B5).
The two inequalities of (iv) increase with $\Lambda$, so they follow from (B5).
Proposition~\ref{prop:confine}(a), (b) give, for $t'\in[T_{j+1},t_1]$,
$D_{j+1}(t')\subset\overline B_{\varrho_j(t')}$ and Equation~\eqref{eq:geometry-moment},
and $\varrho_j(t')\le\Lambda r_{j+1}<\Lconf_js_jr_j\le\frac a8r_j$ by (B5).

\emph{The claim on $D_{j+1}$.} Let ring $j$ be not debris at $t'\in[0,t]$. If
ring $j+1$ is not debris at $t'$, then $S_{j+1}(t')\subset\overline
B_{7r_{j+1}/4}$, $D_{j+2}(t')\subset B_{ar_{j+1}/8}$ by induction, and
$2r_{j+1}=2s_jr_j\le\frac a{16}r_j$; so $D_{j+1}(t')\subset B_{ar_j/8}$.
Otherwise $T_{j+1}<t'\le t$, and $t'\le e_j$ (if $T_j$ is defined, this is
$t'\le T_j$; otherwise $j=k_{\mathrm s}$ and $t'\le\Tend=\tsh_j$); so (a) for
$k=j$ gives the claim. This completes the induction.

\emph{End of \textup{(b)}.} Let $j$ be not debris at $t'\in[0,t]$. Then
$D_{j+1}(t')\subset B_{ar_j/8}$ and the rings $i<j$ do not meet $B_{6r_j}$, so
$h_j$ is defined at $t'$ (Definition~\ref{def:regions}), and in ring-$j$ units
the inner part is supported in $D_{j+1}(t')/r_j\subset B_{a/8}$ and the outer part in
$\R^2\setminus B_6$, as Lemma~\ref{lem:reduction} requires. If $j<n$, the
hypotheses of Lemma~\ref{lem:field-bound} hold at $t'$, since (B1) holds and
$S_i(t')\subset\Nflat_i\subset N_i$ for $i<j$; so $\eta_j(t')\le
\Cvel (\mu\sbar^2+\sbar)\le1$ by (B9).

\emph{\textup{(c)}.} Item (1) of Definition~\ref{def:loc} holds because a ring
$i$ with $t'\le T_i$, or with $T_i$ not defined, is not debris at $t'$, so
$S_i(t')\subset\Nflat_i\subset N_i$. Item (2) concerns a ring $j$ quiet at
$t'$. Its first two conditions, on the rings $i<j$ and on $D_{j+1}(t')$, hold
by (b). Its last clause, the confinement $D_{j+1}(t')\subset\overline
B_{\varrho_j(t')}$ when $T_{j+1}$ is defined and $t'\ge T_{j+1}$, holds by (a)
with $k=j<n$, since a ring $j$ quiet at $t'$ is not debris, so $t'\le e_j$
as above; for $j=n$ it holds as $D_{n+1}=\emptyset$.
\end{proof}

\subsection{Verification of quiet-ring control}
\label{sec:bootstrap-quiet-control}

We verify the hypotheses of the quiet-ring theorem used in
Proposition~\ref{prop:relay-step}. Fix its target time $t$, assume (BA)$(t)$,
and retain the conclusions of its preceding relay steps. Let $b\le t$ be an
endpoint through which Lemma~\ref{lem:geometry} has established
Equations~\eqref{eq:geometry-quiet-support}--\eqref{eq:geometry-field}.
Let $t'\in[0,b]$,
and let $j<n$ be quiet at $t'$. This includes earlier activation times
needed in the hand-off calculation. For Theorem~\ref{thm:quiet}, the
local ring-$i$ coordinate is $s=A_it''$, and its inputs are the pullbacks
\[
\begin{aligned}
\widetilde w_i(s)&=w_i(s/A_i),&\qquad \widetilde h_i(s)&=h_i(s/A_i),\\
\widetilde\eta_i(s)&=\eta_i(s/A_i),&
\widetilde{\maj}_i(s)&=\maj_i(s/A_i).
\end{aligned}
\]
In particular, $s\in[0,A_it']$ corresponds to the physical time interval $[0,t']$.
If $t'=0$, Equation~\eqref{eq:boot-quiet} holds since every ring $i\le j<n$
starts at $\Phi$ in ring-$i$ units, with $\maj_i(0)=0$ and $w_i(0)=0$
(Definition~\ref{def:approx}). Let $t'>0$. Suppose first that $\tsh_j$ is defined;
then the clock is defined up to level $j+1$ and $t'\le\tsh_j$. We apply
Lemma~\ref{lem:O} to $j$ and $t'$: at every $t''\in[0,t']$ every ring $i\le j$
is quiet. The outer and inner support conditions are
Equations~\eqref{eq:geometry-quiet-support} and~\eqref{eq:geometry-quiet-inner},
and $\eta_i(t'')\le1$ is Equation~\eqref{eq:geometry-field}.
Equation~\eqref{eq:BA-quiet} and (B9) give, since $\varepsilon_j=\mu s_j^5$ and
$s_j\le\sbar$,
\begin{equation}
\label{eq:boot-quiet-smallness}
\maj_j(t'')\le2\Abar\varepsilon_j\le2\Abar\mu\sbar^5\le\majst
\qquad(0\le t''\le t');
\end{equation}
and
(B7) gives the three smallness conditions on $\sbar$ in Lemma~\ref{lem:O}.
That lemma gives Equations~\eqref{eq:O-comparison} and~\eqref{eq:O-feedback}
on $[0,t']$ and verifies the hypotheses of Theorem~\ref{thm:quiet} for every
$i\le j$, after rescaling to ring-$i$ time. If $\tsh_j$ is not defined, apply
Lemma~\ref{lem:O} instead to the innermost quiet ring at $t'$. Its activation
time was assigned when the ring just inside it became active, and its
outer rings include ring $j$. The conclusions of Lemma~\ref{lem:O} therefore
give Theorem~\ref{thm:quiet} for ring $j$ on the physical time interval $[0,t']$.
In both cases Theorem~\ref{thm:quiet} gives the bounds~\eqref{eq:boot-quiet}
for every $t''\in[0,t']$.

\subsection{Verification at a hand-off}
\label{sec:bootstrap-handoff-checks}

At a hand-off $\tsh_k\le t$, $k<n$, Equations~\eqref{eq:geometry-confinement}--\eqref{eq:geometry-loc}
and~\eqref{eq:boot-quiet} hold through $\tsh_k$ on their stated intervals.
The preceding relay step has supplied item~(c) of
Proposition~\ref{prop:relay-step} for ring $k+1$ through $T_{k+1}$, and the
hand-offs at the inner levels have already been performed.
The clock is defined up to level $k+1$, and no ring is active on
$(T_{k+1},\tsh_k)$.
By the hand-off argument for ring $k+1$ and
(B6), or by Equation~\eqref{eq:boot-seed-data} if $k+1=n$,
\begin{equation}
\label{eq:boot-kappa}
0<\ksh_{k+1}\le2\kseed .
\end{equation}
We apply Lemma~\ref{lem:kicktrue} with the constants $C_0$ and $K_1$, and check
its hypotheses in the order in which it lists them. First,
$D_{k+1}\subset B_{ar_k/8}$ on $[0,\tsh_k]$ by
Equation~\eqref{eq:geometry-quiet-inner}. Second, item~(c) of
Proposition~\ref{prop:relay-step} for ring $k+1$ supplies the conclusions of
Theorem~\ref{thm:active} and Corollary~\ref{cor:active-moment}, including
Equations~\eqref{eq:relay-active} and~\eqref{eq:relay-active-moment}, on
$\Phase_{k+1}$ with $\eta_{b,k+1}$. They hold also with
$\eta_b\ge\eta_{b,k+1}$, since the bounds increase with $\eta_b$.
Third, the same item verifies hypotheses (i)--(iv) of
Lemma~\ref{lem:drift} for ring $k+1$ with $\bar t=T_{k+1}$ and
$\eta=\eta_{b,k+1}\in[0,\eta_b]$, together with
Equation~\eqref{eq:relay-drift-small}. Fourth, if $k+1<n$: on $[0,\tsh_{k+1}]$
$D_{k+2}\subset B_{ar_{k+1}/8}$ and the rings $i\le k$ lie outside $B_{6r_{k+1}}$
by Equations~\eqref{eq:geometry-quiet-inner} and~\eqref{eq:geometry-quiet-support};
the hypotheses of Theorem~\ref{thm:quiet} for ring $k+1$ on the physical time
interval $[0,\tsh_{k+1}]$ were checked in
Appendix~\ref{sec:bootstrap-quiet-control}; $\maj_{k+1}(\tsh_{k+1})\le K_1\ksh_{k+1}$ is
Equation~\eqref{eq:handoff} for ring $k+1$; the hypotheses of
Proposition~\ref{prop:confine} with $k+1$ in place of $k$ and $t_1=T_{k+1}$ hold
by Lemma~\ref{lem:geometry}(a) for $k+1$; and
$\ksh_{k+1}\ge\Klo s_{k+1}^5$ by Equation~\eqref{eq:relay-amplitudes}. Fifth, the bound
on $q_{k+1}$ for $s\in[T_{k+1},\tsh_k]$ is
Equation~\eqref{eq:geometry-moment} on $\Pin_k$. Finally, conditions (a)--(f) of
Lemma~\ref{lem:kicktrue} are (P1), (P2), ($\kappa$1), ($\kappa$2), ($\kappa$3) at
$\kappa=\ksh_{k+1}$, which Equation~\eqref{eq:boot-kappa} allows, and (B11). So
\begin{equation}
\label{eq:boot-qstar}
|\qs|\in[\tfrac{15}{16},\tfrac{17}{16}]2\pi|\pout_*(G_{k+1})|\delta_0,
\end{equation}
and
\begin{equation}
\label{eq:boot-kick}
\frac{|\beta_k|\,|\qs|}{10|\lamref|(\mu-1)}\le|\Ain(\Tpin)|\le\frac{3|\beta_k|\,|\qs|\,\mu}{|\lamref|(\mu-1)} .
\end{equation}

Next we apply Proposition~\ref{prop:relay} with the constants $C_0$ and $K_1$,
$j=k$ and endpoint $\tsh_k$ in physical time. By
Lemma~\ref{lem:quiet-time} at the time $\tsh_k$ (its hypothesis
$\Kup\sbar^5\le\kseed$ is (B6); its hypotheses on the levels
$i\in(k,n)$ hold by the hand-off argument for these levels, and $\ksh_n\ge\frac12\kseed$ by
Equation~\eqref{eq:boot-seed-data}) and Lemma~\ref{lem:theta-quiet},
\begin{equation}
\label{eq:boot-theta}
\theta_k A_k\tsh_k\le\theta_k\bar\tau^\sharp_k\le\min\Big\{\frac12,\frac1{4\Cls K_1}\Big\}.
\end{equation}
The hypotheses (H1)--(H4) of Proposition~\ref{prop:relay} hold. (H1): (Loc) holds on
$[0,\tsh_k]$ by Equation~\eqref{eq:geometry-loc}. (H2):
Equation~\eqref{eq:O-feedback} holds on $[0,\tsh_k]$ by the application of
Lemma~\ref{lem:O} in Appendix~\ref{sec:bootstrap-quiet-control}, and $\theta_k A_k\tsh_k\le\frac12$ by
Equation~\eqref{eq:boot-theta}. (H3): if $k+1<n$, the
conclusions of Theorem~\ref{thm:quiet} for ring $k+1$ hold on
$[0,\tsh_{k+1}]$, after rescaling, by the verification in
Appendix~\ref{sec:bootstrap-quiet-control}.
(H4): $\maj_{k+1}(\tsh_{k+1})\le K_1\ksh_{k+1}$ by
the hand-off argument for ring $k+1$; if $k+1=n$, then $\tsh_n=0$ and $\maj_n(0)=0$. The
other conditions of (H4), on $\Phase_{k+1}$, hold: $D_{k+2}\subset B_{ar_{k+1}/8}$ by Equation~\eqref{eq:geometry-quiet-inner};
(BA3) for ring $k+1$ holds by (BA)$(t)$, and at $\tsh_{k+1}$ by
Equation~\eqref{eq:boot-switch} for ring $k+1$ or by
Equation~\eqref{eq:boot-start}; and
Equation~\eqref{eq:relay-inner-field} holds on $\Phase_{k+1}$ by the
preceding relay step. Hence, by Proposition~\ref{prop:relay}(ii), Equation~\eqref{eq:relay}
holds at $\tsh_k$.

We bound the right side of Equation~\eqref{eq:relay}. Its last term
$\Xi_{k+1}$ contains the driving before ring $k+1$ activates and the
transient and polynomial inner-field contributions during its active interval. We bound it by ($\kappa$7), after bounding the exponent
$1-\gamma_k$ from below. By
Equation~\eqref{eq:detuning-rates} and the proof of
Lemma~\ref{lem:detuning}, $(1-e_\mu)\sigma\le\sigma_i\le(1+e_\mu)\sigma$
for $i\in\{k,k+1\}$, $e_\mu(\mu+1)\le\frac14(\mu-1)$, $\frac{15}{16}\le1-e_\mu\le1$
and $x_k\le\frac{17}{16}\sigma$. Hence
\begin{gather*}
1-\frac{\sigma_k}{\mu\sigma_{k+1}}\ge1-\frac{1+e_\mu}{\mu(1-e_\mu)}
=\frac{(\mu-1)-e_\mu(\mu+1)}{\mu(1-e_\mu)}\ge\frac34\cdot\frac{\mu-1}\mu,\\
1-\frac{\sigma'}{\sigma_{k+1}}\ge1-\frac{16\sigma'}{15\sigma}=c_\gamma ,
\end{gather*}
and, as $\gamma_k$ is the larger of $\sigma_k/(\mu\sigma_{k+1})$ and
$\sigma'/\sigma_{k+1}$,
\begin{equation}
\label{eq:boot-remainder-exponent}
1-\gamma_k\ge\min\Big(\tfrac34\tfrac{\mu-1}\mu,c_\gamma\Big)\ge\gamma'_\mu.
\end{equation}
As $\ksh_{k+1}<\delta_0$,
$(\ksh_{k+1}/\delta_0)^{1-\gamma_k}\le(\ksh_{k+1}/\delta_0)^{\gamma'_\mu}$, and
($\kappa$7) at $\kappa=\ksh_{k+1}$ gives
\begin{equation}
\label{eq:boot-relay-remainder}
\Xi_{k+1}\le\frac{\delta_0}{\sigb(\mu-1)}.
\end{equation}
So
the bracket in Equation~\eqref{eq:relay} with $j=k$ is at most
$(4C_2+1)\delta_0/(\sigb(\mu-1))+C_D(\Tpin)$; with Equation~\eqref{eq:relay}
this gives the first inequality below, and
Equations~\eqref{eq:detuning-rates} and~\eqref{eq:K1} give the second:
\begin{equation}
\label{eq:boot-quiet-activation}
\maj_k(\tsh_k)\le2\Cvel \varepsilon_ke^{x_k\Tpin}\Big[\frac{(4C_2+1)\delta_0}{\sigb(\mu-1)}+C_D(\Tpin)\Big]\le\Abar\,\varepsilon_k .
\end{equation}
Next we bound from below the coefficient $\cin_k$ produced by the kick. By
Lemma~\ref{lem:coupling}\textup{(2)}, $\cin_k=e^{\zeta_0\Tpin}\Ain(\Tpin)$ with
$\Re\zeta_0=\sigma_k/\mu=x_k$, so that $|\cin_k|=e^{x_k\Tpin}|\Ain(\Tpin)|$. As
$G_k,G_{k+1}\in[0,\Gone]$, Lemma~I.\ref{I-lem:uniform} gives the first bounds below,
and they give the second with $\beta_k=\beta(G_k)s_k^5/(2\pi)$ and the bound
in Equation~\eqref{eq:boot-qstar}:
\begin{equation}
\label{eq:boot-coupling-size}
c_N\le|\beta(G_k)||\pout_*(G_{k+1})|\le C_N,\qquad
\tfrac{15}{16}c_Ns_k^5\delta_0\le|\beta_k||\qs|\le\tfrac{17}{16}C_Ns_k^5\delta_0 .
\end{equation}
With the lower bound in Equation~\eqref{eq:boot-kick}, $\varepsilon_k=\mu s_k^5$
and $\frac{15}{16}\cdot\frac1{10}=\frac3{32}$,
\begin{equation}
\label{eq:boot-driven-lower}
|\cin_k|\ge e^{x_k\Tpin}\frac{|\beta_k||\qs|}{10|\lamref|(\mu-1)}\ge\frac{3c_N}{32\mu}\cdot\frac{\varepsilon_k\delta_0e^{x_k\Tpin}}{|\lamref|(\mu-1)}.
\end{equation}
Multiplying by $\frac12K_1$, Equation~\eqref{eq:K1} and the first inequality of
Equation~\eqref{eq:boot-quiet-activation} give
\begin{equation}
\label{eq:boot-driven-dominates}
\frac12K_1|\cin_k|\ge4\Cvel \varepsilon_ke^{x_k\Tpin}\Big[\frac{(4C_2+1)\delta_0}{\sigb(\mu-1)}+C_D(\Tpin)\Big]\ge2\maj_k(\tsh_k).
\end{equation}

The hypotheses of Lemma~\ref{lem:handoff} hold. For (i), the support
conditions are Equations~\eqref{eq:geometry-quiet-support} and~\eqref{eq:geometry-quiet-inner},
and the hypotheses of Theorem~\ref{thm:quiet} were checked in
Appendix~\ref{sec:bootstrap-quiet-control}. Condition (ii) is given by
Equations~\eqref{eq:O-feedback} and~\eqref{eq:boot-theta}; (iii) is
Equation~\eqref{eq:boot-driven-dominates}; and (iv) follows from
Equation~\eqref{eq:boot-quiet-activation}, $\varepsilon_k\le\mu\sbar^5$ and (B11). So Equation~\eqref{eq:handoff} holds at $\tsh_k$. To bound $\ksh_k$
we use its first bound, the identity
$|\cin_k|=e^{x_k\Tpin}|\Ain(\Tpin)|$ and
Equations~\eqref{eq:boot-kick} and~\eqref{eq:boot-coupling-size}. As $e^{x_k\Tpin}\ge1$ and $\frac{15}{320}\ge\frac1{22}$,
\begin{equation}
\label{eq:boot-amplitude-lower}
\ksh_k\ge\tfrac12|\Ain(\Tpin)|\ge\tfrac12\cdot\tfrac1{10}\cdot\tfrac{15}{16}\cdot\frac{c_N\delta_0s_k^5}{|\lamref|(\mu-1)}\ge\Klo s_k^5;
\end{equation}
and as $x_k\le\frac{17}{16}\sigma$ and
$\frac32\cdot3\cdot\frac{17}{16}=\frac{153}{32}\le5$,
\begin{equation}
\label{eq:boot-amplitude-upper}
\ksh_k\le\frac32e^{x_k\Tpin}|\Ain(\Tpin)|\le\frac32\cdot3\cdot\frac{17}{16}\cdot\frac{C_N\delta_0\mu\,s_k^5}{|\lamref|(\mu-1)}e^{17\sigma\Tpin/16}\le\Kup s_k^5.
\end{equation}
Hence
$0<\ksh_k\le\Kup\sbar^5\le\kseed<\delta_0$ by (B6) and ($\kappa$1), and
$S_k(\tsh_k)\subset\NNm_k\subset N_k$ by Equation~\eqref{eq:boot-quiet}; so the
clock is defined up to level $k$. Since $s\le s_h=M-1$, $C_0\ge K_1C_Q$,
$C_2\ge1$ and $\rhohat_k(\tsh_k)=C_0$, these inequalities also give the
activation estimate~\eqref{eq:boot-switch}.
It gives Equation~\eqref{eq:improved-active} at activation, using no time
after $\tsh_k$.

\subsection{Verification on an active interval}
\label{sec:bootstrap-active-checks}

We retain the target time $t$ and hypothesis (BA)$(t)$ of
Proposition~\ref{prop:relay-step}. Equation~\eqref{eq:relay-seed} or~\eqref{eq:relay-startup} supplies
the initial support, norm and amplitude bounds for ring $k$.

This step verifies the hypotheses
(Init), (Field), (Coef) and (Small) of Theorem~\ref{thm:active} for ring $k$ on
$\Phase_k\cap[0,t]$. Let $t'\in\Phase_k\cap[0,t]$ and
$\bar u=A_k(t'-\tsh_k)\le\hat T_k$. If $t'=\tsh_k$, we need only the two bounds
at the end of this step for $t''=\tsh_k$, and they hold already: $S_k(\tsh_k)\subset\NNm_k\subset\NN_k$ by Equation~\eqref{eq:boot-seed-data}
if $k=n$ and by Equation~\eqref{eq:handoff} if $k<n$, and the bound on
$\|w_k(\tsh_k)\|_{H^s}$ is Equation~\eqref{eq:boot-start}, respectively
\eqref{eq:boot-switch}. So let $t'>\tsh_k$; then $\bar u>0$, and $T_k$ is defined
by Equation~\eqref{eq:boot-seed-data} or~\eqref{eq:relay-startup} and the clock definition. By the invariance of Equation~\eqref{eq:euler} under
translations in time and the scaling of ring-$k$ units, and by
Theorem~\ref{thm:smooth}(5), the solution in ring-$k$ units with time origin
$\tsh_k$ is the solution of Theorem~\ref{thm:smooth} with datum
$A_k^{-1}\omega(\tsh_k,r_k\cdot)$. We apply Lemma~\ref{lem:reduction} and
Theorem~\ref{thm:active} to it, with the parts of the vorticity on $S_k$, on
$D_{k+1}$ and on $S_1\cup\dots\cup S_{k-1}$, with $G=G_k$, $\kappa=\ksh_k$, the
constant $C_0$, $\eta_b=\eta_{b,k}$, $\delta_0$ and $\bar t=\bar u$.
Writing the theorem's time variable as the local coordinate
$u=A_k(t''-\tsh_k)$, its functions are explicitly
\[
w(u)=w_k(\tsh_k+u/A_k),\qquad
h(u)=h_k(\tsh_k+u/A_k),\qquad
c(u)=c_k(\tsh_k+u/A_k).
\]
The support
hypotheses of Lemma~\ref{lem:reduction} on $[0,\bar u]$ are
Equations~\eqref{eq:geometry-quiet-support} and~\eqref{eq:geometry-quiet-inner}. For Theorem~\ref{thm:active}: $G_k\in[0,\Gone]$,
$\delta_0\le\dH$ by (P1), $\eta_{b,k}\le\eta_b\le1$ by (B8), and (Small)
holds: $\ksh_k\le2\kseed\le\kappa_0(C_0)$ by Equation~\eqref{eq:boot-seed-data} if $k=n$, by Equation~\eqref{eq:relay-startup} and (B6) if
$k<n$, and by ($\kappa$4). It remains to check (Init), (Field) and (Coef).

(Init): for $k<n$, Equation~\eqref{eq:handoff} gives
$\|w_k(\tsh_k)\|_{H^{M-1}}\le K_1C_Q\ksh_k\le C_0\ksh_k$ and
$\supp(\Phi+w_k(\tsh_k))\subset\NNm$; for $k=n$, these are the seed estimates
\eqref{eq:boot-seed-data} in ring-$n$ units.

(Field), inner part $h_1(u)=\hin_k(\tsh_k+u/A_k)$. For $k=n$, $h_1=0$, because
$D_{n+1}=\emptyset$ and $G_n=0$. Let $k<n$, $u\in[0,\bar u]$ and
$t''=\tsh_k+u/A_k$. By Equation~\eqref{eq:geometry-confinement}, $\omega_{\ge k+1}(t'')$ is
supported in $\overline B_R$ with $R=\varrho_k(t'')\le\frac a8r_k$,
$A_{k+1}(t''-T_{k+1})=\Tpin+\mu u$, and $D_{k+1}(t'')\subset B_{ar_k/8}$, so
$h_1(u)$ is the field $E_{[k]}$ of Lemma~\ref{lem:inner-field} at $t''$; and
$\|\omega_{\ge k+1}(t'')\|_{L^1}\le\Cstar A_{k+1}r_{k+1}^2$ by
Lemma~\ref{lem:speed}. Lemma~\ref{lem:inner-field}, the inequalities
$2+2\Cstar\Tpin+2\Cstar\mu u\le(2+2\Cstar\Tpin)(1+\Cstar\mu u)$ and
$1+\Cstar\mu u\le(1+\Cstar\mu)\max(1,u)\le(1+2\Cstar)\jap{u}$, the second inequality of (P4), and Equation~\eqref{eq:boot-amplitude-lower} give
\begin{equation}
\label{eq:boot-inner}
\begin{aligned}
\|h_1(u)\|_{C^{M_h}(N')}&\le\Cvel \mu s_k^2(R/r_k)^3=\Cvel \mu s_k^5\big(2+2\Cstar(\Tpin+\mu u)\big)^3\\
&\le\Cvel \mu(2+2\Cstar\Tpin)^3(1+2\Cstar)^3s_k^5\jap{u}^3\\
&\le C_0\Klo s_k^5\jap{u}^3\le C_0\ksh_k\jap{u}^3 .
\end{aligned}
\end{equation}
As $M\le M_h$, this is the first bound of (Field).

(Field), outer part $h_2(u)=\hout_k(\tsh_k+u/A_k)$. We apply
Lemma~\ref{lem:Oprime} with the constants $C_0$, $K_1$ at the time $t'$.
The clock is defined up to level $k$, and (Loc) holds on $[0,t']$ by
Equation~\eqref{eq:geometry-loc}. If $k\ge2$, ring $k-1$ is quiet on $[0,t']$,
so the verification in Appendix~\ref{sec:bootstrap-quiet-control} gives the
hypotheses and conclusions of Lemma~\ref{lem:O} for $j=k-1$ and endpoint
$t'$ in physical time; and $\theta_{k-1}A_{k-1}t'\le\theta_{k-1}A_{k-1}\tsh_{k-1}\le\frac12$
(as ring $k-1$ is quiet at $t'$) by Lemma~\ref{lem:quiet-time} at the time
$\tsh_{k-1}$ and Lemma~\ref{lem:theta-quiet}; the hypotheses of
Lemma~\ref{lem:quiet-time} besides (B6) concern the levels $i\ge k$, whose amplitudes are
bounded below by the hand-off argument at the times $\tsh_i\le\tsh_k$ and, for $i=n$, by
Equation~\eqref{eq:boot-seed-data}. (BA3) for ring $k$ holds on $(\tsh_k,t']$ by (BA)$(t)$ and at
$\tsh_k$ by Equation~\eqref{eq:boot-switch} or \eqref{eq:boot-start}; the inner
bound is Equation~\eqref{eq:boot-inner}. If $k<n$, the conclusions of
Theorem~\ref{thm:quiet} for ring $k$, after rescaling, on the physical time
interval $[0,\tsh_k]$ hold by the verification in
Appendix~\ref{sec:bootstrap-quiet-control},
$\maj_k(\tsh_k)\le K_1\ksh_k$ is Equation~\eqref{eq:relay-startup}, and the
inner support condition is Equation~\eqref{eq:geometry-confinement} on
$[\tsh_k,t']$. Equation~\eqref{eq:etab} therefore gives
\begin{equation}
\label{eq:boot-outer}
\|h_2(u)\|_{C^{M_h}(N')}\le\eta_{b,k}\hat a_k(\tsh_k+u/A_k)
=\eta_{b,k}\ksh_ke^{\sigma_ku}\qquad(0\le u\le\bar u),
\end{equation}
which is the second bound of (Field).

(Coef): we apply Lemma~\ref{lem:drift} to ring $k$ with endpoint in physical time
$\bar t=t'$ and
$\eta=\eta_{b,k}$. Hypotheses (i) and (ii) concern only $k<n$. (i) is
Equation~\eqref{eq:geometry-moment} on $W_k$, since the kick time is
$\tauK=\Tpin+\mu u$.
(ii), with $q'=\qs$ of the hand-off argument: the hypotheses of Lemma~\ref{lem:coupling}(2) hold
by Equations~\eqref{eq:geometry-quiet-support} and~\eqref{eq:geometry-quiet-inner},
so Equation~\eqref{eq:handoff} and
Lemma~\ref{lem:coupling}(2) give
$\ksh_k\ge\frac12|\cin_k|=\frac12e^{x_k\Tpin}|\Ain(\Tpin)|$; Equation~\eqref{eq:boot-kick} gives the lower bound on $|\Ain(\Tpin)|$.
By Equation~\eqref{eq:boot-qstar} and Lemma~I.\ref{I-lem:uniform},
$|\qs|\ge\frac{15}{16}2\pi|\pout_*(G_{k+1})|\delta_0\ge\frac{15}{16}\cout\delta_0>0$,
so $\qs\ne0$, and (P2) gives
$\mu-1\le c'\delta_0\le|\qs|/(10|\lamref|\Cmom\Ip{3}{\sigma/4})$.
Condition (iii) follows from Equation~\eqref{eq:boot-outer}, since the $C^0$
norm is at most the $C^{M_h}$ norm. Condition (iv) is
Equations~\eqref{eq:geometry-quiet-support} and~\eqref{eq:geometry-quiet-inner}. Moreover $\rT(\Tpin)\le\frac1{64}$ by (P3), and
$\Cls\eta_{b,k}\hat T_k\le\frac1{32}$: for $k=1$, $\eta_{b,1}=0$; for $2\le k<n$,
$\hat T_k\le\sigma'^{-1}\big(5\log\frac1{s_k}+\log_+\frac{\delta_0}{\Klo}\big)$,
as $\sigma_k\ge\sigma'$ and $\ksh_k\ge\Klo s_k^5$, and (B8) applies; for
$k=n$, $\hat T_n\le\sigma'^{-1}\log(2\delta_0/\kseed)$ by Equation~\eqref{eq:boot-seed-data} and
$\eta_{b,n}\le\eta_b$, and (B8) applies. These are precisely
Equation~\eqref{eq:drift-smallness} with $\eta=\eta_{b,k}$.
Lemma~\ref{lem:drift} gives (Coef) on $[0,\bar u]$ and, if
$\bar u=\hat T_k$, Equation~\eqref{eq:drift-saturation}.

Theorem~\ref{thm:active} and Corollary~\ref{cor:active-moment} therefore hold on
$[0,\bar u]$, giving Equations~\eqref{eq:relay-active} and~\eqref{eq:relay-active-moment}.
In particular, Equations~\eqref{eq:improved-support} and~\eqref{eq:improved-active}
hold for ring $k$ at every $t''=\tsh_k+u/A_k\in[\tsh_k,t']$.

\subsection{Verification of the quiet-size improvement}
\label{sec:bootstrap-quiet-checks}

At $t'=0$ every quiet-ring Duhamel bound vanishes. Let $t'>0$ belong
to the current interval. On a pinning interval $\Pin_k$,
put $j=k$. On an active interval $\Phase_k$ with $k\ge2$ and $t'>\tsh_k$, put $j=k-1$.
The activation endpoint of ring $k$ is included in $\Pin_k$ for this
argument. There are no outer quiet rings during $\Phase_1$.
Thus ring $j$ is the innermost quiet ring, and $\tsh_j$ is already assigned.
We first bound $\maj_j(t')$ and then use Equation~\eqref{eq:O-comparison}
to control all $i<j$.

Apply Proposition~\ref{prop:relay} to $j$ through $t'$.
Equation~\eqref{eq:geometry-loc} gives (Loc). The application of
Lemma~\ref{lem:O} in Appendix~\ref{sec:bootstrap-quiet-control} gives
Equation~\eqref{eq:O-feedback} on $[0,t']$.
Lemmas~\ref{lem:quiet-time} and~\ref{lem:theta-quiet} give
$\theta_jA_jt'\le\theta_j\bar\tau^\sharp_j\le1/2$.
Their amplitude hypotheses concern only levels $i>j$ that have already
activated. During $\Pin_k$ these are inner levels $i>k$; during $\Phase_k$
they include ring $k$, whose amplitude has just been obtained at the hand-off.
Equation~\eqref{eq:relay-seed} supplies the amplitude of the seeded ring.
The verification of Theorem~\ref{thm:quiet} in
Appendix~\ref{sec:bootstrap-quiet-control} for ring $j+1$ supplies (H3), and
Equation~\eqref{eq:relay-startup} supplies the initial Duhamel bound in (H4).
For the seed the earlier interval is empty and $\maj_n(0)=0$.
The inner support condition in (H4) is Equation~\eqref{eq:geometry-quiet-inner}.
Its inner-field bound is Equation~\eqref{eq:relay-inner-field}, supplied by
the preceding relay step on $\Pin_k$ and by Equation~\eqref{eq:boot-inner}
on $\Phase_k$. Its active norm hypothesis follows from
Equation~\eqref{eq:BA-active} strictly after activation and from
Equation~\eqref{eq:boot-start} or~\eqref{eq:boot-switch} at activation.

In both cases $t'>\tsh_{j+1}$. Proposition~\ref{prop:relay}(ii) or (iii) gives
Equation~\eqref{eq:relay} as an upper bound.
Equations~\eqref{eq:boot-remainder-exponent} and~\eqref{eq:boot-kappa}, with
$j$ in place of $k$, and ($\kappa$7) give Equation~\eqref{eq:boot-relay-remainder}.
Using $x_j\le17\sigma/16$ from Equation~\eqref{eq:detuning-rates} in
Equation~\eqref{eq:relay} yields
\begin{equation}
\label{eq:boot-quiet-innermost}
\maj_j(t')\le2\Cvel\varepsilon_je^{17\sigma\Tpin/16}
\Big[\frac{(4C_2+1)\delta_0}{\sigb(\mu-1)}+C_D(\Tpin)\Big]
=\Abar\varepsilon_j.
\end{equation}
This argument applies at every time up to the current endpoint.

For every outer ring $i<j$, Equations~\eqref{eq:O-comparison} and~\eqref{eq:boot-quiet-innermost} give
\begin{equation}
\label{eq:boot-quiet-outer}
\begin{aligned}
\maj_i(t')&\le\tfrac34\Ptil_{i,j}\maj_j(t')\\
&\le\tfrac34\Abar\varepsilon_i
       \prod_{p=i+1}^{j}\KO\varepsilon_p
\le\tfrac38\Abar\varepsilon_i,
\end{aligned}
\end{equation}
because $\KO\varepsilon_p\le\KO\mu\sbar^5\le1/2$ by (B7).
This includes outer rings whose activation time has not yet been assigned.
Together with Equation~\eqref{eq:boot-quiet-innermost}, this proves
Equation~\eqref{eq:improved-quiet} for every quiet ring.

\section*{Statements and Declarations}
The author has no competing interests to declare.

\end{document}